\documentclass[lfeqn,12pt,oneside]{book}
\usepackage[T1]{fontenc}
\usepackage{titlesec, blindtext, color}
\definecolor{gray75}{gray}{0.75}
\newcommand{\hsp}{\hspace{20pt}}
\titleformat{\chapter}[block]{\Huge\bfseries}{\thechapter\hsp\textcolor{gray75}{|}\hsp}{0pt}{\Huge\bfseries}
\usepackage{fancyhdr}
 
\usepackage[utf8]{inputenc}
\usepackage[english]{babel}
\usepackage{amssymb,mathtools,mathabx,mathrsfs,bbm,pifont
}

\usepackage{lmodern,graphicx}
 
\usepackage[toc,page,header]{appendix}

\usepackage{minitoc}

\usepackage[nottoc]{tocbibind}
 
\usepackage{multirow}
\usepackage[multiple]{footmisc}
\usepackage{xcolor}
\usepackage[colorlinks=true]{hyperref}
\usepackage{chngcntr}
\counterwithout{footnote}{chapter}

\newcommand{\Giu}{{\bigskip\noindent}}
\newcommand{\nl}{{\smallskip\noindent}}
\newcommand{\noi}{{\noindent}}

\newtheorem{theorem}{Theorem}[section]
\newtheorem{definition}[theorem]{Definition}
\newtheorem{proposition}[theorem]{Proposition}
\newtheorem{lemma}[theorem]{Lemma}
\newtheorem{remark}[theorem]{Remark}
\newtheorem{sublemma}[theorem]{Sublemma}
\newtheorem{corollary}[theorem]{Corollary}
\newtheorem{assumption}[theorem]{Assumption}
\newtheorem{notationalrem}[theorem]{Notational Remark}
\newtheorem{tools}[subsection]{$\negsp\negsp$}
\newcommand\asm[1]{ \begin{assumption}\label{#1} }
\newcommand\easm{ \end{assumption} }
\newcommand\dfn[1]{ \begin{definition}\label{#1} }
\newcommand\dfntwo[2]{ \begin{definition}[#2]\label{#1} }
\newcommand\edfn{ \end{definition} }
\newcommand\rem[1]{ \begin{remark}\label{#1} \small \rm}
\newcommand\remtwo[2]{ \begin{remark}[#2]\label{#1} \rm}
\newcommand\erem{ \end{remark} }
\newcommand\thm[1]{ \begin{theorem}\label{#1}}
\newcommand\thmtwo[2]{ \begin{theorem}[#2]\label{#1}}
\newcommand\ethm{ \end{theorem} }
\newcommand\pro[1]{ \begin{proposition}\label{#1}}       
\newcommand\protwo[2]{ \begin{proposition}[#2]\label{#1}}
\newcommand\epro{ \end{proposition} }
\newcommand\lem[1]{ \begin{lemma}\label{#1}}
\newcommand\lemtwo[2]{ \begin{lemma}[#2]\label{#1}}
\newcommand\elem{ \end{lemma} }
\newcommand\sublem[1]{ \begin{sublemma}\label{#1}}
\newcommand\sublemtwo[2]{ \begin{sublemma}[#2]\label{#1}}
\newcommand\esublem{ \end{sublemma} }
\newcommand\cor[1]{ \begin{corollary}\label{#1}}
\newcommand\cortwo[2]{ \begin{corollary}[#2]\label{#1}}
\newcommand\ecor{ \end{corollary} }
\newcommand\notrem[1]{{{ \begin{notationalrem}\label{#1} }\sl}}
\newcommand\enotrem{ \end{notationalrem} }

\newcommand\average[1]{{ \left\langle #1 \right\rangle}}

\newcommand\equ[1]{{\rm (\ref{#1})}}
\newcommand\beq[1]{ \begin{equation}\label{#1} }
\newcommand{\eeq}{ \end{equation} }
\newcommand{\beqno}{ \[ }
\newcommand{\eeqno}{ \] }
\newcommand\beqa[1]{ \begin{eqnarray} \label{#1}}
\newcommand{\eeqa}{ \end{eqnarray} }
\newcommand{\beqano}{ \begin{eqnarray*} }
\newcommand{\eeqano}{ \end{eqnarray*} }

\newcommand{\proof}{\par\medskip\noindent{\bf Proof\ }}

\newcommand{\ie}{{\it i.e.\  }}

\newcommand{\etc}{{\it etc\ }}
\newcommand{\wrt}{{\it w.r.t\ }}

\newcommand{\resp}{{\it resp.\  }}

\newcommand{\dst}{\displaystyle}
\newcommand{\qed}{\hskip.5truecm
            \vrule width 1.7truemm height 3.5truemm depth 0.truemm
            \par\Giu}

\newcommand\ovl[1]{ \overline {#1} }

\newcommand\su[1]{ \frac{1}{ {#1}} }

\newcommand\dist{ {\, \rm dist\, }}
\newcommand\diam{ {\, \rm diam\, }}
\newcommand\minfoc{ {\, \rm minfoc\, }}

\newcommand\conv{ {\, \rm conv\, }}
\newcommand\supp{ {\, \rm supp\, }}

\newcommand\meas{ {\, \rm meas\, }}

\newcommand\sign{ {\, \rm sign\, }}

\newcommand{\diag}{{ \, \rm diag \, }}
\newcommand{\dom}{{ \, \rm dom \, }}

\newcommand{\adj}{ {\rm \, Adj \, }}
\newcommand{\tr}{ {\rm \, tr \, }}
\newcommand{\Id}{ {\rm Id }}
\newcommand{\Iso}{\rm Iso }

\newcommand{\io}{{\infty }}
\newcommand{\ci}{ {C^\infty}   }

\newcommand\igl[2]{{ \int_{
{#1}}^{#2}  }}

\newcommand{\dpr}{ {\partial}   }

\newcommand\eqby[1]{\stackrel{\equ{#1}}{=}}
\newcommand\leby[1]{\stackrel{\equ{#1}}{\le}}
\newcommand\ltby[1]{\stackrel{\equ{#1}}{<}}
\newcommand\geby[1]{\stackrel{\equ{#1}}{\ge}}
\newcommand\gtby[1]{\stackrel{\equ{#1}}{>}}
\renewcommand{\Im}{{\rm \, Im\,}}
\renewcommand{\Re}{{\rm \, Re\,}}

\newcommand{\negsp}{\hspace{-.04truecm}}
\newcommand\tnorm[1]{\left\vvvert #1 \right\vvvert}
\newcommand{\ex}{{\, e}}
\renewcommand{\a }{ {\alpha}   }
\renewcommand{\b}{ {\beta}   }
\newcommand{\g}{ {\gamma}   }
\newcommand{\G}{ {\Gamma}   }
\renewcommand{\d}{ {\delta}   }
\newcommand{\D}{ {\Delta}   }
\newcommand{\vae }{ {\varepsilon}   }

\renewcommand{\th }{ {\theta}   }
\newcommand{\Th }{ {\Theta}   }
\newcommand{\vth }{ {\vartheta}   }
\renewcommand{\k}{ {\kappa}   }
\renewcommand{\l}{ {\lambda}   }
\renewcommand{\L}{ {\Lambda}   }
\newcommand{\m}{ {\mu}   }
\newcommand{\n}{ {\nu}   }
\newcommand{\x }{ {\xi}   }
\newcommand{\X }{ {\Xi}   }
\newcommand{\p}{ {\pi}   }

\renewcommand{\r}{ {\rho}   }
\newcommand{\s}{ {\sigma}   }

\renewcommand{\t}{ {\tau}   }
\newcommand{\f}{ {\varphi}   }
\newcommand{\ph}{ {\phi}   }

\renewcommand{\o}{ {\omega}   }
\renewcommand{\O}{ {\Omega}   }
\newcommand{\torus}{ {\mathbb{ T}}   }
\renewcommand{\natural}{ {\mathbb{ N}}   }
\newcommand{\real}{ {\mathbb{ R}}   }
\newcommand{\integer}{ {\mathbb{ Z}}   }
\newcommand{\complex}{ {\mathbb { C}}   }

\newcommand{\tn}{ {\torus^d} }

\newcommand{\rn}{ {\real^d}   }

\newcommand{\cn}{ {\complex^d }   }
\newcommand{\zn}{ {\integer^d }   }
\newcommand{\nn}{ {\natural^d }   }

\newcommand\ppu{{ (1) }}
\newcommand\ppd{{ (2) }}
\newcommand\ppt{{ (3) }}

\newcommand{\cB}{ {\cal B} }

\font\teneufm=eufm10
\font\seveneufm=eufm7
\font\fiveeufm=eufm5
\newfam\eufmfam
\textfont\eufmfam=\teneufm
\scriptfont\eufmfam=\seveneufm
\scriptscriptfont\eufmfam=\fiveeufm

\newcommand\appA[1]{\subsection{#1}\label{app:A}
\renewcommand{\theequation}{A.\arabic{equation}}
           \setcounter{equation}{0}
\renewcommand{\thetheorem}{A.\arabic{theorem}}
           \setcounter{theorem}{0}
                  }
\newcommand\appB[1]{\subsection{#1}
\renewcommand{\theequation}{B.\arabic{equation}}
           \setcounter{equation}{0}
\renewcommand{\thetheorem}{B.\arabic{theorem}}
           \setcounter{theorem}{0}           
           }
\newcommand\appE[1]{\subsection{#1}
\renewcommand{\theequation}{C.\arabic{equation}}
           \setcounter{equation}{0}
\renewcommand{\thetheorem}{C.\arabic{theorem}}
           \setcounter{theorem}{0}
           }
\newcommand\appC[1]{\subsection{#1}
\renewcommand{\theequation}{D.\arabic{equation}}
           \setcounter{equation}{0}
\renewcommand{\thetheorem}{D.\arabic{theorem}}
           \setcounter{theorem}{0}           
           }
\newcommand\appD[1]{\subsection{#1}
\renewcommand{\theequation}{E.\arabic{equation}}
           \setcounter{equation}{0}
\renewcommand{\thetheorem}{E.\arabic{theorem}}
           \setcounter{theorem}{0}}
\newcommand\appF[1]{\subsection{#1}
\renewcommand{\theequation}{F.\arabic{equation}}
           \setcounter{equation}{0}
\renewcommand{\thetheorem}{F.\arabic{theorem}}
           \setcounter{theorem}{0}}
\newcommand\appG[1]{\subsection{#1}
\renewcommand{\theequation}{G.\arabic{equation}}
           \setcounter{equation}{0}
\renewcommand{\thetheorem}{G.\arabic{theorem}}
           \setcounter{theorem}{0}}

\def\uno{{\mathbbm 1}}

\def\id{{\rm id }}
\newcommand{\wh}{\widehat}
\newcommand{\wt}{\widetilde}

\numberwithin{equation}{section}
\makeatletter
\def\clap#1{\hbox to 0pt{\hss #1\hss}}%
\def\ligne#1{%
\hbox to \hsize{%
\vbox{\centering #1}}}%
\def\haut#1#2#3{%
\hbox to \hsize{%
\rlap{\vtop{\raggedright #1}}%
\hss
\clap{\vtop{\centering #2}}%
\hss
\llap{\vtop{\raggedleft #3}}}}%
\def\bas#1#2#3{%
\hbox to \hsize{%
\rlap{\vbox{\raggedright #1}}%
\hss
\clap{\vbox{\centering #2}}%
\hss
\llap{\vbox{\raggedleft #3}}}}%
\def\maketitle{%
\thispagestyle{empty}\vbox to \vsize{%
\haut{}{\@blurb}{}
\vfill
\hrule height 4pt

\begin{center}
\usefont{OT1}{ptm}{m}{n}
\huge \@title
\end{center}
\par
\hrule height 4pt
\vspace{1cm}
\par
\begin{center}
\usefont{OT1}{phv}{m}{n}
\Large \@author
\par
\end{center}
\vspace{0.5cm}
\vfill
\vfill
\bas{}{ \@date}{}
}
\cleardoublepage
}
\def\date#1{\def\@date{#1}}
\def\author#1{\def\@author{#1}}
\def\title#1{\def\@title{#1}}
\def\location#1{\def\@location{#1}}
\def\blurb#1{\def\@blurb{#1}}
\date{\small \today}
\author{}
\title{}
\location{}\blurb{}
\makeatother
\title{\textcolor{blue}{Quantitative KAM normal forms and sharp measure estimates}}
\author{{\normalsize\emph{Candidate}}\\ Comlan Edmond Koudjinan\\\vspace{0.7cm}
{\normalsize\emph{Supervisor } }\\Prof. Luigi Chierchia\\
\vspace{0.7cm}
{\normalsize\emph{Coordinator }}\\ Prof. Angelo Felice Lopez}
\blurb{
\begin{tabular}{p{0.5\textwidth} p{0.45\textwidth}} 
\raggedright \includegraphics[height=2cm,width=3cm]{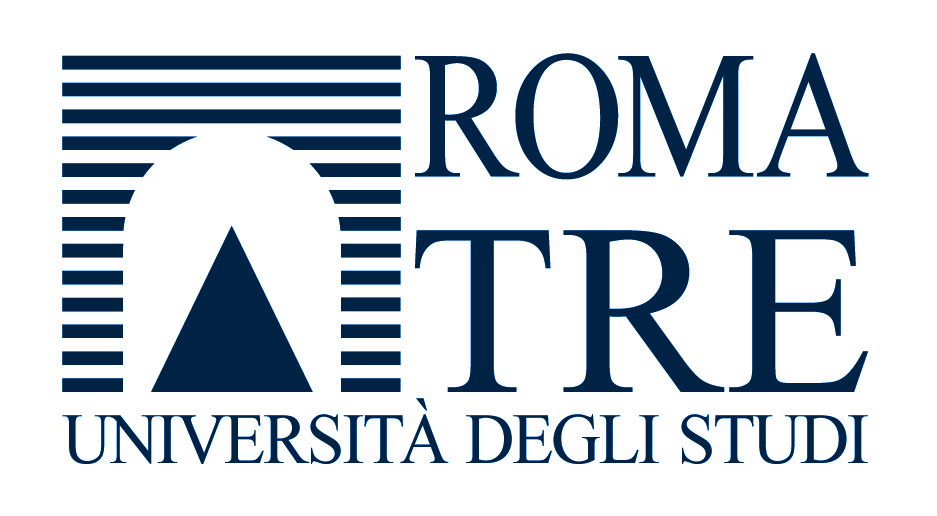}&  \raggedleft\includegraphics[height=2cm,width=3cm]{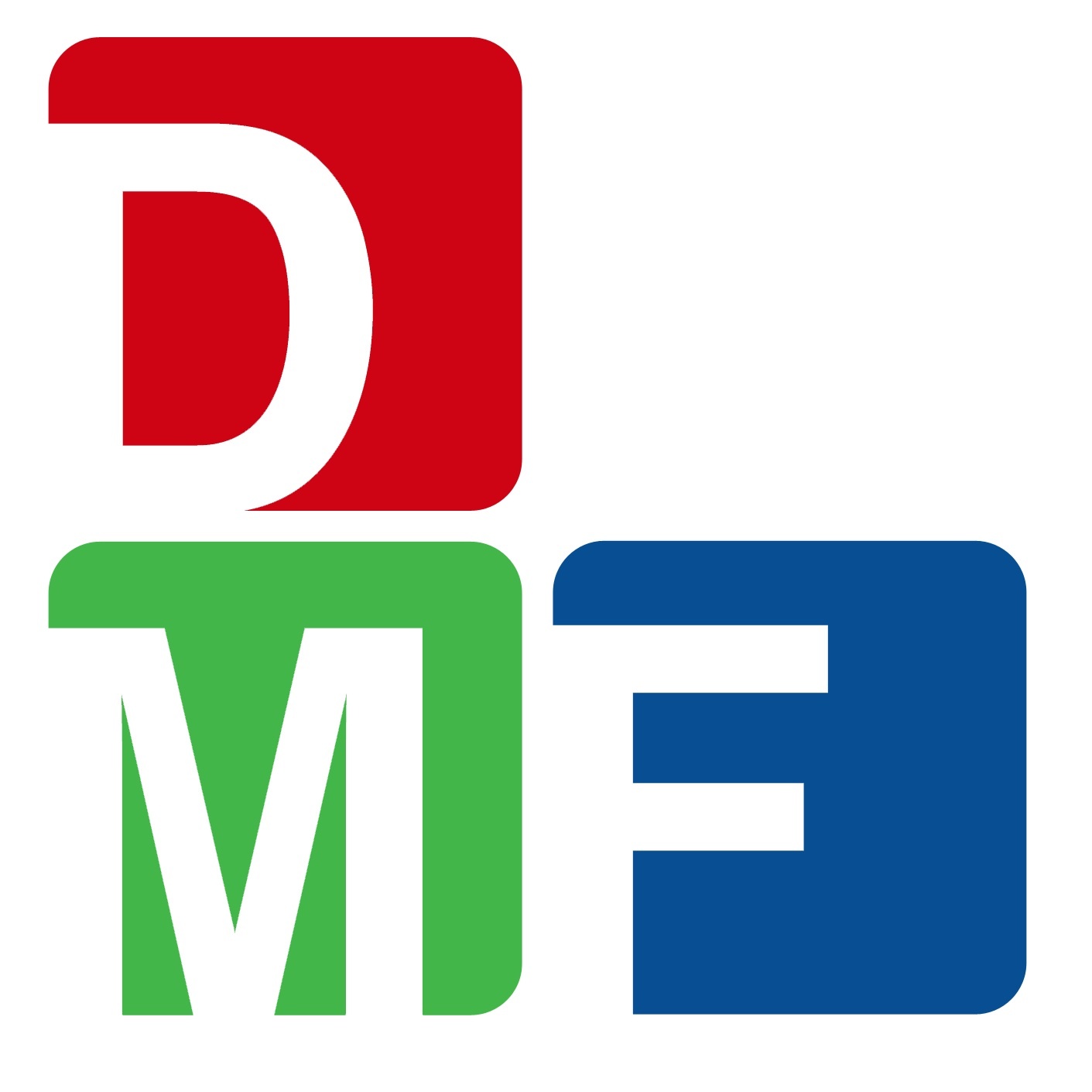}
\end{tabular}
\begin{center}
\bf\large{Universit\`a degli Studi Roma Tre}
\end{center}
\begin{center}
{\bf\large{Dipartimento di Matematica e Fisica}} 
\end{center}
\begin{center}
\bf\large{XXXI Cycle}
\end{center}
\vspace{0.5cm}
\begin{center}
\bf\large{Doctoral Thesis in Mathematics}
\end{center}
}%
\begin{document}
	\dominitoc
	\maketitle
	%
	%
	%
\frontmatter



\chapter{Abstract}
It is widespread since the beginning of KAM Theory that, under ``sufficiently small'' perturbation, of size $\epsilon$, apart a set of measure   $O(\sqrt{\epsilon})$, all the KAM Tori of a non--degenerate integrable Hamiltonian system persist up to a small deformation. However, no explicit, self--contained proof of this fact exists so far. In the present Thesis, we give a detailed proof of how to get rid of a
logarithmic correction (due to a Fourier cut--off) in Arnold’s scheme
and then use it to prove an explicit and ``sharp'' Theorem of  integrability on Cantor--type set. In particular, we give an explicit proof of the above--mentioned measure estimate on the measure of persistent primary KAM tori. We also prove three quantitative KAM normal forms following closely the original ideas of the pioneers Kolmogorov, Arnold and Moser, computing explicitly all the KAM constants involved and fix some ``physical dimension'' issues by means of appropriate rescalings. Finally, we compare those three quantitative KAM normal forms on a simple mechanical system.

\cleardoublepage



%

\tableofcontents


\mainmatter
%
%

\chapter*{Acknowledgements}
I would like to express my deep gratitude to Prof. Luigi Chierchia for introducing to the absolutely fascinating field of {\it KAM theory}, for his guidance, availability, patience,  support which goes much much beyond mathematics and this thesis, and his precious help all along the trip leading to the accomplishment of this thesis.\\

\noi
I would like aslo to express my sincere gratitude to Prof. Francesco Pappalardi for his inestimable support. I am in Roma thanks to him mainly and without him I would not have the chance to meet Prof. Luigi and this thesis would not be possible simply.\\

\noi
I am very grateful to Prof. Luca Biasco and Prof. Michela Procesi for their availability, kindness, the very helpful discussions i had with them.\\

\noi
I am also deeply grateful to Prof. Antonio Siconolfi, with whom I spent almost the half of the three years of my PhD, for his help, generosity, patience and precious advises. I could have completed my PhD with him, but at some some point, due to time constrained, I had to choose and focus on a single subject.\\

\noi
I would like to address a special thanks to Prof. Wilfrid Gangbo whose impact on my career is considerable and this, since 2015 that I met him. It was through him I met Prof. Siconolfi. \\ 

\noi
I would like to address my sincere gratitude to Prof. Théophile Olory for his precious advises and assistance.\\

\noi
I am very grateful to Prof. Bernadin Kpamingan for his help all those years along.\\

\noi
I would like to express my gratitude to Prof. Rafael De La Llave  for his interest in my work, for offering the opportunity to attend a workshop and for drawing my attention on some related papers in the literature.\\

\noi
I would like to express my gratitude to Prof. Jean--Pierre Marco and Prof. Tere M--Seara for their interest in my thesis and their availability to help.\\

\noi
I am also very grateful to Prof. Stefano Luzzatto and Prof. Carlangelo Liverani for their support and precious time.\\

\noi
A big big thanks to my mom, my sisters and brothers, especially Zizi, Igor, Marcelin, Ad\`ele, Jeannette, Justine, Eulalie, ..., for their permanent support.\\

\noi
Finally, I want to thank my queen Cica and my princess Merquela, for their patience and constant encouragement, I am infinitely grateful to you.\\

\noi
This list is of course far from being exhaustive and detailed; I tried to be as short as possible. Otherwise, I would somehow fall into writing my biography and here is certainly not the place for!!!!

\vspace*{2.5cm}
\begin{flushright}
To the memory of my father
\end{flushright}
\newpage
\section*{Notational conventions}
\begin{itemize}
\item $\ex$     denotes the Neper's number \ie $\exp(1)$
\item $\natural=\{1,2,3,\cdots\}$ and $\natural_0=\{1,2,3,\cdots\}$
\item $\real$ and $\complex$ are respectively the set of real and complex numbers
\item $\n!=\n_1!\cdots \n_d!$ and $|\n|_1=\n_1+\cdots +\n_d$, for any $\n=(\n_1,\cdots,\n_d)\in \natural_0^d$
\item $y^\b=y_1^{\b_1}\cdots y_d^{\b_d}$ for any $y,\b\in \rn$
\item $|(y_1,\cdots,y_d)|= \max\{|y_1|,\cdots,|y_d|\}\quad$ and $\quad|(y_1,\cdots,y_d)|_2= \sqrt{y_1^2+\cdots+y_d^2}$
\item $(x_1,\cdots,x_d)\cdot(y_1,\cdots,y_d)=\average{(x_1,\cdots,x_d)\;,(y_1,\cdots,y_d)}=x_1y_1+\cdots+x_dy_d$
\item $\dist$ denotes the distance function
\item $\ovl{A}$  denotes the closure of $A$
\item $\dpr{A}$  denotes the ``boundary'' of $A$
\item $\conv(A)$ denotes the convex--hall of $A$
\item $C^n(A,B)$ (\resp $C^n_c(A,B)$)  denote respectively the set of functions of class $C^n$ (\resp $C^n$ with compact supports) from $A$ into $B$
\item $\meas_d$ denotes the $d$--dimensional Lebesgue--measure
\item $\dom(f)$ denotes the domain of $f$
\item $\supp(f)$ denotes the support of $f$
\item $ B_{r}(p)$ (\resp $D_{r}(p)$) denotes the ball centered at p with radius $r$ in $\rn$ (\resp in $\cn$)
\item $B_{r}(A)$ (\resp $D_{r}(A)$) denotes the r--neighborhood of $A$ in $\rn$ (\resp in $\cn$) 
\item $\D^\t_\a$   denotes the set of $(\a,\t)$--Diophantine vectors
\item $\torus^d_{s}$  denotes the strip of width $s$ around $\tn$ in $\cn$
\item $\average{f}$ denotes the average of $f$ on $\tn$
\item $f_\n=\dpr^\n_y f=\frac{\dpr^\n f}{\dpr y^\n}=\frac{\dpr^{|\n|_1} f}{\dpr y_1^{\n_1}\cdots\dpr y_d^{\n_d}}$ denotes derivative of order $\n$ of $f$
\item $\Iso(V)$  denotes the set of isomorphisms from $V$ onto itself
\item $\mathcal{M}_{n,m}(\cn)$ the set of $n$--by--$m$ matrices with entries in $\cn$ and $\mathcal{S}_n(\cn)\subset \mathcal{M}_n(\cn)\coloneqq \mathcal{M}_{n,n}(\cn)$ the set of symmetric square matrices of order $n$
\item $\adj(A)$ denotes the adjoint of $A$
\item $\det(A)$ denotes the determinant of $A$
\item $A^T$ denotes the transposed of $A$
\item $NM$ and $TM$ are respectively the normal and tangent bundle of the manifold $M$\\
\item $\G(M)$ denotes the space of smooth vector field on $M$\\
\item $\mathfrak{F}(M)$ denotes the space of smooth functions on $M$\\
\item $\minfoc(M)$ denotes the minimal focal distance of the manifold $M$
\end{itemize}
\chapter{Introduction}
In the solar system framework, { Celestical Mechanics}, a branch of { astronomy}, consists ultimately in the study of the $n$--body problem. The $n$--body problem is the dynamical system that governs the motion of $n$ planets interacting according to Newton's gravitation law. A holy--grail question in { Celestical Mechanics} was and remains the { stability} of the solar system, \ie whether the current configuration of the planets will stay unchanged forever under their interaction, or whether some planets will be kicked out of the system or have their trajectories be drastically affected to eventually collapse and give rise to unpredictable behaviors. Across the history of Mathematics, most of the great figures devoted some part of their works to this question. { Laplace }(1773), { Lagrange }(1776), { Poisson }(1809), and { Dirichlet} (1858) used series expansion techniques to study the question of { stability} of the solar system and claimed all to have proved it. Then { Bruns} (1887) proved that, from quantitative point of view, the only method which could solve the $n$--body problem is the series expansions. But,  the works of { Haretu} (1878) and { Poincar\'e}  (1892) (see \cite{poincare189299}) show that all those series expansion techniques fail as the series expansions they use diverge (see \cite{dugas1957history,moulton2012introduction,moser1973stable,dumas2014kam,abraham1978foundations} for more historical details).\\ \ \\
 A new viewpoint is thus undeniably needed to overcome this embarrassing fact. The change of paradigm was made by { Poincar\'e}. Indeed, { Poincar\'e} introduced a completely revolutionary qualitative approach to Mechanics 
 (see \cite{poincare1890probleme,poincare189299}). The point is that, for question such as stability, one needs to study the entire phase portrait, and in particular the asymptotic time behavior of the solutions.\\
  The phase portrait is the family of solutions curves,  which fill up the entire phase space. The phase space is a symplectic manifold (a differentiable manifold together with a symplectic structure). Dynamical system is then just given by a Hamiltonian vector field; this is the Mathematical model for the global study in Mechanics that { Poincar\'e} gave in his qualitative theory. With his new geometrical methods, { Poincar\'e} discovered the non-integrability of the three body problem. In fact, the small divisor problem was well--known to { Poincar\'e}, who was aware that, because of this problem,  nearly--integrable Hamiltonian systems are, in general, not integrable (analytically). { Poincar\'e} and his successors then speculate that most of the classical systems were chaotic, and ergodic. There were even a { gaped proof} of the ergocity of a generic Hamiltonian system by { E. Fermi} in the 1920' (see \cite{fermi1923beweis,fermi1923dimostrazione,fermi1923generalizzazione,fermi1924uber}). This { ergodic hypothesis} was accepted by many, including some of the brightest mind of those times, till the discoveries by Kolmogorov and his followers.\\ 
 At the 1954 International Congress of Mathematician in Amsterdam, against any expectation, { Kolmogorov} (see \cite{Ko1954, kolmogorov1957theorie}) presented a four--pages note where he sketched the proof of the persistence of the majority of { tori} for a nearly--integrable Hamiltonian system. Then, his former student { Arnol'd} (see \cite{ARV63,arnol1963small}) completed the proof in the analytic category, and { Moser} (see \cite{moser1962invariant,moser1966rapidly,moser1966rapidly2}) in the smooth category, whence the acronym { KAM Theory}.\\
The object of { KAM Theory} is the construction of quasiperiodic trajectories, which are sets of perpetual stability,  in Hamiltonian dynamics. A { KAM scheme} is essentially based on the { Newcomb} idea  of successive constructions of change of variable through a { Newton}--like method. Those successive changes of variables are carried out to eliminate, in a { super--exponentianlly} increasing manner, the fast phase variables. A { KAM scheme} of course encounters the small divisor problem that { Poincar\'e} faced. Netherless, the { super--exponentianlly} decay make the whole scheme converge.\\
 Formally, one is given a symplectic manifold $(M,\varpi)$, a (smooth) Hamiltonian $H\colon T^*M \to \real$. To the Hamiltonian $H$, is associated a (unique!) smooth vector field, the Hamiltonian vector field, say $X_H$, given by the equation
 $$
 \varpi(X_H,\cdot)=-dH\;.
 $$
 The smooth vector field $X_H$ then generates a flow, say $\phi^t_H$, by the relation
 \beq{HamFlowInt}
 \frac{d}{dt}\phi^t_H=X_H\circ \phi^t_H\;.
 \eeq
 In particular, if $M=\rn\times\tn$ and $\varpi=dy_1\wedge dx_1+\cdots+dy_d\wedge dx_d$, then $X_H=(-\dpr_x H,\dpr_y H)$ and, therefore,  the equation \equ{HamFlowInt} reads
 \beq{HamFlowIntBis}
 \left\{
 \begin{aligned}
 &\dot{y}=-\dpr_x H\\
 &\dot{x}=\dpr_y H
 \end{aligned}
 \right.
 \,,\qquad
 \phi_H^t=(y(t),\, x(t))\,.
 \eeq
 Then, to construct quasi--periodic trajectories for the Hamiltonian system \equ{HamFlowInt}, one looks for a change of variable
$\phi'\colon (y',x')\mapsto (y,x)=\phi'(y',x')$, with the following properties \\
\begin{itemize}
\item[$(a)$] $\phi$ preserves the Hamiltonian structure of \equ{HamFlowInt}. More precisely, $\phi$ preserves the symplectic form \ie $\phi^*\varpi=\varpi$\,;
\item[$(b)$] $\phi$ conjugates $\phi_H^t$ to a linear flow:
\beq{LinFlowIntR}
\phi^{-1}\circ \phi_H^t\circ\phi(y,x)=(y,\o t+x)\;,\qquad \mbox{with}\qquad \o\coloneqq \dpr_y H(y,x)\;.
\eeq
\end{itemize}
However, one does not solve directly \equ{LinFlowIntR}. Instead, one conjugates the Hamiltonian itself \ie construct $\phi$ in such away that
\beq{ConjHamiInTr}
H\circ \phi(y',x')=H_*(y')\;,
\eeq
as the latter is much easier to carry out than the former. Once \equ{ConjHamiInTr} holds, by the property $(a)$, \equ{LinFlowIntR} follows (see for instance \cite{moser2005notes} for details).\\
The numerical property of the frequency or winding number $\o$ plays a crucial role in the construction of the { invariant tori}. The most common assumption is the Diophantine property. A vector $\o\in \rn$ is said $(\a,\t)$--Diophantine if
\beq{DioIntro}
|\o\cdot k|\ge \frac{\a}{|k|_1^\t}\;,\qquad\forall\; k\in \zn\setminus\{0\},
\eeq
where $|k|_1=|k_1|+\cdots+|k_d|$.\\
{\bf Facts} Let
$$
\D^\t_\a \coloneqq \left\{\o\in \rn\;:|\o\cdot k|\ge \frac{\a}{|k|_1^\t}\;, \ \forall k\in \zn\setminus\{0\}\right\}\;,
$$
and
$$
\D^\t\coloneqq \dst\bigcup_{\a>0} \D^\t_\a\;,
$$
be the set of all $(\a,\t)$--Diophantine vectors. Thus 
\begin{itemize}
\item[$(i)$] If $\t<d-1$, then $\D^\t=\emptyset$ (see \cite{cassells1957introduction});
\item[$(ii)$] If $\t=d-1$, then the set $\D^\t$ has zero Lebesgue--measure, but is of Hausdorff dimension $d$. In particular, the intersection of $\D^\t$ with any open set has the cardinality of $\real$ (see \cite{schmidt1966badly,schmidt1969badly});
\item[$(iii)$] If  $\t>d-1$, then the Lebesgue--measure of $\rn\setminus\D^\t$  is zero. In fact,
 \beq{DiopOutBall}
\meas\left(B_R(0)\setminus \D^\t_\a\right)\le \left(\dst\sum_{k\in\zn\setminus\{0\}}\frac{2^d\sqrt{d}}{|k|_1^{\t+1}}\right)R^{d-1}\;\a\;,\qquad \forall\;R,\a>0. 
 \eeq
\end{itemize}
\ \\
 As soon as the existence of invariant { tori} is established, one can speak of { Kolmogorov sets}, which turn out to be very big.  A { Kolmogorov set} associated to a Hamiltonian $H$ is an union of  its { invariant maximal KAM tori}. A { maximal KAM torus} for $H$ is an embedded, { Lagrangian}, { Kronecker torus} with Diophantine frequency $\o$. A { Kronecker torus} with frequency $\o$ is an embedded torus on which the $H$--flow is conjugated to the linear flow
 $$
\real\times\tn\ni (t,x)\mapsto x+\o t\;.
 $$
\\
{\it In this thesis, we are mainly concerned with ``sharp'' measure estimates of { Kolmogorov sets}, with emphasize on the dependence of those measure estimates upon the geometry of the domain}.\\

\noi
Moser introduced in \cite{moser1967convergent} the original idea of parametrizing a non--degenerate quasi--integrable Hamiltonian by the frequency vectors and then apply the KAM technics 
 (see also \cite{JP, poschel1982integrability}). In \cite{JP}, the author made a short discussion of the measure of the complement of { Kolmogorov set}. Then, very  recently, { Biasco} and { Chierchia} \cite{biasco2018explicit}, give a detailed proof of the measure estimate result in \cite{JP} and show how this measure estimate depends upon the domain. 
{ We revisit the paper} \cite{JP} { in this thesis and our computation in particular fixes a small gap in the statement of}  \cite{JP} (see $\S 4.$ of Remark~\ref{rem001} below).

\ \\
\noi
{ Arnold's scheme} \cite{ARV63,arnold18small} can be summarized as follows. Let $K$ and $P$ be real--analytic in $D_0\coloneqq D_{r_0}(y_0)\times \torus^d_{s_0}$, with $K$ integrable and such that
$$
K_y(y_0)\eqqcolon \o\in \D^\t_\a\qquad\mbox{and}\qquad \det K_{yy}(y_0)\neq 0,
$$
with $\a>0$ and $\t>d-1$. Thus, the torus $\mathcal{T}_{\o,0}\coloneqq \{y_0\}\times \tn$ is a { KAM torus} for $K$ on which its flow $\phi^t_K$ is linear: 
$$\phi^t_K(y_0,x)=(y_0,\o t+x).$$
Then, the idea of { Arnold} is to construct a near--to--identity symplectic change of variables 
$$
\phi_1\colon D_1\coloneqq D_{r_1}(y_1)\times \torus^d_{s_1}\to D_0\;,
$$ 
with $D_1\subset D_0$ such that
\beq{ArnH1v2Intro}
\left\{
\begin{aligned}
& H_1:= H\circ \phi_1=K_1+\vae^2 P_1\ ,\\
& \dpr_{y} K_1(y_1)=\o,\quad \det \dpr^2_{y} K_1(y_1)\neq 0\,.
\end{aligned}
\right.
\eeq
And for $\vae$ small enough, one can iterate the process and build  a sequence of symplectic transformations ($j\ge 1$)
$$
\phi_j\colon D_j\coloneqq D_{r_j}(y_j)\times \torus^d_{s_j}\to D_{j-1}\;,\qquad D_j\subset D_{j-1}\;,
$$ 
and satisfying
\beq{ArnH1v2Intro}
\left\{
\begin{aligned}
& H_j:= H\circ \phi_j=K_1+\vae^{2^j} P_j\ ,\\
& \dpr_{y} K_j(y_j)=\o,\quad \det \dpr^2_{y} K_j(y_j)\neq 0\,.
\end{aligned}
\right.
\eeq
In performing this construction, one is first attempt  to solve the linear PDE 
\beq{GenFunpHiJ}
\dpr_{y} K_j(y)\cdot \dpr_x g_j+P_j(y,x)= \mbox{function of } y \mbox{ exclusively},
\eeq
where $g_j$ is a generating function for $\phi_j$. But \equ{GenFunpHiJ} does not admit solution, because of small divisor problem (see \cite{CL12} for more discussion). The key idea of { Arnold} is then to solve only a truncated version of \equ{GenFunpHiJ}, with the order of truncation large enough so that the error one commits by solving  approximately \equ{GenFunpHiJ} is of order the square of the size of the perturbation. The truncation is the origin of the logarithmic correction in the smallness condition required in order to iterate infinitely many times the { Arnold process}.
In particular, the Lebesgue--measure estimate of the complementary of the { Kolmogorov set} one gets from { Arnold's scheme} is $O(\sqrt{\vae}\;(\log\vae^{-1})^{3(\t+1)/4})$. 
This estimate is not the optimal one, which is $O(\sqrt{\vae})$ (see for instance \cite{biasco2018explicit}, where the constant in front of $\sqrt{\vae}$ in the optimal measure estimate is computed explicitly and the proof uses the KAM Theorem \`a la Moser). The task of getting rid of the logarithmic correction in the { Arnold's scheme} is not obvious.
 The first paper in this direction is the sketchy 7--pages paper  \cite{neishtadt1981estimates}, where { Neishtadt} outlines how to overcome the logarithmic correction. The approach we adopt here is essentially equivalent to the one in \cite{neishtadt1981estimates} though conceptually different. Indeed, in our scheme we fix the frequencies of the { tori} we build up from the beginning once for all. Instead, in \cite{neishtadt1981estimates} as well as in the original paper by { Arnold} \cite{arnold18small}, the { tori} as well as their respective frequencies are constructed iteratively.\\ \ 

\noi
{\it Moreover, in our approach, we focus on the smallest possible $\a$ \ie the situations where the square--root of the sizes of the perturbations are proportional to the { Diophantine} constant $\a$ of the frequency of the { tori}. We then discuss the measure estimate of the { Kolmogorov set} we build up}. The sharpness of the measure of the { Kolmogorov set} is in fact intimately related to the power of the { Diophantine} constant $\a$ in the smallness condition under which one performs the { KAM scheme}. Recently, { Villanueva} \cite{villanueva2008kolmogorov} revisited the classical { Kolmogorov scheme} and succeed to cut down the power of $\a$ in the smallness condition, from $4$ to the optimal which is $2$; but with no measure estimate of { Kolmogorov set} discussion. See also \cite{villanueva2018parameterization} where he got the exponent $2$ for $\a$ in the smallness condition, in the framework of exact symplectic maps in { Euclidean spaces} of even dimensions. 

\ \\
\section{Main results}
As a basic rule in this thesis, we compute explicitly all the { KAM constants}. Investigating the explicit dependence of the ``KAM constants'' upon the parameters in a quasi--integrable Hamiltonian system is of great interest, not only in view of its applications (for instance to the $n$--body problem \cite{celletti2006kam}, to geodesic flows on surfaces, \etc) but also for the discussion of explicit measure estimates of { Kolmogorov sets}. 
The content of this thesis can be described very roughly as follows:
\begin{itemize}
\item[$(i)$]We prove three { quantitative  KAM normal forms} following closely the original ideas of the pioneers { Kolmogorov, Arnold} and { Moser}. We compute in particular explicitly all the { KAM} constants in them and fix physical dimension issues by rescaling conveniently various quantities. 
Then, we compare those three { quantitative  KAM normal forms} on a simple mechanical system.
\item[$(ii)$]  We give detailed proof of how to get rid of the logarithmic correction in the { Arnold's scheme} and then use it to prove an explicit and ``sharp'' Theorem of integrability on {Cantor--type} set.
\item[$(iii)$] {\it We prove three types of { sharp} measure estimate of { Kolmogorov sets}. In the first one, we adopt the global approach which consists in constructing the  { Kolmogorov set} in a given bounded domain and then estimate its measure. In the two others, we slice the domain into relatively small cubes with equilength sides. In each of those cubes, we construct a { Kolmogorov set} associated to the restriction of the Hamiltonian to such a cube and estimate its measure. Then, we sum up the local  { Kolmogorov sets} constructed. \\
One of the local approaches follows the idea in} \cite{biasco2018explicit} {\it and recover its result.\\
In the second local approach, we introduce a { geometric integer constant} of a set which is the minimal number of cubes one needs to cover the set  by cubes with the same side--length, centered on the set and with total ``volume'' not exceeding some fixed amount.  This third approach  is somehow more intrinsic.}    
\end{itemize}
			
\noi

\noi
{$\mathbf{(i)}$} More precisely, we prove in {\bf Theorem~\ref{teo2}} (following Kolmogorov's proof in \cite{Ko1954}, scheme to which a complete proof was given in \cite{CL08,CL12}) that for any small enough perturbation of a non--degenerate Kolmogorov normal form, there exists a symplectic change of variables such that in the new variables, the Hamiltonian reduced to a Kolmogorov normal form.\\

\noi
We prove in {\bf Theorem~\ref{teo4}} (following Arnold \cite{ARV63} and basing on \cite{CL08,CL12}) that, under a sufficiently small perturbation, with size say $\vae$, of a non--degenerate integrable Hamiltonian system, the majority ( $1-O(\sqrt{\vae}\; (\log\vae^{-1})^{3\n/4})$ of the total Lebesgue--measure) of the invariant, Lagrangian, Kronecker tori with Diophantine frequencies 
of the integrable system persist, being only slightly deformed. \\

\noi
In {\bf Theorem~\ref{teo1}}, we prove (following Moser\cite{moser1967convergent} and basing on \cite{JP}) that on a bounded domain, the totality of the invariant maximal KAM tori of the linear normal form, whose frequencies are far enough from the boudary persist under any small enough perturbation. These tori are just a little bit deformed and persist as invariant maximal KAM tori, not of the perturbed Hamiltonian itself, but of the perturbed Hamiltonian plus a small shift of the frequency.\\ 

\noi
In {\bf Chapter~\ref{CQTHMS}}, we compare the explicit KAM mormal forms on a simple mechanical Hamiltonian and compute the numerical values of the thresholds within these Theorems in a concrete case.\\

\noi
In {\bf Theorem~\ref{Extteo4}}, we prove an explicit Theorem of integrability on a Cantor--like set. Namely, for any given sufficiently small perturbation of a non--degenerate integrable Hamiltonian on a bounded domain, we construct 
 a $\ci$--symplectomorphism which conjugates the perturbed Hamiltonian to an integrable Hamiltonian on a Cantor--like set. The Cantor--like set is equipotent to the set of phase points which are at some minimal distance from the boundary and such that their image by the Jacobian of the unperturbed part are Diophantine vectors, with fixed Diophantine parameters. Moreover, the ratio of their respective Lebesgue--measures minus 1 is small with the size of the perturbation.\\

\noi
{$\mathbf{(ii)}$} In {\bf Theorem~\ref{teo4v2}}, we prove a refinement of the Arnold's Theorem by overcoming the {logarithmic correction} looming  from the original scheme. Indeed, we prove that, for any  small enough perturbation of a non--degenerate integrable Hamiltonian system, most of ($1-O(\sqrt{\vae})$ of the total Lebesgue--measure, where $\vae$ is the size of the perturbation) of the  { invariant  maximal KAM tori} 
of the integrable system persist, up to a small deformation. To do so, we isolate the smallness parameter $\vae$ from the super--exponential parameter so that, and this is the whole point, as soon as $\vae$ is chosen conveniently to perform the scheme one time, one can iterate infinitely many times without any other requirement and, in particular, $\vae$ ``disappears'' once for all from the second step on.\\

\noi
{\it In} {\bf Theorem~\ref{Extteo4v2}}, {\it we prove an explicit, intrinsic and sharp integrability Theorem on a Cantor--like set}. Namely, given any small enough, real--analytic perturbation of a non--degenerate integrable Hamiltonian on a bounded domain, we build-up a transformation, $\ci$ in the Whitney sense and symplectic. Actually, the two Cantor--like sets are lipeomorphic\footnote{\ie there exists a bijective Lipschitz continuous function from one onto the other.}. In those new variables,  the nearly--integrable Hamiltonian becomes integrable on a Cantor--like set. The Cantor-like set is  equipotent to the set of phase points which are far enough from the boundary and which images through the  Jacobian of the unperturbed part are $(\a,\t)$--Diophantine, for some $\a>0$ and $\t$ a number larger than the half--dimension minus one. In particular, we get a family of { invariant maximal KAM torus} which complement has a Lebesgue--measure bounded from above by a constant proportional to $\a$.\\

\noi
{$\mathbf{(iii)}$} {\it The novel part of the present thesis consists mainly in {\bf Part~\ref{part2}} and, in particular and more interestingly, the various ``sharp'' and geometric measure estimates of the unstable sets within a Hamiltonian system we provide}.\\

\subsection{Measure estimate I}
Then, we derive in  {\bf Theorem~\ref{Extteo5v2}} the following. Let  $\mathscr{D}\subset\rn$ be a non--empty bounded domain with smooth boundary $\dpr\mathscr D$ and a small enough $\a>0$, depending on the geometry of the hypersurface $\dpr\mathscr D$. Let $H$ be a sufficiently small perturbation, of order $O(\a^2)$, of a non--degenerate integrable Hamiltonian $K$. Then, the set $\mathscr I$ left out of the $H$--{ invariant maximal KAM tori} is bounded in measure by
\begin{align}
\meas(\mathscr I)&\le (3\pi)^d \frac{\mathsf{T}_0}{32d\s_0}\bigg(2\;\mathcal{H}^{d-1}(\dpr \mathscr D)\;\a+ C(d,\s_0,\mathsf{T}_0,\mathbf{R}^{\dpr\mathscr D})\;\a^2+\meas(\mathscr D_\d\setminus \mathscr D_{\d,\a})\bigg)\label{EqIntMes10}\\
   &=O(\a)\,,\label{EqIntMes11}
\end{align}
where\footnote{See Appendix~\ref{appF} for the definitions.
} $\mathcal{H}^{d-1}$ denotes the $(d-1)$--dimensional Hausdorff measure (or equivalently, the $(d-1)$--surface area), $\mathbf{R}^{\dpr\mathscr D}$ denotes the curvature tensor of $\dpr\mathscr D$, $\s_0$ is the loss of analyticity, $\mathsf{T}_0$ is the norm of the inverse of the Hessian $K_{yy}$ of the unperturbed part,
\begin{align*}
&\d\coloneqq\frac{\a\mathsf{T}_0}{32d\s_0}\;,\\
&\mathscr D_\d\coloneqq \{y\in \mathscr D\;:\; \dist(y,\dpr\mathscr D)\ge \d\}\;,\\
&\mathscr D_{\d,\a}\coloneqq \{y\in \mathscr D_\d\ :\ K_{y}(y)\in \D_\a^\t\}\;,\\
&C(d,\s_0,\mathsf{T}_0,\mathbf{R}^{\dpr\mathscr D})\coloneqq \frac{\mathsf{T}_0}{16d\s_0}\dst\sum_{j=1}^{\left[\frac{d-1}{2}\right]}\frac{\mathbf{k}_{2j}(\mathbf{R}^{\dpr\mathscr D})\;}{1\cdot3\cdots (2j+1)}\; \d^{2j-1}\;,
\end{align*}
with $\mathbf{k}_{2j}(\mathbf{R}^{\dpr\mathscr D})$, the $(2j)$--th integrated mean curvature of $\dpr\mathscr D$ in $\rn$.
\rem{DiscRemItr}
{\bf (i)} The first two terms of the {\it r.h.s.} of \equ{EqIntMes10} arise from the estimation of the $\d$--strip around $\dpr\mathscr D$, $\mathscr D\setminus {\mathscr D}_{\d}$, out of which we construct the family of invariant KAM tori. Notice that the last term $\meas(\mathscr D_\d\setminus \mathscr D_{\d,\a})$ is of order $O(\a)$ by \equ{DiopOutBall}, whence \equ{EqIntMes11} holds.\\ 
{\bf (ii)} The estimate \equ{EqIntMes10} might be seen as a ``sharp'' version of the measure estimate of the invariant set in the Two--scale KAM Theorem of \cite{chierchia2010properly}.
\erem
\ \\ 

\noi
The following is proven in {\bf Theorem~\ref{LebLocGen}}. Let $H=K(y)+\vae\; P(y,x)$ be a perturbation of a non--degenerate integrable Hamiltonian $K$, where $\mathfrak{D}$ a non--empty, bounded subset  of $\rn$ and $K, P$ two real analytic function on $\mathfrak{D}\times\tn$ with bounded extension to $D_{\mathsf{r}_0}(\mathfrak{D})\times\torus^d_{s_0}$, for some $\mathsf{r}_0>0$ and $0<s\le 1$. 
We prove that, for a sufficiently small $\vae$ (with explicit upper--bound),  one can construct by ``localization'' argument a family of $H$--{invariant maximal KAM torus}, say $\mathscr K$, which complement has a Lebesgue--measure of order $O(\a)$ and estimated in two ways as follows.
\subsection{Measure estimate II\label{measEst22}}
More specifically, we show one hand that the Kolmogorov set $\mathscr K$ is bounded in measure from above (in the spirit of \cite{biasco2018explicit})  by 
$$
\meas\left((\mathfrak{D}\times\tn)\setminus \mathscr K\right)\le
C\;\mathsf{p}_1\left(\diam\mathfrak{D}+\ell\right)^{d}\a\;,
$$
with 
$C$ a positive universal constant depending only upon the dimension $d$ and $\t$,
\begin{align*}
&\mathsf{p}_1 \coloneqq \frac{\vth\;\eta^2\;\mathsf{T}}{\s_0\;\mathsf{r}_0}\;,\\
&\ell\coloneqq \frac{\mathsf{r}_0}{2^6d\eta^2}\;,\\
&\eta		 \coloneqq \mathsf{T}\mathsf{K}\ge 1\;,\\
&\vth\coloneqq \frac{\mathsf{K}^d}{\varrho}\ge 1\;,\\
&\varrho\coloneqq \dst\inf_{y\in \mathfrak{D}}|\det K_{yy}(y)|> 0\;,\\
&0< \s_0<2^{5-2\t}\;d\;s_0\;,\\
&\mathsf{T} \coloneqq \sup_{y\in\mathfrak{D}}\|K_{yy}(y)^{-1}\|\;, \\
&\mathsf{K} \coloneqq \|K_{yy}\|_{\mathsf{r}_0,\mathfrak{D}}\;.
\end{align*}
This provides an alternative proof to the result in \cite{biasco2018explicit} (alternative in the sense that the proof in \cite{biasco2018explicit} is based upon Moser's idea while here, we use Arnold's scheme) and our proof is somehow more complete as we compute explicitly all the constants while \cite{biasco2018explicit} refers to \cite{JP}, where the constants are left implicit.
%
\subsection{Measure estimate III}
On the other hand, in a more intrinsic way, we build up (under the same basic assumptions as above) a family $\mathscr K$ of $H$--{ invariant maximal KAM torus} such that
$$
\meas\left((\mathfrak{D}\times\tn)\setminus \mathscr K\right)\le C'\; \frac{\vth\;\mathsf{T}}{\s_0}\;n_{\mathfrak{D}}^{\su d}\;\meas(\mathfrak{D})^{\frac{d-1}{d}}\;{\a}\;,
$$
with $\s_0$, $\vth$, $\mathsf{T}$ as in $\S\ref{measEst22}$, $C$ a positive universal constant depending only upon the dimension $d$ and $\t$, 
  and $n_{\mathfrak{D}}\in\natural$ a ``covering number'' of $\mathfrak{D}$ defined morally as follows.\\
Given $R>0$, define the set $\mathscr C_R$ of coverings of $\mathfrak{D}$ by cubes as follows: $F\in \mathscr C_R$ if and only if there exists $n_F\in \natural$ and $n_F$ cubes, say $F_i,$ $1\le i\le n_F$, of equal side--length $2R$, centered at a point  $y_i\in\mathfrak{D}$ and such that 
$$
F=\big\{F_i:1\le i\le n_F\big\} \qquad \mbox{and}\qquad\mathfrak{D}\subset \dst\bigcup_{i=1}^{n_F}F_i\;.
$$
Then define
$$
\mathscr R\coloneqq \bigg\{
R>0: \mathscr C_{R}\neq \emptyset\ \mbox{and}\ \inf_{F\in \mathscr C_{R}}n_F (2R)^d\le 2^d\meas(\mathfrak{D})\bigg\}
$$
and
$$
n_{\mathfrak{D}}\coloneqq \dst\min_{R\in\mathscr R}\min \bigg\{n_F: F\in \mathscr C_{R}\quad\mbox{and}\quad n_F R^d\le \meas(\mathfrak{D})\bigg\}\;.
$$ 
\rem{IntrNd}
In the above definition of ``covering number'', one could replace the coefficient $2^d$ in front of $\meas(\mathfrak{D})$ in the definition of $\mathscr R$ by $\k>1$, leading to
$$
\mathscr R_\k\coloneqq \bigg\{
R>0: \mathscr C_{R}\neq \emptyset\ \mbox{and}\ \inf_{F\in \mathscr C_{R}}n_F R^d\le 2^{-d}\k\cdot \meas(\mathfrak{D})\bigg\}
$$
and
$$
n_{\mathfrak{D},\k}\coloneqq \dst\min_{R\in\mathscr R}\min \bigg\{n_F: F\in \mathscr C_{R}\quad\mbox{and}\quad n_F R^d\le 2^{-d}\k\cdot \meas(\mathfrak{D})\bigg\}\;.
$$
\erem
\section{Notations and set up\label{parassnot}}
Fix\footnote{For us, $\natural\coloneqq \{1,\cdots\}, \quad \natural_0\coloneqq \{0,1,\cdots\}$.} $d\in \natural\setminus\{1\}$. Let $\O\subset \rn$ be non-empty and bounded domain with piecewise smooth boundary and $\tn\coloneqq \rn/2\p\zn$, the $d$--dimensional torus.\\
Given   
$h,\,r,\,s,\,\a,\,\vae_0,\,\t>0,\, n,p\in \natural, \, y_0\in \cn$ and a non-empty $A\subset \cn$, $A'\subset \rn$, we define the following. 
Let\footnote{As usual  $\o\cdot k\coloneqq \o_1k_1+\cdots+\o_d k_d$.}
\beq{dio}\D_\a^\t\coloneqq \left\{\o\in \rn:|\o\cdot k|\geq \frac{\a}{|k|_1^\t}, \forall\ 0\not=k\in \zn\right\},
\eeq
be the set of $(\a,\t)$--Diophantine numbers, where $|k|_1\coloneqq\dst\sum_{j=1}^d |k_j| $ is the $1$--norm on\footnote{And in general on $\cn$ as well as on all its subsets ($\real^n,\,\integer^n,\,\natural^n,$ etc).} $\cn$  
and
\[\O_\a\coloneqq \left\{\o\in \O\cap \D_\a^\t: \dist(\o,\partial\O)\coloneqq \dst\min_{\o_*\in \partial\O}|\o-\o_*| \geq\a\right\},\]
where $|\cdot|$ is some norm on $\cn$; everywhere in this thesis, we shall use 
\[|x|\coloneqq \dst\max_{1\le j\le d}|x_j|,\]
 the sup--norm on $\cn$, except in $\S\ref{algKam}$ where we shall use $|\cdot|=|\cdot|_1$. Let\footnote{We shall nevertheless drop the dimension $d$ as it is fixed once for all, and write $B_r(y_0)$ instead of $B_r^d(y_0)$, \etc.}
\begin{align*}
\mathbb{J}&\coloneqq  \begin{pmatrix}0 & -\uno_d\\
\uno_d & 0\end{pmatrix}\;,\\
\wh{D}^d_r(y_0) &\coloneqq  \left\{y\in \cn: |y-y_0|_1<r\right\}\;,\\
D^d_r(y_0)&\coloneqq  \left\{y\in \cn: |y-y_0|<r\right\}\;,\\
D^d_r(A)&\coloneqq  \dst\bigcup_{y_0\in A}D^d_r(y_0)\;,\\
{B}^d_r(A')&\coloneqq \rn\cap D^d_r(A') \; ,\\
\dst\torus^d_s &\coloneqq  \left\{x\in \cn: |\Im x|<s\right\}/2\p\zn\;,\\
D^d_{r,s}(y_0)&\coloneqq   \left\{y \in \cn: |y-y_0|<r\right\}\times \torus^d_s\; ,\\
 D^d_{r,s}&\coloneqq  D^d_{r,s}(0)\;,\\
 \O_{\a,h} &\coloneqq  \bigcup_{\o\in\O_\a}D_h(w)\;, \\
 W_{r,s,\vae_0} &\coloneqq  D^d_{rs,s}\times \left\{\vae \in \complex: |\vae|<\vae_0\right\}\;,
\end{align*}
where $\uno_d\coloneqq \diag(1)$ is the unit matrice of order $d$.\\
$\cn\times\cn$ will be equipped with the canonical symplectic form $\varpi\coloneqq dy\wedge dx=dy_1\wedge dx_1+\cdots+dy_d\wedge dx_d$.\\
Given a linear operator $\mathcal{L}\colon (V_1,\|\cdot\|)\to (V_2,\|\cdot\|)$, its ``operator--norm'' is given by
\[\|\mathcal{L}\|\coloneqq \sup_{x\in V_1\setminus\{0\}}\frac{\|\mathcal{L}x\|}{\|x\|},\quad \mbox{so that}\quad \|\mathcal{L}x\|\le \|\mathcal{L}\|\, \|x\|\quad \mbox{for any}\quad x\in V_1.\]
Let $\mathcal{A}_{r,s}(y_0)$ (resp. $\mathcal{B}_{r,s}(y_0),\,\mathcal{A}_{r,s,h,d},\,\mathcal{B}_{r,s,\vae_0}$) be the set of real--analytic functions $f$ on $\wh{D}_r(y_0)\times \torus^d_{s}$ (resp. $D_{r,s}\times \O_{\a,h} ,\, W_{r,s,\vae_0}$) with  finite norm  $\tnorm f_{r,s,y_0}$ (resp. $\|f\|_{r,s,y_0},\, \|f\|_{r,s,h,d},\, \|f\|_{r,s,\vae_0}$), defined below.
Let
\[\mathcal{A}^0_{r,s}(y_0)\coloneqq \left\{f\in \mathcal{A}_{r,s}(y_0): \average{f}\coloneqq \frac{1}{(2\pi)^d}\dst\int_{\tn}f(y,x)\, dx=0,\; \forall y\in \wh{D}_r(y_0)\right\},\]
$\mathcal{A}^0_{r,s,h,d}$ and $\mathcal{B}^0_{r,s,\vae_0}$ are defined analogously. Given $\o\in \rn$ and $f\in \mathcal{A}_{r,s}(y_0)\bigcup \mathcal{A}_{r,s,h,d}\bigcup \mathcal{B}_{r,s,\vae_0}$, we define 
\[D_\o f\coloneqq \o\cdot  f_x=\dst\sum_{j=1}^d \o_j \dst f_{{x}_j},\]
write\footnote{As usual, $(y-y_0)^l\coloneqq \dst\prod_{j=1}^d(y_j-y_{0j})^{l_j}$. Here, and henceforth, $\ex\coloneqq \exp(1)$ denotes the Neper number and $i$ a complex--square--root of $-1$: $i^2=-1$.} 
\[f=\dst\sum_{k\in \zn}f_k \ex^{ik\cdot x}=\dst\sum_{k\in \zn,\\l\in\nn}f_{l,k} \ex^{ik\cdot x}(y-y_0)^l,\]
where $\quad f_k:=\dst\frac{1}{(2\pi)^d}\dst\int_{\tn}f(x) \ex^{-ik\cdot x}\, d x \eqqcolon \sum_{l\in \nn}f_{l,k}\; (y-y_0)^l,\quad k\in \zn,\quad $  define
\[T_\k f:=\dst\sum_{|k|_1\leq \k}f_k \ex^{ik\cdot x},\, \k\in \natural 
\]
and define on $\mathcal{A}_{r,s}(y_0)$ (resp. $\mathcal{B}_{r,s}(y_0),\, \mathcal{A}_{r,s,h,d},\,\mathcal{B}_{r,s,\vae_0}$), the norms\footnote{Notice that the above definitions apply also to vector--valued real analytic functions \ie $f\coloneqq (f_1,\cdots,f_n)$ with $\{f_j\}_{j=1}^n\subset \mathcal{A}_{r,s}(y_0),$ etc.}
\begin{align*}
\tnorm f_{r,s,y_0}\coloneqq & \dst\sum_{k,l\in\zn}|f_{l,k}|_1\; \ex^{s|k|_1}r^{|l|_1}\, \left(\mbox{resp.} \quad \|\cdot\|_{r,s,y_0}\coloneqq \dst\sup_{D_{r,s}(y_0)}|f|,\right.\\
    &\left.
\|\cdot\|_{r,s,h,d}\coloneqq \dst\sup_{D_{r,s}\times \O_{\a,h}}|f|,\quad\quad
 \|\cdot\|_{r,s,\vae_0} \coloneqq \dst\sup_{W_{r,s,\vae_0}}|f|\right).
 \end{align*}
 Moreover, for a matrix--valued periodic functions $\mathcal{L}(y,x)$, we 
 define\footnote{With an analogous definitions with the other norms. For instance, $\|\mathcal{L}\|_{r,s,y_0}\coloneqq\dst\sup_{|a|=1}\|\mathcal{L}(\cdot,\cdot)a\|_{r,s,y_0}.$}
\[\tnorm{\mathcal{L}}_{r,s,y_0}\coloneqq\dst\sup_{|a|_1=1}\tnorm{\mathcal{L}(\cdot,\cdot)a}_{r,s,y_0}.\]
and for a given $f\in\mathcal{A}_{r,s,y_0}$, the Fourier's norm of the $3$--tensor $\dpr^3_x f$ is given by
\[\tnorm{\dpr^3_x f}_{r,s,y_0}\coloneqq \dst\sup_{|b|_1=|c|_1=1}\dst\sum_{j=1}^d\tnorm{\dst\sum_{k,l=1}^d\frac{\dpr^3 f}{\dpr x_j \dpr x_k\dpr x_l}b_l c_k}_{r,s,y_0}.\]
Given a map $\f\colon A\subset\complex^n\to \complex^p$, its Lipschitz constant is defined by
\[\|\f\|_{L,A}\coloneqq \sup_{{x,y\in A, x\neq y}}\frac{|\f(x)-\f(y)|}{|x-y|}\leq\infty.\]
\section{General remarks}
\begin{itemize}
\item[$1.$]We have chosen the norms for simplicity but any others\footnote{An algebra norm is a norm satisfying $\|x\cdot y\|\leq \|x\| \|y\|$, for any $x$ and $y$.} ``algebra norms'' maybe be used.
\item[$2.$] As we are going to compare the four theorems on a concrete Hamiltonian and since we use two different norms, we need a kind of equivalence between them. Indeed, we have, for any $r>0,\, 0<\s<s$
\begin{align*}
\|f\|_{r,s-\s,y_0}\leq \tnorm{f}_{r,s-\s,y_0}&= \dst\sum_{n,m\in\zn} |f_{m,n}|_1\ex^{(s-\s)|n|_1}r^{|m|_1}\\
  &\leq \dst\sum_{n,m\in\zn} \frac{d\|f\|_{r,s,y_0}}{r^{|m|_1}}\ex^{-s|n|_1}\ex^{(s-\s)|n|_1}r^{|m|_1}\\
  &= d \|f\|_{r,s,y_0}\dst\sum_{n\in\zn}\ex^{-\s|n|_1}\\
  &= d \|f\|_{r,s,y_0}\left(\dst\sum_{k\in\integer}\ex^{-\s|k|}\right)^d\\
  &= d\|f\|_{r,s,y_0}\left(1+\frac{2}{\ex^{-\s}-1}\right)^d\\
  &=d\|f\|_{r,s,y_0}\tanh^d\left(\frac{\s}{2}\right).
\end{align*}
\item[$3.$] Through the present thesis, we shall denote by $C$ (resp. $c$), at any place (with index or not), a constant depending (eventually) only on $d,\,\t$ and $\bar{\n}$ (see below) and greater (resp. less) or equal than $1$ \ie\ $C=C(d,\t,\bar{\n})\ge 1$ (resp. $c=c(d,\t,\bar{\n})\le 1$).
\end{itemize}
\part{Classical KAM Theorems and Quantitative normal forms\label{part1}}
\chapter{Quantitative KAM normal forms\label{QTHMS}}
\section{Statement of the explicit KAM normal forms theorems\label{statmtExpKAM}}
\subsection{Kolmogorov's normal form \label{SectKolmoNorm}}
\subsubsection{Assumptions\label{AssumpKolmo}}
Let $\quad\a,\,r,\,\vae_0>0,\quad \t\geq d-1,\quad  0<2\sigma<s\le 1\quad$ and
\[s_*\coloneqq s-2\s.\]
 Let's consider  a hamiltonian $H\in \mathcal{B}_{r,s,\vae_0}$  such that $K(y,x):=H(y,x;0)$ has the form\footnote{As usual, $\o\cdot y=\o_1 y_1+\cdots +\o_d y_d$; $Q=O(|y|^2)$ means that  $\dpr^m_yQ(0,x)=0$ for all $m\in\natural^d$ with $|m|_1\le 1$, where
$\dpr^m_y= \frac{\partial^{|m|_1}\phantom{aaaaa}}{\partial y_1^{m_1}\cdots \partial y_d^{m_d}}$ and $|m|_1=m_1+\cdots+m_d$.}
\beq{Knf}
K=\mathrm{K}+ \o\cdot y+ Q(y,x)   \quad{\rm with}
\quad
Q=O(|y|^2) ,\quad \mathrm{K}\in \real, \quad{\rm and}
\eeq
\beq{dioph}
\o\in \D^\t_\a 
\quad\ie\quad \o \quad{\rm is}\quad (\a,\t)\mbox{--diophantine}.
\eeq
\noi
Furthermore, assume that $K$ in \equ{Knf} is non--degenerate in the sense that\footnote{$\average{\cdot}$ being the average over $\torus^d$.}
\beqno
\det \average{\dpr^2_y Q(0,\cdot)}\neq 0\ .
\eeqno
 Write 
 \[\boxed{H\eqqcolon K+\vae P} \]
  and set 
\[M\coloneqq \|P\|_{r,s,\vae_0}\,,\qquad T\coloneqq \average{\partial^2_yQ(0,\cdot)}^{-1}.
\]
Finally define
\begin{align*}
\mathrm{E}&\coloneqq 2 \wh{\mathrm{E}}\coloneqq 2\max\left(r|\o|,\|Q\|_{r,s,\vae_0},|\o|^2\|T\|\right),\\
\mathcal{W} &\coloneqq  \diag(|\o|\, \uno_d,r\sigma|\o|\, \uno_d),\\
r_* &\coloneqq r(s-2\sigma),\\
B_* &\coloneqq   B_{r_*}(0),\\ 
\mathrm{C}_0 &=  2^{d+1-2\t} \sqrt{\G(2\t+1)}\;,\\
\mathrm{C}_1 &\coloneqq  2\cdot 3^{\t} \mathrm{C}_0 , \\
\mathrm{C}_2 &\coloneqq  2 d \mathrm{C}_1 + 2^{-(\t+1)}  ,\\
\mathrm{C}_3 &\coloneqq  d \mathrm{C}_2+2^{-(\t+2)}, \\
\mathrm{C}_4 &\coloneqq   \mathrm{C}_2 + 2^{-2}\mathrm{C}_1  ,\\
\mathrm{C}_5 &\coloneqq  3^{\t}d\mathrm{C}_0 \left( 2d\mathrm{C}_4+2^{-(\t+3)}\right) , \\
\mathrm{C}_6 &\coloneqq  2^{-(\t+2)}\mathrm{C}_4+\mathrm{C}_5,\\
\mathrm{C}_7 &\coloneqq   \frac{3}{2}d \mathrm{C}_5 +81\cdot 2^{-(\t+3)}d^3\mathrm{C}_4 +9\cdot 2^{-(2\t+5)} d^2 ,\\
\mathrm{C}_8 &\coloneqq   18\mathrm{C}_7,  \\
\mathrm{C}_9 &\coloneqq   9 d^2 \mathrm{C}_6^2 +3\cdot 2^{-(2\t+5)}d \mathrm{C}_6 , \\
\bar{\n} &\coloneqq  4\t+10,\\
  \n     &\coloneqq  4\t+12,\\
\bar{\mathrm{C}} &\coloneqq  \max\left(2^{-(2\t+5)}\mathrm{C}_8,\mathrm{C}_9 \right),\\
\wt{\mathrm{C}} &\coloneqq  d\left(3d\bar{\mathrm{C}} +2^{-(2\t+6)}\mathrm{C}_7\right),\\
\mathrm{C}_\sharp &=\mathrm{C}_\sharp(d,\t)\coloneqq  \frac{9d\cdot 2^{4\t+23}}{5}\left(3d\bar{\mathrm{C}} +2^{-(2\t+6)}\mathrm{C}_7\right), \\
\mathrm{c}&=\mathrm{c}(d,\t)\coloneqq  \su{\mathrm{C}_\sharp}\;,\\
\wh{\mathrm{C}}&\coloneqq  \frac{6d}{5}\left(3d\bar{\mathrm{C}}  +2^{-(2\t+6)}\mathrm{C}_7\right) ,\\
\mathrm{C}&=\mathrm{C}(d,\t)\coloneqq  \frac{\wh{\mathrm{C}}}{3\bar{\mathrm{C}}}\ ,\\
\bar L &\coloneqq  \bar{\mathrm{C}} \mathrm{E}^7\sigma^{-\bar\n}r^{-7}\a^{-4}|\o|^{-3}M, \\
\tilde{L}  &\coloneqq \wt{\mathrm{C}} \mathrm{E}^8\sigma^{-\n}r^{-8}\a^{-4}|\o|^{-4}M,\\
L &\coloneqq  \frac{6}{5}r^{-2}|\o|^{-2}\tilde L \mathrm{E}^2=\frac{\mathrm{C}_\sharp}{3\cdot 2^\n} \wh{\mathrm{E}}^{10}\sigma^{-\n}r^{-10}\a^{-4}|\o|^{-6} M,\\
\end{align*}
\subsubsection{Statement of the KAM Theorem}
\thmtwo{teo2}{Komogorov \cite{Ko1954}, pg. 52}
Under the assumptions in $\S\ref{AssumpKolmo}$, the following hold. There exists 
 a real--analytic  symplectomorphism $\phi_*\colon B_*\times\torus^d\overset{into}{\longrightarrow}  B_{r}(0)\times\tn$, depending analytically also on $\vae\in(-\vae_*,\vae_*)$, with 
 \[\boxed{\vae_*\coloneqq \min\left(\vae_0,\mathrm{c}\;\wh{\mathrm{E}}^{-9}\s^{4\t+13}r^{10}\a^4|\o|^6 M^{-1}\right)},\] such that $\phi_*|_{\vae=0}$ is the identity map and, for any $|\vae|<\vae_*$, 
\beqno
H\circ \phi_*(y',x')=K_*(y',x';\vae):=\mathrm{K}_*(\vae)+\o \cdot y'+ Q_*(y',x';\vae), \quad{\rm with}\quad Q_*=O(|y'|^2)\ 
\eeqno
 and
\beq{teo2Est}
\frac{{\mathrm{C}} \mathrm{E}^3}{ r^3|\o|^{3}\sigma^2}\|\mathcal{W}(\phi_*-\id)\|_{r_*,s_*,\vae_*}\;,\quad |\mathrm{K}-\mathrm{K}_*|_{\vae_*},\quad \|Q-Q_*\|_{r_*,s_*,\vae_*}\;,\quad |\o|^2\|T-T_*\| \le \frac{|\vae|L}{3 \sigma}\ .
\eeq
\ethm
\subsection{KAM Theorem apr\`es Arnold\label{AnolKam}}
\subsubsection{Assumptions\label{AssumpArnol}}
Let $\a,r_0>0,\,\t\ge d-1,\, 0<2\s_0<s_0\leq 1,\, y_0\in\rn$ and consider the Hamiltonian parametrized by $\vae\in\real$
\[H(y,x;\vae)\coloneqq K(y)+\vae P(y,x),\]
with 
$$
K,P\in \mathcal{B}_{r_0,s_0}(y_0)\,.
$$
 such that
\beq{ArnoldCond}
\boxed{\o\coloneqq K_y(y_0)\in \D^\t_\a\;,\qquad\quad \det K_{yy}(y_0)\not= 0.}
\eeq
Set
\[
T\coloneqq K_{yy}(y_0)^{-1},\quad M_0\coloneqq  \|P\|_{r_0,s_0,y_0},\quad\mathsf{K}_0\coloneqq \|K_{yy}\|_{r_0,y_0},\quad \mathsf{T}_0\coloneqq \|T\|\,.
\]
Finally, for a given $\vae\not=0$, 
define\footnote{Notice that $\dst \int_{\rn} |y|_1^{\t}\ex^{-|y|_1}dy\ge \int_{\{|y_j|\ge 1:j=1,\cdots,d\}} |y|_1^{\t}\ex^{-|y|_1}dy \ge d^\t \left(\int_{\{|y_1|\ge 1\}} \ex^{-|y_1|}dy_1\right)^d=$\\$=d^\t\left(2\ex^{-1}\right)^d\ge d^{d-1}\left(2\ex^{-1}\right)^d=d^{\frac{d}{2}-1}\left(2\sqrt{d}\ex^{-1}\right)^d>1$ because $\t\ge d-1\ge 1$. Thus, $\mathsf{C}_0>1$ and $\mathsf{C}_1>1$.\label{ftnarc1}}
\begin{align*}
\mathsf{W}_0   &\coloneqq \diag\left(\max\left\{\frac{\mathsf{K}_0}{\a},\frac{1}{r_0}\right\}\uno_d,\uno_d\right)\;,\\
\eta_0	 &\coloneqq \mathsf{T}_0\mathsf{K}_0\;,\\ 
\n	 &\coloneqq \t+1\;,\\ 
\mathsf{C}_0 &\coloneqq 4\sqrt{2}\left(\frac{3}{2}\right)^{2\n+d}\dst\int_{\rn} \left( |y|_1^{\n}+d|y|_1^{2\n}\right)\ex^{-|y|_1}dy\;,\\
\mathsf{C}_1 &\coloneqq 2\left(\frac{3}{2}\right)^{\n+d}\dst\int_{\rn} |y|_1^{\n}\ex^{-|y|_1}dy\;,\\
\mathsf{C}_2 &\coloneqq d^4 3^{8(d-1)}\;,\\
\mathsf{C}_3 &\coloneqq	d^2\mathsf{C}_1^2+6d\mathsf{C}_1 +1\;,\\
\mathsf{C}_4 &\coloneqq \max\left\{\mathsf{C}_0,\,\mathsf{C}_3\right\}\;,\\
\mathsf{C}_5   &\coloneqq 2^{2(\n+d)+11}3^25^{-2}d^2\;,\\
\mathsf{C}_6 &\coloneqq \dst\max\left\{32d\,,\,10^{-\n}\mathsf{C}_7 \right\}\;,\\
\mathsf{C}_7 &\coloneqq \max\{\mathsf{C}_2\;,\,\mathsf{C}_4\}\;,\\
\mathsf{C}_8 &\coloneqq 3\cdot 5^\n\;\mathsf{C}_6\;,\\
\mathsf{C}_9 &\coloneqq \frac{3\cdot 5^{2\n+1}\sqrt{2}}{8}\;\mathsf{C}_6\;,\\
\mathsf{C}_{10} &\coloneqq \max\left\{1\;,\;\left(\frac{3d\cdot 2^{5-d}}{5}\right)^{\su4}\right\}\;,\\
\mathsf{C}_{11} &\coloneqq \frac{\mathsf{C}_5^2\mathsf{C}_9\mathsf{C}_{10}}{3}\;,\\
s_*        &\coloneqq s_0-2\s_0\;,\\
\mathfrak{p}_1&\coloneqq \mathsf{C}_8\;\eta_0\;\s_0^{-(3\n+2d+1)}\;\max\left\{1,\frac{\a}{r_0\mathsf{K}_0}\right\}\;,\\
\mathfrak{p}_2&\coloneqq \mathsf{C}_{11}\;\eta_0^{\frac{17}{4}}\; \s_0^{-(4\n+2d)}\;,\\
\vae_\sharp&\coloneqq \dst\min\left\{\ex^{-1}\,,\,\exp\left(-\frac{\s_0}{5}\left(\frac{12\sqrt{2}}{5}\frac{\a\mathsf{T}_0}{r_0}\right)^{\su{\n}}\right)
\right\}\;,\\
\m_0      &\coloneqq \frac{\mathsf{K}_0|\vae|M_0}{\a^2}\;.
\end{align*}
\noi
	
\subsubsection{Statement of the KAM Theorem}
\thmtwo{teo4}{Arnold \cite{ARV63}}
Under the assumptions in $\S\ref{AssumpArnol}$, the following hold. For any given $\vae$ satisfying
\beq{CondArmm}
\boxed{
\left\{
\begin{aligned}
&\m_0\leq \;\vae_\sharp\,,\\ \ \\
&\mathfrak{p}_1\cdot\max\left\{1\,,\,\mathfrak{p}_2\;\m_0\;\left(\log\m_0^{-1}\right)^{2\n} \right\}\cdot\m_0\;\left(\log\m_0^{-1}\right)^{\n}< 1\,,
\end{aligned}
\right.
}
\eeq
there exist $y_*\in B_{r_0}(y_0)$ and an embedding $\phi_*\colon \tn\to D_{r_0,s_0}(y_0),$ real--analytic on $\torus^d_{s_*}$ and close to the trivial embedding 
\[\phi_0\colon x\in \tn \to (y_*,x)\in D_{r_0,s_0}(y_0),\]
 such that the $d$--torus
\beq{KronTorArn}
\mathcal{T}_{\o,\vae}\coloneqq \phi_*\left(\tn\right)
\eeq
is a non-degenrate invariant Kronecker torus for $H$ \ie
\beq{KronTorArnIE}
\phi^t_H\circ \phi_*(x)=\phi_*(x+\o t).
\eeq
Moreover, 
\beq{estArnTr}
|\mathsf{W}_0(\phi_*-\id)|
\le \s_0^{d+1}  \;,
\eeq
uniformly on  $\{y_*\}\times \torus^d_{s_*}\,.$
\ethm
\rem{ArnolImplyKolmoo}
It is not difficult to see that Theorem~\ref{teo4} is stronger than Theorem~\ref{teo2}. Indeed, let the assumptions in \S\ref{AssumpKolmo} hold. Then, Taylor's expansion yields $H(y,x)=K^\sharp(y)+P^\sharp(y,x)$, where
\begin{align*}
&K^\sharp(y)\coloneqq \mathrm{K}+\o\cdot y+\su2(T^{-1}y)\cdot y\;,\\
&P^\sharp(y,x)\coloneqq \su2\left((\dpr^2_yQ(0,x)-T^{-1})y\right)\cdot y+3\sum_{|\b|_1=3}y^\b\int_0^1\frac{(1-t)^2}{\b!}\dpr^\b_yQ(ty,x)\; dt+\vae P(y,x)\;.
\end{align*}
Thus, 
$$
 K^\sharp_y(0)=\o\in \D^\t_\a\;,\qquad \det K^\sharp_{yy}(0)=\det T^{-1}\not= 0\qquad\mbox{and}\qquad \|P^\sharp\|_{r,s,\vae_0}=O(r^2)+O(\vae)\;.
$$
Consequently, by choosing $r$ proportional to $\sqrt{\vae}$, one can apply Theorem~\ref{teo4} and recover Theorem~\ref{teo2}.
\erem
\subsection{KAM Theorem apr\`es J. Moser (following J. P\"oschel) \label{parstatkampar}}
\subsubsection{Assumptions\label{AssumpPosc}} 
Let $r,h>0,\quad 0<s\leq 1$ and consider the hamiltonian parametrized by $\o$
\[H(y,x,\o)\coloneqq N(I,\o)+P(y,x,\o) \mbox{ where } N(y,\o)\coloneqq K_0(\o)+\o\cdot y,\mbox{ with } K_0,\, P\in \mathcal{A}_{r,s,h,d}\;,\]
 and $X_H$ and $X_N$ the hamiltonian vector fields associate to $H$ and $N$ respectively with respect to the canonical symplectic form $\varpi$. Let $\phi^t_H$ and $\phi^t_N$ be the hamiltonian flow associate to $H$ and $N$ respectively. We have then, $\phi^t_N(y ,x)=(y , x+\o\cdot t)$. 
 Let\footnote{Notice that one could use any $1>\b\geq 1-\frac{1}{{\bar{\n}}}+\frac{1}{\n}$.}
\[\a>0,\quad \n>\bar{\n}=\t+1>d,\quad \epsilon\coloneqq \|P\|_{r,s,h,d}\;,\quad \b\coloneqq 1-\frac{\n-\bar{\n}}{\n{\bar{\n}}},\]
\[\bar{c}\coloneqq \min\left(1,\frac{\n-\bar{\n}}{\n{\bar{\n}}}\ex\right)\;,\quad B\coloneqq B_r(0)\;,\quad W \coloneqq \diag\left( r^{-1}\uno_d,20 s^{-1}\uno_d\right), \]
\[ \wt W\coloneqq \diag\left( 20^{-\bar{\n}}\a s^{\bar{\n}}r^{-1}\uno_d,20^{-\t}\a s^{\t}\uno_d, \uno_d\right),\]
and define\footnote{$\G(z)\coloneqq \int_0^\infty t^{z-1}\ex^{-t}dt, \Re z>0$ is the Euler's gamma function.  Notice that if $\t\ge d-1$ then $C_0\geq 2^{d+1-2\t}\dst \sqrt{4^{[2\t]-3}\cdot 6}\geq \sqrt{3}$, where $[2\t]$ denotes the integer part of $2\t$.}
\beqano
C_0 &\coloneqq & 2^{d+1-2\t} \sqrt{\G(2\t+1)}\;,\\
C_1 &\coloneqq & \frac{4\ex}{3}\;,\\
C_2 &\coloneqq & \dst\sum_{j=0}^{d-1}\frac{(d-1)!}{(d-1-j)!(d-1)^{j+1}}, \\
C_3 &\coloneqq & 4(d+1)C_0\;,\\
C_4 &\coloneqq & 16(d+1)C_0 \;,\\
C_5 &\coloneqq & 16(6C_1+1)C_3\;,\\
C_6 &\coloneqq & 2dC_5\dst\sum_{j=0}^\infty 2^{-2\bar{\n}\left(\left(\frac{3}{2}\right)^j-j-1\right)-j} ,\\
C_7 &\coloneqq & 36d(d+1)(6C_1+1)\dst\exp\left(\frac{36d(d+1)(6C_1+1)}{C_6} \dst\sum_{j=0}^{\infty} 2^{-{\bar{\n}}\left(2\left(\frac{3}{2}\right)^j-j-2\right)}\right),\\
C_8 &\coloneqq & C_7 \dst\sum_{j=0}^{\infty} 2^{-\bar{\n}\left(2\left(\frac{3}{2}\right)^j-j-2\right)-j},\\
C_9 &\coloneqq &\max\left(1,\frac{6C_1+1}{2C_0}\right),\\
C_{10} & \coloneqq  & 9\cdot 4^{\bar{\n}}\left(\max \dst\left\{ 48d(d+1)^2 C_0 , 4^d(d+1)C_2\right\}\right)^2,\\
C_{11} &\coloneqq & \exp\left(\su2\left(\left(\frac{2\bar{\n}}{\bar{c}}\right)^{1/\b}+\frac{1}{20}\right) \right),\\
C_{12}&\coloneqq & \left(\frac{2}{\bar{c}}\right)^{\n}\left( 2\ex^{\frac{1}{40}}C_6\right)^{\n/{\bar{\n}}},\\
C_{13}&\coloneqq & \exp\left(\su2\left(\left(\frac{C_0}{6C_1+1}\right)^{1/\bar{\n}}+\frac{1}{20}\right)\right)
\eeqano
and\footnote{Notice that $C_{12}>C_6$. Moreover, if $\vae<\ex^{\su{40}}\min(20^{\n}C_{11}^{-1},C_{12}^{-1})\s^{\n}$ then $\k_0^{{\bar{\n}}} \s_0^{{\bar{\n}}} \ex^{-\k_0\s_0}\leq \vae  < 1/(2C_6\k_0^{{\bar{\n}}})$ and therefore $\k_0\s>d-1$; compare Appendix~\ref{app:A}, with $\vth=\su2$ and $\vae$ replaced by $\Th\ex^{\su{40}}$.\label{footn5}}
\beq{DefCstPoscSC}
\begin{aligned}
C  &=C(d,\t)  \coloneqq \ex^{\su{40}} C_6,\\
C_*&=C_*(d,\t) \coloneqq  \max\left(C_6C_9,C_8\right),\\
c  &=c(d,\t,\bar{\n})\coloneqq  20^{-\n}\min\left(\frac{1}{4^{\bar{\n}} C_{10}},\frac{20^{\n}\ex^{\su{40}}}{C_{11}},\frac{\ex^{\su{40}}}{C_{12}},\frac{20^{\n}}{C_{13}}\right). 
\end{aligned}
\eeq
Finally, let
\[\Phi_0\colon (x,\o)\in\tn\times \rn\mapsto (0,x)\in B\times \tn,\] be the trivial embedding.
\subsubsection{Statement of the KAM Theorem}
\thmtwo{teo1}{P\"oschel \cite{JP}}
Let $H,\,\Phi_0,\, \epsilon,\, \n,\,\bar{\n},\,C,\, C_*$ and $c$ as in $\S\ref{AssumpPosc}$ and assume that 
\beq{eqAssThm00}
\boxed{\vae\coloneqq \frac{\epsilon}{\a r}\leq c s^{\n}\quad \mbox{and}\quad C\vae\leq \frac{h}{\a}\leq \frac{1}{2\k_0^{{\bar{\n}}}}}\;,
\eeq
where
\[\k_0 \coloneqq  \left[-\frac{40\log\vae-1}{s}\right].\]
Then, there exist a Lipeomorphism $\f\colon\O\to \O $ close to the identity 
and a Lipschitz continuous family of real analytic Lagrangian torus embeddings $\Phi\colon\tn\times \O\to B\times \tn$ closed to the trivial embedding $\Phi_0$ 
such that the following hold. For any $\o\in \O_\a$, $\Phi(\tn,\o)$ is a Lagrangian submanifold and an invariant Kronecker torus for $H_{|\f(\o)}$ with $H_{|\f(\o)}(y ,x,\o)\coloneqq H(y ,x,\f(\o))$, $\ie$
\beq{eq1} \ph^t_{H_{|\f(\o)}}\circ \Phi(x,\o)=\Phi(x+\o t;\o), \forall\, x\in \tn.
\eeq
Moreover, $\Phi(\tn,\o)$ is a Lagrangian submanifold and the maps $x\mapsto \Phi(x,\o)$ is real analytic on $\dst\torus^d_{\frac{s}{2}}$ for each given $\o\in\O$ and one has uniformly on $\dst\dst\torus^d_{\frac{s}{2}}\times \O$ and $\O$ respectively, the following estimates\footnote{Here and in $\S\ref{prteo1}$ as well, we shall denote by $\|f\|_{L,A}$, the uniform Lipschitz' semi-norm of the function $f$ \wrt the $\o$--argument (parameter) varying in the set $A$.  }
\beq{eq2} \left\|W(\Phi-\Phi_0) \right\|,\quad h\left\|W(\Phi-\Phi_0) \right\|_{L,\rn}\leq C_*\frac{\|P\|_{r,s,h}}{ r h}, 
\eeq
\beq{eq3} \left\|\f-\id \right\|,\quad h\left\|\f-\id \right\|_{L,\O}\leq C_*\frac{\|P\|_{r,s,h}}{r}.
\eeq
\ethm

\rem{rem001}

\noi
\begin{itemize}
\item[$1.$] If one chooses $h=\frac{\epsilon C}{r}$, then the assumption~\eqref{eqAssThm00} in Theorem~\ref{teo1} reduced to
\[\epsilon\leq c r\a s^\n.\]
\item[$2.$] Notice that we have some freedom in the choice of $C_6$. Indeed, one just needs to chose
\beq{choC6}
C_6 > \frac{C_5}{2C_0\log 2}\dst\sum_{j=0}^\infty 2^{-\bar{\n}(2\m^j-j-2)}\ .
\eeq
\item[$3.$] To be precise, in Theorem~\ref{teo1},
\[\Phi\colon \{0\}\times\tn\times \O\to B\times \tn.\]
\item[$4.$] Notice that the $\n$ in Theorem~\ref{teo1} is larger than the $\n=\t+1$ P\"oschel uses in Theorem~$A$ and $B$ in \cite{JP}. In fact the Theorem~$A$ and $B$   are not valid for $\n=\t+1$. Indeed\footnote{We are using here the same notations as in \cite{JP}.}, assuming the contrary, then for any $\vae=\frac{\epsilon}{\a r}\leq \g s^\n$, there would exists $\k_0\in \natural$ such that
\beq{NuPoscContr}
\k_0^\n\s^\n\ex^{-\k_0\s}\leq \vae\leq \frac{1}{2\k_0^\n},\quad \mbox{with}\quad \s=\frac{s}{20}.
\eeq
But then\footnote{Because $\k_0\s\geq 1$.}, $\ex^{-\k_0\s}\leq \vae$ \ie $\k_0\s\geq \log\vae^{-1}$. Hence, we would have, if we take $\vae=\g s^\n$ in particular, for any $0<s\leq 1$,
\[\su2\frac{s^\n}{20^\n}\overset{\eqref{NuPoscContr}}{\geq} \vae \k_0^\n\s^\n\geq \vae \left(\log \vae^{-1}\right)^\n=\g s^\n \left(\log (\g s^\n)^{-1}\right)^\n\]
\ie
\[1\geq 2\cdot 20^\n \g \left(\log (\g s^\n)^{-1}\right)^\n,\quad \forall\; 0<s\leq 1;\, \mbox{contradiction}.\]

\end{itemize}
\erem
\subsection{KAM Theorem apr\`es Salamon--Zehnder 
\label{algKam}}
\subsubsection{Assumptions\label{AssumpAlg}} 
Let 
\[0<\hat{s}<s\leq \bar{s}\le 1,\quad \sigma\coloneqq \frac{s-\hat{s}}{2},\quad r,\a,\mathtt{E},E_{j,k}\geq 0,\quad\t\ge d-1,\] where $j,k\in \natural$. Let's consider a hamiltonian $H\in \mathcal{A}_{r,\bar{s}}(y_0)$, for some $y_0\in\rn$ and a pair of real--analytic functions $(u,v)$ on $\torus^d_s$ such that
\[\tnorm{(H_{yy})^{-1}}_{r,{{\bar s}},y_0} \leq \mathtt{E},\quad\tnorm{\partial^j_x\partial^k_y H}_{r,{{\bar s}},y_0}\leq E_{j,k}\;,\]
for any $j,k\in \natural$ and
\[\quad\tnorm{u}_s\le U,\quad \tnorm{v}_s\le V \quad \mbox{and}\quad \r\coloneqq \tnorm{v-y_0}_s<r,\]
 for some 
\[0\leq U\le \bar{s}-s \quad \mbox{and}\quad V>0,\]
\[\mathcal{M}\coloneqq \uno_d+u_\th,\quad \mbox{and}\quad H^0_{yy}(\th)\coloneqq H_{yy}\left(v(\th),\th+u(\th)\right)\]
are invertible for each given $\th\in\tn$ and, defining\footnote{$A^{-T}$ stands for the transpose of the inverse of $A$: $A^{-T}\coloneqq \left(A^{-1}\right)^T$.},
\[\mathcal{T}\coloneqq \mathcal{M}^{-1}H^0_{yy}\mathcal{M}^{-T},\]
$\average{\mathcal{T}}$ is invertible. Let $\o\in \D^\t_\a$ and define $f$ and $g$ by
\beq{DeffEg}
\left\{
\begin{array}{r c l}
\o+D_\o u-H_y(v,\id+u)&=& f\\
D_\o v+ H_x(v,\id+u)&=& g
\end{array}
\right.
\eeq
 Futhermore, assume that
\[\tnorm{\mathcal{M}}_s\le \mathtt{M},\quad \tnorm{\mathcal{M}^{-1}}_s\le \ovl{\mathtt{M}},\quad \tnorm{v_\th}_s\le \wt V,\quad \tnorm{f}_s\le F,\quad \tnorm{g}_s\le G,\quad |\average{\mathcal{T}}^{-1}|\le \wt T, \]
for some $\wt V, F, G\geq 0 $ and $\mathtt{M},\ovl{\mathtt{M}}, \wt T>0$. Finally, define 
\beqano 
\wh V        & \coloneqq & \dst\max\{\wt V, r-\r\},\\
\mathtt{A}_1 & \coloneqq & |\o|\; ,\quad \mathtt{A}_2 \coloneqq \max\left\{\mathtt{A}_1\; ,\, E_{1,1} \right\},\\
\mathtt{A}_3 & \coloneqq & \max\left\{E_{2,1}\;,\, \mathtt{E} E_{1,2}\mathtt{A}_2\; ,\,\mathtt{E}^2 E_{0,3} \mathtt{A}_2^2\right\},\\
\mathtt{A}_4 & \coloneqq & \max\left\{E_{3,0}\; ,\, \mathtt{E} E_{2,1} \mathtt{A}_2\; ,\, \mathtt{E}^2 E_{1,2}\mathtt{A}_1^2 \right\},\\
\mathtt{A}_5 & \coloneqq & \max\left\{\mathtt{A}_4\; , \, \mathtt{E}\mathtt{A}_1 \mathtt{A}_2 \right\},\quad \mathtt{A}_6  \coloneqq  \max\left\{ \mathtt{A}_5\; ,\, \wh V \mathtt{A}_3 \right\},\\
\mathtt{A}_7 & \coloneqq & \max\left\{\mathtt{E} E_{0,2}\; ,\, E_{0,2}\wt T \right\}\cdot \max\left\{\a^{-2} E_{0,2}\mathtt{A}_6\; ,\,\a^{-1}\mathtt{A}_3 \right\},\\
\mathtt{A}_8 & \coloneqq & (s-\hat s)^{2\t}\max\left\{1\; ,\,\frac{\mathtt{E}\mathtt{A}_2}{r-\r}\; ,\, E_{0,2}\wt T\; ,\, E_{1,2}\wt T\; ,\, \mathtt{E} E_{0,3}\mathtt{A}_2\wt T \right\},\\
\mathtt{A}_9 & \coloneqq & \max\left\{\mathtt{A}_7\; ,\,\mathtt{A}_8 \right\}\;,\\
 \mathtt{A}_* & \coloneqq & \max\left\{\mathtt{E} E_{0,2}\; ,\, E_{0,2}\wt T \right\}\cdot \max\left\{\a^{-1} F\; ,\,\a^{-2} E_{0,2}\wh V F\; ,\,\a^{-2}E_{0,2} G \right\}.
\eeqano
\subsubsection{Statement of the KAM Theorem}
\thmtwo{teo3}{Celletti--Chierchia \cite{CACL97}}
Under the assumptions in $\S\ref{AssumpAlg}$, the following holds. There exists a polynomial $\X $ in $(s,\sigma)$ satisfying
\beq{Teo3PolEs}
\frac{5}{4}\le \X(a,b)\le 21+88a\le 109,\quad \forall\;0< a \le 1, \forall\; 0<b<\su 2,
\eeq
such that, if
\beq{CondTeo3}
\boxed{\mathtt{A}_* \mathtt{A}_9 \mathtt{M}^7 \ovl{\mathtt{M}}^9 (s-\hat s)^{-2(2\t+1)} 2^{8\t+13}\t!^4\;\X(s,\sigma)\le 1
},\eeq
then there exists $(\tilde{u},\tilde{v})$, real--analytic on $\torus^d_{\hat{s}}$, $\mathtt{A}_*$--close to $(u,v)$  and solving 
\beq{EmbPartDiffEq}
\left\{
\begin{array}{r c l}
\o+D_\o \tilde{u}-H_y(\tilde{v} ,\id+\tilde{u} )&=& 0\\
D_\o \tilde{v}+ H_x(\tilde{v} ,\id+\tilde{u} )&=& 0
\end{array}
\right.
\eeq
Futhermore, $\average{\tilde{u}}=\average{u}$ and the solution $(\tilde{u},\tilde{v})$ is uniquely determined in the $\mathtt{A}_*$--neighborhood of $(u,v)$ by the condition $\average{\tilde{u}}=\average{u}$.
\ethm

For a proof, see \cite{CACL97}.
\rem{rm1CC}
Notice that instead of the bound $\wt V$ on $v_\th$ used in \cite{CACL97} to define the parameters, here we use $\wh V$. The point is that, with this change, one is then allowed to chose $\wt V=0$ when $v$ is constant.
\erem
\section{Some preliminary facts}
As we are going to use the same idea as in \cite{CL90} to extend maps obtained through the \textit{KAM step} in the proof of \textit{Theorem}~\ref{teo1}, we will need a \textit{cut--off} function.
\lemtwo{cutoff}{Cut--Off}
Then, for any $n\in\natural$, there exists a constant $\mathcal{C}_n >0$ 
such that for any given $R>0$ and a non--empty $\D\subset\rn$, there exists $\chi\in C(\cn)\cap C^{\infty}(\rn)$ with $0\leq \chi\leq 1,\ supp\chi \subset \D_R\coloneqq \dst\bigcup_{\o\in\D}D_R(\o),\ \chi\equiv 1$ on $\D_{\frac{R}{2}}$ 
and for any $m\in \natural^d$ with $|m|_1\leq n$, 
\beq{cutof}
\|\partial_\o^m\chi\|_0 \overset{def}{=} \sup_{\rn}\|\partial_\o^m\chi\|\leq \mathcal{C}_n \frac{(|m|_1+2)!}{R^{|m|_1}}\,.
\eeq
\elem
\proof
Let $a,b>0$ such that $0<a\leq\su4$ and $\su2+a\le b \leq 1-a$. Consider
\[\chi_1\colon t\in\real\mapsto\chi_1(t)= \left\{\begin{aligned}
\ex^{-\frac{1}{1-t^2}} & \quad \mbox{if} \quad |t|<1\\
0                    & \quad \mbox{if} \quad |t|\ge 1
\end{aligned} \right. ,\quad 
\chi_d(\o) \coloneqq R^{-d}N_a\prod_{j=1}^d\chi_1\left(\frac{\o_j}{aR}\right), \, \mbox{with}\]
\[      N_a \coloneqq   \left(a\dst\int_{\real}\chi_1(t)dt \right)^{-d}.\]
Define\footnote{$\dst\chi_{\D_{bR}}$ denotes the characteristic function of the set $\D_{bR}$, \ie $\dst\chi_{\D_{bR}}\equiv 1$ on $\D_{bR}$ and $\dst\chi_{\D_{bR}}\equiv 0$ on $\rn\setminus \D_{bR}$.}
\[\chi(\o)\coloneqq \dst\int_\rn \dst\chi_{\D_{bR}}(y)\chi_d(\o-y)dy.\]
Thus, $\chi\in C(\cn)\cap C^{\infty}(\rn)$ and
\begin{itemize}
\item $supp\chi\subset \D_R$. Indeed, for any $\o\in\rn$, $\chi(\o)>0$ implies that $\dst\chi_{\D_{bR}}(y)>0$ and $\chi_d(\o-y)>0$ for \textit{a.e.} $y\in\rn$ (with respect to the Lebesgue measure on $\rn$); which implies in particular that there exist $y_0\in\rn$ and $\o_*\in\D$ such that $|y_0-\o_*|<bR$ and $|\o-y_0|<aR$, so that $|\o-\o_*|<(a+b)R\leq R$ \ie $supp\chi\subset \D_R$.
\item $\chi\equiv 1$ on $\dst\D_{\frac{R}{2}}$. Indeed, let $\o\in \dst\D_{\frac{R}{2}}$ \ie $|\o-\o_*|<\frac{R}{2}$, for some $\o_*\in \D$. Then, for any $y\in\rn$,
\[|\o-y|<aR\implies |y-\o_*|\leq |y-\o|+|\o-\o_*|<(a+\su2)R\leq bR\implies y\in \D_{bR}.\]
Hence,
\[1\geq \chi(\o)\geq \dst \int_{B^d_{aR}(\o)}\chi_d(\o-y)dy=\dst N_a\left(a\dst\int_{\real}\chi_1(t)dt \right)^{d}=1.\]
\end{itemize}
Moreover, for any $\o\in \D_R$ and for any $n\in\natural$, we have\footnote{In fact, one checks easily that for any $m\in\natural^d$ , $\|\partial_\o^m\chi\|_0\leq \|\partial_{\o_1}^{|m|_1}\chi\|_0=\cdots=\|\partial_{\o_d}^{|m|_1}\chi\|_0$.}
\begin{align*}
\partial^n_{\o_1}\chi(\o)&=\dst\int_\rn \dst\chi_{\D_{bR}}(y)\partial^n_{\o_1}\chi_d(\o-y)dy=\dst\int_{B^d_{aR}(\o)} \dst\chi_{\D_{bR}}(y)\partial^n_{\o_1}\chi_d(\o-y)dy\\
  &= \dst\int_{B^d_{aR}(\o)} \dst\partial^n_{\o_1}\chi_d(\o-y)dy=R^{-d}N_a\dst\int_{\rn} \dst\partial^n_{\o_1}\chi_1\left(\frac{\o_1-y_1}{aR}\right)\prod_{j=2}^d \chi_1\left(\frac{\o_j-y_j}{aR}\right) dy\\
  &=R^{-d}N_a\left((aR)^{-n+1}\dst\int_{\real}\frac{d^n\chi_1(t)}{dt^n}dt \right)\left(aR\dst\int_{\real}\chi_1(t)dt \right)^{d-1}\\
  &=\left((aR)^{-n}\dst\int_{\real}\frac{d^n\chi_1(t)}{dt^n}dt \right)\left(\dst\int_{\real}\chi_1(t)dt \right)^{-1}\left(aR\dst\int_{\real}\chi_1(t)dt \right)^{d}R^{-d}N_a\\
  &=(aR)^{-n}\left(\dst\int_{\real}\chi_1(t)dt \right)^{-1} \dst\int_{\real}\frac{d^n\chi_1}{dt^n}(t)dt.
\end{align*}
One checks easily that for any $t\in\real$ and $n\in\natural$,
\[\frac{d^n\chi_1}{dt^n}(t)=\frac{P_n(t)}{(1-t^2)^{2n}}\chi_1(t),\]
with $\quad P_0\coloneqq 1,\quad P_1(t)\coloneqq -2t\quad$ and for any $n\geq1$
\[ P_{n+1}\coloneqq (-2+4n(1-t^2))tP_n(t)+(1-t^2)^2\frac{dP_n}{dt}(t),\]
and one has, for any $n\geq 1$,
\[\deg(P_n)=3n-2.\]
The existence of the sequence $(\mathcal{C}_n)_n$ then follows easily and in particular, by choosing $a=\frac{1}{4}$ and, then, $b=\frac{3}{4}$, we can take\footnote{We have $\frac{1}{2\ex^2}\leq \dst\int_{\real}\chi_1(t)dt\leq \frac{2\sqrt{6}}{7}$.}
\beq{ConstC1HomEq}
 \mathcal{C}_1\coloneqq \frac{2\ex^2}{3!a}\dst\int_{\real}\left|\frac{d\chi_1}{dt}\right|(t)dt=\frac{4\ex}{3}
\eeq
and
\begin{align*}
 \frac{2\ex^2}{4!a^2}\dst\int_{\real}\left|\frac{d^2\chi_1}{dt^2}\right|(t)dt &\leq\frac{4\ex^2}{3}\cdot 8\dst\int^1_0\frac{1}{(1-t)^4}\ex^{-\frac{1}{2(1-t)}}dt\\
   &= \frac{32\ex^2}{3}\dst\int^\infty_1 s^4 \ex^{-\frac{s}{2}}ds= \frac{832\ex^\frac{3}{2}}{3}\eqqcolon \mathcal{C}_2.
\end{align*}
\qed
The following lemma establishes a bound on the Lipschitz constant of a map obtained as the composition of infinetly many Lipschitz maps and will be used to prove the Lipschitz continuity of the map in \textit{Theorem}~\ref{teo1}, obtained through the infinite iterative KAM scheme.
\lem{LipInf}
Let $(\mathcal{X},\|\cdot\|)$ be a real or complex normed vector space, $\mathcal{L}_j\colon \mathcal{X}\to \mathcal{X}$ be a sequence of invertible linear operators and $(g_j)_{j\geq 0}$ be a sequence of Lipschitz continuous maps from $\mathcal{X}$ to itself. 
Define 
\[l_j\coloneqq \|\mathcal{L}_j(g_j-\id)\mathcal{L}^{-1}_0\|_{L,\mathcal{X}},\quad G_{0}\coloneqq \id,\quad G_{j+1}\coloneqq g_0\circ \cdots \circ g_j,\]
\[ L_{j}\coloneqq \|\mathcal{L}_0(G_j-\id)\mathcal{L}^{-1}_0\|_{L,\mathcal{X}}\quad \mbox{and assume that}\]
\beq{LipLemHyp000}
\d\coloneqq \sup_{j\geq 0} \|\mathcal{L}_j\mathcal{L}^{-1}_{j+1}\|<\infty \quad \mbox{and}\quad \dst\sum_{j=0}^\infty\d^j l_k\leq \su2.
\eeq
Then $(G_j)_{j\geq 0}$ is a sequence of Lipschitz continuous maps and for any $j\geq 0$,
\beq{LipLemEst0077}
\|\mathcal{L}_0(G_{j+1}-\id)\mathcal{L}^{-1}_0\|_{L,\mathcal{X}}\leq 2\dst\sum_{k=0}^j\d^k l_k.
\eeq
\elem
\proof
For any $j\geq 0$, we have
\beqano
L_{j+1}&= & \left\|\mathcal{L}_0(G_j-\id)\mathcal{L}^{-1}_0 \left\{\left(\mathcal{L}_0 \mathcal{L}^{-1}_1\right) \cdots \left(\mathcal{L}_{j-1} \mathcal{L}^{-1}_j\right) \mathcal{L}_j(g_j-\id)\mathcal{L}^{-1}_0+\id\right\}\right.+\\
      &   &+\left.\left(\mathcal{L}_0 \mathcal{L}^{-1}_1\right) \cdots \left(\mathcal{L}_{j-1} \mathcal{L}^{-1}_j\right)\mathcal{L}_j(g_j-\id)\mathcal{L}^{-1}_0\right\|_{L,\mathcal{X}}\\
      &\leq & (1+\d^jl_j)L_j+\d^jl_j.
\eeqano
Hence, we get, inductively, that for any $j\geq 0$, 
$L_{j+1}\leq L_j+2\d^j l_j\leq 2\dst\sum_{k=0}^j\d^k l_k$. 
\qed

We recall the very famous \textit{Cauchy estimate} used to control the derivatives of an analytic function.
\lemtwo{Cau}{Cauchy's estimate}
Let $p\in \natural,\,f\in \mathcal{A}_{r,s,h,d}$.
Then,
 for any multi--index $(l,k)\in \natural^d\times\natural^d$ with $|l|_1+|k|_1\le p$ and for any $0<r'<r,\, 0<s'<s$,\footnote{As usual, $\dpr_y^l\coloneqq \frac{\dpr^{|l|_1}}{\dpr y_1^{l_1}\cdots\dpr y_d^{l_d}},\, \forall\, y\in\rn,\, l\in\zn $.\label{notDevPart}}
\[\|\partial_{y}^l \partial_{x}^k f\|_{r',s',h,n}\leq p!\; \|f\|_{r,s,h,n}(r-r')^{|l|_1}(s-s')^{|k|_1}.\]
\elem
In the next lemma, we recall some properties of the Fourier's coefficients of an analytic function.
\lem{fce}
Let $\k\in \natural,\,f\in \mathcal{A}_{r,s,h,d},\, 0<\s < s$ with $\k>\frac{d-1}{\s}$. Then
\beqano
&&(i)\qquad\  |f_k(y,\o)|\leq \ex^{-|k|_1 s}\|f\|_{r,s,h,d}\; ,\quad \forall\; k\in \zn,\,y\in D_r(0),\,\o\in \O_{\a,h},\\
&&(ii)\qquad \|f-T_\k f\|_{r,s-\s,h,d}\leq 4^d C_2 \k^d \ex^{-\k\s}\|f\|_{r,s,h,d}. 
\eeqano
\elem

\proof
\begin{itemize}
\item[$(i)$]Let $k\in \zn, \,y\in D_r(0),\,\o\in \O_{\a,h} $. Then
\[f_k(y,\o) = \dst\frac{1}{(2\pi)^d}\dst\int_{\tn}f(y,x,\o)\ex^{-ik\cdot x}\, dx.\]
But, for any given $\b\in\rn\ \mbox{ such that } |\b|<s$, by periodicity of $v$ in each argument and Cauchy's theorem, we get
 \[f_k(y,\o) = \dst\frac{1}{(2\pi)^d}\dst\int_{\tn}f(y,x-i\b,\o)\ex^{-ik\cdot (x-i\b)}\, dx.\]
Now, we choose $\b\coloneqq (s-\s)\left( \sign(k_1),\cdots,\sign(k_d)\right)$ with $0<\s<s$. Thus, we get
\[|f_k(y,\o)|\leq  \ex^{-|k|_1(s-\s)}\|f\|_{r,s,h,d}\, ,\]
and letting $\s \to 0^+$, we get the desired inequality.
\item[$(ii)$] We have
\begin{eqnarray*}
\|f-T_\k f\|_{r,s-\s,h,d}
			   & \leq & \|f\|_{r,s,h,d} \dst\sum_{|k|_1>\k} \ex^{-|k|_1\s}\\
			   & =   & \|f\|_{r,s,h,d} \dst\sum_{ l=\k+1}^\infty\sum_{\substack{k\in \zn \\ |k|_1=l}} \ex^{-|k|_1\s}\\
			   & \leq & \|f\|_{r,s,h,d} \dst\sum_{ l=\k+1}^\infty 4^d l^{d-1}\ex^{-l\s}\\
			   & \leq & \|f\|_{r,s,h,d} \dst\int_\k^\infty 4^d t^{d-1} \ex^{-t\s}\, dt \\
			   & \leq & 4^d \k^d \ex^{-\k\s} \|f\|_{r,s,h,d} \dst\int_0^\infty (t+1)^{d-1} \ex^{-t(d-1)}\, dt\\
			   &= & 4^d C_2 \k^d \ex^{-\k\s} \|f\|_{r,s,h,d}.
\eeqano
\end{itemize}
\qed
In the following lemma, we recall some facts about the homological equation.
\lemtwo{sde}{\cite{CC95}}
Let $p\in \natural,\,\o\in \dst\D_\a^\t$ and $f\in \mathcal{A}^0_{r,s,h,d}$. Then, for any $0<\s<s$, the equation
\[D_\o g=f\]
has a unique solution in $\mathcal{A}^0_{r,s-\s,h,d}$ and there exist constants $\bar B_p=\bar B_p(d,\t)\ge 1$ and $k_p=k_p(d,\t)\ge 1$ such that for any multi--index $k\in \natural^d$ with $|k|_1\le p$
\[\|\partial_\o^k g\|_{r,s-\s,h,d}\leq \bar B_p\frac{\|f\|_{r,s,h,d}}{\a} \s^{-k_p} .\]
In particular, one can take $\bar B_1=C_0$ 
 (see \cite{RH75,CC95}). 
\elem
\ \\ 
Now, we recall the classical implicit function theorem, in a quantitative framework.
\lemtwo{IFTLem}{Implicit Function Theorem I\cite{CL12}}
Let $ r,s>0,\, n,m\in \natural,\, (y_0,x_0)\in \complex^n\times\complex^m$ and \footnote{Let us point out that any other norm (different!) may be used on $\complex^n,\, \complex^m$ and $\complex^{n+m}$.}
\[F\colon (y,x)\in D^n_r(y_0)\times D^m_s(x_0)\subset \complex^{n+m}\mapsto F(y,x)\in\complex^n\]
be continuous with continuous Jacobian matrix $F_y$. Assume that $F_y(y_0,x_0)$ is invertible with inverse $T\coloneqq F_y(y_0,x_0)^{-1}$ 
 such that
\beq{HypIFT}
\sup_{D^n_r(y_0)\times D^m_s(x_0)}\|\uno_n-TF_y(y,x)\|\leq c<1 \quad \mbox{and}\quad \sup_{ D^m_s(x_0)}|F(y_0,\cdot)|\leq \frac{(1-c) r}{\|T\|}.
\eeq
Then, there exists a unique continuous function $ g\colon D^m_s(x_0)\to D^n_r(y_0)$ such that the following are equivalent
\begin{itemize}
\item[$(i)$] $(y,x)\in D^n_r(y_0)\times D^m_s(x_0)$ and $F(y,x)=0$;
\item[$(ii)$] $x\in D^m_s(x_0)$ and $y=g(x)$.
\end{itemize}
Moreover, $g$ satisfies
\beq{EstIFT}
\sup_{D^m_s(x_0)}|g-y_0|\leq \frac{\|T\|}{1-c}\sup_{D^m_s(x_0)}|F(y_0,\cdot)|.
\eeq
\elem

\noi
Finally, we recall some inversion function theorems.\\
Taking $n=m,\ c=c'=\su2$ and $F(y,x)=f(y)-x$ in Lemma~\ref{IFTLem}, for a given $f\in C^1(D^n_r(y_0),\complex^n)$, then the following holds.
\lemtwo{inv1}{Inversion Function Theorem I}
Let $f\colon D^n_r(y_0)\to \complex^n$ be a $C^1$ function with invertible Jacobian $f_y(y_0)$ and assume that
\[\dst\sup_{D^n_r(y_0)}\|\uno_n-Tf_y\|\leq \su2,\quad T\coloneqq f_y(y_0)^{-1}.\]
Then, there exists a unique $C^1$ function 
\[g\colon D^n_s(x_0)\to D^n_r(y_0),\quad x_0\coloneqq f(y_0),\quad s\coloneqq \frac{r}{2\|T\|},\]
such that
\[f\circ g(x)=x,\quad g\circ f(y)=y,\quad \forall\; x\in D^n_s(x_0),\quad \forall\; y\in g\left(D^n_s(x_0)\right).\]
Moreover,
\beq{estInvg}
\dst\sup_{D^n_s(x_0)}\|g_x\|\le 2\|T\|\;.
\eeq
\elem
\rem{remIFTLm}
\begin{itemize}
\item[$(i)$] Notice that in Lemmata~\ref{IFTLem} and \ref{inv1} if, in addition, $F$ is periodic in $x$ (resp. analytic, real on reals) then so is $g$.
\item[$(ii)$] Notice that Lemmata~\ref{IFTLem} and \ref{inv1} still hold if, everywhere therein, open balls are replaced by closed balls or 
complex--balls by real--balls. 
\end{itemize}
\erem
\noi
Another consequence of the Implicit Function Theorem is the following version of Inversion Function Theorem.
\lemtwo{inv}{Inversion Function Theorem II \cite{JP}}
Assume that $f$ is a real analytic function from $\O_{\a,h}$ into $\cn$ such  that 
\[\left\|f-\id \right\|\leq \d \leq \frac{h}{4}.\]
Then $f$ has a real analytic inverse $g$ defined on $\dst \O_{\a,\frac{h}{4}}$ and on which it satisfies
\[\|g-\id\|,\quad \frac{h}{4}\|Dg-\Id\|\leq \d.\]
\elem
\rem{rem2}
Notice that, all the Lemmata above are valid if one replace $(\mathcal{A}_{r,s,h,d},\|\cdot\|_{r,s,h,d})$ by $(\mathcal{B}_{r,s,\r_0},\|\cdot\|_{r,s,\r_0})$ or $(\mathcal{B}_{r,s}(y_0),\|\cdot\|_{r,s,y_0})$.
\erem
\section{Proofs}
\subsection{Proof of Theorem~\ref{teo2}}
The proof is essentially the one in \cite{CL08} though one needs to re--scale various quantities; therefore we shall skip some details. First of all, notice that $K,P\in \cB_{s,\vae_0}$. For simplicity, sometimes, the explicit dependence on $r$ or $\vae_0,\,\vae_*$ will not be denoted in the norm $\|\cdot\|_{r,s,\vae_0}$ or in the $\cB$--spaces, etc, as $r,\,\vae_0$ and $\vae_*$ will not be changed during the iteration. We begin by describing completely one step of the scheme, namely the KAM step, which will be then iterated infinitely many time to compute the symplectic change of variable.

\Giu
{\bf KAM step}
Kolmogorov's idea is to construct a near--to--the--identity symplectic transformation $\phi_1$, such that 
\beq{H1}
H_1:= H\circ \phi_1=K_1+\vae^2 P_1\ ,\quad
K_1= K_1^* +\o\cdot y'+
Q_1(y',x')\ , Q_1=O(|y'|^2)\ ;
\eeq
if this is achieved, the Hamiltonian $K_1$ has the same basic properties of $K$
(the linear part in $y$ is the same and, being $\phi_1$ close to the identity,  $K_1$ is non--degenerate) and the procedure can be iterated.

\nl
For, Kolmogorov considers the generating function of $\phi_1$ of the form\footnote{Compare \cite{DAS01,AKN06} for generalities on symplectic transformations and their generating functions.
For simplicity, we do not report in the notation the dependence of 
various functions on $\vae$, but, in fact, $P=P(y,x;\vae)$, $\mathrm{s}=\mathrm{s}(x;\vae)$, $a=a(x;\vae)$, etc.}
\beq{generating}
g(y',x):=y'\cdot x+\vae\big(b\cdot x+ \mathrm{s}(x) +  y' \cdot a(x)\big)\ ,
\eeq
where, $\mathrm{s}$  and $a$ are (respectively, scalar and vector--valued, $\vae$--dependent)  real--analytic functions on
$\torus^d$ with zero average and $b\in \real^d$.
Define\footnote{As usual, we denote $\mathrm{s}_x=\dpr_x \mathrm{s}=(\mathrm{s}_{x_1},...,\mathrm{s}_{x_d})$ and  
$a_x$ denotes the  matrix   $ (a_x)_{ij}:=\frac{\dpr a_j}{\dpr x_i}$; as above, we often do not report in the notation the dependence upon $\vae$ (but $u_0$, $A$ and $u$ do depend also on $\vae$).}
\beqno u_0=u_0(x):= b + \mathrm{s}_x \ ,\qquad
A=A(x):= a_x\quad {\rm and} \quad u=u(y',x):=u_0+A \, y'\ . 
\eeqno
Then $\phi_1$ is implicitely defined by
\beqno
\left\{
\begin{array}{l}
y=y'+\vae u(y',x):=y'+\vae(u_0(x)+A(x)y')\\
\ \\
x'=x+\vae a(x) \ .
\end{array}\right.
\eeqno
Moreover, for $\vae$ small, $x\in\torus^d\to x+\vae a(x)\in\torus^d$ defines a diffeomorphism of $\torus^d$ with inverse
\beqno
x=\f(x'):=x'+\vae \wt{\f}(x';\vae)\ ,
\eeqno
for a suitable real--analytic function $\wt{\f}$.
Thus $\phi_1$ is explicitly given by
\beq{phi1}
\phi_1:(y',x') \to 
\left\{
\begin{array}{l}y=y'+\vae\; u\big(y',\f(x')\big)\\\ \\
x=\f(x') \ .
\end{array}\right.
\eeq
To determine $b$, $s$ and $a$, observe that by Taylor's formula
\beq{f1}
H(y'+\vae\; u,x)=\mathrm{K}+\o\cdot y'+Q(y',x)+\vae \Big[
\o\cdot u + Q_y(y',x)\cdot u+P(y',x)\Big] 
+ \vae^2 P'(y',x)
\eeq
where $P':=P'(y',x;\vae):=P^\ppu+ P^\ppd$ with
\beq{P'}
\left\{
\begin{array}{l}
P^\ppu:= \su{\vae^2}[ Q(y'+\vae u,x)- Q(y',x)-\vae Q_y(y',x)\cdot u]=
{\dst\igl01}(1-t)Q_{yy}(y'+t\vae u,x)\, u\cdot u\, dt\\  
P^\ppd:=\su{\vae}[P(y'+\vae u,x)-P(y',x)]=
{\dst\igl01}P_{y}(y'+t\vae u,x)\cdot u\, dt\ .
\end{array}
\right.
\eeq
Note that 
\beq{Q1'}
Q_y(y',x)\cdot (a_x y')=: Q^\ppu(y',x)=O(|y'|^2)\ ,
\eeq
and that (again by Taylor's formula)
\beq{Q23'}
\left\{
\begin{array}{l}
Q_y(y',x)\cdot u_0=Q_{yy}(0,x) y'\cdot u_0 + Q^\ppd(y',x)\ ,\quad
Q^\ppd:= {\dst\igl01}(1-t)Q_{yyy}(ty',x)y'\cdot y'\cdot u_0\, dt\\ \ \\
P(y',x)=P(0,x)+P_y(0,x)\cdot y' + Q^{(3)}(y',x)\ ,\quad
Q^{(3)}:={\dst\igl01}(1-t)P_{yy}(ty',x)\, y'\cdot y'\, dt\ .
\end{array}
\right.
\eeq
Thus, since\footnote{Recall that $\o\cdot \mathrm{s}_x=D_\o \mathrm{s}$ and $\o\cdot (a_x y')=(D_\o a)\cdot y'$.} $\o\cdot u=\o\cdot b+D_\o \mathrm{s}+ D_\o a\cdot y'$, we find 
\beq{HH}
H(y'+\vae u,x)=\mathrm{K}+\o\cdot y' + Q(y',x)+\vae Q'(y',x)+\vae F'(y',x)+\vae^2P'(y',x)
\eeq
with $P'$ as in \equ{f1}--\equ{P'} and
\beq{F'Q'}
\left\{
\begin{array}{l}
Q'(y',x):= Q^\ppu+Q^\ppd+Q^{(3)}=O(|y'|^2)\ \\ \ \\
F'(y',x):=\o\cdot b + D_\o \mathrm{s}+P(0,x)+\Big\{ D_\o a+ Q_{yy}(0,x) b + Q_{yy}(0,x) \mathrm{s}_x + P_y(0,y')\Big\}\cdot y'
\end{array}
\right.
\eeq
By Lemma~\ref{sde}, {\sl there exist a unique constant $b$ and unique functions  $s$ and $a$} (with zero average)
{\sl such that $F'$ is constant}. In fact, if
\beqno
\left\{\begin{array}{l}
\mathrm{s}:=- D_\o^{-1} \Big(P(0,x)-P_0(0)\Big)= - \sum_{n\in\integer^d \atop n\neq 0}{\dst \frac{P_n(0)}{ i \o\cdot n} \ex^{i n\cdot x}}\\
b:= - \average{Q_{yy}(0,\cdot)}^{-1}\Big( \average{Q_{yy}(0,\cdot) \mathrm{s}_x}+ \average{P_y(0,\cdot)}\Big) \\
a:=-D_\o^{-1}\Big( Q_{yy}(0,x)(b+\mathrm{s}_x) + P_y(0,x)\Big)
\end{array}\right.
\eeqno
then $F'=\o\cdot b+ P_0(0)$. Thus, with this determination of $g$ in \equ{generating}, recalling \eqref{phi1}, we find that  \equ{H1} holds with
\beqno
\left\{
\begin{array}{l}
\mathrm{K}_1\coloneqq \mathrm{K}+\vae \wt{\mathrm{K}}\ , \phantom{\ AAAAAAAAAAAAA}\ \wt{\mathrm{K}}\coloneqq \o\cdot b + P_0(0)\\  
Q_1(y',x')\coloneqq Q(y',x')+\vae \wt Q(y',x')\ ,\qquad \wt Q\coloneqq {\dst \igl01}Q_x(y',x'+t\vae\a(x'))\cdot\a\, dt\ 
+Q'(y',\f(x'))\\
P_1(y',x')\coloneqq P'(y',\f(x'))\ .
\end{array}
\right.
\eeqno
Clearly, for $\vae$ small enough $\average{\dpr^2_{y'} Q_1(0,\cdot)}$ is invertible and,
if  $T:=\average{Q_{yy}(0,\cdot)}^{-1}$,  we may write
\beq{defT1}
T_1:=\average{\dpr^2_{y'} Q_1(0,\cdot)}^{-1}=:
T+\vae \wt T\ .
\eeq

\nl
Next, we provide the KAM step with carefull estimates. Actually, we shall do the estimates in term of a lower bound of $|\o|$ instead of $|\o|$, so that, by taking this lower bound equal to $|\o|$, we shall get the estimates in Theorem~\ref{teo2}\footnote{The point is that, if we replace $|\o|$ by $\underline{\o}$ everywhere in \S\ref{SectKolmoNorm}, except in the expressions of $\mathrm{E}$ and $\wh{\mathrm{E}}$, then, Theorem~\ref{teo2} holds for any $\o$ such that $|\o|\geq \underline{\o}$.}. Thus, we fix, for the remainder of the proof,
\beq{LowBNo}
0<\underline{\o}\leq |\o|\ .
\eeq
 Recall the definition in $\S\ref{SectKolmoNorm}$; in particular\footnote{ The notation in Eq. \equ{C} means that each term on the l.h.s. is bounded by the r.h.s.
}
\beq{C}
 r|\o|\, ,\, \|Q\|_{s,\vae_0} ,  |\o|^2\|T\|< \mathrm{E}\ .
\eeq
Finally, fix\footnote{The parameter $s'$ will be the size of the domain of analyticity of the new symplectic variables $(y',x')$, domain on which we shall bound the Hamiltonian $H_1=H\circ\phi_1$, while ${{\bar s}}$ is an intermediate domain where we shall bound various functions of $y'$ and $x$. Note that $\sigma<\frac{1}{2}$.
} 
\beqno
0<\sigma<\frac{s}{2} \quad {\rm and\ define}\quad{{\bar s}}:=s-\frac23\sigma\ ,\quad s':=s-\sigma\ .
\eeqno
\lem{lem:1}
Then\footnote{Here $\|\cdot\|_{{{\bar s}}}=\|\cdot\|_{{{\bar s}},\vae_0}$.}
\beq{es1}
\begin{array}{l}
\underline{\o}\|\mathrm{s}_x\|_{{{\bar s}}},
\underline{\o}|b|,|\wt{\mathrm{K}}|, r\sigma\underline{\o}\|a\|_{{{\bar s}}},r\underline{\o}\|a_x\|_{{{\bar s}}},
\underline{\o}\|u_0\|_{{{\bar s}}},\underline{\o}\|u\|_{{{\bar s}}}, \|Q'\|_{{{\bar s}}}, r^2\sigma^2\|\dpr^2_{y'}  Q'(0,\cdot)\|_0\le \bar L\
\end{array}
\eeq
Furthermore, if $\vae_*\le \vae_0$ satisfies 
\beq{cond1}
\vae_* \, r^{-1}\sigma^{-1}\underline{\o}^{-1} \bar L\le \frac{\sigma}{3}\ ,
\eeq
then
\beq{es1prim}
\|P'\|_{\bar s} \le r^{-1}\sigma^{-1}\underline{\o}^{-1}\bar L M \ .
\eeq
and the following hold. For $|\vae|<\vae_*$,    the map $\psi_\vae(x):=  x+\vae a(x)$ has an analytic inverse 
 $\f(x')=x'+\vae \wt{\f}(x';\vae)$  such that, for all $|\vae|<\vae_*$, 
\beq{boundal}
\|\wt{\f}\|_{s',\vae_*}\le  r^{-1}\sigma^{-1}\underline{\o}^{-1} \bar L \qquad {\rm and}\ , \forall\ |\vae|<\vae_*\ ,
\quad
\f=\id + \vae \wt{\f} : \torus^d_{s'}\to\torus^d_{{\bar s}}
\ ;
\eeq
for any $(y',x,\vae)\in W_{{{\bar s}}_,\vae_*}$, $|y'+\vae u(y',x)|<rs$; the map $\phi_1$ is a symplectic diffeomorphism and
\beq{phiok}
\phi_1=\big( y'+\vae u(y', \f(x')),\f(x')\big): W_{s',\vae_*}\to D_{rs}\times\torus^d_s,\quad{\rm and}
\quad \|\mathcal{W}\,\tilde \phi\|_{s',\vae_*}\le \bar L\ ,
\eeq
where $\tilde \phi$ is defined by the relation $\phi_1=:\id + \vae \tilde \phi$.

\noi
Finally,  if\;\footnote{Notice that $L\ge (r^{-1}|\o|^{-1}\mathrm{E})^2\tilde L\ge \tilde L\ge r^{-1}\sigma^{-1}|\o|^{-1}\bar LC\ge \bar L $ since $r|\o|\le \mathrm{E}$, so that \equ{cond2} implies \equ{cond1}.\label{footnote 20}}
\beq{cond2}
\vae_* \, \frac{L}{\mathrm{E}} \le \frac{\sigma}{3}\ ,
\eeq
then 
\beq{tesit}
\left\{
\begin{array}{l}
|\wt{\mathrm{K}}|\ ,\ \|\wt Q\|_{s',\vae_*}\ ,\ |\o|^2\|\wt T\|\ ,\ \|\mathcal{W}\,\tilde \phi\|_{s',\vae_*}\le L
\\ \ \\
\|P_1\|_{s',\vae_*}\le \frac{L M}{\mathrm{E}}\ .
\end{array}
\right.
\eeq
\elem

\proof 
We begin by estimating $\|\mathrm{s}_x\|_{{{\bar s}}}$. Actually these estimate will be given on a larger intermediate domain, namely, $W_{s-\frac\sigma3,\vae_0}$, allowing to give the remaining bounds on the smaller domain $W_{{{\bar s}}_,\vae_0}$. Let $f(x):=P(0,x)-\average{P(0,\cdot)}$. By definition of $\|\cdot\|_s$ and $M$, it follows that $\|f\|_s\le \|P(0,x)\|_s+\|\average{P(0,\cdot)}\|\le  2 M$. 
By Lemma~\ref{sde} with $p=1$ and $s'=s-\frac\sigma3$, one gets $$\|\mathrm{s}_x\|_{s-\frac\sigma3}\le C_0 \frac{2M}\a \, 3^{\t} \sigma^{-\t}=\mathrm{C}_1 \sigma^{-\t}\a^{-1}M\le \bar L,$$
so that
$$\underline{\o}\|\mathrm{s}_x\|\le 2^{-(3\t+10)}\mathrm{C}_1 \mathrm{E}^7\sigma^{-\bar{\n}} r^{-7}\a^{-4}\underline{\o}^{-3}M\le \bar{L}$$
\noi
Next, we estimate $b$. By definitions and Lemma~\ref{Cau}, we have\footnote{Remember that $\sigma<1/2 ,\,\|T\|\le |\o|^{-2}\mathrm{E}\le \underline{\o}^{-2}\mathrm{E}$ and $r\a\le r|\o|\le \mathrm{E}$.}
\beqano
|b| &\le & \|T\|\left(\dst\max_{1\le l \le d}\sum_{j=1}^d \|Q_{y_l y_j}\|_{s-\sigma}\|\mathrm{s}_{x_j}\|_{s-\frac\sigma3}+\dst\max_{1\le j \le d} \|P_{y_j}\|_{s-\sigma}\right)\\
   &\le & \underline{\o}^{-2}\mathrm{E}\left(2 d \mathrm{E} r^{-2}\sigma^{-2}\mathrm{C}_1 \sigma^{-\t}\a^{-1}M+ M r^{-1}\sigma^{-1}\right)\\
   &\le &  \left(2d \mathrm{C}_1 \mathrm{E} + \sigma^{\t+1}r\a\right)\underline{\o}^{-2}\mathrm{E} \sigma^{-(\t+2)}r^{-2}\a^{-1}M\\
   &\le &\left(2d \mathrm{C}_1 + 2^{-(\t+1)}\right)\underline{\o}^{-2}\mathrm{E}^2 \sigma^{-(\t+2)}r^{-2}\a^{-1}M\\
   &= & \mathrm{C}_2 \mathrm{E}^2 \sigma^{-(\t+2)}r^{-2}\a^{-1}\underline{\o}^{-2}M,
\eeqano
so that
$$\underline{\o}|b| \le 2^{-(3\t+8)}\mathrm{C}_2 \mathrm{E}^7\sigma^{-\bar{\n}} r^{-7}\a^{-4}\underline{\o}^{-3}M\le \bar{L}$$
\noi
Next, we estimate $\wt{\mathrm{K}}$. We have
\beqano
|\wt{\mathrm{K}}|&\le & d|\o|\cdot|b|+\|P\|_{s}\\
	   &\le & d\underline{\o}\mathrm{C}_2 \mathrm{E}^2 \sigma^{-(\t+2)}r^{-2}\a^{-1}\underline{\o}^{-2}M+M\\
       &\le &(d\underline{\o}^{-1}\mathrm{C}_2\mathrm{E}^2+\sigma^{\t+2}r^2\a)\sigma^{-(\t+2)}r^{-2}\a^{-1}M \\
       &\le & (d\underline{\o}^{-1}\mathrm{C}_2\mathrm{E}^2+2^{-(\t+2)}\underline{\o}^{-1}r^2|\o|^2)\sigma^{-(\t+2)}r^{-2}\a^{-1}M\\
       &= & \mathrm{C}_3 \mathrm{E}^2\sigma^{-(\t+2)}r^{-2}\a^{-1}\underline{\o}^{-1}M\\
       &\le & 2^{-(3\t+8)}\mathrm{C}_3 \mathrm{E}^7\sigma^{-\bar{\n}} r^{-7}\a^{-4}\underline{\o}^{-2}M \le \bar L.
\eeqano
\noi
Next, we estimate $u_0$. We have
\beqano
\|u_0\|_{{{\bar s}}} &\le & |b|+ \|s_x\|_{s-\frac\sigma3}\le \mathrm{C}_2 \mathrm{E}^2 \sigma^{-(\t+2)}r^{-2}\a^{-1}\underline{\o}^{-2}M + \mathrm{C}_1 \sigma^{-\t}\a^{-1}M \\
		&\le & (\mathrm{C}_2\mathrm{E}^2\underline{\o}^{-2}+ \mathrm{C}_1r^{2}\sigma^{2}) \sigma^{-(\t+2)}r^{-2}\a^{-1}M\\
		&\le & (\mathrm{C}_2\mathrm{E}^2\underline{\o}^{-2}+ 2^{-2}\mathrm{C}_1r^{2}|\o|^2\underline{\o}^{-2}) \sigma^{-(\t+2)}r^{-2}\a^{-1}M\\
		&\le & (\mathrm{C}_2+ 2^{-2}\mathrm{C}_1)\mathrm{E}^2 \sigma^{-(\t+2)}r^{-2}\a^{-1}\underline{\o}^{-2}M\\
		& = & \mathrm{C}_4 \mathrm{E}^2 \sigma^{-(\t+2)}r^{-2}\a^{-1}\underline{\o}^{-2}M,
\eeqano
so that
$$\underline{\o}\|u_0\|_{{{\bar s}}}\le  2^{-(3\t+8)}\mathrm{C}_4 \mathrm{E}^7\sigma^{-\bar{\n}} r^{-7}\a^{-4}\underline{\o}^{-3}M \le \bar L$$
\noi
Next, we estimate $a$ and $a_x$. Let $f(x)\coloneqq Q_{yy}(0,x)(b+\mathrm{s}_x)+P_y(0,x)$. Then, by Lemma~\ref{Cau}\ref{sde}, we have
\beqano
\|f\|_{s-\frac\sigma3}&\le &\dst\max_{1\le l \le d}\sum_{j=1}^d \|Q_{y_l y_j}\|_{s-\sigma}\left(|b_j|+\|\mathrm{s}_{x_j}\|_{s-\frac\sigma3}\right)+\dst\max_{1\le j \le d} \|P_{y_j}\|_{s-\sigma}\\
	&\le &\dst\max_{1\le l \le d}\sum_{j=1}^d \|Q_{y_l y_j}\|_{s-\sigma}\left(|b|+\|\mathrm{s}_{x}\|_{s-\frac\sigma3}\right)+\dst\max_{1\le j \le d} \|P_{y_j}\|_{s-\sigma}\\
	&\le & 2d \mathrm{E}\sigma^{-2}r^{-2}\cdot \mathrm{C}_4 \mathrm{E}^2 \sigma^{-(\t+2)}r^{-2}\a^{-1}\underline{\o}^{-2}M+ \sigma^{-1}r^{-1}M\\
	&\le & \left( 2d \mathrm{C}_4\mathrm{E}^3\underline{\o}^{-2}+\sigma^{\t+3}r^3\a\right)  \sigma^{-(\t+4)}r^{-4}\a^{-1}M\\
	&\le &\left( 2d \mathrm{C}_4\mathrm{E}^3\underline{\o}^{-2}+\sigma^{\t+3}r^3|\o|^3\underline{\o}^{-2}\right)  \sigma^{-(\t+4)}r^{-4}\a^{-1}M\\
	&\le &\left( 2d \mathrm{C}_4+2^{-(\t+3)}\right)  \mathrm{E}^3\sigma^{-(\t+4)}r^{-4}\a^{-1}\underline{\o}^{-2}M
\eeqano
Thus, by Lemma~\ref{sde}, we obtain\footnote{The factor $d$ comes from the fact \quad $\|a_x\|_{{{\bar s}}}=\|A\|_{{{\bar s}}}\le \dst\max_{1\le l \le d}\sum_{j=1}^d \|A_{lj}\|_{{{\bar s}}}.$}
\beqano
\|a\|_{{{\bar s}}},\|a_x\|_{{{\bar s}}}&\le &dC_0\frac{\|f\|_{s-\frac\sigma3}}{\a}3^{\t}\sigma^{-\t}\\
	&\le & 3^{\t}dC_0 \left( 2d \mathrm{C}_4 +2^{-(\t+3)}\right)  \mathrm{E}^3\sigma^{-(2\t+4)}r^{-4}\a^{-2}\underline{\o}^{-2}M\\
	&= & \mathrm{C}_5\mathrm{E}^3\sigma^{-(2\t+4)}r^{-4}\a^{-2}\underline{\o}^{-2}M,
\eeqano
so that
$$r\sigma\underline{\o}\|a\|_{{{\bar s}}},r\underline{\o}\|a_x\|_{{{\bar s}}} \le 2^{-(2\t+6)}\mathrm{C}_5 \mathrm{E}^7\sigma^{-\bar{\n}} r^{-7}\a^{-4}\underline{\o}^{-3}M \le \bar L $$
\noi
Next, we estimate $u$. We have
\beqano
\|u\|_{{{\bar s}}} &\le &\|u_0\|_{{{\bar s}}}+\|A y'\|_{{{\bar s}}}\le \|u_0\|_{{{\bar s}}}+ \dst\max_{1\le l \le d}\sum_{j=1}^d  \|A_{lj}\|_{{{\bar s}}}r{{\bar s}}\\
		&\le &  \mathrm{C}_4 \mathrm{E}^2 \sigma^{-(\t+2)}r^{-2}\a^{-1}\underline{\o}^{-2}M +d\frac{\mathrm{C}_5}{d} \mathrm{E}^3\sigma^{-(2\t+4)}r^{-3}\a^{-2}\underline{\o}^{-2}M\\
		&\le &  (2^{-(\t+2)}\mathrm{C}_4r\a+\mathrm{C}_5\mathrm{E})\mathrm{E}^2\sigma^{-(2\t+4)}r^{-3}\a^{-2}\underline{\o}^{-2}M \\
		&\le &  (2^{-(\t+2)}\mathrm{C}_4+\mathrm{C}_5)\mathrm{E}^3\sigma^{-(2\t+4)}r^{-3}\a^{-2}\underline{\o}^{-2}M \\
		&= &    \mathrm{C}_6 \mathrm{E}^3\sigma^{-(2\t+4)}r^{-3}\a^{-2}\underline{\o}^{-2}M,
\eeqano
so that
$$\underline{\o}\|u\|_{{{\bar s}}}\le  2^{-(2\t+6)}\mathrm{C}_6 \mathrm{E}^7\sigma^{-\bar{\n}} r^{-7}\a^{-4}\underline{\o}^{-3}M \le \bar L$$
\noi
Next, we estimate $Q'$. To do this, we need to estimate $Q^\ppu,Q^\ppd$ and $Q^\ppt$. By definitions and Lemma~\ref{Cau}, we have
\beqano
\|Q^\ppu\|_{{{\bar s}}} &\le &\dst\sum_{1\le l,j\le d}\|Q_{y_l}\|_{{{\bar s}}}\|A_{lj}\|_{{{\bar s}}}r{{\bar s}}\\ 
	&\le & d^2 \mathrm{E}\frac{3}{2}\sigma^{-1}r^{-1}\cdot \frac{\mathrm{C}_5}{d} \mathrm{E}^3\sigma^{-(2\t+4)}r^{-3}\a^{-2}\underline{\o}^{-2}M\\
	&=& \frac{3}{2}d \mathrm{C}_5 \mathrm{E}^4\sigma^{-(2\t+5)}r^{-4}\a^{-2}\underline{\o}^{-2}M 
\eeqano
\beqano
\|Q^\ppd\|_{{{\bar s}}} &\le & \dst\int_0^1 (1-t)\dst\sum_{1\le j,l,k\le d }\|Q_{y_jy_ly_k}\|_{{{\bar s}}}\| y_j'\|_{{{\bar s}}}\| y_l'\|_{{{\bar s}}}\| ({u_0})_k\|_{{{\bar s}}}dt\\
	&\le & \frac{d^3}{2}6\mathrm{E}\frac{27}{8}\sigma^{-3}r^{-3}\cdot r^{2}\bar{s}^2\cdot  \mathrm{C}_4 \mathrm{E}^2 \sigma^{-(\t+2)}r^{-2}\a^{-1}\underline{\o}^{-2}M \\
	&=  & \frac{81d^3\mathrm{C}_4}{8} \mathrm{E}^3 \sigma^{-(\t+5)}r^{-3}\a^{-1}\underline{\o}^{-2}M\\
	&\le & \frac{81d^3\mathrm{C}_4}{8} \mathrm{E}^4 \sigma^{-(\t+5)}r^{-4}\a^{-2}\underline{\o}^{-2}M\\
	&\le & 81\cdot 2^{-(\t+3)}d^3\mathrm{C}_4 \mathrm{E}^4 \sigma^{-(2\t+5)}r^{-4}\a^{-2}\underline{\o}^{-2}M
\eeqano
and
\beqano
\|Q^\ppt\|_{{{\bar s}}} &\le & \dst\int_0^1 (1-t)\dst\sum_{1\le j,l\le d }\|P_{y_ly_j}\|_{{{\bar s}}}\| y_j'\|_{{{\bar s}}}\| y_l'\|_{{{\bar s}}}dt\\
	&\le & \frac{d^2}{2}2M\frac{9}{4}\sigma^{-2}r^{-2}\cdot r^{2}\bar{s}^2\ \\
	&\le & \frac{9d^2}{4} \sigma^{-2}M
\eeqano
Thus
$$\|Q'\|_{{{\bar s}}} \le \|Q^\ppu\|_{{{\bar s}}}+\|Q^\ppd\|_{{{\bar s}}}+\|Q^\ppt\|_{{{\bar s}}}\le \mathrm{C}_7 \mathrm{E}^4\sigma^{-(2\t+5)}r^{-4}\a^{-2}\underline{\o}^{-2}M,$$
so that
$$\|Q'\|_{{{\bar s}}}\le  2^{-(2\t+5)}\mathrm{C}_7 \mathrm{E}^7\sigma^{-\bar{\n}} r^{-7}\a^{-4}\underline{\o}^{-3}M \le \bar L.$$
Finally, we estimate $\dpr_{y'}^2Q'(0,\cdot)$. We have, once again by Lemma~\ref{Cau},
$$\|\dpr_{y'}^2 Q'(0,\cdot)\|_0\le \|\dpr_{y'}^2Q'(0,\cdot)\|_{s-\sigma}\le 2\mathrm{C}_7 \mathrm{E}^4\sigma^{-(2\t+5)}r^{-4}\a^{-2}\underline{\o}^{-2}M\cdot 9\sigma^{-2}r^{-2}=\mathrm{C}_{8}\mathrm{E}^4\sigma^{-(2\t+7)}r^{-6}\a^{-2}M, $$
so that
$$r^2\sigma^2\|\dpr_{y'}^2Q'(0,\cdot)\|_0\le  2^{-(2\t+5)}\mathrm{C}_8 \mathrm{E}^7\sigma^{-\bar{\n}} r^{-7}\a^{-4}\underline{\o}^{-3}M \le \bar L.$$
\noi
Now, under the assumption \equ{cond1}, we prove \equ{es1prim}. For $(y',x;\vae)\in W_{{{\bar s}},\vae_0}$ and $0\le t \le 1$, by \equ{es1} one has \beq{ineq1}
|y'+t\vae u(x)|\le r{{\bar s}} + \vae \|u\|_{{{\bar s}}} \le r{{\bar s}} + \vae_* \underline{\o}^{-1}\bar L \le r{{\bar s}} + r\frac{\sigma}{3} =rs-r\frac{\sigma}{3}<rs,
\eeq
so that
\beqano
\|P^\ppu\|_{{{\bar s}}}&\le & \dst\int_0^1 (1-t)\dst\sum_{1\le j,l\le d }\|Q_{y_ly_j}\|_{s-\frac\sigma3}\| u_j\|_{{{\bar s}}}\| u_l\|_{{{\bar s}}}dt\\
	&\le & \frac{d^2}{2}2 \mathrm{E}\cdot 9\sigma^{-2}r^{-2}\cdot \left(\mathrm{C}_6 \mathrm{E}^3\sigma^{-(2\t+4)}r^{-3}\a^{-2}\underline{\o}^{-2}M \right)^2\\
	&=& 9d^2\mathrm{C}_6^2 \mathrm{E}^7 \sigma^{-(4\t+10)}r^{-8}\a^{-4}\underline{\o}^{-4} M^2 
\eeqano
 and
\beqano
 \|P^\ppd\|_{{{\bar s}}}&\le & \dst\int_0^1 \dst\sum_{1\le j\le d }\|P_{y_j}\|_{s-\frac\sigma3}\| u_j\|_{{{\bar s}}}dt\\
        &\le & d  M\cdot 3\sigma^{-1}r^{-1}\cdot \mathrm{C}_6 \mathrm{E}^3\sigma^{-(2\t+4)}r^{-3}\a^{-2}\underline{\o}^{-2}M\\
        &=   & 3d \mathrm{C}_6 \mathrm{E}^3\sigma^{-(2\t+5)}r^{-4}\a^{-2}\underline{\o}^{-2}M^2\\
        &\le & 3\cdot 2^{-(2\t+5)}d \mathrm{C}_6 \mathrm{E}^3\underline{\o}^{-2}r^4\a^2|\o|^2 \sigma^{-(4\t+10)}r^{-8}\a^{-4}\underline{\o}^{-2} M^2\\
        &\le & 3\cdot 2^{-(2\t+5)}d \mathrm{C}_6 \mathrm{E}^7 \sigma^{-(4\t+10)}r^{-8}\a^{-4}\underline{\o}^{-4} M^2.
\eeqano
Thus
$$\|P'\|_{{{\bar s}}}\le \|P^\ppu\|_{{{\bar s}}}+\|P^\ppd\|_{{{\bar s}}}\le \mathrm{C}_9 \mathrm{E}^7 \sigma^{-(4\t+10)}r^{-8}\a^{-4}\underline{\o}^{-4} M^2\le r^{-1}\sigma^{-1}\underline{\o}^{-1}\bar L M.$$

\noi
Next, we show how \equ{cond1} implies the existence of the inverse of $\psi_\vae$ satisfying \equ{boundal}. 
The defining relation $\psi_\vae\circ \f=\id$ implies that  $\wt{\f}(x')=-a(x'+\vae \wt{\f}(x'))$,  where $\wt{\f}(x')$ is short for $\wt{\f}(x';\vae)$ and  such relation   is a fixed point equation for   the non--linear operator $f: u \to f(v):=-a(\id + \vae v)$. To find a fixed point for this equation one can use a standard contraction Lemma (see  \cite{Ko.man}). Let $Y$ denote the closed ball  (with respect to the sup--norm) of continuos functions $v:\torus^d_{s'}\times\{|\vae|< \vae_*\}\to\complex^d$ such that $\|v\|_{s',\vae_*}\le r^{-1}\sigma^{-1}\underline{\o}^{-1}\bar L$. By \equ{cond1}, $|\Im(x'+\vae v(x'))|<s'+\vae_* r^{-1}\sigma^{-1}\underline{\o}^{-1}\bar L<s'+\frac\sigma3={{\bar s}}$, for any $v\in Y$, and any $x'\in \torus^d_{s'}$;
thus, $\|f(v)\|_{r,s',\vae_*}\le \|a\|_{{{\bar s}}}\le r^{-1}\sigma^{-1}\underline{\o}^{-1}\bar L$ by \equ{es1},
so that $f:Y\to Y$; notice that, in particular, this means that $f$ sends $x$--periodic functions into $x$--periodic functions. Moreover, \equ{cond1} implies also that $f$ is  a contraction: if $v_1,v_2\in Y$, then, by the mean value theorem and \equ{es1}, $|f(v_1)-f(v_2)|\le 
\|a_x\|_{{{\bar s}}}$ $|\vae|$ $|v_1-v_2|\le r^{-1}\sigma^{-1}\underline{\o}^{-1}\bar L |\vae|\ |v_1-v_2|$, so that, by taking the sup--norm, one has $\|f(v_1)-f(v_2)\|_{s'}\le \vae_* r^{-1}\sigma^{-1}\underline{\o}^{-1}\bar L \|v_1-v_2\|_{s'}<\frac{1}{6} \|v_1-v_2\|_{s'}$ showing that $f$ is a contraction. Thus, there exists a unique $\wt{\f}\in Y$ such that $f(\wt{\f})=\wt{\f}$. Furthermore, recalling  that the fixed point is achieved as the uniform limit $\lim_{n\to\io} f^n(0)$ ($0\in Y$) and since $f(0)=-a$ is analytic, so is $f^n(0)$ for any $n$ and, hence,  by Weierstrass Theorem on the uniform limit of analytic function, the limit $\wt{\f}$ itself is analytic. In conclusion, $\f\in\cB_{s',\vae_*}$ and  \equ{boundal} holds. 
\noi 
Next, \equ{boundal} and \equ{ineq1} imply \equ{phiok} and therefore, $\phi_1$ defines a symplectic diffeomorphism\footnote{Notice, in particular that the matrix $\uno_d + \vae a_x$ is, for any $x\in \torus^d_{{{\bar s}}}$, invertible  with inverse 
$\uno_d + \vae S(x;\vae)$; in fact,  since $\|\vae a_x\|_{{{\bar s}}}<\vae_* L/\mathrm{E}\le \sigma/3\le 1/6$ the matrix $\uno_d + \vae a_x$ is invertible  with inverse given by the ``Neumann series''  $(\uno_d+\vae a_x)^{-1}=\uno_d + \sum_{k=1}^\io (-1)^k (\vae a_x)^k =: \uno_d + \vae S(x;\vae) $,  so that
$\|S\|_{{{\bar s}},\vae_*}\le (\|a_x\|_{{{\bar s}},\vae_*})/(1-|\vae| \|a_x\|_{{{\bar s}},\vae_*})< \frac{6}{5}r^{-1}\sigma^{-1}\underline{\o}^{-1}\bar L$.
} satisfying \equ{phiok} and the fourth inequality in the first line of  \equ{tesit}.
\nl
It remains to show the other estimates in \equ{tesit}.
Since $L\ge \bar L$, the bound on $|\wt E|$ follows \equ{es1}. By 
\equ{es1prim}, \equ{phiok} and \equ{ineq1}, one has $\|P_1\|_{s',\vae_*}\le \|P'\|_{{{\bar s}},\vae_*}\le r^{-1}\sigma^{-1}\underline{\o}^{-1}\bar L M\le  LM/\mathrm{E}$.
Now, by Cauchy estimates, \equ{es1}, \equ{cond1} and \equ{phiok}, it follows  that
\beqano
\|\wt Q\|_{s',\vae_*}&\le & \dst\int_0^1\dst\sum_{j=1}^d \| Q_{x_j}\|_{{{\bar s}}}\|\wt{\f}_j\|_{s'}dt+\| Q'\|_{{{\bar s}}}\\
    &\le & d \mathrm{E}\frac{3}{2}\sigma^{-1}r^{-1}\sigma^{-1}\underline{\o}^{-1}\bar{L}+\| Q'\|_{{{\bar s}}}\\
    &\le & \frac{3}{2}d \bar{\mathrm{C}} \mathrm{E}^8\sigma^{-(\bar\n+2)}r^{-8}\a^{-4}\underline{\o}^{-4}M+\mathrm{C}_7 \mathrm{E}^4\sigma^{-(2\t+5)}r^{-4}\a^{-2}\underline{\o}^{-2}M\\
    &\le & \left(\frac{3}{2}d\bar{\mathrm{C}} \mathrm{E}^4\underline{\o}^{-2} +2^{-(2\t+7)}r^4\a^2\mathrm{C}_7\right)\mathrm{E}^4\sigma^{-(\bar\n+2)}r^{-8}\a^{-4}\underline{\o}^{-2}M\\
    &\le & \left(\frac{3}{2}d\bar{\mathrm{C}}\underline{\o}^{-2}  +2^{-(2\t+7)}\mathrm{C}_7\underline{\o}^{-2}\right)\mathrm{E}^8\sigma^{-(\bar\n+2)}r^{-8}\a^{-4}\underline{\o}^{-2}M\\
    &=  & \frac{\wt{\mathrm{C}}}{2d} \mathrm{E}^8\sigma^{-(\bar\n+2)}r^{-8}\a^{-4}M\le \wt L
\eeqano
and\footnote{Recall that $0<2\sigma<s$ so that $s'-\sigma=s-2\sigma>0$.}
\beqano
\|\dpr_{y'}^2\wt Q(0,\cdot)\|_{0,\vae_*}&\le & \dst\max_{1\le l\le d}\sum_{j=1}^d\|\wt Q_{y_l'y_j'}\|_{s'-\sigma,\vae_*}\\
    &\le &  2d \frac{\wt{\mathrm{C}}}{2d} \mathrm{E}^8\sigma^{-(\bar\n+2)}r^{-8}\a^{-4}M\cdot \sigma^{-2}r^{-2}\\
    &= & \wt{\mathrm{C}} \mathrm{E}^8\sigma^{-(\bar\n+2)}r^{-10}\a^{-4}M=r^{-2}\wt L.
\eeqano
so that\footnote{It is only here that a constant $L>r^{-1}\sigma^{-1}|\o|^{-1}\bar L\mathrm{E}$ is needed;
the (irrelevant) factor $6\underline{\o}^{-2}\mathrm{E}^2/5$ has been introduced for later convenience.
}
\beq{es1*}
\|\wt Q\|_{s',\vae_*}\ ,\  \frac{6}{5}\underline{\o}^{-2}\mathrm{E}^2 \, \| \dpr^2_{y'} \wt Q(0,\cdot)\|_{0,\vae_*}\le 6r^{-2}\underline{\o}^{-2}\wt L \mathrm{E}^2/5= L\ ,
\eeq
Thus, 
\beqa{Q*}
\average{\dpr^2_{y'} Q_1(0,\cdot)}
&=& \average{\dpr^2_y Q(0,\cdot)} + \vae \average{\dpr^2_{y'} \wt Q(0,\cdot)}
=T^{-1} \left(\uno_d+\vae T \average{\dpr^2_{y'}\wt Q(0,\cdot)}\right) \nonumber\\
&=:& T^{-1}(\uno_d+\vae R) \ ,
\eeqa
and, in view of \equ{C} and \equ{es1*}, we have
\beqano
\|R\| &\le & \|T\|\, \left\| \average{\dpr^2_{y'}\wt Q(0,\cdot)} \right\|\\
	  &\le & \underline{\o}^{-2} \mathrm{E}\| \dpr^2_{y'} \wt Q(0,\cdot)\|_{0,\vae_*}\le \frac{5L}{6\mathrm{E}}
\eeqano
Therefore, by \equ{cond2},
$\vae_*\|R\|\le \sigma/3\le 1/6<1$, implying that $(1+\vae R)$ is invertible and 
$$
(\uno_d+\vae R)^{-1}= \uno_d +\sum_{k=1}^\io (-1)^k \vae^k R^k=: 1+\vae D$$
with $\|D\|\le \|R\|/(1-|\vae|\, \|R\|)<L/\mathrm{E}$. In conclusion, by \equ{Q*}, and the estimate on $\|D\|$,
$$
T_1=(1+\vae R)^{-1} T = T +\vae DT=: T +\vae \wt T\ , \qquad |\o|^2
\|\wt T\|\le \|D\|\, |\o|^2
\|T\|\le \|D\|\mathrm{E}\le \frac{L}{\mathrm{E}}\, \mathrm{E}=L\ ,
$$
proving last estimate in \equ{tesit} and, hence,   Lemma~\ref{lem:1}. \qed
\nl
Next Lemma shows that, for $|\vae|$ small enough, Kolmogorov's construction can be iterated and convergence proved.

\lem{lem:2} 
Fix $0<s_*<s$  and, for $j\ge 0$,  let\footnote{Notice that $s_{j}\downarrow s_*$.}
\beqno
\left\{\begin{array}{l}
s_0:=s\\ \ \\
\dst\sigma_0:=\frac{s-s_*}2
\end{array}\right.\qquad
\left\{
\begin{array}{l}
\dst\sigma_j:=\frac{\sigma_0}{2^j}\\ \ \\
s_{j+1}:=s_j-\sigma_j=s_*+\frac{\sigma_0}{2^j}
\end{array}\right.\;.
\eeqno
 Let also   $H_0:=H$, $\mathrm{K}_0:=\mathrm{K}$, $Q_0:=Q$, $K_0:=K$, $P_0:=P$, with $\mathcal{W}$, $\wh{\mathrm{C}}$, $\wt{\mathrm{C}}$, $\mathrm{E}$, $L$ and $\n$  as in $\S\ref{SectKolmoNorm}$ and 
assume that $\vae_*\le \vae_0$ satisfies
\beq{cond}
\vae_* \, \mathrm{e}_*\,{\,d_*}\, \|P\|_{s,\vae_0}\le 1\quad {\rm where}\quad 
\mathrm{e}_*:=  3\wh{\mathrm{C}}\,\sigma_0^{-(\n+1)}\,   \mathrm{E}^{9}r^{-10}\a^{-4}\underline{\o}^{-6}\ ,\ {\,d_*}:=2^{\n+1}\ .
\eeq
Then,
one can construct a sequence of symplectic transformations 
\beq{phij}
\phi_j:W_{rs_j,s_{j},\vae_*}\to D_{rs_{j-1}}\times\torus^d_{s_{j-1}}\ ,
\eeq
so that
\beq{Hj}
H_j:=H_{j-1}\circ\phi_j=: K_j + \vae^{2^j} P_j\ ,
\eeq
converges uniformly to a Kolmogorov's normal form. More precisely,
$\vae^{2^j} P_j$, 
$\Phi_j:=\phi_1\circ\phi_2\circ \cdots\circ \phi_j$, 
$\mathrm{K}_j$, $K_j$, $Q_j$ converge uniformly on 
$W_{s_*,\vae_*}$ to, respectively, $0$, $\phi_*$, $\mathrm{K}_*$, $K_*$, $Q_*$, which are real--analytic on $W_{s_*,\vae_*}$ and  $H\circ\phi_*=K_*=\mathrm{K}_*+\o\cdot y+ Q_*$ with $Q_*=O(|y|^2)$. 
Finally, the following estimates hold for any $|\vae|< \vae_*$ and for any $i\ge 0$:
\begin{align}
& |\vae|^{2^i}M_i:=|\vae|^{2^i}\|P_i\|_{s_i,\vae_*}\le \frac{(\, |\vae| \mathrm{e}_*{\,d_*}M)^{2^i}}{\mathrm{e}_* {\,d_*}^{i+1}}\ ,\label{estfin1}\\
&
\frac{{\mathrm{C}} \mathrm{E}^3}{ r^3\underline{\o}^{3}\sigma^2}\|\mathcal{W}(\phi_*-\id)\|_{s_*},\ |\mathrm{K}-\mathrm{K}_*|,\ \|Q-Q_*\|_{s_*},\  |\o|^2\|T-T_*\| \le \frac{|\vae|L}{3 \sigma}\ ,\label{estfin2}
\end{align}
where $T_*:=\average{\dpr_y^2 Q_*(0,\cdot)}^{-1}$. 
\elem
\proof
Notice that 
 \equ{cond} implies \equ{cond2} (and, hence, \equ{cond1}). For $i\ge 0$, define
 $$\mathcal{W}_i\coloneqq \diag(\underline{\o}\,\uno_d,r\sigma_i\underline{\o}\,\uno_d)\quad\mbox{and}\quad \bar L_i\coloneqq \bar{\mathrm{C}} \mathrm{E}^7\sigma_i^{-\bar\n}r^{-7}\a^{-4}\underline{\o}^{-3}M_i\, .$$
Let us assume ({\sl inductive hypothesis}) that we can iterate $j\ge 1$ times Kolmogorov transformation obtaining $j$ symplectic transformations $\phi_{i+1}:W_{rs_{i+1},s_{i+1},\vae_*}\to D_{rs_i}\times\torus^d_{s_i}$, for $0\le i\le j-1$, and $j$  Hamiltonians $H_{i+1}=H_i\circ\phi_{i+1}=K_{i+1}+\vae^{2^{i+1}} P_{i+1}$ real--analytic on $W_{s_{i+1},\vae_*}$ such that, for any $0\le i\le j-1$,
\beq{bbb}
\left\{
\begin{array}{l}
r|\o|,\|Q_i\|_{s_i},|\o|^2\|T_i\|\le \mathrm{E}\ \\  \ \\ 
|\vae|^{2^i} L_i:=|\vae|^{2^i}\wh{\mathrm{C}} \mathrm{E}^{10} \sigma_0^{-\n}\, 2^{\n i}r^{-10}\a^{-4} M_i  \le \mathrm{E}\frac{\sigma_i}{3}\ . 
\end{array}\right.
\eeq
Observe that for $j=1$, it is $i=0$ and \equ{bbb} is implied by the definition of $\mathrm{E}$  and by condition \equ{cond}.

\noi
Because of \equ{bbb}, \equ{cond2} holds 
for $H_i$ and Lemma~\ref{lem:1}  can be applied to $H_i$ and  one has, for $0\le i\le j-1$ and for any $|\vae|< \vae_*$ (compare \equ{tesit}): 
\beqa{C.1}
&& |\mathrm{K}_{i+1}|\le |\mathrm{K}_i|+|\vae|^{2^i} L_i\ ,\quad \|Q_{i+1}\|_{s_{i+1}}\le \|Q_i\|_{s_i}+|\vae|^{2^i} L_i\ ,\quad |\o|^2\|T_{i+1}\|\le |\o|^2\|T_i\|+|\vae|^{2^i} L_i\ ,
\nonumber\\
&&\|\mathcal{W}_i(\phi_{i+1}-\id)\|_{s_{i+1}}\le |\vae|^{2^i}\bar L_i\ ,\quad M_{i+1}\le M_iL_i\mathrm{E}^{-1}\ .
\eeqa 
Observe that, by definition of $\mathrm{e}_*$, ${\,d_*}$ in \equ{cond} and of $L_i$ in \equ{bbb}, one has $|\vae|^{2^j} L_j (3 \sigma_j^{-1}\mathrm{E}^{-1})=\mathrm{e}_*{\,d_*}^j |\vae|^{2^j} M_j=:\theta_j/d_*$, so that $L_i\mathrm{E}^{-1}<\mathrm{e}_* {\,d_*}^i M_i$, thus by last relation in \equ{C.1}, for any $0\le i\le j-1$,  $|\vae|^{2^{i+1}}M_{i+1}<\mathrm{e}_* {\,d_*}^i (M_i|\vae|^{2^i})^2$ $.i.e.$ $\theta_{i+1}<\theta_i^2$, which iterated, yields \equ{estfin1}  $.i.e.$ $\theta_i\le \theta_0^{2^i}$ for $0\le i\le j$. 

\nl
Next, we  show that, thanks to \equ{cond}, \equ{bbb} holds also for $i=j$. In fact, by \equ{bbb} and the definition of $\mathrm{E}$ in $\S\ref{SectKolmoNorm}$, we have 
$$\|Q_i\|_{s_j}\le \|Q\|_{s}+\sum_{i=0}^{j-1} \vae_*^{2^i} L_i\le \|Q\|_{s}+ \su{3}\mathrm{E} \sum_{i\ge 0}\sigma_i<\|Q\|_{s}+\su{3}\mathrm{E}\sum_{i\ge 0} 2^{-(i+1)}=\|Q\|_{s}+\su{3}\mathrm{E}< \mathrm{E}\ .$$
The  bound for  $\|T_i\|$ is proven in an identical manner. 
Now, by $\equ{estfin1}_{i=j}$ and \equ{cond}, 
$$\theta_j/d_*=|\vae|^{2^j} L_j (3 \sigma_j^{-1}\mathrm{E}^{-1})=\mathrm{e}_*{\,d_*}^j |\vae|^{2^j} M_j\le \mathrm{e}_*{\,d_*}^j (\mathrm{e}_*{\,d_*}\vae_* M)^{2^j}/(\mathrm{e}_*{\,d_*}^{j+1})\le 1/{\,d_*}<1\ ,$$
which implies the second inequality in \equ{bbb} with $i=j$; the proof of the induction is finished and one can construct an {\sl infinite sequence} of Kolmogorov transformations satisfying \equ{bbb}, \equ{C.1} and \equ{estfin1} {\sl for all $i\ge 0$}. 

\nl
To check \equ{estfin2}, we observe that 
$$|\vae|^{2^i} L_i\mathrm{E}^{-1}=\frac{\sigma_0}{3 \cdot 2^i}\ \mathrm{e}_*{\,d_*}^i |\vae|^{2^i} M_i < \su{2^{i+1}d_*}(|\vae| \mathrm{e}_*{\,d_*}M)^{2^i}\le \su{d_*} \Big(\frac{|\vae|\mathrm{e}_*{\,d_*}M}2\Big)^{i+1}$$
and therefore 
$$\sum_{i\ge 0}  |\vae|^{2^i} L_i\le \frac{\mathrm{E}}{d_*}\sum_{i\ge 1} \Big(\frac{|\vae|\mathrm{e}_*{\,d_*}M}2\Big)^{i}\le |\vae|\mathrm{e}_*\mathrm{E}M =\frac{|\vae|L}{3 \sigma_0}\ .$$ 
Thus,  
$$\|Q-Q_*\|_{s_*}\le \sum_{i\ge 0}|\vae|^{2^i}\|\tilde Q_i\|_{s_i}\le \sum_{i\ge 0}|\vae|^{2^i} L_i\le\frac{|\vae|L}{3 \sigma_0}\ ;$$
and analogously for $|\mathrm{K}-\mathrm{K}|_*$ and $\|T-T_*\|$. 

\nl
Next, we prove that $\Phi_j$ is convergent by proving that it is Cauchy. For any $j\geq 1$, we have\footnote{Notice that $\Phi_0=\id$, for any $j\ge 0,\, L_j\ge r^{-1}\sigma_{j}^{-2}\underline{\o}^{-1}\bar L_{j}\mathrm{E}$ and, by \eqref{es1}, \eqref{cond1},  \eqref{boundal} and \eqref{cond}, we have, $$\wt\phi_j\coloneqq \phi_j-\id\colon W_{s_*,\vae_*}\to W_{s_*+\sigma_j/3,\vae_*}\quad \mbox{and}\quad s_*+\sigma_j/3= s_{j-1}-\frac{11}{6}\sigma_{j-1}.$$ }
\beqano
\|\mathcal{W}_0(\Phi_j-\Phi_{j-1})\|_{s_*,\vae_*}&=&\|\mathcal{W}_0\Phi_{j-1}\circ\phi_j-\mathcal{W}_0\Phi_{j-1}\|_{s_{*},\vae_*}\\
           &\leq &\|\mathcal{W}_0d\Phi_{j-1}\mathcal{W}_j^{-1}\|_{\vae_*\bar L_{j-1}}\, \|\mathcal{W}_j(\phi_j-\id)\|_{s_*,\vae_*}\\
           &\leq &\|\mathcal{W}_0d\Phi_{j-1}\|_{s_{*}+\sigma_j/3,\vae_*}\, \|\mathcal{W}_j(\phi_j-\id)\|_{s_*,\vae_*}\\
           &\le & \|\mathcal{W}_0d\Phi_{j-1}\|_{s_{j-1}-\frac{11}{6}\sigma_{j-1},\vae_*}\, \vae^{2^{j-1}}\bar L_{j-1}\\
           &\leq &\frac{6}{11}\|\mathcal{W}_0\Phi_{j-1}\|_{s_{j-1},\vae_*}\sigma_{j-1}^{-1}\max\left(r^{-1}\underline{\o}^{-1}, r^{-1}\sigma_{j-1}^{-1}\underline{\o}^{-1}   \right)\vae^{2^{j-1}}\bar L_{j-1}\\
           &\le &\frac{6}{11}\|\mathcal{W}_0\Phi_{0}\|_{s_{0},\vae_*}\cdot\vae^{2^{j-1}}\cdot r^{-1}\sigma_{j-1}^{-2}\underline{\o}^{-1}\bar L_{j-1}\\
           &\leq &\frac{6}{11}\max\left(rs_0\underline{\o},r\sigma_{0}\underline{\o}s_0\right)\vae^{2^{j-1}}L_{j-1}\mathrm{E}^{-1}\\
           &\leq & \frac{6}{11} \cdot r s_0\underline{\o} \cdot \su3 \sigma_{j-1}.
\eeqano
Therefore, for any $n\geq 0,\, j\geq 1$,
$$\|\Phi_{n+j}-\Phi_n\|_{s_{n+j},\vae_*}\leq  \sum_{i=n}^{n+j}\|\Phi_{i+1}-\Phi_i\|_{s_i,\vae_*}\leq \frac{2}{11} r s_0\underline{\o}\sum_{i=n}^{n+j}\sigma_{i}
.$$
Hence $\Phi_j$ converges uniformly on $W_{s_*,\vae_*}$ to some $\phi_*$, which is then real--analytic function on $W_{s_*,\vae_*}$.

\nl
To estimate $\|\mathcal{W}_0(\phi_*-\id)\|_{s_*}$, observe that  
$$\|\mathcal{W}_0(\Phi_i-\id)\|_{s_i}\le \|\mathcal{W}_0(\Phi_{i-1}\circ\phi_i-\phi_i)\|_{s_i}+\|\mathcal{W}_0(\phi_i-\id)\|_{s_i}\le \|\mathcal{W}_0(\Phi_{i-1}-\id)\|_{s_{i-1}}+|\vae|^{2^i}\bar L_i\ ,$$
which iterated yields 
\begin{align*}
\|\mathcal{W}_0(\Phi_i-\id)\|_{s_i}&\le \sum_{k=0}^i |\vae|^{2^k}\bar L_k=\frac{\bar{\mathrm{C}}}{\wh{\mathrm{C}}}\mathrm{E}^{-3}r^3\underline{\o}^{3} \sum_{k=0}^i|\vae|^{2^k} L_k\sigma_k^2\le \frac{\bar{\mathrm{C}}}{\wh{\mathrm{C}}}\mathrm{E}^{-3}r^3\underline{\o}^{3} |\vae|\mathrm{e}_*M \sigma_0^2\\
     &=\frac{\bar{\mathrm{C}}}{\wh{\mathrm{C}}}\mathrm{E}^{-3}r^3\underline{\o}^{3} |\vae|L\sigma_0=\frac{r^3\underline{\o}^{3}}{3\mathrm{C}\mathrm{E}^3}|\vae|L\sigma.
\end{align*}
Therefore, taking the limit over $i$ completes the proof of \equ{estfin2}, Lemma~\ref{lem:2} and, whence, of Kolmogorov's Theorem.  \qed
%
\subsection{Proof of Theorem~\ref{teo4}}
\lemtwo{lem:1bis}{KAM step}
Let $r>0,\,0<2\s<s\leq 1$ and consider the hamiltonian parametrized by $\vae\in\real$
$$
H(y,x;\vae)\coloneqq K(y)+\vae P(y,x),
$$
with 
$$
K,P\in \mathcal{B}_{r,s}(\mathsf{y})\,.
$$
Assume that\footnote{In the sequel, $K$ and $P$  stand for  generic real analytic hamiltonians which, later on, will respectively play the roles of $K_j$ and $P_j$,  and $\mathsf{y},\,r$, the roles of $y_j,\,r_j$ in the iterative step.}\textsuperscript{,}\footnote{Notice that $\mathsf{T}\mathsf{K}\ge \mathsf{T}\|K_{yy}(\mathsf{y})\|\ge \|T\|\|K_{yy}(\mathsf{y})\|=\|T\|\|T^{-1}\|\ge 1 $.\label{ftnTK1}}
\beq{RecHypArn}
\begin{aligned}
&\det K_{yy}(\mathsf{y})\neq 0\;, \qquad\qquad\;\, T\coloneqq K_{yy}(\mathsf{y})^{-1}\;,\\
&\|K_{yy}\|_{r,\mathsf{y}}\le \mathsf{K}\;,\qquad\qquad\quad\ \; \|T\|\le \mathsf{T}\;,\\
& \|P\|_{r,s,\mathsf{y}}\le M \;,\qquad\quad\qquad\,\,\  \o\coloneqq K_{yy}(\mathsf{y})\in \D^\t_\a\;. 
\end{aligned}
\eeq
Fix $\vae\not=0$ and let
\beq{DefNArn}
\begin{aligned}
&\l\ge \frac{4}{5}\log\left(\s^{2\n+d}\frac{\a^2}{|\vae|MK}\right),\quad\k\coloneqq 5\s^{-1}\l, \quad 
\bar{r}\le  \dst\min\left\{\frac{\a}{2d\mathsf{K}\k^{\n}}\,,\, \frac{5}{24d}\frac{r}{\mathsf{T}\mathsf{K}} \right\},\\
& \ \bar{s}\coloneqq s-\frac{2}{3}\s,\quad s'\coloneqq s-\s \,,
\end{aligned}
\eeq
Finally, define\footnote{Notice that $\mathsf{L}\ge \s^{-d}\ovl{\mathsf{L}}\ge \ovl{\mathsf{L}}$ since $\s\le 1$.
}
\begin{align*}
\ovl{\mathsf{L}}&\coloneqq\frac{\mathsf{C}_0}{\sqrt{2}} \max\left\{1,\frac{\a}{r\mathsf{K}}\right\}\frac{M \mathsf{K}}{\a^2 }\s^{-(2\n+d)}\;,\\
\mathsf{L}&\coloneqq M\dst\max\left\{\frac{8\mathsf{T}  }{r\bar{r}}\s^{-(\n+d)}\,,\, \frac{\mathsf{C}_7}{\sqrt{2}} \max\left\{1,\frac{\a}{r\mathsf{K}}\right\}\frac{\mathsf{K}}{\a^2}\s^{-2(\n+d)}\right\}\\
		  &= M\dst\max\left\{\frac{8\mathsf{T}  }{r\bar{r}}\s^{-(\n+d)}\,,\,\frac{4}{\mathsf{K} r^2}\,,\,\frac{\mathsf{C}_7}{\sqrt{2}} \max\left\{1,\frac{\a}{r\mathsf{K}}\right\}\frac{\mathsf{K}}{\a^2}\s^{-2(\n+d)}\right\}
\;.
\end{align*}
Then, there exists a generating function $g\in \mathcal{B}_{\bar r,\bar s}(\mathsf{y})$
 with the following properties:
\beq{Est1Lem1b}
\left\{
\begin{aligned}
&\|g_x\|_{\bar{r},\bar{s},\mathsf{y}}\le  \mathsf{C}_1 \frac{M}{\a} \s^{-(\n+d)}\,,\\
& \|g_{y'}\|_{\bar{r},\bar{s},\mathsf{y}},\, \|\dpr_{y'x}^2 g\|_{\bar{r},\bar{s},\mathsf{y}}\le \ovl{\mathsf{L}}\,,\\
&\|\dpr_{y'}^2\wt K\|_{\bar{r},\mathsf{y}}\le 
\mathsf{K}\mathsf{L}\,,
\end{aligned}
\right.
\eeq
where 
$$\wt K(y')\coloneqq \average{P(y',\cdot)}\;.$$
 If, in addition,  
\beq{cond1Bis}
|\vae|\leq \vae_\sharp\quad \mbox{and}\quad {|\vae| }{\mathsf{L}}\le \frac{\sigma}{3}
\ ,
\eeq
then, there exists $\mathsf{y}'\in\rn$ such that 
\beq{convEst}
\left\{
\begin{aligned}
&\dpr_{y'} K'(\mathsf{y}')=\o \,,\qquad \qquad\quad\qquad\quad \det\dpr^2_{y'} K'(\mathsf{y}')\neq 0\,,\\
&|\vae|\|g_x\|_{\bar{r},\bar{s},\mathsf{y}}\le \frac{r}{3}\,,\qquad \qquad\quad\qquad \ 
|\mathsf{y}'-\mathsf{y}| 
\le \frac{8|\vae|\mathsf{T} M}{r}\,,\\
&|\vae|\|\wt T\|\le \mathsf{T}|\vae|\mathsf{L}
\,, \qquad\qquad\quad\quad\quad\ \;\|P_+\|_{\bar{r},\bar s,\mathsf{y}} \le  \mathsf{L}M\,,
\end{aligned}
\right. 
\eeq
where 
$$
K'\coloneqq K+\vae\wt K\;,\qquad \left(\dpr^2_{y'} K'(\mathsf{y}')\right)^{-1}\eqqcolon T+\vae\;\wt T\;,\qquad P_+(y',x)\coloneqq P(y'+\vae g_x(y',x),x)\;.
$$  
and the following hold. For $y'\in D_{\bar{r}}(\mathsf{y})$,    the map $\psi_\vae(x):=  x+\vae g_{y'}(y',x)$ has an analytic inverse 
 $\f(x')=x'+\vae \wt{\f}(y',x';\vae)$  such that
\beq{boundalBis}
\|\wt{\f}\|_{\bar{r}, s',\mathsf{y}}\le   \ovl{\mathsf{L}} \qquad {\rm and}
\quad
\f=\id + \vae \wt{\f} : D_{\bar{r}/2,s'}(\mathsf{y}')\to \torus^d_{\bar{s}} \ ;
\eeq
for any $(y',x)\in D_{\bar{r},\bar s}(\mathsf{y})$, $|y'+\vae g_x(y',x)-\mathsf{y}|<\frac{2}{3} r$; the map $\phi'$ is a symplectic diffeomorphism and
\beq{phiokBis0}
\phi'=\big( y'+\vae g_x(y', \f(y',x')),\f(y',x')\big): D_{\bar{r}/2,s'}(\mathsf{y}')\to D_{2r/3, \bar{s}}(\mathsf{y}),
\eeq
with
\beq{phiokBis1}
\|\mathsf{W}\,\tilde \phi\|_{\bar{r}/2,s',\mathsf{y}'}\le \s^d{\mathsf{L}}\,,
\eeq
where $\tilde \phi$ is defined by the relation $\phi'=:\id + \vae \tilde \phi$,
$$
\mathsf{W}\coloneqq \begin{pmatrix}
\max\{\frac{\mathsf{K}}{{\a}}\;,\frac{1}r\}\;\uno_d & 0\\ \ \\
0			& \uno_d 
\end{pmatrix}
$$
and
\beq{tesitBis}
\|P'\|_{\bar{r}/2, s',\mathsf{y}'}\le  \mathsf{L}M\;,
\eeq
with
$$
P'(y',x')\coloneqq P_+(y',\f(x'))=P\circ \phi'(y',x')\;.
$$
\elem
\proof\\
\Giu
{\bf Step 1: Construction of the Arnold's transformation }  We seek for a near--to--the--identity symplectic transformation 
\[\phi'\colon D_{r_1,s_1}(\mathsf{y}')\to D_{r,s}(\mathsf{y}),\]
with $D_{r_1,s_1}(\mathsf{y}')\subset D_{r,s}(\mathsf{y})$,   generated by a function of the form $y'\cdot x+\vae g(y',x)$, so that
\beq{ArnTraKam}
\phi'\colon \left\{\begin{aligned}
y  &=y'+\vae g_x(y',x)\\
x' &=x+\vae g_{y'}(y',x)\, ,
\end{aligned}
\right.
\eeq
such that
\beq{ArnH1}
\left\{
\begin{aligned}
& H':= H\circ \phi'=K'+\vae^2 P'\ ,\\
& \dpr_{y'} K'(\mathsf{y}')=\o,\quad \det \dpr^2_{y'} K'(\mathsf{y}')\neq 0\,.
\end{aligned}
\right.
\eeq
By Taylor's formula, we get\footnote{Recall that $\average{\cdot}$ stands for the average over $\tn$.}
\beq{Arneq11}
\begin{aligned}
H(y'+\vae g_x(y',x),x)=&K(y')+\vae \wt K(y') +\vae \left[K_y(y')\cdot g_x +T_{\k} P(y',\cdot)-\wt K(y') \right]+\\
						&+\vae^2 \left( P^\ppu+P^\ppd+ P^\ppt\right)(y',x) \\
			= & K'(y')+\vae \left[K_y(y')\cdot g_x +T_{\k} P(y',\cdot)-\wt K(y') \right]+ \vae^2 P_+(y',x),
\end{aligned}
\eeq
with $\k\in\natural$, which will be chosen large enough so that $P^\ppt=O(\vae)$ 
and 
\beq{ArnDefPs}
\left\{
\begin{aligned}
P_+&\coloneqq P^\ppu+P^\ppd+ P^\ppt\\
P^\ppu &\coloneqq \su{\vae^2}\left[K(y'+\vae g_x)-K(y')-\vae K_y(y')\cdot g_x \right]=\dst\int^1_0(1-t)K_{yy}(\vae t g_x)\cdot g_x\cdot g_x dt\\
P^\ppd &\coloneqq \su\vae \left[P(y'+\vae g_x,x)-P(y',x)\right]=\dst\int_0^1P_y(y'+\vae t g_x,x)\cdot g_x dt\\
P^\ppt &\coloneqq \su\vae \left[ P(y',x)-T_{\k} P(y',\cdot)\right]=\su\vae \dst\sum_{|n|_1>\k} P_n(y')\ex^{in\cdot x}\; .
\end{aligned}
\right.
\eeq
By the non--degeneracy condition in \eqref{RecHypArn}, for $\vae$ small enough (to be made precised below), $\det\dpr_{y'}^2 K'(\mathsf{y})\neq0$
 and, therefore, by Lemma~\ref{IFTLem}, there exists a unique $\mathsf{y}'\in D_r(\mathsf{y})$ such that the second part of \eqref{ArnH1} holds. 
In view of \eqref{Arneq11}, in order to get the first part of \eqref{ArnH1}, we need to find $g$ such that  $K_y(y')\cdot g_x +T_{\k} P(y',\cdot)-\wt K(y')$ vanishes; such a $g$ is indeed given by
 \beq{HomEqArn}
 g\coloneqq \dst\sum_{0<|n|_1\leq \k} \frac{-P_n(y')}{iK_y(y')\cdot n}\ex^{in\cdot x},
 \eeq
provided that 
\beq{CondHomEqArn}
K_y(y')\cdot n\neq 0, \quad \forall\; 0<|n|_1\leq \k,\quad \forall\; y'\in D_{r_1}(\mathsf{y}')\quad  \left(\subset D_{r}(\mathsf{y})\right).
\eeq
But, in fact, since $K_y(\mathsf{y})$ is rationally independent, then, given any $\k\in\natural$, there exists $\bar{r}\leq r$ such that
\beq{CondHomEqArnBis}
K_y(y')\cdot n\neq0,\quad \forall\; 0<|n|_1\leq \k, \quad\forall\; y'\in D_{\bar{r}}(\mathsf{y}).
\eeq
The last step is to invert the function $x\mapsto x+\vae g_{y'}(y',x)$ in order to define $P'$. But, by Lemma~\ref{IFTLem}, for $\vae$ small enough, the map $x\mapsto x+\vae g_{y'}(y',x)$ admits an real--analytic inverse of the form
\beq{InvComp2Fi}
\f(y',x';\vae)\coloneqq x'+\vae \wt{\f}(y',x';\vae),
\eeq
so that the Arnod's symplectic transformation is given by
\beq{ArnTrans0}
\phi'\colon (y',x')\mapsto \left\{
\begin{aligned}
y &= y'+\vae g_x(y',\f(y',x'))\\
x &= \f(y',x';\vae)= x'+\vae \wt{\f}(y',x';\vae) .
\end{aligned}
\right.
\eeq
Hence, \eqref{ArnH1} holds with
\beq{DefP1Ar}
P'(y',x')\coloneqq P_+(y', \f(y',x')).
\eeq
{\bf Step 2: Quantitative estimates}\\
First of all, notice that\footnote{Recall footnote \textsuperscript{\ref{ftnTK1}}.}
\beq{rrbarAs}
\bar{r}\le \frac{5r}{24d\mathsf{T}\mathsf{K}}<\frac{r}{2}\;.
\eeq
\noi
We begin by extending the ``diophantine condition w.r.t. $K_y$'' uniformly to $D_{\bar{r}}(\mathsf{y})$ up to the order $\k$. Indeed, by the Mean Value Inequality and $K_y(\mathsf{y})=\o\in\D^\t_\a$,  we get, for any $0<|n|_1\leq \k$ and any $y'\in D_{\bar{r}}(\mathsf{y})$,
\begin{align}
|K_y(y')\cdot n|&=|\o\cdot n +(K_y(y')-K_y(\mathsf{y}))\cdot n|\geq |\o\cdot n|\left(1-d\frac{\|K_{yy}\|_{\bar{r},\mathsf{y}}}{|\o\cdot n|}|n|_1\bar{r}\right) \nonumber\\
         &\geq \frac{\a}{|n|_1^\t}\left(1-\frac{d\mathsf{K}}{\a }|n|_1^{\t+1}\bar{r} \right)\geq \frac{\a}{|n|_1^\t}\left(1-\frac{d\mathsf{K}}{\a }\k^{\t+1}\bar{r} \right)\ge \frac{\a}{2|n|_1^\t},\label{ArnExtDiopCond}
\end{align}
so that, by Lemma~\ref{fce}--$(i)$, we have
\begin{align*}
\|g_x\|_{\bar{r},\bar{s},\mathsf{y}} &\overset{def}{=}\dst\sup_{D_{\bar{r},\bar{s}}(\mathsf{y})}\left|\dst\sum_{0<|n|_1\leq \k}\frac{nP_n(y')}{K_y(y')\cdot n}\ex^{in\cdot x} \right|\leq \dst\sum_{0<|n|_1\leq \k}\frac{\|P_n\|_{\bar{r},\bar{s}, \mathsf{y}}}{|K_y(y')\cdot n|}|n|_1\ex^{\left(s-\frac{2}{3}\s\right)|n|_1}\\
   &\leq \dst\sum_{0<|n|_1\leq \k} M\ex^{-s|n|_1}\frac{2|n|_1^{\n}}{\a}\ex^{\left(s-\frac{2}{3}\s\right)|n|_1}\leq \frac{2M}{\a}\dst\sum_{n\in\zn} |n|_1^{\n}\ex^{-\frac{2}{3}\s|n|_1}\\
   &\leq \frac{2M}{\a}\dst\int_{\rn} |y|_1^{\n}\ex^{-\frac{2}{3}\s|y|_1}dy\\
   &= \left(\frac{3}{2\s}\right)^{\n+d}\frac{2M}{\a}\dst\int_{\rn} |y|_1^{\n}\ex^{-|y|_1}dy\\
   &= \mathsf{C}_1 \frac{M}{\a} \s^{-(\n+d)}\,,
\end{align*}
\begin{align*}
\|\dpr_{y'}g\|_{\bar{r},\bar{s},\mathsf{y}} &\overset{def}{=}\dst\sup_{D_{\bar{r},\bar{s}}(\mathsf{y})}\left|\dst\sum_{0<|n|_1\leq \k}\left(\frac{ \dpr_yP_n(y')}{K_y(y')\cdot n}-P_n(y')\frac{ K_{yy}(y')n}{(K_y(y')\cdot n)^2}\right)\ex^{in\cdot x} \right|\\
   &\leq \dst\sum_{0<|n|_1\leq \k}\dst\sup_{D_{\bar{r}}(\mathsf{y})}\left(\frac{\|(P_y)_n\|_{\bar{r},s, \mathsf{y}}}{|K_y(y')\cdot n|}+d\|P_n\|_{r,s, \mathsf{y}}\frac{\|K_{yy}\|_{r,\mathsf{y}}|n|_1}{|K_y(y')\cdot n|^2}\right)\ex^{\left(s-\frac{2}{3}\s\right)|n|_1}\\
   &\stackrel{\equ{RecHypArn}+\equ{ArnExtDiopCond}}{\le} \dst\sum_{0<|n|_1\leq \k}\left( \frac{M}{r-\bar{r}}\ex^{-s|n|_1}\frac{2|n|_1^{\t}}{\a}+dM\ex^{-s|n|_1}\mathsf{K}|n|_1\left(\frac{2|n|_1^{\t}}{\a}\right)^2\right)\ex^{\left(s-\frac{2}{3}\s\right)|n|_1}\\
   &\leby{rrbarAs} \frac{4M}{\a^2 r}\dst\sum_{0<|n|_1\leq \k}\left( |n|_1^{\t}\a +dr\mathsf{K}|n|_1^{2\t+1}\right)\ex^{-\frac{2}{3}\s|n|_1}\\
   &\le \max\left\{\a,r\mathsf{K}\right\}\frac{4M}{\a^2 r}\dst\sum_{0<|n|_1\leq \k}\left( |n|_1^{\t}+d|n|_1^{2\t+1}\right)\ex^{-\frac{2}{3}\s|n|_1}\\
   &\leq \max\left\{1,\frac{\a}{r\mathsf{K}}\right\}\frac{4M \mathsf{K}}{\a^2 }\dst\int_{\rn} \left( |y|_1^{\t}+d|y|_1^{2\t+1}\right)\ex^{-\frac{2}{3}\s|y|_1}dy \\
   &= \left(\frac{3}{2\s}\right)^{2\t+d+1}\max\left\{1,\frac{\a}{r\mathsf{K}}\right\}\frac{4M \mathsf{K}}{\a^2 }\dst\int_{\rn} \left( |y|_1^{\t}+d|y|_1^{2\t+1}\right)\ex^{-|y|_1}dy\\
   &\le \frac{\mathsf{C}_0}{\sqrt{2}} \max\left\{1,\frac{\a}{r\mathsf{K}}\right\}\frac{M \mathsf{K}}{\a^2 }\s^{-(2\t+d+1)}\\
   &\le \ovl{\mathsf{L}} \;,
\end{align*}
and, analogously,
\begin{align*}
\|\dpr^2_{y'x}g\|_{\bar{r},\bar{s},\mathsf{y}} &\overset{def}{=}\dst\sup_{D_{\bar{r},\bar{s}}(\mathsf{y})}\left|\dst\sum_{0<|n|_1\leq \k}\left(\frac{ \dpr_yP_n(y')}{K_y(y')\cdot n}-P_n(y')\frac{ K_{yy}(y')n}{(K_y(y')\cdot n)^2}\right)\cdot n\ex^{in\cdot x} \right|\\
   &\leq \dst\sum_{0<|n|_1\leq \k}\dst\sup_{D_{\bar{r}}(\mathsf{y})}\left(\frac{\|(P_y)_n\|_{\bar{r},s, \mathsf{y}}}{|K_y(y')\cdot n|}+d\|P_n\|_{r,s, \mathsf{y}}\frac{\|K_{yy}\|_{r,\mathsf{y}}|n|_1}{|K_y(y')\cdot n|^2}\right)|n|_1\ex^{\left(s-\frac{2}{3}\s\right)|n|_1}\\
   &\le \max\{\a,r\mathsf{K}\}\frac{4M}{\a^2 r}\dst\sum_{0<|n|_1\leq \k}\left( |n|_1^{\t}+d|n|_1^{2\t+1}\right)|n|_1\ex^{-\frac{2}{3}\s|n|_1}\\
   &\leq \max\left\{1,\frac{\a}{r\mathsf{K}}\right\}\frac{4M \mathsf{K}}{\a^2 }\dst\int_{\rn} \left( |y|_1^{\t}+d|y|_1^{2\t+1}\right)|y|_1\ex^{-\frac{2}{3}\s|y|_1}dy \\
   &= \left(\frac{3}{2\s}\right)^{2\t+d+2}\max\left\{1,\frac{\a}{r\mathsf{K}}\right\}\frac{4M \mathsf{K}}{\a^2 }\dst\int_{\rn} \left( |y|_1^{\t+1}+d|y|_1^{2\t+2}\right)\ex^{-|y|_1}dy\\
   &= \frac{\mathsf{C}_0}{\sqrt{2}} \max\left\{1,\frac{\a}{r\mathsf{K}}\right\}\frac{M \mathsf{K}}{\a^2 }\s^{-(2\n+d)} \\
   &=\ovl{\mathsf{L}}\;,
\end{align*}
and, for $|\vae|< {\vae_*}$,
\[\|\wt K_y\|_{r/2,\mathsf{y}}=\| \left[P_y\right]\|_{r/2,\mathsf{y}}\leq \|P_y\|_{r/2,\bar{s}, \mathsf{y}}\leq  \frac{M}{r-\frac{r}{2}}\leq \frac{2M}{r} \;,\]
\[\|\dpr_{y'}^2\wt K\|_{r/2,\mathsf{y}}=\| \left[P_{yy}\right]\|_{r/2,\mathsf{y}}\leq \|P_{yy}\|_{r/2,\bar{s}, \mathsf{y}}\leq  \frac{M}{(r-\frac{r}{2})^2}\leq \frac{4M}{r^2}\le \mathsf{K}\mathsf{L} 
\;.\]
Next, we prove the existence and uniqueness of $\mathsf{y}'$ in \eqref{ArnH1}. Consider then
\begin{align*}
F\colon D_{\bar{r}}(\mathsf{y})\times D^1_{2|\vae|}(0) &\longrightarrow \qquad \cn\\
		(y,\eta)\quad &\longmapsto K_y(y)+\eta \wt K_{y'}(y)-K_y(\mathsf{y}).
\end{align*}
Then
\begin{itemize}
\item $F(\mathsf{y},0)=0,\quad F_y(\mathsf{y},0)^{-1}=K_{yy}(\mathsf{y})^{-1}=T$;
\item For any $(y,\eta)\in D_{\bar{r}}(\mathsf{y})\times D^1_{2|\vae|}(0)$,
\begin{align*}
\|\uno_d-TF_y(y,\eta)\|&\leq \|\uno_d-TK_{yy}\|+|\eta|\;\|T\|\;\|\dpr_{y'}^2\wt K\|_{r/2,\mathsf{y}}\\
	  &\leq d\|T\|\|K_{yyy}\|_{\bar{r},\mathsf{y}}\bar{r}+ 2|\vae|\mathsf{T}\frac{4M}{r^2}\\
      &\leq d\mathsf{T} \mathsf{K}\frac{\bar{r}}{r-\bar{r}}+8\mathsf{T}\frac{|\vae| M}{ r^2}\\
      &\leby{rrbarAs}d\mathsf{T} \mathsf{K} \frac{2\bar{r}}{ r}+|\vae|\frac{8\mathsf{T} M}{r^2}\\
      &\le 2d\mathsf{T} \mathsf{K}\frac{\bar{r}}{ r}+\su2{|\vae|}\mathsf{L}\\
      &\overset{\equ{rrbarAs}+\equ{cond1Bis}}{\leq}\frac{5}{12}+\frac{\s}{6}\\
      &\le \frac{5}{12}+\su{12}=\su2\;;
\end{align*}
\item Recalling $\s\le\su2$, we have
\begin{align}
2\|T\|\|F(\mathsf{y},\cdot)\|_{2|\vae|,0}&=2\|T\|\dst\sup_{B^1_{2|\vae|}(0)}|\eta \wt K_{y'}(\mathsf{y})|\nonumber\\
		&\leq 2\mathsf{T} \frac{4|\vae| M}{r}\nonumber\\
		&< \bar{r}\s^d|\vae|\mathsf{L}\label{distyy1I0}\\
		&\ltby{cond1Bis} \bar{r}\;\frac{\s}{3}\nonumber\\
		&< \frac{\bar{r}}{2}\;.\nonumber
\end{align}
\end{itemize}
Therefore, Lemma~\ref{IFTLem} applies. Hence, there exists a function $g\colon D^1_{2|\vae|}(0)\to D_{\bar{r}}(\mathsf{y})$ such that its graph coincides with $F^{-1}(\{0\})$. In particular, $\mathsf{y}'\coloneqq g(\vae)$ is the unique $y\in D_{\bar{r}}(\mathsf{y})$ satisfying $0=F(y,\vae)=\dpr_y K'(y)-\o$ \ie the second part of \eqref{ArnH1}. Moreover, 
\beq{EcarY1Y0}
|\mathsf{y}'-\mathsf{y}|\leq 2\|T\|\|F(\mathsf{y},\cdot)\|_{2|\vae|,0}\leq \frac{8|\vae|\mathsf{T} M}{r}\leby{distyy1I0} \bar{r}\;\s^{d}{|\vae|}\mathsf{L}< \frac{\bar{r}}{2}\;,
\eeq
so that
\beq{NextSetArn}
D_{\frac{\bar{r}}{2}}(\mathsf{y}')\subset D_{\bar{r}}(\mathsf{y}).
\eeq
Next, we prove that $\dpr^2_y K'(\mathsf{y}')$ is invertible. Indeed, by Taylor' formula, we have
\begin{align*}
\dpr^2_y K'(\mathsf{y}')&= K_{yy}(\mathsf{y})+ \dst\int_0^1 K_{yyy}(\mathsf{y}+t\vae \wt y)\cdot\vae\wt y dt+\vae \wt K_{yy}(\mathsf{y}')\\
           &= T^{-1}\left(\uno_d+\vae T\left(\dst\int_0^1 K_{yyy}(\mathsf{y}+t\vae \wt y)\cdot\wt y dt+ \wt K_{yy}(\mathsf{y}')\right)\right)\\
           &\eqqcolon T^{-1}(\uno_d+\vae A),
\end{align*}
and, by Cauchy's estimate, 
\begin{align*}
|\vae|\|A\|&\leq \|T\|\left(d\|K_{yyy}\|_{r/2,\mathsf{y}}|\vae||\mathsf{y}'-\mathsf{y}|+ |\vae|\|\dpr_{y'}^2\wt K\|_{r/2,\mathsf{y}}\right)\\
     &\leq \|T\|\left(\frac{d\|K_{yy}\|_{r,\mathsf{y}}}{r-\frac{r}{2}}|\vae||\mathsf{y}'-\mathsf{y}|+|\vae|\|\wt K_{yy}\|_{r/2,\mathsf{y}}\right)\\
	 &\leby{EcarY1Y0} \mathsf{T}\left(\frac{2d\mathsf{K}}{r}\frac{8|\vae|\mathsf{T} M}{r}+\frac{4|\vae|M}{r^2} \right)\\
	 &\leq \frac{4|\vae|\mathsf{T}M}{r^2}(4d\mathsf{T}\mathsf{K}+1)\\
	 &\leq\frac{20d|\vae|\mathsf{T}^2\mathsf{K} M}{r^2}\\
	 &\leby{rrbarAs} \frac{25}{6d}|\vae|\frac{\mathsf{T}  M}{r\bar{r}}\\
	 &< \su 2|\vae|\mathsf{L}\\
	 &\leby{cond1Bis}\frac{\s}{6}\\
	 &\le\su2.
\end{align*}
Hence $\dpr_{y'}^2 K'(\mathsf{y}')$ is invertible with
\[\dpr_{y'}^2 K'(\mathsf{y}')^{-1}=(\uno_d+\vae A)^{-1}T=T+\dst\sum_{k\geq 1}(-\vae)^k A^k T\eqqcolon T+\vae \wt T,\]
and
\[|\vae|\|\wt T\|\leq |\vae|\frac{\|A\|}{1-|\vae|\|A\|}\|T\|\leq 2|\vae|\|A\| \|T\|
\le |\vae|\mathsf{L}\mathsf{T}
\le 2\frac{\s}{6}\mathsf{T}
= \mathsf{T}\frac{\s}{3}\,.\]
Next, we prove estimate on $P_+$. We have,
\[|\vae|\|g_x\|_{\bar{r},\bar{s},\mathsf{y}}\leq |\vae|\mathsf{C}_1 \frac{M}{\a} \s^{-(\t+d+1)} \le |\vae| \frac{r}{3}\mathsf{L}\leby{cond1Bis}\frac{r}{3}\frac{\s}{3}\le \frac{r}{3}\]
so that, for any $(y',x)\in D_{\bar{r},\bar{s}}(\mathsf{y})$,
\[ |y'+\vae g_x(y',x)-\mathsf{y}|\leq \bar{r}+\frac{r}{3}< \frac{r}{8d}+\frac{r}{3}<\frac{2r}{3}<r\,,\]
and thus
\begin{align*}
\|P^\ppu\|_{\bar{r},\bar{s},\mathsf{y}}&\leq d^2 \|K_{yy}\|_{r,\mathsf{y}}\|g_x\|_{\bar{r},\bar{s},\mathsf{y}}^2\leq d^2 \mathsf{K}\left( \mathsf{C}_1 \frac{M}{\a} \s^{-(\n+d)}\right)^2\\
   &=d^2\mathsf{C}_1^2 \frac{\mathsf{K}M^2}{\a^2} \s^{-2(\n+d)}, 
\end{align*}
\begin{align*}
\|P^\ppd\|_{\bar{r},\bar{s},\mathsf{y}}&\leq d\|P_y\|_{\frac{5r}{6},\bar{s},\mathsf{y}}\|g_x\|_{\bar{r},\bar{s},\mathsf{y}}\leq d\frac{6M}{r}\mathsf{C}_1 \frac{M}{\a} \s^{-(\n+d)}\\
     &= 6d\mathsf{C}_1 \frac{M^2}{\a r}\s^{-(\n+d)}
\end{align*}
and by Lemma~\ref{fce}--$(i)$, we have,
\begin{align*}
|\vae|\|P^\ppt\|_{\bar{r},s-\frac{\s}{2},\mathsf{y}}&\leq \dst\sum_{|n|_1>\k}\|P_n\|_{\bar{r},\mathsf{y}}\ex^{(s-\frac{\s}{2})|n|_1}\leq M\dst\sum_{|n|_1>\k}\ex^{-\frac{\s |n|_1}{2}}\\
  &\leq M\ex^{-\frac{ \k\s}{4}}\dst\sum_{|n|_1>\k}\ex^{-\frac{\s |n|_1}{4}}\leq M\ex^{-\frac{ \k\s}{4}}\dst\sum_{|n|_1>0}\ex^{-\frac{\s |n|_1}{4}}\\
  &= M\ex^{-\frac{ \k\s}{4}} \left(\left(\dst\sum_{k\in \integer}\ex^{-\frac{\s |k|}{4}}\right)^d-1\right)=M\ex^{-\frac{ \k\s}{4}}\left(\left(1+\frac{2\ex^{-\frac{\s }{4}}}{1-\ex^{-\frac{\s }{4}}} \right)^d-1\right)\\
  &= M\ex^{-\frac{ \k\s}{4}}\left(\left(1+\frac{2}{\ex^{\frac{\s }{4}}-1} \right)^d-1\right)\leq M\ex^{-\frac{ \k\s}{4}}\left(\left(1+\frac{2}{\frac{\s }{4}} \right)^d-1\right)\\
  &\leq \s^{-d} M\ex^{-\frac{ \k\s}{4}}\left(\left(\s +8 \right)^d-\s^d\right)\leq d 8^{d}\s^{-d} M\ex^{-\frac{ \k\s}{4}}\\
  &= \mathsf{C}_2\s^{-d} M\ex^{-\frac{5}{4}\l}\\
  &\leby{DefNArn} \mathsf{C}_2\s^{-d} M\s^{-(2\n+d)}\frac{|\vae|{M}\mathsf{K}}{{\a}^2}\\
  &= \mathsf{C}_2 M\frac{|\vae|{M}\mathsf{K}}{{\a}^2}\s^{-2(\n+d)}\,.
\end{align*}
Hence\footnote{Recall that $\s<1$.},
\begin{align*}
\|P_+\|_{\bar{r},\bar{s},\mathsf{y}}&\leq \|P^\ppu\|_{\bar{r},\bar{s},\mathsf{y}}+\|P^\ppd\|_{\bar{r},\bar{s},\mathsf{y}}+\|P^\ppt\|_{\bar{r},\bar{s},\mathsf{y}}\\
  &\leq d^2\mathsf{C}_1^2 \frac{\mathsf{K}M^2}{\a^2} \s^{-2(\n+d)}+6d\mathsf{C}_1 \frac{M^2}{\a r}\s^{-(\n+d)}+\mathsf{C}_2 M\frac{|\vae|{M}\mathsf{K}}{{\a}^2}\s^{-2(\n+d)}\\
  &= \left(d^2\mathsf{C}_1^2 r\mathsf{K}+6d\mathsf{C}_1 \a \s^{\n+d}+\mathsf{C}_2 r\mathsf{K}\right)\frac{M^2}{\a^2 r}\s^{-2(\t+d+1)}\\
  &\le \left(d^2\mathsf{C}_1^2+6d\mathsf{C}_1 +\mathsf{C}_2\right)\max\left\{\a,r\mathsf{K}\right\}\frac{M^2}{\a^2 r}\s^{-2(\t+d+1)}\\
  &\leby{RecHypArn} \frac{\mathsf{C}_3}{\sqrt{2}} \max\left\{1,\frac{\a}{r\mathsf{K}} \right\}\frac{M^2\mathsf{K}}{\a^2 }\s^{-2(\n+d)}\\
  &\le \mathsf{L}M\;.
\end{align*}
The proof of the claims on $\phi'$ and $P'$ are proven in a similar way as in Lemma~\ref{lem:1}.
\qed

\noi
Finally, we prove the convergence of the scheme by mimicking Lemma~\ref{lem:2}.
\lem{lem:2Bis} 
Let   $H_0\coloneqq H$, $K_0\coloneqq K$, $P_0\coloneqq P$, $\phi^0=\phi_0\coloneqq\id $, and $r_0$, $s_0$, $s_*$, $\s_0$, $\m_0$, $\mathsf{W}_0$, $M_0$, $\mathsf{K}_0$, $\mathsf{T}_0$ and $\vae_\sharp$ be  as in $\S\ref{AnolKam}$. For a given $\vae\not=0$, 
 define\footnote{Notice that $s_{j}\downarrow s_*$ and $r_{j}\downarrow 0$.}
\begin{align*}
 \dst\sigma_j&\coloneqq \frac{\sigma_0}{2^j}\,,\\
  s_{j+1}&\coloneqq s_j-\sigma_j=s_*+\frac{\sigma_0}{2^j}\,,\\
  \bar s_{j}&\coloneqq s_j-\frac{2\s_i}{3}\,,\\
 \mathsf{K}_{j+1}&\coloneqq \mathsf K_0\dst\prod_{k=0}^{j}(1+\frac{\s_k}{3})\le \mathsf K_0\ex^{\frac{2\s_0}{3}}<\mathsf{K}_0\sqrt{2}\,,  \\
\mathsf{T}_{j+1}&\coloneqq \mathsf T_0\dst\prod_{k=0}^{j}(1+\frac{\s_k}{3})\le \mathsf T_0\ex^{\frac{2\s_0}{3}}<\mathsf{T}_0\sqrt{2}\,,\\
\l_0    &\coloneqq \log\m_0^{-1}\,,\\
 \dst\l_j&\coloneqq 2^j\l_0\,,\\
       \k_j&\coloneqq 5\s_j^{-1}\l_j\,,\\
 r_{j+1}&\coloneqq  \dst\min\left\{\frac{\a}{4d\sqrt{2}\mathsf{K}_0\k_j^{\n}}\,,\, \frac{5}{96d}\frac{r_j}{\eta_0} \right\} \,,\\
  \mathsf{d}_* &\coloneqq \mathsf{C}_5\;\eta_0^2\,,\\
 \mathsf{e}_* &\coloneqq \mathsf{C}_9\frac{\mathsf{K}_0}{\a^2} \s_0^{-(4\n+2d+1)}\;\l_0^{2\n}\;,\\
 \mathsf{f}_*&\coloneqq \mathsf{C}_8\;\max\left\{1,\frac{\a}{r_0\mathsf{K}_0}\right\}\eta_0\;\s_0^{-(3\n+2d+1)}\;\m_0\;\l_0^{\n}\,.
\end{align*}
Assume that $\vae$ is such that
\beq{condBis}
\m_0\leq \vae_\sharp 
\quad\mbox{and}\quad  \;\mathsf{f}_*\max\left\{1\,,\,\frac{\mathsf{C}_{10}}{3}\;\s_0\;\;\eta_0^{\su4}\;\mathsf{e}_*{\;d_*^2}\;|\vae|\; M_0\right\}< 1\,. 
\eeq
Then, one can construct a sequence of symplectic transformations 
\beq{phijBis}
\phi_j:D_{r_j,s_{j}}(y_j)\to D_{r_{j-1},s_{j-1}}(y_{j-1})\ ,
\eeq
so that
\beq{HjBis}
H_j:=H_{j-1}\circ\phi_j=: K_j + \vae^{2^j} P_j\ ,
\eeq
converges uniformly. More precisely, 
$\vae^{2^j} P_j$, 
$\phi^j\coloneqq \phi_0\circ\phi_1\circ\phi_2\circ \cdots\circ \phi_j$, 
 $K_j$, $y_j$ converge uniformly on 
$\{y_*\}\times\dst\torus^d_{s_*}$ to, respectively, $0$, $\phi_*$, $K_*$, $y_*$ which are real--analytic on $\dst\torus^d_{s_*}$ and  $H\circ\phi_*=K_*$ with $\det\dpr^2_yK_*(y_*)\not=0$. 
Finally, the following estimates hold for any $i\ge 1$:
\beqa{estfin1Bis}
&& 
|\vae|^{2^i}\|P_i\|_{r_i,s_i,y_i}\le\frac{\left(\frac{1}{3}\s_0\mathsf{f}_*\mathsf{e}_*{d_*^2}|\vae| M_0\right)^{2^{i-1}}}{\mathsf{e}_* {\,d_*}^{i+1}}\ ,\\
&&
|\mathsf{W}(\phi_*-\id)|
\le \s_0^{d+1}
\quad \mbox{on}\quad \{y_*\}\times\torus^d_{s_*}
\label{estfin2Bis}\ .
\eeqa
\elem 
\proof
For $i\ge 0$, define
\begin{align*}
 \mathsf{W}_i&\coloneqq \diag\left(\max\left\{\frac{\mathsf{K}_i}{\a},\su{r_i}\right\}\uno_d\;,\uno_d\right)\,,\\ 
 \ovl{\mathsf{L}}_i&\coloneqq \mathsf{C}_0 \max\left\{1,\frac{\a}{r_i\mathsf{K}_i}\right\}\frac{M_i \mathsf{K}_0}{\a^2 }\s_i^{-(2\n+d)}\;,\\
\mathsf{L}_i&\coloneqq M_i\dst\max\left\{\frac{4\sqrt{2}\mathsf{T}_0  }{r_ir_{i+1}}\s_i^{-(\n+d)}\,,\, \mathsf{C}_7 \max\left\{1,\frac{\a}{r_i\mathsf{K}_i}\right\}\frac{\mathsf{K}_0}{\a^2}\s_i^{-2(\n+d)}\right\}\\
          &\ge M_i\dst\max\left\{\frac{4\mathsf{T}_i  }{r_ir_{i+1}}\s_i^{-(\n+d)}\,,\,\frac{4}{\mathsf{K_i} r_i^2}\,,\, \mathsf{C}_7 \max\left\{1,\frac{\a}{r_i\mathsf{K}_i}\right\}\frac{\mathsf{K}_0}{\a^2}\s_i^{-2(\n+d)}\right\}
 \, .
 \end{align*}
Let us assume ({\sl inductive hypothesis}) that we can iterate $j\ge 1$ times the KAM step obtaining $j$ symplectic transformations\footnote{Compare \equ{phiokBis0}.} 
\beq{bes06}
\phi_{i+1}:D_{r_{i+1},s_{i+1}}(y_{i+1})\to D_{2r_i/3, s_i}(y_i),\quad \mbox{for}\quad 0\le i\le j-1,
\eeq
 and $j$  Hamiltonians $H_{i+1}=H_i\circ\phi_{i+1}=K_{i+1}+\vae^{2^{i+1}} P_{i+1}$ real--analytic on $D_{r_{i+1},s_{i+1}}(y_{i+1})$ such that, for any $0\le i\le j-1$,
\beq{bbbBis}
\left\{
\begin{array}{l}
\|\dpr_y^2 K_i\|_{r_i,y_i}\le \mathsf{K}_i\,,\ \\  \ \\
\|T_i\|\le \mathsf{T}_i\,,\ \\  \ \\ 
 \|P_{i}\|_{r_{i},s_{i},y_{i}}\le M_{i}\,,\ \\  \ \\ 
 \l_i\ge \frac{4}{5}\log\left(\s_i^{2\n+d}\frac{\a^2}{|\vae|^{2^i}M_iK_i}\right)\,,\ \\  \ \\ 
|\vae|^{2^i} \mathsf{L}_i \le \frac{\sigma_i}{3}
\ . 
\end{array}\right.
\eeq
Observe that for $j=1$, it is $i=0$ and \equ{bbbBis} is implied by the definitions of $\mathsf{K}_0,\,\mathsf{T}_0,\,\l_0,\, M_0$  and by condition \equ{condBis}.

\noi
Because of \equ{condBis} and \equ{bbbBis}, \equ{cond1Bis} holds 
for $H_i$ and Lemma~$\ref{lem:1bis}$  can be applied to $H_i$ and  one has, for $0\le i\le j-1$  (see \equ{Est1Lem1b}, \equ{convEst}, \equ{phiokBis1} and \equ{tesitBis}): 
\beqa{C.1Bis}
\left\{
\begin{aligned}
&|y_{i+1}-y_i|\le 2\s_i^{\n+d}r_{i+1}|\vae|^{2^i}\mathsf L_i\;,\\ 
 &\|T_{i+1}\|\le \|T_i\|+\mathsf T_i|\vae|^{2^i}\mathsf L_i\;,  \\
&\|\mathsf{K}_{i+1}\|_{r_{i+1},y_{i+1}}\le \|\mathsf{K}_i\|_{r_i,y_{i}}+|\vae|^{2^i}M_i \;,\\
&\|\dpr_y^2\mathsf{K}_{i+1}\|_{r_{i+1},y_{i+1}}\le \|\dpr_y^2\mathsf{K}_i\|_{r_i,y_{i}}+\mathsf K_i|\vae|^{2^i}\mathsf L_i \;,\\
&\|\mathsf{W}_i(\phi_{i+1}-\id)\|_{r_{i+1},s_{i+1},y_{i+1}}\le \s_i^{d}\;|\vae|^{2^i}{\mathsf L}_i \;,\\
&\|P_{i+1}\|_{r_{i+1},s_{i+1},y_{i+1}}\le M_{i+1}\coloneqq M_i \mathsf L_i\;.
\end{aligned}
\right.
\eeqa
Let $0\le i\le j-1$. Since
$$
\m_0\le \vae_\sharp \implies \frac{\a}{4d\sqrt{2}\mathsf{K}_0\k_0^{\n}}\le \frac{5}{96d}\frac{r_0}{\eta_0} \,,  
$$
then
$$
r_1=\frac{\a}{4d\sqrt{2}\mathsf{K}_0\k_0^{\n}}\;,
$$
and therefore
\begin{align*}
r_{i+1}&= \dst\min\left\{\frac{\a}{4d\sqrt{2}\mathsf{K}_0\k_i^{\n}}\,,\, \frac{5}{96d}\frac{r_i}{\eta_0} \right\}\\
	   &= \dst\min\left\{\frac{\a}{4d\sqrt{2}\mathsf{K}_0\k_i^{\n}}\,,\, \frac{5}{96d\eta_0} \frac{\a}{4d\sqrt{2}\mathsf{K}_0\k_{i-1}^{\n}}\,,\, \left(\frac{5}{96d\eta_0}\right)^2r_{i-1}\right\}\\
	   &\,\ \vdots \\
	   &= \dst\min\left\{\frac{\a}{4d\sqrt{2}\mathsf{K}_0\k_i^{\n}}\,,\, \frac{5}{96d\eta_0} \frac{\a}{4d\sqrt{2}\mathsf{K}_0\k_{i-1}^{\n}}\,,\,\cdots\,,\, \left(\frac{5}{96d\eta_0}\right)^ir_{1}\right\}\\
	   &= \frac{r_1}{4^i}\dst\min\left\{1\,,\, \frac{5}{24d\eta_0} \,,\,\cdots\,,\, \left(\frac{5}{24d\eta_0}\right)^i\right\}\\
	   &=  \left(\frac{5}{96d\eta_0}\right)^ir_1\,.
\end{align*}
Thus, since
\beq{k0sup10}
\m_0\le \vae_\sharp \implies \m_0\le \ex^{-1}\implies \k_0\ge 5\s_0^{-1}\ge 10\,,
\eeq
we have
\begin{align*}
|\vae| \mathsf L_0 (3 \sigma_0^{-1})&= 3|\vae|M_0\dst\max\left\{\frac{4\sqrt{2}\mathsf{T}_0  }{r_0r_{1}}\s_0^{-(\n+d)}\,,\, \mathsf{C}_7 \max\left\{1,\frac{\a}{r_0\mathsf{K}_0}\right\}\frac{\mathsf{K}_0}{\a^2}\s_0^{-2(\n+d)}\right\}\\
&\le 3\dst\max\left\{4\sqrt{2}\mathsf{T}_0\frac{ \a }{{r}_{1}}\frac{ \a }{r_0 \mathsf{K}_0}\,,\,\mathsf{C}_7 \max\left\{1,\frac{\a}{r_0\mathsf{K}_0}\right\}\right\}\s_0^{-2(\n+d)-1}\frac{ \mathsf{K}_0|\vae| M_0}{\a^2}\\
&= 3\dst\max\left\{32d\eta_0\k_0^{\n}\frac{ \a }{r_0 \mathsf{K}_0}\,,\,\mathsf{C}_7 \max\left\{1,\frac{\a}{r_0\mathsf{K}_0}\right\}\right\}\s_0^{-2(\n+d)-1}\frac{ \mathsf{K}_0|\vae| M_0}{\a^2}\\
&\le 3\dst\max\left\{32d\,,\,10^{-\n}\mathsf{C}_7 \right\}\cdot\max\left\{1,\frac{\a}{r_0\mathsf{K}_0}\right\}\eta_0\;\s_0^{-2(\n+d)-1}\frac{ \mathsf{K}_0|\vae| M_0}{\a^2}\k_0^{\n}\\
&= \mathsf{C}_8\;\max\left\{1,\frac{\a}{r_0\mathsf{K}_0}\right\}\eta_0\;\s_0^{-(3\n+2d+1)}\m_0\l_0^{\n}\\
&=\mathsf{f}_*\leby{condBis}1\;.
\end{align*}
Now, fix $i\ge 1$. We have
\beq{rikial3}
r_i\mathsf{K}_i\le r_1\mathsf{K}_0\sqrt{2}\leby{k0sup10}\frac{\a }{4d\cdot 10^\n}<\a
\eeq
so that
\begin{align*}
|\vae|^{2^i} \mathsf L_i (3 \sigma_i^{-1})&= 3|\vae|^{2^i} M_i\dst\max\left\{\frac{4\sqrt{2}\mathsf{T}_0  }{r_ir_{i+1}}\s_i^{-(\n+d)}\,,\, \mathsf{C}_7 \max\left\{1,\frac{\a}{r_i\mathsf{K}_i}\right\}\frac{\mathsf{K}_0}{\a^2}\s_i^{-2(\n+d)}\right\}\s_i^{-1}\\
&=3|\vae|^{2^i} M_i\dst\max\left\{\frac{4\sqrt{2}\mathsf{T}_0  }{r_ir_{i+1}}\s_i^{-(\n+d)}\,,\,\mathsf{C}_7 \frac{ 1}{\a r_i}\s_i^{-2(\n+d)}\right\}\s_i^{-1}\\
&\le 3\dst\max\left\{\frac{4\sqrt{2}\a\mathsf{T}_0  }{r_{i+1}}\,,\,\mathsf{C}_7 \right\}\s_i^{-2(\n+d)-1}\frac{|\vae|^{2^i} M_i}{\a r_i}\\
&= 3\dst\max\left\{32d\eta_0\k_0^\n\left(\frac{96d\eta_0}{5}\right)^i\,,\,\mathsf{C}_7 \right\}\s_i^{-2(\n+d)-1}\frac{|\vae|^{2^i} M_i}{\a r_i}\\
&\le 3\dst\max\left\{32d\,,\,10^{-\n}\mathsf{C}_7 \right\}\left(\frac{96d\eta_0}{5}\right)^i\eta_0\;\k_0^\n\;\s_i^{-2(\n+d)-1}\frac{|\vae|^{2^i} M_i}{\a r_i}\\
&\le 12d\sqrt{2}\dst\max\left\{32d\,,\,10^{-\n}\mathsf{C}_7 \right\}\left(\frac{96d\eta_0}{5}\right)^{2i-1}\eta_0\;\k_0^{2\n}\;\s_i^{-2(\n+d)-1}\frac{\mathsf{K}_0|\vae|^{2^i} M_i}{\a^2}\\
&=\mathsf{C}_9\frac{\mathsf{K}_0}{\a^2} \s_0^{-(4\n+2d+1)}\;\l_0^{2\n}\; \left(2^{2(\n+d)+11}3^25^{-2}d^2\eta_0^2\right)^i\;{|\vae|^{2^i} M_i}\\
&=\mathsf{e}_*\mathsf{d}_*^i |\vae|^{2^i}M_i\eqqcolon \frac{\theta_i}{\mathsf d_*}\;,
\end{align*}
 so that 
 $$
 \mathsf L_i<\mathsf{e}_*\,\mathsf d_*^i  M_i\;,
 $$ 
 thus by last relation in \equ{C.1Bis}, for any $1\le i\le j-1$,  
 $$
 |\vae|^{2^{i+1}}M_{i+1}<\mathsf{e}_*\mathsf d_*^i |\vae|^{2^{i+1}} M_i^2
 $$ 
 \ie $\theta_{i+1}<\theta_i^2$ , which iterated, yields 
 $\theta_i\le \theta_1^{2^{i-1}}$ for $1\le i\le j$. 
Next, we  show that, thanks to \equ{condBis}, \equ{bbbBis} holds also for $i=j$. In fact, by \equ{bbbBis} and \equ{C.1Bis},  
 we have 
$$
\|T_{i+1}\|\le \|T_i\|+\mathsf T_i|\vae|^{2^i}\mathsf L_i\le \mathsf T_i+\mathsf T_i\frac{\s_i}{3}=\mathsf T_{i+1}\,,
$$
and similarly for  $\|\dpr_y^2 K_{i+1}\|_{r_{i+1},y_{i+1}}$.
\nl 
Now, by $\equ{estfin1Bis}_{i=j}$\;,
$$|\vae|^{2^j}\mathsf L_j (3 \sigma_j^{-1})\le\frac{\theta_j}{\mathsf{d}_*}\le\su{\mathsf{d}_*}(\mathsf{e}_*{\,d_*^2}\vae^2 M_1)^{2^{j-1}}\le\su{\mathsf{d}_*}\left(\frac{\s_0}{3}\mathsf{f}_*\mathsf{e}_*{d_*^2}|\vae| M_0\right)^{2^{j-1}}
\leby{condBis} \su{\mathsf{d}_*}<1\ ,$$
which implies the last inequality in \equ{bbbBis} with $i=j$.\\
Next, we check the fourth inequality in \equ{bbbBis} for $i=j$. We have\footnote{Notice that $\mathsf{L}_i\ge M_i\mathsf{C}_7\frac{\mathsf{K}_0}{\a^2}\s_{i}^{-2(\n+d)},\, \forall\; i\ge 0$.\label{ftn600}}
\begin{align*}
\l_j&=2\l_{j-1}\\
	&\overset{\equ{bbbBis}_{i=j-1}}{\ge}\frac{4}{5}\log\left(\s_{j-1}^{2(2\n+d)}\frac{\a^4}{|\vae|^{2^{j}}M_{j-1}^2K_{j-1}^2}\right)\\
	&\ge \frac{4}{5}\log\left(\s_{j-1}^{2(2\n+d)}\frac{\a^4}{|\vae|^{2^{j}}M_{j-1}K_{j-1}^2}\cdot\s_{j-1}^{-2(\n+d)}\frac{\mathsf{C}_7\mathsf{K}_0}{\a^2\mathsf{L}_{j-1}}\right)\\
	&= \frac{4}{5}\log\left(\s_{j-1}^{2\n}\frac{\a^2}{|\vae|^{2^{j}}M_{j}K_{j-1}}\cdot\frac{\mathsf{C}_7\mathsf{K}_0}{\mathsf{K}_{j-1}}\right)\\
	&> \frac{4}{5}\log\left(\s_{j}^{2\n+d}\frac{\a^2}{|\vae|^{2^{j}}M_{j}K_{j}}\right)\,.
\end{align*}
The proof of the induction is then finished and one can construct an {\sl infinite sequence} of Arnold's transformations satisfying \equ{bbbBis}, \equ{C.1Bis} and \equ{estfin1Bis} {\sl for all $i\ge 0$}. \\
\nl
Next, we prove that $\phi^j$ is convergent by proving that it is Cauchy. For any $j\ge 3$, we have, using again Cauchy's estimate,\footnote{Notice that $\equ{rikial3}\implies W_i=\diag(\su{r_i}\;\uno_d\;,\uno_d)\;,\,\forall\; i\ge 1$ and recall that $2^{i-1}\ge i,\,\forall\; i\ge 0$.}
\beqano
\|\mathsf{W}_{j-1}(\phi^{j-1}-\phi^{j-2})\|_{r_j,s_j,y_j}&=&\|\mathsf{W}_{j-1}\phi^{j-2}\circ\phi_{j-1}-\mathsf{W}_{j-1}\phi^{j-2}\|_{r_j,s_j,y_j}\\
           &\leby{bes06}& \|\mathsf{W}_{j-1}D\phi^{j-2}\mathsf{W}_{j-1}^{-1}\|_{2r_{j-1}/3, s_{j-1},y_{j-1}}\, \|\mathsf{W}_{j-1}(\phi_{j-1}-\id)\|_{r_j,s_j,y_j}\\
           &\leby{C.1Bis}&  \max\left(r_{j-1}\frac{3}{r_{j-1}},\frac{3}{2\s_{j-1}}\right)    \|\mathsf{W}_{j-1}\phi^{j-2}\|_{r_{j-1}, s_{j-1},y_{j-1}} \times \\
           &&\qquad \times\|\mathsf{W}_{j-1}(\phi_{j-1}-\id)\|_{r_j,s_j,y_j}\\
           &=&  \frac{3}{2\s_{j-1}}   \|\mathsf{W}_{j-1}\phi^{j-2}\|_{r_{j-1}, s_{j-1},y_{j-1}} \, \|\mathsf{W}_{j-1}(\phi_{j-1}-\id)\|_{r_j,s_j,y_j}\\
           &\le & \frac{1}{2}    \|\mathsf{W}_{j-1}\phi^{j-2}\|_{r_{j-1}, s_{j-1},y_{j-1}} \cdot \s_{j-1}^d\left(|\vae|^{2^{j-1}}{\mathsf{L}}_{j-1}3\s_{i-1}^{-1}\right)\\
           &\le & \frac{1}{2}    \|\mathsf{W}_{j-1}\phi_1\|_{r_{2}, s_2,y_{2}} \cdot \s_{j-1}^d\;\th_{j-1}\\
           &\le &  \frac{1}{2}\left(\dst\prod_{i=1}^{j-2}\|\mathsf{W}_{i+1}\mathsf{W}_{i}^{-1}\| \right)\|\mathsf{W}_{1}\phi_1\|_{r_{2}, s_2,y_{2}} \cdot \s_{j-1}^d\;\th_{j-1}\\
           &\eqby{rikial3}& \frac{1}{2}\left(\dst\prod_{i=1}^{j-2}\frac{r_i}{r_{i+1}} \right)\|\mathsf{W}_{1}\phi_1\|_{r_{2}, s_2,y_{2}} \cdot \s_{j-1}^d\;\th_{j-1}\\
           &=& \frac{r_1}{2r_{j-1}}\|\mathsf{W}_{1}\phi_1\|_{r_{2}, s_2,y_{2}} \cdot \s_{j-1}^d\;\th_{j-1}\\
           &=& \frac{48d}{5}\s_2^d\;\eta_0\;\|\mathsf{W}_{1}\phi_1\|_{r_{2}, s_2,y_{2}} \cdot \left(\frac{3d\cdot 2^{5-d}\eta_0}{5}\right)^{j-3}\cdot \;\th_1^{2^{j-2}}\\
           &\le& \frac{48d}{5}\s_2^d\;\eta_0\;\|\mathsf{W}_{1}\phi_1\|_{r_{2}, s_2,y_{2}} \cdot \left(\frac{3d\cdot 2^{5-d}\eta_0}{5}\right)^{2^{j-4}}\cdot \;\th_1^{2^{j-2}}\\
           &=& \frac{48d}{5}\s_2^d\;\eta_0\;\|\mathsf{W}_{1}\phi_1\|_{r_{2}, s_2,y_{2}} \cdot \left(\left(\frac{3d\cdot 2^{5-d}}{5}\right)^{\su4}\eta_0^{\su4}\; \th_{1}\right)^{2^{j-2}}\\
           &=& \frac{48d}{5}\s_2^d\;\eta_0\;\|\mathsf{W}_{1}\phi_1\|_{r_{2}, s_2,y_{2}} \cdot \left(\mathsf{C}_{10}\;\eta_0^{\su4}\; \th_{1}\right)^{2^{j-2}}  \;.
\eeqano
\noi
Therefore, for any $n\ge 1,\, j\geq 0$,
\begin{align*}
\|\mathsf{W}_{1}(\phi^{n+j+1}-\phi^n)\|_{r_{n+j+2},s_{n+j+2},y_{n+j+2}}&\leq  \sum_{i=n}^{n+j}\|\mathsf{W}_{1}(\phi^{i+1}-\phi^i)\|_{r_{i+2},s_{i+2},y_{i+2}}\\
&\le \sum_{i=n}^{n+j}\left(\dst\prod_{k=1}^{i}\|\mathsf{W}_{k}\mathsf{W}_{k+1}^{-1}\| \right)\|\mathsf{W}_{i+1}(\phi^{i+1}-\phi^i)\|_{r_{i+2},s_{i+2},y_{i+2}}\\
&\eqby{alfhtrikpi} \sum_{i=n}^{n+j}\dst\prod_{k=1}^{i}\max\left\{1\;,\frac{r_{k+1}}{r_k} \right\}\|\mathsf{W}_{i+1}(\phi^{i+1}-\phi^i)\|_{r_{i+2},s_{i+2},y_{i+2}}\\
&= \sum_{i=n}^{n+j}\|\mathsf{W}_{i+1}(\phi^{i+1}-\phi^i)\|_{r_{i+2},s_{i+2},y_{i+2}}\\
&\le \frac{48d}{5}\s_2^d\;\eta_0\;\|\mathsf{W}_{1}\phi_1\|_{r_{2}, s_2,y_{2}}\dst\sum_{i=n}^{n+j} \left(\mathsf{C}_{10}\;\eta_0^{\su4}\; \th_{1}\right)^{2^{i}}
\end{align*}
and
$$
\mathsf{C}_{10}\;\eta_0^{\su4}\; \th_{1}\ltby{condBis}1\;.
$$
Hence, 
$\phi^j$ converges uniformly on $\{y_*\}\times\torus^d_{s_*}$ to some $\phi^*$, which is then real--analytic map in $x\in\torus^d_{s_*}$.

\nl
To estimate $|\mathsf{W}_0(\phi^*-\id)|$ on $\{y_*\}\times\torus^d_{s_*}$, observe that
, for $i\ge 1$,\footnote{Recall that $2^{i}\ge i+1,\, \forall\, i\ge 0$.}
$$\s_{i}^d\;|\vae|^{2^i}\mathsf L_i\le \frac{\s_0^{d+1}}{3 \cdot 2^{i(d+1)}}\ \frac{\th_1^{2^{i-1}}}{\mathsf{d}_*} \le \frac{\s_0^{d+1}}{3 \cdot 2^{i(d+1)} \mathsf{d}_*}\th_1^{{i}}= \frac{\s_0^{d+1}}{3\mathsf{d}_*} \Big(\frac{\th_1}{2^{d+1}}\Big)^{i}$$
and therefore 
$$\dst\sum_{i\ge 1}  \s_{i}^d\;|\vae|^{2^i}\mathsf L_i\le \frac{\s_0^{d+1}}{3\mathsf{d}_*}\sum_{i\ge 1}\Big(\frac{\th_1}{2^{d+1}}\Big)^{i}\le \frac{\s_0^{d+1}\;\th_1}{3\cdot {2^{d}}\;\mathsf{d}_*}\ltby{condBis}\su2\;\s_0^{d+1}
\ .$$ 
Moreover, for any $i\ge 1$,
\begin{align*}
\|\mathsf{W}_1(\phi^i-\id)\|_{r_{i+1},s_{i+1},y_{i+1}}&\le \|\mathsf{W}_1(\phi^{i-1}\circ\phi_i-\phi_i)\|_{r_{i+1},s_{i+1},y_{i+1}}+\|\mathsf{W}_1(\phi_i-\id)\|_{r_{i+1},s_{i+1},y_{i+1}}\\
&\le \|\mathsf{W}_1(\phi^{i-1}-\id)\|_{r_{i},s_{i},y_{i}}+ \left(\dst\prod_{j=0}^{i-1}\|\mathsf{W}_{j}\mathsf{W}_{j+1}^{-1}\| \right) \|\mathsf{W}_{i}(\phi_i-\id)\|_{r_{i+1},s_{i+1},y_{i+1}}\\
&= \|\mathsf{W}_1(\phi^{i-1}-\id)\|_{r_{i},s_{i},y_{i}}+ \|\mathsf{W}_{i}(\phi_i-\id)\|_{r_{i+1},s_{i+1},y_{i+1}}\\
&= \|\mathsf{W}_1(\phi^{i-1}-\id)\|_{r_{i},s_{i},y_{i}}+ \|\mathsf{W}_{i}(\phi_i-\id)\|_{r_{i+1},s_{i+1},y_{i+1}}\\
&\le \|\mathsf{W}_1(\phi^{i-1}-\id)\|_{r_{i},s_{i},y_{i}}+\s_{i}^d\;|\vae|^{2^{i}}{\mathsf{L}}_{i}\ ,
\end{align*}
which iterated yields
\begin{align*}
\|\mathsf{W}_1(\phi^i-\id)\|_{r_i,s_i,y_i}&\le \dst\sum_{k=1}^{i-1}\s_{k}^d\; |\vae|^{2^k}{\mathsf{L}}_k\\
&\le  \dst\sum_{k\ge 1}\s_{k}^d\;|\vae|^{2^k}{\mathsf{L}}_k\\
&\le \su2\;\s_0^{d+1}
\,.
\end{align*}
Therefore, taking the limit over $i$ completes yields, uniformly on $\{y_*\}\times \torus^d_{s_*}$,
$$
|\mathsf{W}_1(\phi^*-\id)|\le \su2\;\s_0^{d+1}\;.
$$
Now, to complete the proof of the Lemma and, consequently, of the Theorem, just set $\phi_*\coloneqq \phi_0\circ \phi^*$ and observe that, uniformly on $\{y_*\}\times \torus^d_{s_*}$,
\beqano
|\mathsf{W}_0(\phi_*-\id)|&\le& |\mathsf{W}_0(\phi_0\circ \phi^*-\phi^*)|+|\mathsf{W}_0(\phi^*-\id)|\\
&\le& \|\mathsf{W}_0(\phi_0-\id)\|_{r_1,s_1,y_1}+\|\mathsf{W}_0\mathsf{W}_1^{-1}\|\;|\mathsf{W}_1(\phi^*-\id)|\\
&=& \|\mathsf{W}_0(\phi_0-\id)\|_{r_1,s_1,y_1}+\max\left\{\frac{r_1\mathsf{K}_0}{\a}\;,\,\frac{r_1}{r_0}\;,\,1 \right\}|\mathsf{W}_1(\phi^*-\id)|\\
&\eqby{rikial3}& \|\mathsf{W}_0(\phi_0-\id)\|_{r_1,s_1,y_1}+|\mathsf{W}_1(\phi^*-\id)|\\
&\le& \s_{0}^d\;|\vae|{\mathsf{L}}_0+\su2\;\s_0^{d+1}\\
&\leby{condBis}& \su3\;\s_0^{d+1}+\su2\;\s_0^{d+1}\\
&<& \s_0^{d+1}
\,.
\eeqano
\qed
\subsection{Proof of Theorem~\ref{teo1}\label{prteo1}}
As usual, the proof is inductive: at each step $j\in \natural$, a small perturbation of some normal form $N_j=e_j(\o)+\o\cdot y$,
\[H_j=N_j+P_j\]
is considered. Then, a coordinates and parameter transformation $\mathcal{F}_j$ is constructed so that
\[H_j\circ \mathcal{F}_j=N_{j+1}+P_{j+1}\]
with another normal form $N_{j+1}$, some much smaller error term $P_{j+1}$ satisfying
\[\|P_{j+1}\|\leq C\|P_j\|^\frac{3}{2}\]
for some constant $C>0$ and the sequence $\mathcal{F}^{j+1}\coloneqq \mathcal{F}_0\circ \cdots \circ \mathcal{F}_j $ converges to an embedding of an invariant Kronecker torus.\\
The first step, called \textbf{KAM step}, will be then to describe one cycle of this iterative scheme in which, for readability, we drop the subscribe $j$ and consider a generic hamiltonian $H=N+P$. First of all, instead of $H$, we consider the hamiltonian $\bar{H}$ obtained from $H$ by first linearizing the perturbation $P$ in $y$ and then truncating its Fourier series in $x$ at some suitable high order $\k$.\\
The transformation $\mathcal{F}$ is of the form
\[\mathcal{F}\coloneqq(\Phi,\f)\coloneqq\left(\wt\Phi\circ(\pi_1,\pi_2;\f\circ\pi_3),\f\circ\pi_3\right)\colon (y,x;\o)\mapsto \left(\wt\Phi(y,x;\f(\o)),\f(\o)\right),\]
where $\wt\Phi$ is obtained as the time--1--map of the flow $\Phi_{F}^t$ of some hamiltonian $F$ and
\beq{proJectn}
\pi_j\colon \cn\times\cn\times\cn \to \cn ,\ j=1,2,3 : \pi_1(y,x,\o)=y,\ \pi_2(y,x,\o)=x,\ \pi_3(y,x,\o)=\o\ .
\eeq
 In particular, $\Phi$ is then symplectic for $\o$ fixed\footnote{Indeed, denoting the Lie derivative by $\mathcal{L}$ and the contraction operator by $\iota$, we have
\[\dst{\frac{d}{dt}\left(\phi^t_{F}\right)^*\varpi = \left(\phi^t_{F}\right)^*\dst\mathcal{L}_{X_{F}}\varpi = \left(\phi^t_{F}\right)^*\left(\iota_{X_{F}}\underbrace{d\varpi}_{0}+d\underbrace{\iota_{X_{F}}\varpi}_{-dF}\right) = 0 \implies \left(\phi^t_{F}\right)^*\varpi=\left(\phi^0_{F}\right)^*\varpi=\id^*\varpi=\varpi}.\]
}. Then, we iterate this cycle and prove the convergences.

In all this section, the sup--norm on $D_{r,s}\times \O_{\a,h}$ will be denoted by $\|\cdot\|_{r,s,h}\coloneqq \|\cdot\|_{r,s,h,d}$, while on $D_{r,s}\times \cn$ (resp. $D_{r,s}\times\rn$), it will be denoted by $\|\cdot\|_{r,s,\infty}$ (resp. $\|\cdot\|_{r,s,0}$).
\subsubsection{KAM step}
\lemtwo{KamStpLem}{KAM step}
Assume that $\|P\|_{r,s,h}\leq \epsilon$ with 
\begin{itemize}
\item[$(a)$] $\epsilon\leq \frac{1}{C_4} \a\eta r \s^\n$,
\item[$(b)$] $\epsilon\leq \frac{1}{C_6} hr$,
\item[$(c)$] $h \leq  \frac{\a}{2\k^{\bar{\n}}}$,
\end{itemize}
for some $0<\eta<1/8,\, 0<\s<s/10$ and sufficiently large $\k> \frac{d-1}{\s}$. 
Then there exist $\f\colon \rn\to \rn$ a $C^\infty$--diffeomorphism with $\f(\O)=\O$ and $\f\equiv \id$ on $\rn\setminus\O$ , $\wt\Phi \colon D_{\eta r,s-5\s}\times \cn \to D_{r,s}$ a ($\o$--) family of symplectic transformations parametrized over $\cn$, each being real--analytic with holomorphic extention to $D_{\eta r,s-5\s}$ and  $C^\infty$ in $\o$ on $\rn$ and such that, if $\mathcal{F}\coloneqq(\Phi,\f)\coloneqq\left(\wt\Phi\circ(\pi_1,\pi_2;\f\circ\pi_3),\f\circ\pi_3\right)$ 
, the following hold: its restriction map
\[\mathcal{F}_h = (\Phi,\f)\colon D_{\eta r,s-5\s}\times \O_{\a,\frac{h}{4}}\to D_{r,s}\times \O_{\a,h}\]
is well--defined, real--analytic (in all arguments),  $H\circ \mathcal{F}_h=N_+ + P_+$ with another normal form $N_+\coloneqq e_+(\o)+\o\cdot y$ and 
\[\|P_+\|_{\eta r,s-5\s,\frac{h}{4}} \leq \frac{\sqrt{C_{10}}}{3\cdot 2^{\bar\n}}\left(\frac{\epsilon^2}{\a r \s^{\bar\n}}+(\eta^2+\k^n\ex^{-\k\s})\epsilon \right).\]
Moreover\footnote{We denote by $D\Phi$ and $D\f$, respectively, the jacobian of $\Phi$ with respect to $(y,x,\o)$ and of $\f$ with respect to $\o$.},
\beq{kamest1} 
\|\bar W(\mathcal{F}-\id)\|_{\frac{r}{4},s-4\s,0}\le \left(\frac{\s_0}{\s}\right)^{2\bar{\n}-1}\cdot \dst\max\left(4C_3\frac{\epsilon}{\a r \s^{\bar{\n}}}, \frac{C_3}{C_0} \frac{\epsilon}{rh} \right)\ ,
\eeq
\beq{kamest2} 
\left\|\bar W(D\mathcal{F}-\Id)\bar W^{-1} \right\|_{\frac{r}{8} ,s-5\s,0}\le \left(\frac{\s_0}{\s}\right)^{2\bar{\n}-1} \cdot\dst\max\left(dC_5 \frac{\epsilon}{\a r \s^{\bar{\n}}},\frac{C_5}{C_0}\frac{\epsilon}{rh} \right)\ ,
\eeq
\beq{kamest3}
\left\|\f-\id \right\|_\infty,\, h\left\|D\f-\Id \right\|_0\leq \frac{C_5}{2C_0}\frac{\epsilon}{r}\ ,
\eeq
 for a given $\s_0\ge \s$, with 
$$
\bar{W}\coloneqq \diag\left(\su r\uno_d,\ \su{\s}\left(\frac{\s_0}{\s}\right)^{2\bar{\n}-1}\uno_d,\ \su h\uno_d\right).
$$
\elem
\proof
For convenience, we will follow the scheme of the proof in \cite{JP} and add two more steps allowing us, later, as we said, to estimate the Lipschitz' semi--norm of the symplectic transformation we are going to build--up without invoking the Whitney's extension theorem.\\
{\bf 1. Truncation.} Let $Q(x,\o)\coloneqq P(0,x,\o)+ P_y(0,x,\o)\cdot y$, the linerization of $P$ and $R\coloneqq T_\k Q$. Then by Cauchy's estimate we get
\[\|Q\|_{r,s,h}\leq \|P\|_{r,s,h}+d\frac{\|P\|_{r,s,h}}{r}r\leq (d+1)\epsilon,\]
And
\beqano
\|P-Q\|_{2\eta r,s,h}& =  & \left\|\dst\int_0^1 (1-t) P_{yy}(ty ,x,\o)(y,y) dt \right\|_{2\eta r,s,h}\\
				 & \leq & \dst \sum_{1\leq j,l\leq d}\, \dst\sup_{(y,x,\o)\in D_{2\eta r,s}\times \O_{\a,h} }\dst\int_0^1 (1-t)\|\partial_{y_j y_l}^2 P(ty,x,\o)\|\|y_j\|\|y_l\|dt \\
				 &\leq &  \dst \sum_{1\leq j,l\leq d} \dst\int_0^1 (1-t) \frac{\|P\|_{r,s,h}}{(r-2\eta r)^2} (2\eta r)^2 dt \\
				 & = & \frac{2\eta^2 d^2}{(1-2\eta)^2}\epsilon\\
				 & \leq & \frac{32 d^2}{9}\eta^2\epsilon.
\eeqano
By Lemma~\ref{fce}, we have 
\[\|R-Q\|_{r,s-\s,h}\leq 4^d C_2 \k^d \ex^{-\k\s}\|Q\|_{r,s,h} \leq 4^d(d+1) C_2 \k^d \ex^{-\k\s}\epsilon,\]
and therefore
\beqano
\|R\|_{r,s-\s,h}&\leq & \|R-Q\|_{r,s-\s,h}+\|Q\|_{r,s-\s,h}\\
	&\leq & \left(4^d C_2 \k^d \ex^{-\k\s}+1 \right)(d+1)\epsilon  \leq 2(d+1)\epsilon=\frac{C_3}{2C_0}\epsilon,
\eeqano
because $C_2\leq C_{10}$ and, later, $\k$ will be chosen so that $\k^d \ex^{-\k\s}\leq (4^{\bar{\n}} C_{10})^{-1}$.\\
{\bf 2. Extending the Diophantine condition.} The Diophantine condition (compare \eqref{dio}) is assumed to hold only on $\O_\a$. Nevertheless, given $\o\in \O_{\a,h}$, there exits $\o_*\in \O_\a$ such that $|\o-\o_*|<h$, so that, for any $|k|_1\leq \k$
\[|k\cdot(\o-\o_*)|\leq |k|_1\cdot |\o-\o_*|\leq \k h \overset{(c)}{\leq}\frac{\a}{2\k^\t}\leq \frac{\a}{2|k|_1^\t}\]
and thanks to \eqref{dio}, we get, for any $\o\in \O_{\a,h}$
\beq{edc}
|k\cdot \o|\geq \frac{\a}{2|k|_1^\t},\quad \forall \, 0\neq |k|_1\leq \k.
\eeq
{\bf 3. Finding the hamiltonian F by solving a homological equation.} We have
\[H=\bar{H}+ (P-Q)+(Q-R) \mbox{ with } \bar{H}=N+R\]
Let's remind that we are looking at for a hamiltonian $F$ such that its flow $\phi^t_F$
 satisfies\footnote{In fact, rigorously, one should write $H\circ (\phi^1_F,\id)$ and so on.}
 \[H\circ \phi^1_F=N_+^1 + P_+^1,\]
for some hamiltonian $N_+^1$ closed to a normal form and much smaller error term $P_+^1$. We have
\[H\circ \phi^1_F= \bar{H}\circ \phi^1_F+ (P-Q)\circ \phi^1_F + (Q-R)\circ \phi^1_F.\]
Next we expend $\bar{H}\circ \phi^t_F$ around $t=0$ to get\footnote{Given a function K and denoting the Poisson bracket by $\{\cdot,\cdot\}$, we have
\[\frac{d}{dt}K\circ \phi^t_F= dK(\phi^t_F)\cdot \frac{d}{dt}\phi^t_F= dK(\phi^t_F)\cdot \mathbb{J}dF(\phi^t_F)=\{K,F\}\circ \phi^t_F.\]}
\beqano
\bar{H}\circ \phi^1_F & = & N\circ\phi^1_F + R\circ \phi^1_F \\
					  &= & N+\frac{d}{ds}N\circ \phi^s_F{\big|}_{s=0}+ \dst\int_0^1 (1-t) \frac{d^2}{ds^2}N\circ \phi^s_F{\big|}_{s=t}\,dt +\\
					  & &+ R+\dst\int_0^1 \frac{d}{ds}R\circ \phi^s_F{\big|}_{s=t}\, dt \\
					  & = & N+\{N,F\}+ \dst\int_0^1 (1-t)\{\{N,F\},F\}\circ \phi^t_F\, dt+R+\int_0^1 \{R,F\}\circ \phi^t_F\, dt\\
					  &= & \underbrace{N+[R]}_{\coloneqq N_+^1}+ \{N,F\}+ (R-[R]) +\underbrace{ \dst\int_0^1 \{(1-t)\{N,F\}+R,F\}\circ \phi^t_F\, dt}_{\coloneqq P_+^0}
\eeqano
Since $Q$ is affine in the variable $y$ , then so is $R$ and a fortiori $[R]$; moreover $[R]$ does not depend on $x$. Therefore, there exist analytic functions $e_+^0\colon \o \to e_+^0(\o)$ and $v\colon \o\to v(\o)$ such that $[R](y,\o)=e_+^0(\o)+v(\o)\cdot y$ (in fact $v=[R_y]$) so that
\beq{nplus}
N_+^1=\underbrace{e(\o)+e_+^0(\o)}_{\coloneqq e_+^1(\o)}+ \underbrace{(\o+v(\o))}_{\coloneqq \rho(\o)}\cdot y 					
\eeq
Let
\beq{pplus}
P_+^1\coloneqq P_+^0 + (P-Q)\circ \phi^1_F + (Q-R)\circ \phi^1_F
\eeq
The main point is then to determine $F$ by solving the homological equation
\[\{N,F\}+ (R-[R])=0\]
$i.e.$ (recall that $N=e(\o)+\o\cdot y$)
\beq{homeq}
D_\o F=\{F,N\}=R-[R]
\eeq
so that we have
\[(1-t)\{N,F\}+R= (1-t)([R]-R)+R=(1-t)[R]+tR,\]
\beq{poplus}
P_+^0= \dst\int_0^1 \{(1-t)[R]+tR,F\}\circ \phi^t_F\, dt,
\eeq
and
\beq{newfrm}
H\circ \phi^1_F=N_+^1 + P_+^1.
\eeq
Since $[R-[R]]=0,\, \|[R]\|_{r,h}\leq \|R\|_{r,s-\s,h}$ and $\|R-[R]\|_{r,s-\s,h}\leq 2\|R\|_{r,s-\s,h}<\infty$, Lemma~\ref{sde} applies to \eqref{homeq} and we find $F$ with
\beq{fb}\|F\|_{r,s-2\s,h}\leq \frac{2C_0}{\a\s^\t}\|R-[R]\|_{r,s-\s,h}\leq \frac{4C_0}{\a\s^\t}\|R\|_{r,s-\s,h}\leq 2C_3\frac{\epsilon}{\a\s^\t}
\eeq
Then by Cauchy's estimate we get
\[\|F_x\|_{r,s-3\s,h}\leq \frac{\|F\|_{r,s-2\s,h}}{\s}\leq 2C_3\frac{\epsilon}{\a\s^{\bar{\n}}}, \]
\[ \|F_y\|_{\frac{r}{2},s-2\s,h}\leq \frac{2\|F\|_{r,s-2\s,h}}{r}\leq 4C_3\frac{\epsilon}{\a r \s^\t}\]
and 
\beq{fob}
\|F_\o\|_{r,s-2\s,\frac{h}{2}}\leq \frac{2\|F\|_{r,s-2\s,h}}{h}\leq 4C_3\frac{\epsilon}{\a \s^\t h}
\eeq
so that
\[\frac{1}{r}\|F_x\|_{\frac{r}{2},s-3\s,h},\, \frac{1}{\s}\|F_y\|_{\frac{r}{2},s-3\s,h}\leq 4C_3\frac{\epsilon}{\a r \s^{\bar{\n}}}\]
and by using assumption $(a)$, we get
\beq{fthb}
\|F_x\|_{r,s-3\s,h}\leq \frac{2C_3}{C_4} \eta r\leq \eta r\leq \frac{r}{8}
\eeq
\beq{fib}
\|F_y\|_{\frac{r}{2},s-2\s,h}\leq \frac{4C_3}{C_4} \eta \s\leq \s
\eeq
{\bf 4. Extending the hamiltonian $F$.} Thanks to Lemma~\ref{cutoff}, there exists a cut--off $\chi_1\in C(\cn)\cap C^{\infty}(\rn)$ with $0\leq \chi_1\leq 1,\ supp\chi_1 \subset \O_{\a,\frac{h}{2}}$ and $\chi_1\equiv 1$ on $\O_{\a,\frac{h}{4}}$. Now we extend $F$, witch we call $\tilde{F}$, as follows: $\tilde{F}\equiv 0$ on 
 $D_{r,s-2\s}\times\left( \cn\setminus \O_{\a,\frac{h}{2}}\right)$ and $\tilde{F}(y,x,\o)=\chi_1(\o)F(y,x,\o)$ on $D_{r,s-2\s}\times \O_{\a,\frac{h}{2}}$. Thus 
\begin{itemize}
\item[$(i)$] $\tilde{F}$ concide with $F$ on $D_{r,s-2\s}\times \O_{\a,\frac{h}{4}}$, is continuous on $D_{r,s-2\s}\times\cn$, $C^{\infty}$ on $D_{r,s-2\s}\times\rn$ and for any $\o\in \cn$ given, the map $(y,x)\mapsto \tilde{F}(y,x,\o)$ is real--analytic with holomorphic extention to $D_{r,s-2\s}$.
\item[$(ii)$] 
\begin{align}
\|\tilde{F}\|_{r,s-2\s,\infty}&\overset{def}{=} \sup_{D_{r,s-2\s}\times\cn} \|\tilde{F}\|\leq \|F\|_{r,s-2\s,h}\leq 2C_3\frac{\epsilon}{\a \s^\t} \, ,\label{ftildest1}\\
\|\tilde{F}_y\|_{\frac{r}{2},s-2\s,\infty} &\leq \|F_y\|_{\frac{r}{2},s-2\s,h}\leq 4C_3\frac{\epsilon}{\a r \s^\t}\leq \s \, ,\label{ftildest2}\\
\|\tilde{F}_x\|_{r,s-3\s,\infty}&\leq \|F_x\|_{r,s-3\s,h}\leq 2C_3\frac{\epsilon}{\a \s^{\bar{\n}}}\leq \frac{r}{8} \label{ftildest3}
\end{align}
and by using \eqref{cutof}, \eqref{fb} and \eqref{fob}, we get
\beqa{ftildest4}\|\tilde{F}_\o\|_{r,s-2\s,0} &\leq & 
\|\partial_\o \dst\chi_1\|_0 \cdot \|F\|_{r,s-2\s,h}+ \|F_\o\|_{r,s-\s,\frac{h}{2}}\\
&\leq & 24 C_1C_3\frac{\epsilon}{\a \s^\t h}+4C_3\frac{\epsilon}{\a \s^\t h} \le \frac{C_5}{4} \frac{\epsilon}{\a \s^\t h}\, .\eeqa
\end{itemize}
{\bf 5. Transforming coordinates.} As we said, the coordinates transformation $\wt\Phi$ is obtained as the time--1--map of the flow $\dst\phi^t_{\tilde{F}}$ of the hamiltonian $\tilde{F}$ with equations of motion
\[\dot{y}=-\tilde{F}_x,\quad \dot{x}=\tilde{F}_y \mbox{ or equivalently } \frac{d}{dt}\dst\phi^t_{\tilde{F}}=\mathbb{J}d\tilde{F}(\dst\phi^t_{\tilde{F}}).\]
By using \eqref{ftildest2} and \eqref{ftildest3}, we deduce that, given $\o\in \cn$, the flow  $\dst\phi^t_{\tilde{F}}$ is well--defined, real--analytic with holomorphic extention to $D_{\frac{r}{4},s-4\s}$ and $C^\infty$ in $\o$ on $\rn$ for any $0\leq t\leq 1$, with
\beq{defPh}
\phi^t_{\tilde{F}}\colon D_{\frac{r}{4},s-4\s} \to D_{\frac{r}{2},s-3\s}
\eeq
and, setting $\wt\Phi\coloneqq\phi^1_{\tilde{F}}\eqqcolon (U,V)$, we have
\beq{PhUe}
\|U-\id\|_{\frac{r}{4},s-4\s,0}\leq \|\tilde{F}_x\|_{\frac{r}{2},s-3\s,\infty}\leq 2C_3\frac{\epsilon}{\a\s^{\bar{\n}}}
\eeq
\beq{PhVe}
\|V-\id\|_{\frac{r}{4},s-4\s,0}\leq \|\tilde{F}_y\|_{\frac{r}{2},s-3\s,\infty}\leq 4C_3\frac{\epsilon}{\a r \s^\t}\ .
\eeq
Moreover, since $R$ is affine in the variable $y$, then so is $F$ and then $\tilde{F}$ so that $\tilde{F}_y$ and $V$ do not depend on $y$, therefore the jacobian of $\Phi$ is of the form
\beq{jcbPh}
D\wt\Phi=\begin{pmatrix}
U_y & U_x\\
0 & V_x
\end{pmatrix},
\eeq
with, by using Cauchy's estimate, the following bounds
\begin{align}
\|U_y-\Id\|_{\frac{r}{8},s-5\s,0}&\leq  d\frac{8\|U-\id\|_{\frac{r}{4},s-4\s,0}}{r}\leq 16d C_3\frac{\epsilon}{\a r \s^{\bar{\n}}}\, ,\label{Uimp00}\\
\|U_x\|_{\frac{r}{8},s-5\s,0}&\leq  d\frac{\|U-\id\|_{\frac{r}{4},s-4\s,0}}{\s}\leq 2dC_3\frac{\epsilon}{\a \s^{{\bar{\n}}+1}}\, ,\label{Uimp01}\\
\|V_x-\Id\|_{\frac{r}{8},s-5\s,0}&\leq  d\frac{\|V-\id\|_{\frac{r}{4},s-4\s,0}}{\s}\leq 4dC_3\frac{\epsilon}{\a r \s^{\bar{\n}}}\, ,\label{Uimp02}\\
\|U_\o\|_{\frac{r}{4},s-4\s,0}&\leq \|\tilde{F}_{\o\th}\|_{r,s-3\s,0}\leq d\frac{\|\tilde{F}_{\o}\|_{r,s-2\s,0}}{\s}\leq \frac{dC_5}{4} \frac{\epsilon}{\a \s^{\bar{\n}} h}\, ,\label{Uimp03}\\
\|V_\o\|_{\frac{r}{4},s-4\s,0}&\leq  \|\tilde{F}_{\o y}\|_{\frac{r}{2},s-2\s,0}\leq d\frac{2\|\tilde{F}_{\o}\|_{r,s-2\s,0}}{r}\leq \frac{dC_5}{2} \frac{\epsilon}{\a r \s^\t h}\ .\label{Uimp04}
\end{align}
{\bf 6. New error term estimate.} To estimate $P_+^1$ (compare \eqref{pplus}), we need to estimate $\{R,F\}$. By Cauchy's estimate, we have
\beqano
\|\{R,F\}\|_{\frac{r}{2},s-3\s,\frac{h}{2}}& \leq & \dst\sum_{j=1}^d \|\dst R_{x_j}\|_{\frac{r}{2},s-3\s,h} \|\dst F_{y_j}\|_{\frac{r}{2},s-3\s,\frac{h}{2}}+\|\dst R_{y_j}\|_{\frac{r}{2},s-3\s,h} \|\dst F_{x_j}\|_{\frac{r}{2},s-3\s,\frac{h}{2}}\\
							   & \leq &\dst\sum_{j=1}^d \frac{\|\dst R\|_{\frac{r}{2},s-2\s,h}}{\s} \|\dst F_{y}\|_{\frac{r}{2},s-3\s,\frac{h}{2}}+ \frac{2\|\dst R\|_{r,s-3\s,h}}{r} \|\dst F_{x}\|_{\frac{r}{2},s-3\s,\frac{h}{2}}\\
							  & \leq & \dst\sum_{j=1}^d \frac{C_3}{2C_0}\frac{\epsilon}{\s}\cdot 4C_3\frac{\epsilon}{\a r \s^\t}+ \frac{2C_3}{2C_0}\frac{\epsilon}{r}\cdot 2C_3\frac{\epsilon}{\a \s^{\bar{\n}}}\\
							  &  =  & \dst \frac{4dC_3^2}{C_0}\frac{\epsilon^2}{\a r\s^{\bar{\n}}},
\eeqano
\beqano
\left\|\left\{[R],F\right\}\right\|_{\frac{r}{2},s-3\s,\frac{h}{2}}& = & \left\|[R]_y\cdot F_x\right\|_{\frac{r}{2},s-3\s,\frac{h}{2}}= \left\|[R_y]\cdot F_x\right\|_{\frac{r}{2},s-3\s,\frac{h}{2}}\\
&\leq & d \left\|[R_y]\right\|_{\frac{r}{2},s-3\s,h}\cdot \| F_x\|_{\frac{r}{2},s-3\s,\frac{h}{2}} \\
								 &\leq & d \|R_y\|_{\frac{r}{2},s-3\s,h}\cdot \| F_x\|_{\frac{r}{2},s-3\s,\frac{h}{2}}\\
								 &\leq & d\frac{2\|R\|_{r,s-3\s,h}}{r} \| F_x\|_{\frac{r}{2},s-3\s,\frac{h}{2}}\leq  \frac{2dC_3^2}{C_0}\frac{\epsilon^2}{\a r\s^{\bar{\n}}}
\eeqano
and therefore 
\beqano
\|P_+^0\|_{\eta r,s-5\s,\frac{h}{2}}&\overset{\eqref{poplus}}{=}& \left\|\dst\int_0^1 \{(1-t)[R]+tR,F\}\circ \phi^t_F\, dt \right\|_{\eta r,s-5\s,\frac{h}{2}}\\
						&\overset{\eqref{fthb}+\eqref{fib}}{\leq}&\dst\int_0^1 \left\| \{(1-t)[R]+tR,F\} \right\|_{2\eta r,s-4\s,\frac{h}{2}}\, dt \\
						&\leq & \dst\int_0^1 \left\| \{(1-t)[R]+tR,F\} \right\|_{\frac{r}{2},s-3\s,\frac{h}{2}}\, dt\\
						&\leq & \frac{3 dC_3^2}{C_0}\frac{\epsilon^2}{\a r\s^{\bar{\n}}}\, .
\eeqano
Hence
\begin{eqnarray}
\|P_+^1\|_{\eta r,s-5\s,\frac{h}{2}} &\overset{\eqref{pplus}}{=} & \left\|P_+^0 + (P-Q)\circ \phi^1_F + (Q-R)\circ \phi^1_F\right\|_{\eta r,s-5\s,\frac{h}{2}}\nonumber\\
                       &\overset{\eqref{fthb}+\eqref{fib}}{\leq}& \left\|P_+^0\right\|_{\eta r,s-5\s,\frac{h}{2}} + \|P-Q\|_{2\eta r,s-4\s,\frac{h}{2}} + \|Q-R\|_{2\eta r,s-4\s,\frac{h}{2}}\nonumber\\
                       &\leq & \frac{3dC_3^2}{C_0}\frac{\epsilon^2}{\a r\s^{\bar{\n}}}+\frac{32d^2}{9}\eta^2 \epsilon+ 4^d(d+1)C_2 \k^d \ex^{-\k\s}\epsilon\nonumber\\
                       &\leq & \max\left(\frac{3dC_3^2}{C_0}, 4^d(d+1)C_2\right)\left(\frac{\epsilon^2}{\a r\s^{\bar{\n}}}+\left(\eta^2 + \k^d \ex^{-\k\s}\right)\epsilon\right)\nonumber\\
                       & =   &\frac{\sqrt{C_{10}}}{3\cdot 2^{{\bar{\n}}}}\left(\frac{\epsilon^2}{\a r\s^{\bar{\n}}}+\left(\eta^2 + \k^d \ex^{-\k\s}\right)\epsilon\right) .\label{ppluse}
\end{eqnarray}
{\bf 7. Transforming frequencies.} In view of \eqref{nplus}, we need to invert the map 
\[\rho\colon \o \mapsto \o+v(\o)=\o+[R_y].\]
But we have
\[\|\rho-\id\|_{h}=\|\dst[R_y]\|_h\leq \dst\|R_y\|_{\frac{r}{2},s-\s,h}\leq \frac{C_3}{C_0}\frac{\epsilon}{r}\overset{(b)}{\leq}\frac{C_3}{C_0C_6}h\leq \frac{h}{64\sqrt{3}}\leq \frac{h}{4}\]
Therefore, we apply Lemma~\ref{inv}  and get a real analytic map $\tilde{\f}\colon \O_{\a,\frac{h}{4}}\to \O_{\a,\frac{h}{2}}$, inverse of $\rho$ and satisfies
\[\|\tilde{\f}-\id\|_{\frac{h}{4}},\, \frac{h}{4}\|D\tilde{\f}-\Id\|_{\frac{h}{4}}\leq \frac{C_3}{C_0}\frac{\epsilon}{r}\]
Now we extend $\tilde{\f}$: by Lemma~\ref{cutoff}, there exists a cut--off $\chi_2\in C(\cn)\cap C^{\infty}(\rn)$ with $0\leq \chi_2\leq 1,\ supp\chi_2 \subset \O_{\a,\frac{h}{4}},\ \chi_2\equiv 1$ on $\O_{\a,\frac{h}{8}}$. Then let $\f\coloneqq \id+(\tilde{\f}-\id)\chi_2$ on $\O_{\a,\frac{h}{8}}$ and $\f\coloneqq \id$ on $\cn\setminus  \O_{\a,\frac{h}{4}}$. Thus, setting
\begin{align*}
N_+(y,\o_+)&\coloneqq N_+^1(y,\tilde{\f}(\o_+))=e^1_+\circ\tilde{\f}(\o_+)+ \o_+\cdot y,\\
\mathcal{F}&\coloneqq (\Phi,\f)\coloneqq \left(\wt\Phi\circ(\pi_1,\pi_2,\f\circ\pi_3),\f\circ\pi_3\right),\\
P_+&\coloneqq P_+^1\circ \mathcal{F}_{|D_{\frac{r}{4},s-4\s}\times \O_{\a,\frac{h}{8}}}\ ,
\end{align*}
 we have $H\circ \mathcal{F}=N_+ + P_+$ on $D_{\frac{r}{4},s-4\s}\times \O_{\a,\frac{h}{8}}$, with all the required properties.
Moreover,
\beq{tataf0}
\|\f-\id\|_\infty \leq \|\tilde{\f}-\id\|_{\frac{h}{4}} \leq \frac{C_3}{C_0} \frac{\epsilon}{r}\overset{(b)}{\leq} \frac{h}{64\sqrt{3}}
\eeq
and
\begin{align}
\|D\f-\Id\|_0 & \leq  \|D\tilde{\f}-\Id\|_{\frac{h}{4}}+ \|\tilde{\f}-\id\|_{\frac{h}{4}} \|D_\o \chi_2\|_0 \overset{\eqref{cutof}}{\leq} \frac{4C_3}{C_0} \frac{\epsilon}{hr}+\frac{C_3}{C_0} \frac{\epsilon}{r}\cdot \frac{24C_1}{h}\nonumber\\
           &=  \frac{C_5}{4C_0}\frac{\epsilon}{hr}\overset{(b)}{\leq} \frac{ C_5}{4C_0C_6}<\frac{1}{2}\ .\label{jacobfi}
\end{align}
So, in particular, $\f_{|\rn}$ is $C^\infty$--diffeomorphism from $\rn$ onto itself and since\footnote{because $\O_\a\subset \O,\  \dist(\O_\a,\partial\O)\geq \a$ and $h$ will be chosen (just below) in such away that $h\leq \frac{\a}{2}$ so that $\O_{\a,\frac{h}{4}}\cap\rn\subset \O$.} $\f=\id$ outside of $\O$ then $\f_{|\O}$ is $C^\infty$--diffeomorphism from $\O$ onto itself.\\
{\bf 8. Estimating $\Phi$.} By \eqref{PhUe}, \eqref{PhVe} and \eqref{tataf0}, we have\footnote{Recall that $\s_0\ge \s$.}
\begin{align*}
\|\bar W(\mathcal{F}-\id)\|_{\frac{r}{4},s-4\s,0}&\le \dst\max\left(6C_3 \frac{\epsilon}{\a r \s^{\bar{\n}}}\left(\frac{\s_0}{\s}\right)^{2\bar{\n}-1},\  \frac{C_3}{C_0} \frac{\epsilon}{rh} \right)\\
    &\le \left(\frac{\s_0}{\s}\right)^{2\bar{\n}-1}\cdot\dst\max\left(6C_3 \frac{\epsilon}{\a r \s^{\bar{\n}}},\  \frac{C_3}{C_0} \frac{\epsilon}{rh} \right)\ .
\end{align*}
Now\footnote{Recall that $D\Phi$ denotes the Jacobian of $\Phi$ \wrt $(y,x)$, $D\f$ the Jacobian of $\f$ \wrt $\o$ and $\Phi-\id$ means $\Phi-(\pi_1,\pi_2)$.}, since
$$
\dpr_y\Phi=\dpr_y\wt\Phi|_{(\pi_1,\pi_2;\f\circ\pi_3)},\quad\dpr_x\Phi=\dpr_x\wt\Phi|_{(\pi_1,\pi_2;\f\circ\pi_3)},\quad \dpr_\o(\Phi-\id)=\dpr_\o\wt\Phi\big|_{(\pi_1,\pi_2;\f\circ\pi_3)}\circ D\f|_{\pi_3},
$$
$$\bar W(D\mathcal{F}-\Id)\bar W^{-1}= \begin{pmatrix}
\wt U_y-\Id & \left(\frac{\s}{\s_0}\right)^{2\bar{\n}-1}\frac{\s}{r}\wt U_x & \frac{h}{r}\wt U_\o \\ 
	0 			    & \wt V_x-\Id & \left(\frac{\s_0}{\s}\right)^{2\bar{\n}-1}\frac{h}{\s}\wt V_\o \\
  0					&  0		  & \f_\o-\Id
\end{pmatrix}
\ ,
$$
with $\Phi\eqqcolon (\wt U,\wt V)$, then by \eqref{Uimp00}--\eqref{Uimp04} and \eqref{jacobfi}, 
 we have
\begin{align*}
\left\|\bar W(D\mathcal{F}-\Id)\bar W^{-1} \right\|_{\frac{r}{8} ,s-5\s,0}&\le \dst\max\left\{ \left(16dC_3+2dC_3+\frac{dC_5}{4}\right)\frac{\epsilon}{\a r \s^{\bar{\n}}},\right.\\
			&\qquad \qquad \left.\left(4dC_3+\frac{dC_5}{2} \right) \frac{\epsilon}{\a r \s^{\bar{\n}}}\left(\frac{\s_0}{\s}\right)^{2\bar{\n}-1},\ \frac{C_5}{C_0}\frac{\epsilon}{rh} \right\}\\
			&\leq \dst\max\left(dC_5 \frac{\epsilon}{\a r \s^{\bar{\n}}}\left(\frac{\s_0}{\s}\right)^{2\bar{\n}-1},\ \frac{C_5}{C_0}\frac{\epsilon}{rh} \right)\\
			&\le \left(\frac{\s_0}{\s}\right)^{2\bar{\n}-1}\cdot\dst\max\left(dC_5 \frac{\epsilon}{\a r \s^{\bar{\n}}},\ \frac{C_5}{C_0}\frac{\epsilon}{rh} \right)\ .
\end{align*}
 Since $D_{\frac{r}{8},s-5\s}\supseteq D_{\eta r,s-5\s}$, then the estimates on $\Phi$ are proven.
\qed

\subsubsection{Iteration of the KAM step}
Since we are going to iterate the \emph{KAM step} infinitely many times, we need to choose the sequences $r_j,\, s_j,\, h_j,\, \k_j,\, \s_j,\, \eta_j$ conveniently so that at each step, all the assumptions in the \emph{KAM step} hold. See \cite{JP}, for details on how those sequences are choosen. First, we set up the sequences, then we prove that at each step they meet all the assumptions in \emph{KAM step} and then we prove the iterative lemma.\\
Let then $\m\coloneqq \frac{3}{2}$ and (recall \textsuperscript{\ref{footn5}})\\
$
\left\{
\begin{array}{ll}
0< s_0 \leq 1 \\ \\
s_{j+1} = s_j - 5\s_j
\end{array}
\right.,\quad
$
$\left\{
\begin{array}{rl}
\s_0 &= \frac{s_0}{20}\\ \\
\s_{j+1} &=\frac{\s_j}{2}
\end{array}
\right.,\quad
$
$\left\{
\begin{array}{rl}
E_0&\leq  20^{\n}c \s_0^{\n-{\bar{\n}}} \\ \\
E_{j+1} & =C_{10}^{\m-1}E_j^\m
\end{array}
\right.,\vspace{0.1cm}\\ \ \\
$
$\left\{
\begin{array}{rl}
\k_0 & =\left[-\frac{40\log\vae-1}{s_0}\right]\\ \\
\k_{j+1} & =4 \k_j
\end{array}
\right.,\quad
$
$\left\{
\begin{array}{rl}
\a C\vae& \leq h_0\leq \frac{\a}{2\k_0^{{\bar{\n}}}}\\ \\
h_{j+1} & =\frac{h_j}{4^{\bar{\n}}}
\end{array}
\right. \mbox{ and }
$\vspace{0.1cm}
\ \\
$\left\{
\begin{array}{rl}
0<r_0\leq 1\\ \\
\eta_j^2=E_j,\, r_{j+1}=\eta_j r_j \mbox{ and } \epsilon_j=\a E_j r_j \s_j^{{\bar{\n}}}
\end{array}
\right..\\ \ \\
$
Thus the following hold
\lem{kamass}
For any $j\in \natural$,\\
$\begin{array}{rrl}
(i) & \quad  \epsilon_j & \leq \frac{1}{C_4} \a \eta_j r_j \s_j^{{\bar{\n}}}\\
(ii) & \quad \epsilon_j &\leq \frac{1}{C_6} h_j r_j\\
(iii) & \quad h_j &\leq \frac{\a}{2 \k_j^{{\bar{\n}}}}\vspace{0.1cm}\\
(iv) & \quad \epsilon_{j+1} &\geq \frac{\sqrt{C_{10}}}{3\cdot 2^{{\bar{\n}}}}\left(\epsilon_j E_j+(\eta_j^2+\k_j^d \ex^{-\k_j \s_j})\epsilon_j\right)\vspace{0.05cm}\\
(v) & \quad \k_j\s_j & > d-1,\quad 0<\eta_j <1/8,\quad \mbox{and} \quad 0<\s_j < s_j/10
\end{array}$
\elem
\proof
$(i)$ As $E_j$ is decreasing (super--exponentially) and $(i)\Leftrightarrow C_4^2 E_j\leq  1$ then it is enough to check it for $j=0$. But by definitions, we have
\[E_0\leq \frac{1}{4^{\bar{\n}} C_{10}}\leq \frac{1}{C_4^2}.\]
$(ii)+(iii)$  By definitions, it follows\footnote{See Appendix~\ref{app:A}.} 
\[\k_0^{{\bar{\n}}} \s_0^{{\bar{\n}}} \ex^{-\k_0\s_0}\leq E_0 \s_0^{{\bar{\n}}}= \frac{1}{\a} \frac{\epsilon_0}{ r_0} \leq \frac{h_0}{\a C_6}\leq \frac{1}{2C_6  \k_0^{{\bar{\n}}}}.\]
Now let $j\in \natural$ and assume
\[\k_j^{\bar{\n}} \s_j^{{\bar{\n}}} \ex^{-\k_j\s_j}\leq E_j \s_j^{{\bar{\n}}}= \frac{1}{\a} \frac{\epsilon_j}{ r_j} \leq \frac{h_j}{\a C_6}\leq \frac{1}{2 C_6 \k_j^{\bar{\n}}}.\]
Then, by using the above definitions we get
\beqano
\k_{j+1}^{\bar{\n}} \s_{j+1}^{\bar{\n}} \ex^{-\k_{j+1}\s_{j+1}}&=& 4^{\bar{\n}} \k_j^{\bar{\n}} \s_{j+1}^{\bar{\n}} \ex^{-2 \k_j\s_j}\\
         &\leq & \left(\k_j^{\bar{\n}}  \ex^{-\k_j\s_j}\right)^2\s_{j+1}^{\bar{\n}}\leq E_j^2 \s_{j+1}^{\bar{\n}}\\
         &\leq & C_{10}^{\m-1}E_j^\m \s_{j+1}^{\bar{\n}} = E_{j+1} \s_{j+1}^{\bar{\n}}=\frac{1}{\a}\frac{\epsilon_{j+1}}{r_{j+1}}\\
         & =   & \sqrt{C_{10} E_j}E_j \frac{\s_j^{\bar{\n}}}{2^{\bar{\n}}}\leq \frac{1}{4^{\bar{\n}}}E_j \s_j^{\bar{\n}}\\
         &\leq & \frac{h_j}{\a 4^{\bar{\n}}C_6}= \frac{h_{j+1}}{\a C_6}\\
         &\leq & \frac{1}{2C_6\cdot 4^{\bar{\n}} \k_j^{\bar{\n}}}= \frac{1}{2 C_6 \k_{j+1}^{\bar{\n}}}
\eeqano
which ends the proof of $(ii)$ and $(iii)$.\\
$(iv)$ We have 
\beqano
\frac{\sqrt{C_{10}}}{3\cdot 2^{\bar{\n}}}\left(\epsilon_j E_j+(\eta_j^2+\k_j^d \ex^{-\k_j \s_j})\epsilon_j\right) &\leq  &\frac{\sqrt{C_{10}}}{3\cdot 2^{\bar{\n}}}\left(\epsilon_j E_j+(E_j+E_j)\epsilon_j\right)\\
& = &\frac{\sqrt{C_{10}}}{ 2^{\bar{\n}}}\epsilon_j E_j =\frac{\sqrt{C_{10}}}{2^{\bar{\n}}} \a E_j^2 r_j \s_j^{\bar{\n}} \\
& = &\a C_{10}^{\m-1}E_j^\m \sqrt{E_j}\frac{r_{j+1}}{\eta_j} \left(\frac{\s_j}{2}\right)^{\bar{\n}}\\
& = & \a r_{j+1} \s_{j+1}^{\bar{\n}} E_{j+1} = \epsilon_{j+1}.
\eeqano
$(v)$ We have
\[k_j\s_j=2^j \k_0 \s_0\geq \k_0 \s_0 > d-1,\]
\[0 < \eta_j^2 = E_j = \frac{1}{C_{10}}\dst (C_{10} E_0)^{\m^j}\leq \frac{1}{C_{10} 4^{\bar{\n}\m^j}}< \frac{1}{ 4^3} ,\]
\[s_j-10\s_j = \frac{s_0}{2}>0.\]
\qed
Now we arrive to the iterative lemma. Given $j\in \natural$, let\footnote{Notice that $\s_0/\s_j=2^j$.} 
\[D_j\coloneqq D_{r_j,s_j},\quad O_j\coloneqq \O_{\a,h_j},\quad
 \bar{W}_j\coloneqq \diag\left(r_j^{-1}\uno_d,\ 2^{j(2\bar{\n}-1)}\s_j^{-1}\uno_d,\ h_j^{-1}\uno_d\right),\] 
\[
u_j\coloneqq \frac{\epsilon_j}{r_j h_j},\quad v_j \coloneqq 
C_{10} E_j.\]
Then, by the definitions, we have $u_{j+1}=2^{\bar{\n}} u_j\sqrt{v_j} $ and $v_{j+1}=v_j^\m$. Thus, 
\[v_j=\dst v_0^{\m^j},\quad v_0=C_{10} E_0\leq 4^{-\bar{\n}}, \quad u_j=2^{j\bar{\n}}u_0 \dst v_0^{\m^j-1},\quad u_0\leq  \frac{1}{C_6}\quad \mbox{and} \quad \frac{u_j}{h_j}=\frac{u_0}{h_0}2^{3j\bar{\n}}\dst v_0^{\m^j-1}.\]
In particular for any $j\geq 0$,
\beq{eqUjSuHj}u_j\leq 2^{j\bar{\n}}\cdot u_0\cdot 4^{-{\bar{\n}}(\m^j-1)},\quad \frac{u_j}{h_j}\leq 2^{3j\bar{\n}}\cdot \frac{u_0}{h_0}\cdot 4^{-\bar{\n}(\m^j-1)}.
\eeq 
\lem{IterLemmTeo1}
Suppose $H_0\coloneqq N+ P_0$ is real analytic on $D_0\times O_0$ with
\[\|P_0\|_{r_0,s_0,h_0}\leq \epsilon_0=\a E_0 r_0 \s_0^{\bar{\n}}.\]
Then for each $j\in \natural$, there exist a normal form $N_j$ and a transformation \\
$\mathcal{F}^j=(\Phi^j;\f^j)\coloneqq\mathcal{F}_0\circ \cdots \circ \mathcal{F}_{j-1}\colon D_j\times \rn\to D_0\times \rn$ such that
\begin{itemize}
\item[(i)] $\f^j\colon \rn\to \rn$ is a $C^\infty$--diffeomorphism with $\f^j\equiv \id$ on $\O$ and $\Phi^j\colon D_j\times \rn\to D_0$ is a ($\o$--) family of real--anatylitic, 
symplectic transformations parametrized over $\rn$ and $C^\infty$ in $\o$;
\item[(ii)] $\mathcal{F}^j$ is Lipschitz--continuous in $\o$ with
\beq{iterlip}
\|\bar{W}_0\left(\mathcal{F}^j-\id\right)\|_{L,\rn}\leq C_9 \frac{\epsilon_0}{r_0 h_0^2},
\eeq
uniformly on $D_j\times \rn$.
\item[(iii)] The restriction $\tilde{\mathcal{F}}^j\coloneqq \mathcal{F}^j_{|D_j\times O_j}$ is real--analytic with holomorphic extension to $D_j\times O_j$ for each given $\o\in O_j$ 
and satisfies $\tilde{\mathcal{F}}^j\colon D_j\times O_j \to D_0\times O_0$, $H\circ \tilde{\mathcal{F}}^j=N_j+P_j$ with 
\[\|P_j\|_{r_j,s_j,h_j}\leq \epsilon_j=\a E_j r_j \s_j^{\bar{\n}}.\]
\end{itemize}
Furthermore
\beq{iterCau0}
\|\bar{W}_0(\mathcal{F}^{j+1}-\mathcal{F}^j)\|_{r_{j+1},s_{j+1},0}\leq C_7 \frac{\epsilon_0}{r_0 h_0}\cdot 2^{-{\bar{\n}}(2\m^j+3j-2)+j}\ .
\eeq 

\elem
\proof
For $j=0$ we take $\mathcal{F}^1=\mathcal{F}_0=\id,\, N_1=N, P_1=P_0$ and we are done.\\
Next we pick $j\geq 0$ and we assume that it holds at the step $j$. Then we have to check it for the step $j+1$. But, thanks to lemma~\ref{kamass}, we can apply the \emph{KAM step} to $H_j$ to get a transformation $\mathcal{F}_j=(\Phi_j;\f_j)\colon D_{j+1}\times\rn\to D_j\times \rn$ for which every properties in \emph{KAM step} hold. So, its restriction 
\[\tilde{\mathcal{F}}_j\coloneqq\mathcal{F}_{j|D_{j+1}\times O_{j+1}}\colon D_{j+1}\times O_{j+1} \to D_j\times O_j\]
and there exists a normal form $N_{j+1}$ such that $H_{j+1}=H_j \circ \tilde{\mathcal{F}}_j=N_{j+1}+ P_{j+1}$ with
\[\|P_{j+1}\|_{r_{j+1},s_{j+1},h_{j+1}} \leq \frac{\sqrt{C_{10}}}{3\cdot 2^{{\bar{\n}}}}\left(\epsilon_j E_j+\left(\eta_j^2+\k_j^d \ex^{-\k_j \s_j}\right)\epsilon_j\right).\]
Then we apply $(iv)$ of lemma~\ref{kamass} to obtain
\[\|P_{j+1}\|_{r_{j+1},s_{j+1},h_{j+1}} \leq \epsilon_{j+1}=\a E_{j+1} r_{j+1} \s_{j+1}^{\bar{\n}}.\]
Therefore 
\[\mathcal{F}^{j+1}=\mathcal{F}_0\circ \cdots \circ \mathcal{F}_{j}=\left(\Phi^j\circ(\Phi_j,\f_j);\f^j\circ\f_j\right) \colon D_{j+1}\times \rn \to  D_0\times \rn\]
is a transformation such that $H\circ \tilde{\mathcal{F}}^{j+1}=H_j\circ \tilde{\mathcal{F}}_{j}=N_{j+1}+ P_{j+1}$ with all the required properties in $(i)$ and $(iii)$.\\
It remains the estimates on $\mathcal{F}^j$. By \eqref{kamest1} and \eqref{kamest2} 
we have\footnote{Notice that $\frac{h_j}{\a \s_j^{\bar{\n}}}=\frac{h_0}{\a \s_0^{\bar{\n}}}\cdot\frac{1}{2^{j{\bar{\n}}}}\leq \frac{h_0}{\a \s_0^{\bar{\n}}}\leq \frac{C_5}{4C_0C_4}$.\label{footn13}}
\beqa{}
& &\|\bar{W}_j(\mathcal{F}_j-\id)\|_{r_{j+1},s_{j+1},0}\ ,\quad\|\bar{W}_j(D\mathcal{F}_j-\Id)\bar{W}^{-1}_j\|_{r_{j+1},s_{j+1},0}\leq  \nonumber \\
 &\le &\left(\frac{\s_0}{\s_j}\right)^{2\bar{\n}-1}\cdot \dst\max\left(dC_5 \frac{\epsilon_j}{\a r_j \s_j^{\bar{\n}}},\frac{C_5}{C_0}\frac{\epsilon_j}{r_jh_j} \right)\nonumber\\
 &\le &
 2^{-j(2\bar{\n}-1)}\frac{dC_5^2}{4C_0C_4}u_j\label{proofCau0}\\
 &\leby{eqUjSuHj}&  \frac{dC_5^2}{4C_0C_4} \frac{\epsilon_0}{r_0 h_0}\cdot 2^{-{\bar{\n}}(2\m^j+j-2)+j}\ . \label{proofCau}
\eeqa
Thus
\beqano
\|\bar{W}_0(\mathcal{F}^{j+1}-\mathcal{F}^j)\|_{r_{j+1},s_{j+1},0}&= & \|\bar{W}_0(\mathcal{F}^j\circ \mathcal{F}_j-\mathcal{F}^j)\|_{r_{j+1},s_{j+1},0}\\
				&\leq & \|\bar{W}_0 D\mathcal{F}^j\bar{W}_j^{-1}\|_{r_j,s_j,0}\cdot \|\bar{W}_j(\mathcal{F}_j-\id)\|_{r_{j+1},s_{j+1},0}\\
				&\overset{\eqref{proofCau}}{\leq} & \frac{dC_5^2}{4C_0C_4} \frac{\epsilon_0}{r_0 h_0}\cdot 2^{-{\bar{\n}}(2\m^j+j-2)+j}\cdot \|\bar{W}_0 D\mathcal{F}^j\bar{W}_j^{-1}\|_{r_j,s_j,0}.
\eeqano
Next, we need to bound $\|\bar{W}_0 D\mathcal{F}^j\bar{W}_j^{-1}\|$ uniformly on $D_j\times \rn$. But for any $j\geq 0$, we have\footnote{Recall that any $j\geq 0$, $\eta_j\le \sqrt{c}\le 1/\sqrt{4^{\bar{\n}}C_{10}}<1/(3\cdot 4^{\bar{\n}}) $.}
\begin{align}
\|\bar{W}_{j} \bar{W}_{j+1}^{-1}\|&=\left\|diag\left(\frac{r_{j+1}}{r_{j}}\uno_d,\su{2^{2\bar{\n}-1}}\frac{\s_{j+1}}{\s_{j}}\uno_d,\frac{h_{j+1}}{h_{j}}\uno_d\right)\right\|\nonumber\\
    &=\max\left(\frac{r_{j+1}}{r_{j}},\su{2^{2\bar{\n}-1}}\frac{\s_{j+1}}{\s_{j}},\frac{h_{j+1}}{h_{j}}\right)\nonumber\\
    &=\max\left(\eta_j,\su{4^{\bar{\n}}}\right)\nonumber\\
    &= \frac{1}{4^{\bar{\n}}}.\label{proofNormInf1}
\end{align}
and
\beqano
\bar{W}_0 D\mathcal{F}^j \bar{W}_j^{-1} &=&\bar{W}_0 D\mathcal{F}_0\circ \cdots \circ D\mathcal{F}_{j-1}\bar{W}_j^{-1}\\
			   &=& \left(\bar{W}_0 D\mathcal{F}_0\bar{W}_0^{-1}\right)\left(\bar{W}_0 \bar{W}_1^{-1}\right) \cdots   \left(\bar{W}_{j-1} D\mathcal{F}_{j-1}\bar{W}_{j-1}^{-1}\right)\left(\bar{W}_{j-1} \bar{W}_j^{-1}\right)\ ,
\eeqano
so that,
\beqano
\|\bar{W}_0 D\mathcal{F}^j\bar{W}_j^{-1}\|&\leq & \|\bar{W}_0 D\mathcal{F}_0\bar{W}_0^{-1}\| \|\bar{W}_0 \bar{W}_1^{-1}\| \cdots   \|\bar{W}_{j-1} D\mathcal{F}_{j-1}\bar{W}_{j-1}^{-1}\| \|\bar{W}_{j-1} \bar{W}_j^{-1}\|\\
				&\overset{\eqref{proofNormInf1}}{\leq} & \|\bar{W}_0 D\mathcal{F}_0\bar{W}_0^{-1}\|  \cdots   \|\bar{W}_{j-1} D\mathcal{F}_{j-1}\bar{W}_{j-1}^{-1}\|\\
				&\overset{\eqref{proofCau}}{\leq} & 4^{-j\bar{\n}}\dst\prod_{k=0}^{j-1}\left(1+\frac{dC_5^2}{4C_0C_4} \frac{\epsilon_0}{r_0 h_0}\cdot 2^{-{\bar{\n}}(2\m^k+k-2)+k} \right)\\
				&\leq & 4^{-j\bar{\n}}\dst\prod_{k=0}^{\infty}\left(1+\frac{dC_5^2}{4C_0C_4C_6} \cdot 2^{-{\bar{\n}}(2\m^k+k-2)+k} \right)\\ 
				&\leq &4^{-j\bar{\n}}\exp\left(\frac{dC_5^2}{4C_0C_4C_6} \dst\sum_{k=0}^{\infty}  2^{-{\bar{\n}}(2\m^k+k-2)+k} \right)\\
				& =  & 4^{-j\bar{\n}}\frac{4C_0C_4C_7}{dC_5^2}\ . 
\eeqano
$$
C_6\coloneqq \frac{C_5}{C_0}
$$
$$
C_7\coloneqq \frac{dC_5^2}{4C_0C_4}\dst\exp\left(\frac{dC_5^2}{4C_0C_4C_6} \sum_{j=0}^{\infty}  2^{-{\bar{\n}}(2\m^j+j-2)+j} \right)
$$
Therefore,
\[\|\bar{W}_0(\mathcal{F}^{j+1}-\mathcal{F}^j)\|_{r_{j+1},s_{j+1},0}\leq C_7 \frac{\epsilon_0}{r_0 h_0}\cdot 2^{-{\bar{\n}}(2\m^j+3j-2)+j}\ .\]
Finally, using again \eqref{proofCau0} and \eqref{proofNormInf1}, we get
\begin{align*}
\|\bar{W}_0 \left(D\mathcal{F}^{j+1}-\Id\right)\bar{W}_{j+1}^{-1}\|&= \|\bar{W}_0 \left(D\mathcal{F}^j\circ D\mathcal{F}_{j}-\Id\right)\bar{W}_{j+1}^{-1}\|\\
		&= \|\bar{W}_0 \left(D\mathcal{F}^j-\Id\right)\bar{W}_j^{-1}\circ \bar{W}_{j} D\mathcal{F}_{j}\bar{W}_{j}^{-1}\circ\bar{W}_{j}\bar{W}_{j+1}^{-1}+\\
		&\qquad\bar{W}_0\left(D\mathcal{F}_{j}-\Id\right)\bar{W}_{j+1}^{-1}\|\\
		&\le \|\bar{W}_0 \left(D\mathcal{F}^j-\Id\right)\bar{W}_j^{-1}\| \|\bar{W}_{j} D\mathcal{F}_{j}\bar{W}_{j}^{-1}\| \|\bar{W}_{j}\bar{W}_{j+1}^{-1}\|+\\
		&\quad\   \|\left(\bar{W}_0\bar{W}_{1}^{-1}\right)\cdots\left(\bar{W}_{j-1}\bar{W}_{j}^{-1}\right)\bar{W}_{j}\left(D\mathcal{F}_{j}-\Id\right)\bar{W}_{j}^{-1}\circ \bar{W}_{j}\bar{W}_{j+1}^{-1}\|\\
		&\le\su{4^{\bar{\n}}}\|\bar{W}_0 \left(D\mathcal{F}^j-\Id\right)\bar{W}_j^{-1}\|\left(1+\su{2^{j(2\bar{\n}-1)}}\frac{dC_5^2}{4C_0C_4}u_{j} \right)+\\
		&\qquad \su{4^{(j+1)\bar{\n}}}\su{2^{j(2\bar{\n}-1)}}\frac{dC_5^2}{4C_0C_4} u_{j}\ .
\end{align*}
Therefore, letting
$$
w_j\coloneqq \log\left(\|\bar{W}_0 \left(D\mathcal{F}^j-\Id\right)\bar{W}_j^{-1}\|+\su{4^{j\bar{\n}}}\right),\quad z_j\coloneqq \su{2^{j(2\bar{\n}-1)}}\frac{dC_5^2}{4C_0C_4} u_j,\quad j\ge 0,
$$
we get, for any $j\ge 1$,
$$
w_{j+1}\le w_j+\log\left(\su{4^{\bar{\n}}}\left(1+ z_{j}\right)\right)\quad\mbox{and}\quad w_1=\log\left(\su{4^{\bar{\n}}}\right)\ ,
$$
so that, for any $j\ge 1$, 
\begin{align*}
 w_{j} & \le \dst\sum_{k=0}^{j-1}\log\left(\su{4^{\bar{\n}}}\left(1+ z_{k}\right)\right)\\
	&\le \log\left(\su{4^{j\bar{\n}}}\dst\prod_{k=0}^{j-1}\left(1+ z_{k}\right)\right)\\
	&\le \log\left(\su{4^{j\bar{\n}}}\right)+\log\left(\dst\prod_{k=0}^\infty\left(1+ z_{k}\right)\right)\\
	&\le \log\left(\su{4^{j\bar{\n}}}\right)+\dst\sum_{k=0}^\infty  z_{k}\\
	&\le \log\left(\su{4^{j\bar{\n}}}\right)+C_8 u_0\ ,
\end{align*}
	%
\ie
$$
\|\bar{W}_0 \left(D\mathcal{F}^{j}-\Id\right)\bar{W}_{j}^{-1}\|\le \su{4^{j\bar{\n}}}\left(\ex^{C_8 u_0}-1\right).
$$
In particular, for any $j\ge 1$,
\begin{align*}
h_{j}\|\bar{W}_0 \left(\dpr_\o\mathcal{F}^{j}-\Id\right)\|&\le\|\bar{W}_0 \left(D\mathcal{F}^{j}-\Id\right)\bar{W}_{j}^{-1}\|\\
       &\le \su{4^{j\bar{\n}}}\left(\ex^{C_8 u_0}-1\right)\ ,
\end{align*}
\ie\footnote{Recall that $u_0\le 1/C_6$ and $\ex^a-1\le a\ex^a,\ \forall a\ge 0 $.}
$$
h_0\|\bar{W}_0 \left(\dpr_\o\mathcal{F}^{j}-\Id\right)\|\le \ex^{C_8 u_0}-1\le C_8 u_0\ex^{C_8 u_0}\le C_8 \ex^{C_8/C_6}u_0=C_9 u_0\ ,
$$
\ie
$$
\|\bar{W}_0 \left(\mathcal{F}^{j}-\id\right)\|_{L,\rn}\le  C_9\frac{\epsilon_0}{r_0 h_0^2}\quad \mbox{uniformly on}\quad D_j\times\rn\ .
$$
\qed
\subsubsection{Deduction of Theorem~\ref{teo1}}
We set $P_0=P,\quad s_0=s,\quad r_0=r,\quad h_0=h,\quad \mbox{and} \quad \epsilon_0=\epsilon=\|P\|_{r,s,h}$;
 thus Lemma~\ref{IterLemmTeo1} applies. Hence, by \eqref{iterCau0}, $(\mathcal{F}_j)_j$ is a Cauchy sequence and therefore converges uniformly to some $\mathcal{F}=(\Phi,\f)$ on
\[\dst\bigcap_{j\geq 0}D_j\times \rn=T_*\times \rn,\mbox{ where } T_*\coloneqq \{0\}\times \torus^d_{\frac{s}{2}},\]
with the map $x \mapsto \Phi(0,x, \o) $ real analytic on $\dst\torus^d_{\frac{s}{2}}$ for each given $\o \in \rn$ (by Weierstrass's theorem) and for any $j\geq 1$,
\beqano
\|\bar{W}_0(\mathcal{F}^j-\id)\|_{r_j,s_j,h_j} &\leq & \|\bar{W}_0(\mathcal{F}^j-\mathcal{F}^{j-1})\|_{r_j,s_j,h_j}+\cdots + \|\bar{W}_0(\mathcal{F}^2-\mathcal{F}^0)\|_{r_0,s_0,h_0}\\
							  &\leq &   C_{10} u_0,
\eeqano
$$
C_{10}\coloneqq  C_7\dst\sum_{j=0}^\infty 2^{-{\bar{\n}}(2\m^j+3j-2)+j}
$$
and letting $j\to \infty$, we get, uniformly on $T_*\times \rn$,
\beq{estthma1gr}
\|\bar{W}_0(\mathcal{F}-\id)\| \leq  C_{10} \frac{\epsilon}{ r h}. 
\eeq
Moreover, by letting $j\to \infty$ in \eqref{iterlip}, we get, uniformly on $T_*\times \rn$,
\beq{estthma2}
\|\bar{W}_0(\mathcal{F}-\id)\|_{L,\rn} \leq  C_9 \frac{\epsilon}{r h}.
\eeq
Let's prove that $\f$ is a lipeomorphism from $\O$ onto itself. 
Indeed, for any $j\ge 0$,
\begin{align*}
\|D\f^{j+1}-\Id\|_0 +1&= \|(D\f^{j}-\Id)D\f_j+ (D\f_j-\Id)|_0+1\\
    &\le \|D\f^{j}-\Id\|_0\left( \|D\f_j-\Id\|_0+1\right)+ \|D\f_j-\Id)\|_0+1\\
    &= \left(\|D\f^{j}-\Id\|_0+1\right)\left( \|D\f_j-\Id\|_0+1\right)\ ,
\end{align*}
so that\footnote{Recall that $\f^1=\f_0=\id$, so that $\|D\f^{1}-\Id\|_0=0$.}
\begin{align*}
\|D\f^{j+1}-\Id\|_0&\le -1+ \left(\|D\f^{1}-\Id\|_0+1\right)\dst\prod_{k=0}^j\left( \|D\f_k-\Id\|_0+1\right)\\
    &\leby{kamest3} -1+ \dst\prod_{k=0}^j\left( 1+\frac{C_5}{2C_0}u_k\right)\\
    &\le \exp\left( \frac{C_5}{2C_0}\dst\sum_{k=0}^\infty u_k\right)-1\\
    &\leby{eqUjSuHj}  \exp\left( \frac{C_5}{2C_0}u_0\dst\sum_{k=0}^\infty 2^{-\bar{\n}(2\m^k-k-2)}\right)-1\\
    &\le \exp\left( \frac{C_5}{2C_0C_6}\dst\sum_{k=0}^\infty 2^{-\bar{\n}(2\m^k-k-2)}\right)-1\\
    &\leq \ex^{\log\left(\frac{3}{2}\right)}-1=\su2<1
\end{align*}
Hence, $\f$ is a lipeomorphism (Lipschitz continuous bijection with inverse Lipschitz continuous as well) from $\rn$ onto itself closed to the identity. Furthermore, $\f\equiv \id$ outside of $\O$ since each $\f_j$ is so, so that $\f$ restricted to $\O$ is a lipeomorphism from $\O$ onto itself.\\
Next, we prove that for each $\o\in\O_\a$, $\Phi(0,x,\o)$ is an invariant Kronecker torus for $H_{|\f(\o)}(y,x)\coloneqq H(y,x,\f(\o))$. Indeed, by letting $j\to\infty$ in \emph{Iterative Lemma}, $(iii)$, we obtain
\[H_{|\f(\o)}\circ \Phi(y,x,\o)=H\circ \mathcal{F}(y,x,\o)\eqqcolon e_\infty(\o)+\o\cdot y,\]
\[\mbox{on}\quad \dst\bigcap_{j\geq 0} D_j\times O_j= T_{*}\times \O_\a.\]
Thus,
\[\Phi^{-1}\circ \dst\phi_{H_{|\f(\o)}}\circ \Phi(y,x;\o)=(y,\o t+x)\quad \ie\quad \phi_{H_{|\f(\o)}}\circ \Phi(y,x;\o)=\Phi(y,\o t+x;\o), \]
on $T_*\times\O_\a$.\\
It remains just to prove that the tori are Lagragian. Indeed, since $T_{(0,x)}T_* \cong \{0\}\times \cn$ for any $x\in\torus^d_{\frac{s}{2}}$, each $\Phi^j$ is symplectic and $\varpi$ is smooth, then we have, for any $\o\in\rn$,
\[\Phi^*\varpi\|_{T^*} = \dst\lim_{j\to \infty}(\Phi^j)^*\varpi\|_{T^*}=\varpi\|_{T^*}=0.\]
\qed
\rem{remPfPo}
Notice that one could apply Lemma~\ref{LipInf} as well to prove that $\f$ is a lipeomorphism, provided that $C_6$ is chosen a little bigger. In fact, by \eqref{kamest3}, we have\footnote{Recall the notations in the proof of \emph{Iterative Lemma.}}
\[
\|\f_j-\id\|_{L,\rn}\leq \frac{C_5}{2C_0} u_j.
\]
Thus, for 
$$
C_6\ge \frac{C_5}{C_0} \dst\sum_{j=0}^\infty 2^{-\bar{\n}(2\m^j-j-2)}\ ,
$$
by taking With $\mathcal{L}_j\equiv \id,\, g_j=\f_j,\, \d=1,\, l_j=\frac{C_5}{2C_0} u_j$, since  
$\sum_{j=0}^{\infty}  l_j\leq  \frac{C_5}{2C_0C_6} \dst\sum_{j=0}^\infty 2^{-\bar{\n}(2\m^j-j-2)}\leq \su2$, for any $j\geq 0$ so that we apply again Lemma~\ref{LipInf} to get
\[
\|\f^j-\id\|_{L,\rn}\leq \frac{C_5}{2C_0}\dst\sum_{k= 0}^\infty u_k\leq \frac{1}{2}.
\]
Therefore
\[
\|\f-\id\|_{L,\rn}\leq \frac{1}{2}<1.
\]
\erem
\chapter{Comparison of the KAM theorems on a mechanical Hamiltonian\label{CQTHMS}}
We consider the simple mechanical Hamiltonian\footnote{As usual, $y^2=y\cdot y=y_1^2+\cdots+y_d^2$.}
\[H_0(y,x;\vae)\coloneqq \frac{y^2}{2} +\vae P_0(x)\coloneqq \frac{y^2}{2} +\vae\left(\cos x_1+ \dst\sum_{j=1}^{d-1}\cos(x_{j+1}-x_j)\right),\]
and we choose
\[\boxed{
s={\bar s}= \frac{10 s_*}{9}=\frac{10\hat s}{9}=20\sigma=1,\quad 
\vae_0=\infty,\quad \t\ge d-1,\quad\a>0,\quad y_0=\o\in \D^\t_\a}.\]
Moreover, we have
\beq{AplDevHamTayl}
H_0(y+y_0,x;\vae)=\frac{\o^2}{2}+\o\cdot y+\frac{y^2}{2}+\vae P_0(x)\eqqcolon H(y,x;\o;\vae).
\eeq
\section{Application of Theorem~\ref{teo2}}
By \eqref{AplDevHamTayl}, we have
\[K_0=\frac{\o^2}{2},\quad Q(y,x)=\frac{y^2}{2},\quad P(y,x;\vae)= P_0(x),\quad T=\uno_d . \]
Hence\footnote{See $\S\ref{AppTeo3}$ for an idea to compute $\|P\|_{r,s,\vae_0}$.}
\[M=\|P\|_{r,s,\vae_0}= \cosh s+(d-1)\cosh (2s)
 ,\quad \wh{\mathrm{E}}=\max\left\{\frac{\o^2}{2},r|\o|,\frac{d r^2}{2},|\o|^2\right\},\]
 \[ \quad L=\frac{\mathrm{C}_\sharp}{3} \wh{\mathrm{E}}^{10}r^{-10}\a^{-4}|\o|^{-6} M,\]
\[\vae_*\coloneqq\mathrm{c}\;\wh{\mathrm{E}}^{-9}\s^{4\t+13}r^{10}\a^4|\o|^6 M^{-1}= \frac{\mathrm{c}\;\wh{\mathrm{E}}^{-9}\s^{4\t+13}r^{10}\a^4|\o|^6}{ \cosh s+(d-1)\cosh (2s)}\; .\]
Therefore, Theorem~\ref{teo2} holds for
\beq{SmalCondAppli22}
\boxed{|\vae|< \vae_*\coloneqq\frac{\mathrm{c}\;\wh{\mathrm{E}}^{-9}\s^{4\t+13} r^{10}\a^4|\o|^6}{\cosh s+(d-1)\cosh (2s)}\;.
}\eeq
\section{Application of Theorem~\ref{teo4}}
By the very definition of $H_0$, we have
\[P=P_0,\quad K(y)=\frac{y^2}{2}, \quad K_y(y)=y,\quad T=K_{yy}(y)=\uno_d,\]
so that\footnote{See $\S\ref{AppTeo3}$ for an idea to compute $\|P_0\|_{r,s,y_0}$ and $\|K_y\|_{r,y_0}$; for the later, writing $y_0=(y_{01},\cdots,y_{0d})$, one can just choose the family $y_a\coloneqq (y_{01}+a\sign(y_{01}),\cdots, y_{0d}+ a\sign(y_{0d})),\ 0<a<r$. }
\[M=\|P_0\|_{r,s,y_0}=\cosh s+(d-1)\cosh (2s),\qquad \|K_y\|_{r,y_0}\overset{def}{=}\dst\sup_{|y-y_0|<r}|y|= r+|y_0|=r+|\o|,\]
\[ \|K_{yy}\|_{r,y_0}\overset{def}{=}\dst\sup_{|a|=1}\|K_{yy}(\cdot)a\|_{r,y_0}=\dst\sup_{|a|=1}\sup_{|y-y_0|<r}|K_{yy}(y)a|=\dst\sup_{|a|=1}\sup_{|y-y_0|<r} |a|= 1,\]
\[ \|T\|=1,\qquad \wh{\mathsf{E}}= \dst\max\left\{ r(r+|\o|),\, |\o|^2\right\}.\]
Therefore, Theorem~\ref{teo4} holds for
\beq{SmalCondAppli44}
\boxed{|\vae|< \vae_*\coloneqq \frac{\a^2}{ M_0}\m_*\;,
}\eeq
where
\begin{align*}
&\m_*\coloneqq\dst\max\left\{0<\m\le \vae_\sharp\ : \ \mathfrak{p}_1\cdot\max\left\{1\,,\,\mathfrak{p}_2\;\m\;\left(\log\m^{-1}\right)^{2\n} \right\}\cdot\m\;\left(\log\m^{-1}\right)^{\n}< 1\right\}\;,\\
&\vae_\sharp\coloneqq \min\left\{\ex^{-1}\,,\,\exp\left(-\frac{\s}{5}\left(\frac{12\sqrt{2}}{5}\frac{\a}{r}\right)^{\su{\n}}\right)
\right\}\;,\\
&\mathfrak{p}_1 \coloneqq \mathsf{C}_8\;\s_0^{-(3\n+2d+1)}\;\max\left\{1,\frac{\a}{r}\right\}\;,\\
&\mathfrak{p}_2 \coloneqq \mathsf{C}_{11}\; \s_0^{-(4\n+2d)}\;,
\end{align*}
\section{Application of Theorem~\ref{teo1}\label{AppTeo3}}
By \eqref{AplDevHamTayl}, we have
\[K_0(\o)=\frac{\o^2}{2},\quad P(y,x;\o)=\frac{y^2}{2}+\vae P_0(x).
\]
We choose
\[ 0<r< \frac{2c\a s^\n}{d} ,\quad h\coloneqq\frac{\epsilon C}{r},\]
with $c,\ C$ as in \eqref{DefCstPoscSC} and $\epsilon  =\|P\|_{r,s,h}$. Next, we compute $\epsilon$.
\begin{align*}
\epsilon=\dst\sup_{D^d_{r,s}\times \O^d_{\a,h}}|P| & \leq \frac{dr^2}{2}+\dst\sup_{x\in \torus^d_1}|\vae|\left(|\cos x_1|+ \dst\sum_{j=1}^{d-1}|\cos(x_{j+1}-x_j)|\right)\\
		 &\le \frac{d r^2}{2}+|\vae|(\cosh s+(d-1)\cosh (2s)),
\end{align*}
 Now, choosing
\[y_a\coloneqq \left(a\; i^{\frac{1-\sign\vae}{2}},a\;i^{\frac{1-\sign\vae}{2}},\cdots,a\;i^{\frac{1-\sign\vae}{2}} \right), \quad 0<a<r\]
and
\[ x_b\coloneqq \left( b\; i,-b\;i,\cdots,\underbrace{(-1)^{j+1}b\;i}_{j^{th}\mbox{ term}},\cdots ,(-1)^{d+1}b\;i\right),\quad 0<b<s,\]
we get
\begin{align*}
\left|P(y_a,x_b;\o)\right| &=\left|\frac{d a^2}{2}\sign(\vae)+\vae\left(\cosh b+ \dst\sum_{j=1}^{d-1}\cosh (2b)\right)\right|\\
				    &= \frac{d a^2}{2}+|\vae|\left(\cosh b+ (d-1)\cosh (2b)\right).
\end{align*}
Therefore,
\begin{align*}
\epsilon &\ge \dst\sup_{\substack{0<a<r\\0<b<s}}\left|P(y_a,x_b;\o)\right|\\
	     &= \dst\sup_{\substack{0<a<r\\0<b<s}}\frac{d a^2}{2}+|\vae|\left(\cosh b+ (d-1)\cosh (2b)\right)\\
	     &= \frac{d r^2}{2}+|\vae|\left(\cosh s+ (d-1)\cosh (2s)\right).
\end{align*}
Thus
\begin{align*}
\epsilon &= \frac{d r^2}{2}+|\vae|(\cosh s+(d-1)\cosh (2s))\\
		  &= \frac{1}{2}c\a r+|\vae|(\cosh s+(d-1)\cosh (2s)).
\end{align*}
 Consequently, 
 if 
\beq{SmalCondAppli1}
\boxed{|\vae|\leq \vae_*\coloneqq \frac{2c\a s^{\n}-d r^2}{2(\cosh s+(d-1)\cosh (2s))}},
\eeq
then Theorem~\ref{teo1} holds.
\section{Application of Theorem~\ref{teo3}}
We choose
\[u\equiv 0,\quad v\equiv y_0=\o.\]
Thus, 
\[H_{yy}\equiv\uno_d,\quad\mathcal{M}=\uno_d,\quad\mathcal{T}=\uno_d,\quad f\equiv 0,\quad g=\vae\nabla_x P_0.\]
Therefore, we can take
\[\mathtt{E}=1,\quad E_{j,k}=0\quad\mbox{if}\quad jk>0 \quad\mbox{or}\quad k\geq 3,\quad E_{0,2}=1,\quad \r= U=\wt V=0,, \]
\[\wt T=1,\quad V=|\o|_1,\quad \mathtt{M}=\ovl{\mathtt{M}}=1,\quad F=0.\]
Next, we compute $G=|\vae|\tnorm{\nabla_x P_0}_{r,s,y_0}$. We have, 
\begin{align*}
 G &= |\vae|\dst\sup_{|a|_1=1}\tnorm{-a_1\frac{\ex^{ix_1}-\ex^{ix_2}}{2i}+\dst\sum_{j=1}^{d-1}(a_j-a_{j+1})\frac{\ex^{i(x_{d+1}-x_d)}-\ex^{-i(x_{d+1}-x_d)}}{2i} }_{r,s,y_0}\\
   &=|\vae|	\dst\sup_{|a|_1=1}\left(|a_1|\frac{\ex^{s}+\ex^{s}}{2}+\dst\sum_{j=1}^{d-1}|a_j-a_{j+1}|\frac{\ex^{2s}+\ex^{2s}}{2}\right)\\
   &=|\vae|	\dst\sup_{|a|_1=1}\left(|a_1|\ex^s+\dst\sum_{j=1}^{d-1}(|a_j|+|a_{j+1}|)\ex^{2s}\right)\\
   &\le |\vae|\min(\ex^s+(d-1)\ex^{2s},\ 2\ex^{2s})\;.
\end{align*}
But then, taking $a_0=(1,0)$ if $d=2$ and
\[a_0=\left(0,1,0,\cdots,0\right)\quad \mbox{if}\quad d\ge 3,\]
we obtain
\[G\ge |\vae|\tnorm{\nabla_x P_0(\cdot)a_0}_{r,s,y_0}= |\vae|\min(\ex^s+(d-1)\ex^{2s},\ 2\ex^{2s}).\]
Hence
\[G=|\vae|\min(\ex^s+(d-1)\ex^{2s},\ 2\ex^{2s})\eqqcolon |\vae|\wh G\;.\]
It remains the choice of $E_{3,0}$. Writting 
\beq{P0Four0}
P_0=\su2(\ex^{ix_1}+\ex^{-ix_1})+\su2\dst\sum_{j=1}^{d-1}\ex^{i(x_{j+1}-x_j)}+e^{-i(x_{j+1}-x_j)}\eqqcolon\dst\sum_{m\in \zn}P_{0m}\ex^{i m\cdot x},
\eeq
we have, for any $j,k,l=1,\cdots,d$,
\beq{3dprP0}
\frac{\dpr^3 P_0}{\dpr x_j\dpr x_k\dpr x_l}=-i\dst\sum_{m\in \zn} m_j m_k m_l P_{0m}\ex^{i m\cdot x},
\eeq
so that
\newpage
\begin{align*}
\tnorm{\dpr^3_x P_0}_{r,s,y_0}&\overset{def}{=} \dst\sup_{|b|_1=|c|_1=1}\dst\sum_{j=1}^d\sup_{|a|_1=1}\tnorm{\dst\sum_{k,l=1}^d\frac{\dpr^3 P_0}{\dpr x_j \dpr x_k\dpr x_l}(\cdot)a_j b_k c_l}_{r,s,y_0}\\
  &= \dst\sup_{\substack{|b|_1=1\\|c|_1=1}}\dst\sum_{j=1}^d\sup_{|a|_1=1}\left(\tnorm{\dst\sum_{\substack{l=1\\ \;}}^d\frac{\dpr^3 P_0}{\dpr x_j^2\dpr x_l}(\cdot)a_j b_j c_l}_{r,s,y_0}+ \right.\\	  
  				  &\qquad\left.  \tnorm{\dst\sum_{\substack{1\le k\le d\\ k\neq j}}\frac{\dpr^3 P_0}{\dpr x_j^2 \dpr x_k}(\cdot)a_j b_k c_j}_{r,s,y_0} +\tnorm{\dst\sum_{\substack{1\le k\le d\\ k\neq j}}\frac{\dpr^3 P_0}{\dpr x_j \dpr x_k^2}(\cdot)a_j b_k c_k}_{r,s,y_0} \right)\\
  				  &\overset{\eqref{3dprP0}}{=} \dst\sup_{|b|_1=|c|_1=1}\dst\sum_{j=1}^d\sup_{|a|_1=1}\left(\dst|a_j| |b_j|\sum_{m\in \zn} |m_j|^2\left|\sum_{\substack{l=1\\ \;}}^d  c_l m_l\right| |P_{0m}|\ex^{ s|m|_1}+ \right.\\
  				  &\hspace*{4cm}\dst|a_j| |c_j|\sum_{m\in \zn} |m_j|^2\left|\sum_{\substack{1\le k\le d\\ k\neq j}}  b_k m_k\right| |P_{0m}|\ex^{ s|m|_1}+ \\
  				  & \hspace*{4cm}\left.\dst|a_j| \sum_{m\in \zn} |m_j|\left|\sum_{\substack{1\le k\le d\\ k\neq j}} b_k c_k m_k^2\right| |P_{0m}|\ex^{ s|m|_1}\right)\\
  				  &\le \dst\sum_{m\in \zn}|P_{0m}|\ex^{ s|m|_1}\sup_{\substack{|b|_1=1\\|c|_1=1}}\left(\dst \sum_{j=1}^d|b_j| |m_j|^2\left|\sum_{\substack{l=1\\ \;}}^d  c_l m_l\right| + \right.\\
  				  &\hspace*{3.5cm}\left.\dst\sum_{j=1}^d |c_j| |m_j|^2\left|\sum_{\substack{1\le k\le d\\ k\neq j}}  b_k m_k\right| + \dst\sum_{j=1}^d |m_j|\left|\sum_{\substack{1\le k\le d\\ k\neq j}} b_k c_k m_k^2\right| \right)\\
  				  &\le \dst\sum_{m\in \zn}|P_{0m}|\ex^{ s|m|_1}(|m|^3+|m|^3+|m|_1|m|^2 )\\
  				  &\overset{\eqref{P0Four0}}{=} 2\su2\left( 3\ex^s+4(d-1)\ex^{2s}\right)\\
  				  &=3\ex^s+4(d-1)\ex^{2s}\eqqcolon E_{3,0}.
\end{align*}
Therefore
\begin{align*}
\wh V        &=  r,\qquad \mathtt{A}_1 = \mathtt{A}_2 =|\o|,\qquad \mathtt{A}_3 = 0,\qquad \mathtt{A}_4 = 3\ex^s+4(d-1)\ex^{2s} ,\\
\mathtt{A}_5 &= \mathtt{A}_6  =  \max\left\{ 3\ex^s+4(d-1)\ex^{2s}\; ,|\o|^2 \right\},\qquad \mathtt{A}_7 = \a^{-2} \max\left\{ 3\ex^s+4(d-1)\ex^{2s}\; ,|\o|^2 \right\},\\
\mathtt{A}_8 &= (s-\hat s)^{2\t}\max\left\{1\; ,\,\frac{|\o|}{r} \right\},\qquad \mathtt{A}_9 = \max\left\{\mathtt{A}_7\; ,\,\mathtt{A}_8 \right\},\qquad  \mathtt{A}_* =  \a^{-2}|\vae|\wh G.
\end{align*}
Therefore, Theorem~\ref{teo3} holds for
\beq{SmalCondAppli033}
\boxed{|\vae|\le \vae_*\coloneqq \frac{\a^2 (s-\hat s)^{2(2\t+1)}}{109\cdot 2^{8\t+13}\t!^4\mathtt{A}_9\wh G}\;.
}\eeq

\noi
In particular, for $d=2$, we have the following.
\cor{cor1}
Consider the hamiltonian $H(y_1,y_2,x_1,x_2;\vae)\coloneqq \frac{y_1^2+y_2^2}{2}+\vae(\cos x_1+\cos(x_2-x_1))$ and $\o\coloneqq \left(\frac{\sqrt{5}-1}{2},1\right)$. Then, for any $|\vae|<\vae_*$, there exists a Kronecker's invariant torus $\mathscr{T}_{y_0,\o} $ for $H$ \ie 
\[\]
\begin{table}[h]
   \centering
\begin{tabular}[c]{|l|c|c|}
\hline
KAM theorem  & Parameters & $\vae_*$ \\
\hline
%
Kolmogorov					       & $r=1,\ \s=1/{20}$ & $9.18337\times 10^{-30}$\\
					       \hline
Arnold  & $r=1,\ \s=1/{20}$ &   $2.02258\times 10^{-49}$ \\
\hline
\multirow{2}{*}{Moser}   & $r=1.73502\times 10^{-15}, \quad \s=1/{20}$  &  \multirow{2}{*}{$6.12208\times 10^{-37} $}       \\
     	  & $h=2.53148\times 10^{-10}+4.46141\times 10^{20} |\vae|$ & \\
\hline
Salamon--Zehnder  & $r=1 ,\ s= 1, \ \hat s=1/10 $ &  $7.38385\times 10^{-27}  $\\
\hline   
\end{tabular}
 \caption{\label{TableVaeStar}  Values of $\vae_*$ according to the KAM theorem}
\end{table}
\ecor
%
\chapter{Global symplectic extension of Arnold's theorem\label{ExtArn013}}
\section{Assumptions\label{AssumpExtArnol}}
Let $\a,r_0>0,\,\t\ge d-1,\, 0<s_0\leq 1,\, \mathscr D\subset\rn$ be  a non--empty, bounded domain\footnote{\ie open and connected.} and consider the Hamiltonian parametrized by $\vae\in\real$
\[H(y,x;\vae)\coloneqq K(y)+\vae P(y,x),\]
where $K,P$ are real--analytic functions with bounded holomorphic extensions to\footnote{Recall the notations in $\S\ref{parassnot}$}
$$
D_{r_0,s_0}(\mathscr D)\coloneqq \dst\bigcup_{y_0\in\mathscr D}D_{r_0,s_0}(y_0)\,,
$$
the norm being
$$
\|\cdot\|_{r_0,s_0,\mathscr D}\coloneqq \dst\sup_{D_{r_0,s_0}(\mathscr D)}|\cdot|\,.
$$
Assume that 
\beq{ArnoldCondExt}
|\det K_{yy}(y)|\not=0\;,\quad\forall\;y\in \mathscr D\;.
\eeq
Define 
\begin{align*}
\D_\a^\t &\coloneqq \left\{\o\in\rn:\quad |\o\cdot k|\ge \frac{\a}{|k|_1^\t}\,,\quad\forall\;k\in\zn\setminus\{0\} \right\}\,,\\
\mathscr D_{r_0,\a} &\coloneqq 
\left\{y_0\in\mathscr D:\ \dist(y_0,\dpr\mathscr D)\ge \frac{r_0}{32d}\ \mbox{ and }\ K_y(y_0)\in \D_\a^\t
\right\}\,,\\
T&\colon \mathscr D_{r_0,\a}\ni y_0 \mapsto K_{yy}(y_0)^{-1}\in \Iso(\rn)\,.
\end{align*}
\noi
Finally, for $\vae\not=0$ given, let\footnote{Recall from footnote\textsuperscript{\ref{ftnarc1}} that $\mathsf{C}_0,\mathsf{C}_1>1$.}
\begin{align*}
M_0 &\coloneqq 
\|P\|_{r_0,s_0,\mathscr D_{r_0,\a}}\;,\\
\mathsf{K}_0 &\coloneqq \|K_{yy}\|_{r_0,\mathscr D_{r_0,\a}}\;,\\
\epsilon   &\coloneqq \frac{\mathsf{K}_0|\vae|M_0}{\a^2 }\;,\\
\mathsf{T}_0 &\coloneqq \|T\|_{\mathscr D_{r_0,\a}}=\sup_{y_0\in\mathscr D_{r_0,\a}}\|T(y_0)\|\;,\\
\mathsf{K}_\infty &\coloneqq \mathsf{K}_0\ex^{\frac{1}{3}}\;,\\
\mathsf{T}_\infty &\coloneqq \mathsf{T}_0\ex^{\frac{1}{3}}\;,\\
\mathsf{C}_0 &\coloneqq 4\left(\frac{3}{2}\right)^{2\t+d+2}\dst\int_{\rn} \left( |y|_1^{\t+1}+d|y|_1^{2\t+2}\right)\ex^{-|y|_1}dy\;,\\
\mathsf{C}_1 &\coloneqq 2\left(\frac{3}{2}\right)^{\t+d+1}\dst\int_{\rn} |y|_1^{\t+1}\ex^{-|y|_1}dy\;,\\
\mathsf{C}_2 &\coloneqq  2^{3d}d\;,\\ 
\mathsf{C}_3 &\coloneqq	d^2\mathsf{C}_1^2+6d\mathsf{C}_1 +\mathsf{C}_2\;,\\
\mathsf{C}_4 &\coloneqq \mathsf{C}_0+\frac{24\ex}{3}\mathsf{C}_1\;,\\
\mathsf{C}_5   &\coloneqq \max\left\{ 4d\mathsf{T}_\infty\mathsf{K}_\infty\,,\,\frac{3}{ 2^{3\t+7}}\mathsf{C}_1\right\}\;,\\
0&<\s_0< \min\left\{\frac{s_0}{2},\,2^{-2(\t-1)}\mathsf{C}_5\sqrt{2} \right\}\;,\\
s_*			&\coloneqq s_0-2\s_0\;,\\
\mathsf{C}_6 &\coloneqq 
\frac{16\mathsf{C}_5\sqrt{2}}{\s_0}\;,\\\
\mathsf{C}_7 &\coloneqq \max\left\{\mathsf{C}_3,\,\mathsf{C}_4\right\}\;,\\
\mathsf{C}_8 &\coloneqq \left(1+\su3\ex^{\su3}\right)^d-1\;,\\
\mathsf{R} &\coloneqq 4\a\mathsf{T}_\infty \;,\\
\l_0	   &\coloneqq \log\epsilon^{-1}\;,\\
r_1		     &\coloneqq \frac{\mathsf{R}}{16\mathsf{C}_5\left(4\s_0^{-1}\l_0\right)^{\t+1}}\;,\\
\mathsf{a}_* &\coloneqq 3\cdot 2^{2\t+\frac{5}{2}}\s_0^{-(3\t+2d+4)}\mathsf{C}_5\max\left\{2^{2(\t+2)}\s_0^{d+1}\,,\, \mathsf{C}_7 \l_0^{-(\t+1)}\sqrt{2}\right\}\;,\\
\mathsf{b}_* &\coloneqq  3\cdot 2^{2(\t+1)}\mathsf{C}_5\;\s_0^{-(3\t+2d+4)}\max\left\{\frac{16\a\mathsf{T}_\infty}{r_0}\s_0^{\t+d+1}\,,\,  \mathsf{C}_7 \max\left\{1,\,\frac{r_0\mathsf{K}_\infty}{\a} \right\}\right\}\;,\\
\mathsf{c}_* &\coloneqq \exp\left(-\frac{\s_0}{4}\left(\frac{\mathsf{R}}{r_0\s_0}\right)^{\su{\t+1}}\right)\;,\\
\mathsf{d}_* &\coloneqq  2^{2\t+2d+3}\mathsf{C}_6^2\,,
\end{align*}
\begin{align*}
\mathsf{e}_* &\coloneqq  
		      \frac{\l_0^{2(\t+1)}}{\a^2 \mathsf{T}_\infty}\cdot \mathsf{a}_* \;,\\
\mathsf{f}_* &\coloneqq  
             \frac{\l_0^{\t+1}}{\a^2 \mathsf{T}_\infty}\cdot \mathsf{b}_*\;,\\ 
\mathsf{g}_*&\coloneqq \frac{\mathsf{b}_*}{\mathsf{T}_\infty\mathsf{K}_0}\;,\\
\mathsf{h}_*&\coloneqq \frac{4\mathsf{a}_*\;\mathsf{b}_*\;\mathsf{d}_*\s_0}{3\mathsf{T}_\infty^2\mathsf{K}_0^2}\;,\\
\vae_1      &\coloneqq \su3 \mathsf{f}_*\s_0^{d+1}|\vae|M_0+\frac{4}{9} \mathsf{f}_*\mathsf{e}_*\mathsf{d}_* \s_0 |\vae|^2 M_0^2= \su3 \mathsf{g}_*\; \s_0^{d+1} \;\epsilon \left(\log\epsilon^{-1}\right)^{\t+1}+\su3 \mathsf{h}_*\; \epsilon^2 \left(\log\epsilon^{-1}\right)^{3(\t+1)} \;,\\
\vae_2      &\coloneqq \su3 \mathsf{f}_*\s_0 |\vae|M_0 +\frac{4}{9} \mathsf{f}_*\mathsf{e}_*\mathsf{d}_* \s_0 |\vae|^2 M_0^2=\su3 \mathsf{g}_*\; \s_0 \;\epsilon \left(\log\epsilon^{-1}\right)^{\t+1}+\su3 \mathsf{h}_* \;\epsilon^2 \left(\log\epsilon^{-1}\right)^{3(\t+1)}\;,\\
\vae_\sharp &\coloneqq 
\mathsf{c}_* \frac{\a^2 }{\mathsf{K}_0 M_0}
\;.
\end{align*}
\section{Statement of the extension Theorem}
\thm{Extteo4}
Under the assumptions in $\S\ref{AssumpExtArnol}$, we have the following. For any given $\vae$ such that\footnote{Notice that $\equ{smcEAr0}$ equivals to: $|\vae|\leq \vae_\sharp\,,\quad |\vae|\,\mathsf{f}_*\,\|P\|_{r_0,s_0,\mathscr D_{r_0,\a}}\le 1\quad\mbox{and}\quad 4|\vae|^2\, \mathsf{f}_*\, \mathsf{e}_*\,\mathsf{d}_*^2\, \s_0\, \|P\|_{r_0,s_0,\mathscr D_{r_0,\a}}^2\le 3$.}
\beq{smcEAr0}
\boxed{\epsilon\le \mathsf{c}_*\;,\qquad \mathsf{g}_*\; \epsilon \left(\log\epsilon^{-1}\right)^{\t+1}\le 1\;,\qquad \mathsf{h}_* \;\epsilon^2 \left(\log\epsilon^{-1}\right)^{3(\t+1)}\le 1\;,}
\eeq
there exist $\mathscr D_*\subset {\mathscr D}$ having the same cardinality  as $\mathscr D_{r_0,\a}$, a lipeomorphism $G_*\colon \mathscr D_{r_0,\a}\overset{onto}{\longrightarrow}\mathscr D_*$,  
a $\ci$ map $K_*\colon \mathscr D\to \real$ and  a $\ci$--symplectomorphism $\phi_*\colon \mathscr D\times \tn\righttoleftarrow,$ real--analytic in $x\in\torus^d_{s_*}$ 
and such that the following hold.
\begin{align}
\dpr_{y_*}K_*\circ G_*&=\dpr_{y}K  \qquad\quad\quad\qquad\ \quad \mbox{on} \quad \mathscr D_{r_0,\a}\;,\label{conjCaneq00}\\
\dpr^\b_{y_*}H\circ \phi_*(y_*,x)&=\dpr^\b_{y_*}K_*(y_*),\qquad \forall\;(y_*,x)\in \mathscr D_*\times\tn,\quad \forall\; \b\in\natural_0^d \label{conjCaneq0}
\end{align}
and 
\begin{align}
|\meas(\mathscr D_*)-\meas(\mathscr D_{r_0,\a})|&\le \mathsf{C}_8\; \vae_2\; \ex^{\vae_2}\meas(\mathscr D_{r_0,\a})\,,\label{measKKst}\\
|\mathsf{W}_0(\phi_*-\id)|
&\le \vae_1
\qquad\qquad\qquad\ \mbox{on}\quad \mathscr D_*\times\torus^d_{s_*}
\,,\label{estArnTrExt}
\end{align}
where
$$
\mathsf{W}_0  \coloneqq \diag\left(\su{4r_1}\uno_d,\uno_d\right)\;.
$$
\ethm
\rem{extArnRem1}
From \equ{conjCaneq0}, on deduces that the $d$--tori
\beq{KronTorArnExt}
\mathcal{T}_{\o_*,\vae}\coloneqq \phi_*\left(y_*,\tn\right),\qquad y_*\in \mathscr D_*\,,\quad \o_*\coloneqq \dpr_{y_*}K_*(y_*)\,,
\eeq
are non-degenrate invariant Kronecker tori for $H$ \ie
\beq{KronTorArnIEExt}
\phi^t_H\circ \phi_*(y_*,x)=\phi_*(y_*,x+\o_* t)\,, \qquad \forall\; x\in\tn.
\eeq
\erem
\section{Proof of Theorem~\ref{Extteo4}}
\noi
\Giu
{\bf KAM step} Given $r,s,K,P,\mathscr D,\mathscr D_\sharp$ satisfying the assumptions $\ref{AssumpExtArnol}$, 
 we seek for $r_1<r,\,s_1<s$, a set $\mathscr D_\sharp'\subset D_{r_1}(\mathscr D_\sharp)$ having the same cardinality as $\mathscr D_\sharp$ and  a near--to--the--identity real--analytic symplectic transformation $\phi_1: \mathscr D\times\tn\righttoleftarrow$ satisfying
\[\phi_1\colon D_{r_1,s_1}(\mathscr D_\sharp')\to D_{r,s}(\mathscr D_\sharp),\]
with $D_{r_1,s_1}(\mathscr D_\sharp')\subset D_{r,s}(\mathscr D_\sharp)$ and $\phi_1$  generated by an extension $y'\cdot x+\vae\hat g(y',x)$ of a function of the form $y'\cdot x+\vae g(y',x)$ \ie
\beq{ArnTraKamExt}
\phi_1\colon \left\{\begin{aligned}
y  &=y'+\vae\hat g_x(y',x)\\
x' &=x+\vae\hat g_{y'}(y',x)\, ,
\end{aligned}
\right.
\eeq
such that
\beq{ArnH1Ext}
\left\{
\begin{aligned}
& H_1:= H\circ \phi_1=K_1+\vae^2 P_1\ ,\quad K_1= K_1(y'),\quad\mbox{on }D_{r_1,s_1}(\mathscr D_\sharp')\,,\\
& \det\dpr_{y'}^2 K_1(y_1)\not=0\,,\quad\forall\; y_1\in \mathscr D_\sharp',\\
& \dpr_{y'} K_1(\mathscr D_\sharp')= \dpr_{y} K(\mathscr D_\sharp)\,.
\end{aligned}
\right.
\eeq
By Taylor's formula, we get\footnote{Recall that $\average{\cdot}$ stands for the average over $\tn$.}
\beq{Arneq11Ext}
\begin{aligned}
H(y'+\vae g_x(y',x),x)=&K(y')+\vae P_0(y') +\vae \left[K_y(y')\cdot g_x +T_{\k} P(y',\cdot)-P_0(y') \right]+\\
						&+\vae^2 \left( P^\ppu+P^\ppd+ P^\ppt\right)(y',x) \\
			= & K_1(y')+\vae \left[K_y(y')\cdot g_x +T_{\k} P(y',\cdot)-P_0(y') \right]+ \vae^2 P'(y',x),
\end{aligned}
\eeq
with $\k\in\natural$, which will be chosen large enough so that $P^\ppt=O(\vae)$, $P_0(y')\coloneqq \average{P(y',\cdot)}$ and 
\beq{ArnDefPsExt}
\left\{
\begin{aligned}
K_1    &\coloneqq K(y')+\vae P_0(y')\eqqcolon K(y')+\vae \wt K(y')\\
P'&\coloneqq P^\ppu+P^\ppd+ P^\ppt\\
P^\ppu &\coloneqq \su{\vae^2}\left[K(y'+\vae g_x)-K(y')-\vae K_y(y')\cdot g_x \right]=\dst\int^1_0(1-t)K_{yy}(\vae t g_x)\cdot g_x\cdot g_x dt\\
P^\ppd &\coloneqq \su\vae \left[P(y'+\vae g_x,x)-P(y',x)\right]=\dst\int_0^1P_y(y'+\vae t g_x,x)\cdot g_x dt\\
P^\ppt &\coloneqq \su\vae \left[ P(y',x)-T_{\k} P(y',\cdot)\right]=\su\vae \dst\sum_{|n|_1>\k} P_n(y')\ex^{in\cdot x}\; .
\end{aligned}
\right.
\eeq
By the non--degeneracy condition in \eqref{ArnoldCondExt} and Lemma~\ref{inv1}, 
 for $\vae$ small enough (to be made precised below), there exists $\bar{r}\le r$ such that 
for each $y_0\in \mathscr D_\sharp$, there exists a unique $y_1\in D_{\bar{r}}(y_0)$ satisfying $\dpr_{y'}K_1(y_1)=\dpr_y K(y_0)$ and $\det\dpr_{y'}^2 K_{1}(y_1)\not=0$; $\mathscr D_\sharp'$ is precisely the set of those $y_1$ when $y_0$ runs in $\mathscr D_\sharp$.  
More precisely, $\mathscr D_\sharp'$ and $\mathscr D_\sharp$ are 
``diffeomorphic''\footnote{\ie there a exits a bijection from $\mathscr D_\sharp$ onto $\mathscr D_\sharp'$ which extends to a diffeomorphism on some neighborhood of $\mathscr D_\sharp$.}, say via $G$,
and, for each $y_1\in\mathscr D_\sharp'$, the matrix $\dpr^2_{y'}K_1(y_1)$ is invertible with inverse of the form
$$
T_1(y_1)\coloneqq \dpr^2_{y'}K_1(y_1)^{-1}\eqqcolon T(y_0)+\vae\wt T(y_1),\quad y_1=G(y_0).
$$
Write
\beq{DefwtYExt}
y_1\eqqcolon y_0+\vae \wt y \,,\quad y_1=G(y_0)\,,\quad\forall\; y_1\in \mathscr D_\sharp'.
\eeq 
In view of \eqref{Arneq11Ext}, in order to get the first part of \eqref{ArnH1Ext}, we need to find $g$ such that  $K_y(y')\cdot g_x +T_{\k} P(y',\cdot)-P_0(y')$ vanishes; such a $g$ is indeed given by
 \beq{HomEqArn1}
 g\coloneqq \dst\sum_{0<|n|_1\leq \k} \frac{-P_n(y')}{iK_y(y')\cdot n}\ex^{in\cdot x},
 \eeq
provided that 
\beq{CondHomEqArn1}
K_y(y')\cdot n\not= 0, \quad \forall\; 0<|n|_1\leq \k,\quad \forall\; y'\in D_{r_1}(\mathscr D_\sharp')\quad  \left(\subset D_{r}(\mathscr D_\sharp)\right).
\eeq
But, in fact, since $K_y(y_0)$ is rationally independent, for each $y_0\in\mathscr D_\sharp$, then, given any $\k\in\natural$, there exists $r'\leq r$ such that
\beq{CondHomEqArnExt}
K_y(y')\cdot n\not=0,\quad \forall\; 0<|n|_1\leq \k, \quad\forall\; y'\in D_{{r}'}(\mathscr D_\sharp).
\eeq
Then we invert the function $x\mapsto x+\vae\hat g_{y'}(y',x)$ in order to define $P_1$. But, by Lemma~\ref{IFTLem}, for $\vae$ small enough, the map $x\mapsto x+\vae \hat g_{y'}(y',x)$ admits an real--analytic inverse of the form
\beq{InvComp2Fi}
\f(y',x';\vae)\coloneqq x'+\vae \wt{\f}(y',x';\vae),
\eeq
so that the Arnod's symplectic transformation is given by
\beq{ArnTrans0}
\phi_1\colon (y',x')\mapsto \left\{
\begin{aligned}
y &= y'+\vae\hat g_x(y',\f(y',x'))\\
x &= \f(y',x';\vae)= x'+\vae \wt{\f}(y',x';\vae) .
\end{aligned}
\right.
\eeq
Hence, \eqref{ArnH1Ext} holds with
\beq{DefP1ArExt}
P_1(y',x')\coloneqq P'(y', \f(y',x')).
\eeq
\noi
Finally, we extend  $K_1$.\\

\noi
Next, we make a quantitative evaluation of the above construction. 
Assume that\footnote{In the sequel, $K$ and $P$  stand for  generic real analytic Hamiltonians which, later on, will respectively play the roles of $K_j$ and $P_j$,  and $y_0,\,r$, the roles of $y_j,\,r_j$ in the iterative step.} 
 $H(y,x;\vae)\coloneqq K(y)+\vae P(y,x)$, where $K,P$ are real--analytic functions with bounded holomorphic extensions to $D_{r,s}(\mathscr D)$ and
\beq{RecHypArnExt}
\begin{aligned}
&\mathscr D_\sharp\subset \left\{y_0\in\mathscr D:\ \dist(y_0,\dpr\mathscr D)\ge \frac{r}{32d}\ \mbox{ and }\ K_y(y_0)\in \D_\a^\t
\right\}\;,\\ 
&  
\det K_{yy}(y)\not=0
\;,\qquad\qquad\qquad T(y)\coloneqq K_{yy}(y)^{-1}\;, \forall\;y\in\mathscr D_\sharp\\
&\|K_{yy}\|_{r,\mathscr D_\sharp}\le \mathsf{K}<\mathsf{K}_\infty\;,\qquad\qquad \|T\|_{\mathscr D_\sharp}\le \mathsf{T}<\mathsf{T}_\infty\;,\\
& \|P\|_{r,s,\mathscr D_\sharp}\le M \;,\qquad\qquad\,\,\, \o\in \D^\t_\a\;,\qquad\qquad r\le r_0 \;,
\end{aligned}
\eeq
\noi
Fix $0<2\s<s\leq 1$ and 
 fix $\vae\not=0$ in such away that,
\beq{DefNArnExt}
\begin{aligned}
&\l\ge \log\left(\frac{\a^2 }{\mathsf{K}|\vae|M}\right)>1,\quad\k\coloneqq 4\s^{-1}\l, \quad 
\bar{r}\coloneqq 
\su{4\mathsf{C}_{5}}\dst\min\left(\frac{\mathsf{R}}{\k^{\t+1}}\,,\, r\s \right),\\
&\tilde r\coloneqq \frac{\bar{r}}{4\mathsf{C}_5},\quad \bar{s}\coloneqq s-\frac{2}{3}\s,\quad s'\coloneqq s-\s \,,
\end{aligned}
\eeq
so that
\footnote{Recall footnote \textsuperscript{\ref{ftnTK1}}.}
\beq{rrbarAsExt}
\bar{r}\le\frac{r\s}{4\mathsf{C}_{5}}= \frac{\s}{16d\mathsf{T}_\infty\mathsf{K}_\infty}r\le \frac{r}{32d}< \frac{r}{2}\quad \mbox{and}\quad \k> 8\,.\\
\eeq
\lem{lem:1Ext}
Let\footnote{Notice that $\mathsf{L}\ge \s^{-d}\ovl{\mathsf{L}}\ge \ovl{\mathsf{L}}$ since $\s\le 1$. Notice also that $\mathsf{T}\mathsf{K}\ge 1$, so that $\frac{16\mathsf{T}_\infty  }{r\bar{r}}\s^{-(\t+d+1)}\ge \frac{16\mathsf{T} }{r^2}\ge \frac{4}{\mathsf{K}r^2}$.
}
\begin{align*}
\ovl{\mathsf{L}}&\coloneqq\mathsf{C}_4 \max\{\a,r\mathsf{K}\}\frac{ M}{\a^2 \bar{r}}\s^{-(2\t+d+2)}\;,\\
\mathsf{L}&\coloneqq M\dst\max\left\{\frac{16\mathsf{T}_\infty  }{r\bar{r}}\s^{-(\t+d+1)}\,,\,
\mathsf{C}_7 \max\{\a,r\mathsf{K}\}\frac{ 1}{\a^2 \bar{r}}\s^{-2(\t+d+1)}\right\}\\
          &=M\dst\max\left\{\frac{16\mathsf{T}_\infty  }{r\bar{r}}\s^{-(\t+d+1)}
\,,\,\frac{4}{\mathsf{K} r^2}\,,\,\mathsf{C}_7 \max(\a,r\mathsf{K})\frac{ 1}{\a^2 \bar{r}}\s^{-2(\t+d+1)}\right\}
\;.
\end{align*}
Then
\beq{Est1Lem1bExt}
\left\{
\begin{aligned}
&\|g_x\|_{\bar{r},\bar{s},\mathscr D_\sharp}\le  \mathsf{C}_1 \frac{M}{\a} \s^{-(\t+d+1)}\,,\\
& \|g_{y'}\|_{\bar{r},\bar{s},\mathscr D_\sharp},\, \|\dpr_{y'x}^2 g\|_{\bar{r},\bar{s},\mathscr D_\sharp}\le \ovl{\mathsf{L}}\,,\\
&\|\dpr_{y'}^2\wt K\|_{\bar{r},\mathscr D_\sharp}\le 
\mathsf{K}\mathsf{L}\,.
\end{aligned}
\right.
\eeq
 If $\vae_*>0$ satisfies 
\beq{cond1ExtExt}
\vae_*\leq \vae_\sharp\quad \mbox{and}\quad {\vae_* }{\mathsf{L}}\le \frac{\sigma}{3}
\ ,
\eeq
then, for any $|\vae|\leq \vae_*$, there exists 
a  diffeomorphism $G\colon  D_{\tilde{r}}(\mathscr D_\sharp){\to}G( D_{\tilde{r}}(\mathscr D_\sharp))$, 
$\dpr_{y'}K_1\circ G=\dpr_{y}K $ and such that $\mathscr D_\sharp'\coloneqq G(\mathscr D_\sharp)\subset B_{\bar{r}}(\mathscr D_\sharp)$,
\beq{convEstExtRe}
\left\{
\begin{aligned}
&|\vae|\|g_x\|_{\bar{r},\bar{s},\mathscr D_\sharp}\le \frac{r}{3}\,,\qquad \qquad\quad\,\,
\|G-\id\|_{\tilde{r},\mathscr D_\sharp}\le \frac{\bar{r}}{2}\;,\\
&|\vae|\|\wt T\|_{\mathscr D_\sharp'}\le \mathsf{T}|\vae|\mathsf{L}
\,, \qquad\quad\quad\quad \|\dpr_z G-\uno_d\|_{\tilde{r},\mathscr D_\sharp}\le |\vae|\mathsf{L}\,,\\
&\|P'\|_{\bar{r},\bar s,\mathscr D_\sharp} \le  \mathsf{L}M\,, \qquad\qquad\quad\,\,\, B_{\bar{r}/4}(\mathscr D_\sharp')\subset B_{\bar{r}/2}(\mathscr D_\sharp)\subset\mathscr D
\end{aligned}
\right. 
\eeq
and the following hold. $g$ has and extension $\hat g\colon \rn\times\tn\to\real$ and, for any $|\vae|\leq\vae_*$ and $y'\in D_{\bar{r}/2}(\mathscr D_\sharp)$, the map $\psi_\vae(x):=  x+\vae\hat g_{y'}(y',x)$ has an analytic inverse 
 $\f(x')=x'+\vae \wt{\f}(y',x';\vae)$  such that, for all $|\vae|\le\vae_*$, 
\beq{boundalExtExtRe}
\|\wt{\f}\|_{\bar{r}/4, s',y_0}\le   \ovl{\mathsf{L}} \qquad {\rm and}
\quad
\f=\id + \vae \wt{\f} : D_{\bar{r}/4,s'}(\mathscr D_\sharp')\to \torus^d_{\bar{s}} \ ;
\eeq
for any $y_0\in\mathscr D_\sharp$ and $(y',x,\vae)\in D_{\bar{r},\bar s}(y_0)\times D^1_{\vae_*}(0)$, $|y'+\vae g_x(y',x)-y_0|<\frac{5}{6} r$; the map $\phi_1$ is a symplectic diffeomorphism and
\beq{phiokExt0Ext}
\phi_1=\big( y'+\vae g_x(y', \f(y',x')),\f(y',x')\big): D_{\bar{r}/4,s'}(\mathscr D_\sharp')\to D_{2r/3, \bar{s}}(\mathscr D_\sharp),
\eeq
with
\beq{phiokExt1Ext}
\|\mathsf{W}\,\tilde \phi\|_{\bar{r}/4,s',\mathscr D_\sharp'}\le \ovl{\mathsf{L}}
\,,
\eeq
where $\tilde \phi$ is defined by the relation $\phi_1=:\id + \vae \tilde \phi$,
$$
\mathsf{W}\coloneqq \begin{pmatrix}
\su{\bar r}\uno_d & 0\\
0			& \uno_d 
\end{pmatrix}
$$
and
\beq{tesitExtExtRe}
\begin{aligned}
&\|P_1\|_{\bar{r}/4, s',\mathscr D_\sharp'}\le  \mathsf{L}M\;.
\end{aligned}
\eeq
Moreover, $K_1$ possesses a $\ci$--extensions $\wh K_1\colon \rn\to \real$
such that for any $n\in\natural_0$, there exists $C_n\in \natural$ and for any $\b_1,\b_2\in\natural_0^d$ with $|\b_1|_1+|\b_2|_1\le n$,
\begin{align}
\bar{r}^{|\b_1|_1}\s^{|\b_2|_1}\|\dpr_{y'}^{\b_1}\dpr_{x'}^{\b_2}\mathsf W(\phi_1-\id)\|_{0,0}&\le C_n|\vae|\ovl{\mathsf L}\,,\label{fi1estextExt}\\
\bar{r}^{|\b_1|_1}\|\dpr_{y'}^{\b_1}(\wh K_1- K)\|_0&\le C_n |\vae| M\;. \label{KestextExt}
\end{align}
\elem
\proof
We begin by extending the ``diophantine condition w.r.t. $K_y$'' uniformly to $D_{\bar{r}}(\mathscr D_\sharp)$ up to the order $\k$. Indeed, 
for any $y_0\in\mathscr D_\sharp$, $0<|n|_1\leq \k$ and $y'\in D_{\bar{r}}(y_0)$,
\begin{align}
|K_y(y')\cdot n|&=|\o\cdot n +(K_y(y')-K_y(y_0))\cdot n|\geq |\o\cdot n|\left(1-d\frac{\|K_{yy}\|_{\bar{r},\mathscr D_\sharp}}{|\o\cdot n|}|n|_1\bar{r}\right) \nonumber\\
         &\geq \frac{\a}{|n|_1^\t}\left(1-\frac{d\mathsf{K}}{\a }|n|_1^{\t+1}\bar{r} \right)\geq \frac{\a}{|n|_1^\t}\left(1-\frac{d\mathsf{K}}{\a }\k^{\t+1}\bar{r} \right)\geq \frac{\a}{|n|_1^\t}\left(1-\frac{d\mathsf{K}_\infty}{\a }\k^{\t+1}\bar{r} \right)\nonumber\\
         &\ge \frac{\a}{2|n|_1^\t},\label{ArnExtDiopCondExt}
\end{align}
so that, by Lemma~\ref{fce}--$(i)$, we have 

\begin{align*}
\|g\|_{\bar{r},\bar{s},\mathscr D_\sharp} &\overset{def}{=}\dst\sup_{D_{\bar{r},\bar{s}}(\mathscr D_\sharp)}\left|\dst\sum_{0<|n|_1\leq \k}\frac{P_n(y')}{K_y(y')\cdot n}\ex^{in\cdot x} \right|\leq \dst\sum_{0<|n|_1\leq \k}\frac{\|P_n\|_{\bar{r},\bar{s}, \mathscr D_\sharp}}{|K_y(y')\cdot n|}\ex^{\left(s-\frac{2}{3}\s\right)|n|_1}\\
   &\leq \dst\sum_{0<|n|_1\leq \k} M\ex^{-s|n|_1}\frac{2|n|_1^{\t}}{\a}\ex^{\left(s-\frac{2}{3}\s\right)|n|_1}\leq \frac{2M}{\a}\dst\sum_{n\in\zn} |n|_1^{\t}\ex^{-\frac{2}{3}\s|n|_1}\\
   &\leq \frac{2M}{\a}\dst\int_{\rn} |y|_1^{\t}\ex^{-\frac{2}{3}\s|y|_1}dy\\
   &= \left(\frac{3}{2\s}\right)^{\t+d}\frac{2M}{\a}\dst\int_{\rn} |y|_1^{\t}\ex^{-|y|_1}dy\\
   &\le \mathsf{C}_1 \frac{M}{\a} \s^{-(\t+d)}
\end{align*}
and analogously,
\begin{align*}
\|g_x\|_{\bar{r},\bar{s},\mathscr D_\sharp} &\overset{def}{=}\dst\sup_{D_{\bar{r},\bar{s}}(\mathscr D_\sharp)}\left|\dst\sum_{0<|n|_1\leq \k}\frac{nP_n(y')}{K_y(y')\cdot n}\ex^{in\cdot x} \right|\leq \dst\sum_{0<|n|_1\leq \k}\frac{\|P_n\|_{\bar{r},\bar{s}, \mathscr D_\sharp}}{|K_y(y')\cdot n|}|n|_1\ex^{\left(s-\frac{2}{3}\s\right)|n|_1}\\
   &\leq \dst\sum_{0<|n|_1\leq \k} M\ex^{-s|n|_1}\frac{2|n|_1^{\t+1}}{\a}\ex^{\left(s-\frac{2}{3}\s\right)|n|_1}\leq \frac{2M}{\a}\dst\sum_{n\in\zn} |n|_1^{\t+1}\ex^{-\frac{2}{3}\s|n|_1}\\
   &\leq \frac{2M}{\a}\dst\int_{\rn} |y|_1^{\t+1}\ex^{-\frac{2}{3}\s|y|_1}dy\\
   &= \left(\frac{3}{2\s}\right)^{\t+d+1}\frac{2M}{\a}\dst\int_{\rn} |y|_1^{\t+1}\ex^{-|y|_1}dy\\
   &\le \mathsf{C}_1 \frac{M}{\a} \s^{-(\t+d+1)}\,,
\end{align*}
\begin{align*}
\|\dpr_{y'}g\|_{\bar{r},\bar{s},\mathscr D_\sharp} &\overset{def}{=}\dst\sup_{D_{\bar{r},\bar{s}}(\mathscr D_\sharp)}\left|\dst\sum_{0<|n|_1\leq \k}\left(\frac{ \dpr_yP_n(y')}{K_y(y')\cdot n}-P_n(y')\frac{ K_{yy}(y')n}{(K_y(y')\cdot n)^2}\right)\ex^{in\cdot x} \right|\\
   &\leq \dst\sum_{0<|n|_1\leq \k}\dst\sup_{D_{\bar{r}}(\mathscr D_\sharp)}\left(\frac{\|(P_y)_n\|_{\bar{r},s,\mathscr D_\sharp}}{|K_y(y')\cdot n|}+d\|P_n\|_{r,s,\mathscr D_\sharp}\frac{\|K_{yy}\|_{r,\mathscr D_\sharp}|n|_1}{|K_y(y')\cdot n|^2}\right)\ex^{\left(s-\frac{2}{3}\s\right)|n|_1}\\
   &\stackrel{\equ{RecHypArnExt}+\equ{ArnExtDiopCondExt}}{\le} \dst\sum_{0<|n|_1\leq \k}\left( \frac{M}{r-\bar{r}}\ex^{-s|n|_1}\frac{2|n|_1^{\t}}{\a}+dM\ex^{-s|n|_1}\mathsf{K}|n|_1\left(\frac{2|n|_1^{\t}}{\a}\right)^2\right)\ex^{\left(s-\frac{2}{3}\s\right)|n|_1}\\
   &\leby{rrbarAsExt} \frac{4M}{\a^2 r}\dst\sum_{0<|n|_1\leq \k}\left( |n|_1^{\t}\a +dr\mathsf{K}|n|_1^{2\t+1}\right)\ex^{-\frac{2}{3}\s|n|_1}\\
   &\le \max(\a,r\mathsf{K})\frac{4M}{\a^2 r}\dst\sum_{0<|n|_1\leq \k}\left( |n|_1^{\t}+d|n|_1^{2\t+1}\right)\ex^{-\frac{2}{3}\s|n|_1}\\
   &\leq \max(\a,r\mathsf{K})\frac{4M}{\a^2 r}\dst\int_{\rn} \left( |y|_1^{\t}+d|y|_1^{2\t+1}\right)\ex^{-\frac{2}{3}\s|y|_1}dy \\
   &= \left(\frac{3}{2\s}\right)^{2\t+d+1}\max(\a,r\mathsf{K})\frac{4M}{\a^2 r}\dst\int_{\rn} \left( |y|_1^{\t}+d|y|_1^{2\t+1}\right)\ex^{-|y|_1}dy\\
   &\le \mathsf{C}_0 \max(\a,r\mathsf{K})\frac{M}{\a^2 r}\s^{-(2\t+d+1)}\\
   &\le \ovl{\mathsf{L}} \;,
\end{align*}
\begin{align*}
\|\dpr^2_{y'x}g\|_{\bar{r},\bar{s},\mathscr D_\sharp} &\overset{def}{=}\dst\sup_{D_{\bar{r},\bar{s}}(\mathscr D_\sharp)}\left|\dst\sum_{0<|n|_1\leq \k}\left(\frac{ \dpr_yP_n(y')}{K_y(y')\cdot n}-P_n(y')\frac{ K_{yy}(y')n}{(K_y(y')\cdot n)^2}\right)\cdot n\ex^{in\cdot x} \right|\\
   &\leq \dst\sum_{0<|n|_1\leq \k}\dst\sup_{D_{\bar{r}}(\mathscr D_\sharp)}\left(\frac{\|(P_y)_n\|_{\bar{r},s, \mathscr D_\sharp}}{|K_y(y')\cdot n|}+d\|P_n\|_{r,s,\mathscr D_\sharp}\frac{\|K_{yy}\|_{r,\mathscr D_\sharp}|n|_1}{|K_y(y')\cdot n|^2}\right)|n|_1\ex^{\left(s-\frac{2}{3}\s\right)|n|_1}\\
   &\le \max(\a,r\mathsf{K})\frac{4M}{\a^2 r}\dst\sum_{0<|n|_1\leq \k}\left( |n|_1^{\t}+d|n|_1^{2\t+1}\right)|n|_1\ex^{-\frac{2}{3}\s|n|_1}\\
   &\leq \max(\a,r\mathsf{K})\frac{4M}{\a^2 r}\dst\int_{\rn} \left( |y|_1^{\t}+d|y|_1^{2\t+1}\right)|y|_1\ex^{-\frac{2}{3}\s|y|_1}dy \\
   &= \left(\frac{3}{2\s}\right)^{2\t+d+2}\max(\a,r\mathsf{K})\frac{4M}{\a^2 r}\dst\int_{\rn} \left( |y|_1^{\t+1}+d|y|_1^{2\t+2}\right)\ex^{-|y|_1}dy\\
   &= \mathsf{C}_0 \max(\a,r\mathsf{K})\frac{M}{\a^2 r}\s^{-(2\t+d+2)} \\
   &=\ovl{\mathsf{L}}\;,
\end{align*}
and, for $|\vae|< {\vae_*}$,
\[\|\dpr_{y'}\wt K\|_{\bar{r},\mathscr D_\sharp}=\| \left[P_y\right]\|_{\bar{r},\mathscr D_\sharp}\leq \|P_y\|_{\bar{r},\bar{s}, \mathscr D_\sharp}\leq  \frac{M}{r-\bar{r}}\leq \frac{2M}{r} \;,\]
\[\|\dpr_{y'}^2\wt K\|_{\bar{r},\mathscr D_\sharp}=\| \left[P_{yy}\right]\|_{\bar{r},\mathscr D_\sharp}\leq \|P_{yy}\|_{\bar{r},\bar{s}, \mathscr D_\sharp}\leq  \frac{M}{(r-\bar{r})^2}\leq \frac{4M}{r^2}\le \mathsf{K}\mathsf{L} 
\;.\]
\noi
Now, we extend the generating function and $K_1$ to $\rn\times\tn$ and $\rn$ respectively, by making use of a cut--off function. Let then $\chi_1\in C(\cn)\cap\ci(\rn)$, with $0\le \chi_1\le 1,\, \supp\chi_1\subset D_{\bar{r}}(\mathscr D_\sharp),\, \chi_1\equiv 1$ on $D_{\bar{r}/2}(\mathscr D_\sharp)$ and satisfying \equ{cutof}. Thus, given $x\in \torus^d_{\bar{s}}$, set $\hat{g}(y',x')=\chi_1(y'){g}(y',x'),\ \hat{K}_1=K+\chi_1\cdot( K_1- K)$ if $y'\in D_{\bar{r}}(\mathscr D_\sharp)$, $\hat{g}(y',x')=0,\ \hat{K}_1= K$ if $y'\in \left(\cn\setminus D_{\bar{r}}(\mathscr D_\sharp)\right)\cap \mathscr D$ and $\hat{K}_1= K\equiv 0$ on $\cn\setminus\left( D_{\bar{r}}(\mathscr D_\sharp)\cup \mathscr D\right)$.
\begin{align}
\|\hat g\|_{0,0}&\le \|g\|_{\bar{r},\bar{s},\mathscr D_\sharp}\le \mathsf{C}_1 \frac{M}{\a} \s^{-(\t+d)}\;,\label{gext001}\\
\|\hat g_x\|_{0,0}&\le \|g_x\|_{\bar{r},\bar{s},\mathscr D_\sharp}\le \mathsf{C}_1 \frac{M}{\a} \s^{-(\t+d+1)}\;,\label{gext002}\\
\|\hat g_{y'}\|_{0,0}&\le \|\dpr_{y'}\chi_1\|_0\|g\|_{\bar{r},\bar{s},\mathscr D_\sharp} +\|g_{y'}\|_{\bar{r},\bar{s},\mathscr D_\sharp}\leby{cutof} \mathsf{C}_4 \max(\a,r\mathsf{K})\frac{M}{\a^2 \bar{r}}\s^{-(2\t+d+1)}\le \ovl{\mathsf{L}}\;,\label{gext003}\\
\|\hat g_{y'x}\|_{0,0}&\le \|\dpr_{y'}\chi_1\|_0\|g_x\|_{\bar{r},\bar{s},\mathscr D_\sharp} +\|g_{y'x}\|_{\bar{r},\bar{s},\mathscr D_\sharp}\leby{cutof}\mathsf{C}_4 \max(\a,r\mathsf{K})\frac{M}{\a^2 \bar{r}}\s^{-(2\t+d+2)}\le \ovl{\mathsf{L}}\;.\label{gext004}
\end{align}
And generally, for any $n\in\natural_0$, there exists $C_n\in \natural$ and for any $\b_1,\b_2\in\natural_0^d$ with $|\b_1|_1+|\b_2|_1\le n$,
\beq{gextciest}
\bar{r}^{|\b_1|_1}\s^{|\b_2|_1-1}\|\dpr_{y'}^{\b_1}\dpr_{x'}^{\b_2}\hat{g}_{y'}\|_{0,0}\le C_n\ovl{\mathsf L}\; .
\eeq
\equ{KestextExt} is a straightforward consequence of Leibniz's rule.\\
\noi
Next, we construct $\mathscr D_\sharp'$ 
 in \eqref{ArnH1Ext} for $|\vae|\le {\vae_*}$. 
 For, fix $|\vae|\leq {\vae_*}$, $y_0\in\mathscr D_\sharp$ 
  and consider
\begin{align*}
F\colon D_{\bar{r}}(y_0)\times D_{\tilde{r}}(y_0) &\longrightarrow \qquad \cn\\
		(y,z)\quad &\longmapsto K_y(y)+\vae \wt K_{y'}(y)-K_y(z).
\end{align*}
Then
\begin{itemize}
\item $F_y(y_0,y_0)=\dpr^2_{y}K(y_0)+\vae\dpr^2_{y'}\wt K(y_0) = T(y_0)^{-1}\left(\uno_d+\vae T(y_0)\dpr^2_{y'}\wt K(y_0) \right)\eqqcolon T(y_0)^{-1}(\uno_d+\vae A_0)$ and
\begin{align*}
\|\vae A_0\|&\le \|T(y_0)\|\|\vae \dpr^2_{y'}\wt K(y_0)\|\le \mathsf T \frac{4|\vae|  M}{r^2}\leby{rrbarAsExt} |\vae|\frac{2\mathsf T_\infty M}{r\bar{r}}
\le \su2\vae_*\mathsf L\leby{cond1ExtExt} \frac{\s}{6}<\su2.
\end{align*}
Hence, $F_y(y_0,y_0)$ is invertible, with inverse 
$$
T_0\coloneqq (\uno_d+\vae A_0)^{-1}T(y_0)=\left(\uno_d+\dst\sum_{k\geq 1}(-\vae)^k A_0^k\right) T(y_0)
$$
 satisfying
\beq{T1estExt}
\|T_0\|\le \frac{\|T(y_0)\|}{1-\|\vae A_0\|}\le 2\mathsf T.
\eeq
\item For any $(y,z)\in D_{\bar{r}}(y_0)\times D_{\tilde{r}}(y_0)$,
\begin{align*}
	\|\uno_d-T_0F_y(y,z)\|&\leq \|T_0\|\|\dpr^2_{y}K(y_0)- K_{yy}\|_{\bar{r},\mathscr D_\sharp}+|\vae|\;\|T_0\|\;|\dpr^2_{y'}\wt K(y_0)|+|\vae|\;\|T_0\|\;\|\dpr_{y'}^2\wt K\|_{\bar{r},\mathscr D_\sharp}\\
		  &\leq d\cdot 2\mathsf T\|K_{yyy}\|_{\bar{r},\mathscr D_\sharp}\cdot\bar{r}+ 4|\vae|\mathsf{T}\frac{4M}{r^2}\\
	      &\leq 2d\mathsf{T}_\infty \mathsf{K}_\infty\frac{\bar{r}}{r-\bar{r}}+16\mathsf{T}_\infty\frac{|\vae| M}{ r^2}\\
	      &\leby{rrbarAsExt}\frac{\mathsf{C}_5}{2} \frac{2\bar{r}}{ r}+|\vae|\frac{16\mathsf{T}_\infty M}{r\bar{r}}\\
	      &\le \mathsf{C}_5\frac{\bar{r}}{ r}+{\vae_*}\mathsf{L}\\
	      &\overset{\equ{DefNArnExt}+\equ{cond1ExtExt}}{\leq}\su4+\frac{\s}{3}\\
	      &\le \su4+\su4=\su2\;;
\end{align*}
\item For any $z\in D_{\tilde{r}}(y_0)$,
\begin{align*}
2\|T_0\||F(y_0,z)|&\le 4\mathsf{T}|K_y(z)-K_y(y_0)|+ 4\mathsf{T}|\vae|\| \wt K_{y'}\|_{\bar{r},\mathscr D_\sharp}\\
&\leq 4\mathsf{T}\|K_{yy}\|_{\bar{r},\mathscr D_\sharp}\cdot \tilde{r}+ 4\mathsf{T} \frac{2|\vae| M}{r}\\
&\le 4\mathsf{T}_\infty\mathsf{K}_\infty\tilde{r}+\frac{8|\vae|\mathsf{T}_\infty M}{r}\\
&\le \frac{\mathsf{C}_5}{d}\frac{\bar{r}}{4\mathsf{C}_5}+{\vae_*}\frac{\bar{r}}{2}\mathsf{L}\\
&\leby{cond1ExtExt}\frac{\bar{r}}{8}+ \frac{\bar{r}}{12}\\
&<\frac{\bar{r}}{4}\;,
\end{align*}
\ie 
$$
2\|T_0\|\|F(y_0,\cdot)\|_{\tilde{r},y_0}\le \frac{\bar{r}}{4}\;.
$$
\end{itemize}

\noi
Therefore, Lemma~\ref{IFTLem} applies. Hence, there exists a real--analytic map $G^{y_0}\colon D_{\tilde{r}}(y_0)\to D_{\bar{r}}(y_0)$ such that its graph coincides with $F^{-1}(\{0\})$ \ie  $y_1=y_1(z,y_0,\vae)\coloneqq G^{y_0}(z)$ is the unique $y\in D_{\bar{r}}(y_0)$ satisfying $0=F(y,z)=\dpr_y K_1(y)-K_y(z)$, for any $z\in D_{\tilde{r}}(y_0)$
. Moreover, $\forall\;z\in D_{\tilde{r}}(y_0)$,
\beq{EcarY1Y0ReExt}
|G^{y_0}(z)-y_0|\leq 2\|T_0\|\|F(y_0,\cdot)\|_{\tilde{r},y_0}
\leq \frac{\bar{r}}{4}\;,
\eeq
\beq{gy0zRe}
|G^{y_0}(z)-z|\le |G^{y_0}(z)-y_0|+|y_0-z|\le \frac{\bar{r}}{4}+\tilde{r}< \frac{\bar{r}}{2},
\eeq
so that
\beq{NextSetArnExt}
D_{{\bar{r}}/{4}}(G^{y_0}(z))\subset D_{\bar{r}/2}(y_0).
\eeq
\noi
Next, we prove that $\dpr^2_{y'}K_1(y_1)$ is invertible, where $y_1=G^{y_0}(z)$ for some given $z\in D_{\tilde{r}}(y_0)$. Indeed, by Taylor's formula, we have,
\begin{align*}
\dpr^2_{y'} K_1(y_1)&= K_{yy}(y_0)+ \dst\int_0^1 K_{yyy}(y_0+t(y_1-y_0))(y_1-y_0) dt+\vae \dpr^2_{y'}\wt K(y_1)\\
           &= T(y_0)^{-1}\left(\uno_d+T(y_0)\left(\dst\int_0^1 K_{yyy}(y_0+t(y_1-y_0))(y_1-y_0) dt+ \dpr^2_{y'}\wt K(y_1)\right)\right)\\
           &\eqqcolon T(y_0)^{-1}(\uno_d+\vae A),
\end{align*}
and, by Cauchy's estimate, for any\footnote{Recall footnote \textsuperscript{\ref{ftnTK1}} 
.} $|\vae|\leq {\vae_*}$,
\begin{align*}
|\vae|\|A\|&\leq \|T(y_0)\|\left(d\|K_{yyy}\|_{\bar{r},\mathscr D_\sharp}|y_1-y_0|+ |\vae|\|\dpr_{y'}^2\wt K\|_{\bar{r},\mathscr D_\sharp}\right)\\
     &\leq \|T\|_{\mathscr D_\sharp}\left(\frac{d\|K_{yy}\|_{r,\mathscr D_\sharp}}{r-\bar{r}}|y_1-y_0|+|\vae|\|\dpr^2_{y'}\wt K\|_{\bar{r},\mathscr D_\sharp}\right)\\
     &\leq \mathsf{T}\left(\frac{2d\mathsf{K}}{r}\frac{\bar{r}}{2}+\frac{4M}{r^2} \right)\\
	 &\leby{rrbarAsExt} \mathsf{T}\left(\frac{\s}{16\mathsf{T}_\infty}+\su{4\mathsf{T}_\infty}\vae_*\mathsf{L} \right)\\
	 &\leby{cond1ExtExt} \mathsf{T}\left(\frac{\s}{16\mathsf{T}_\infty}+\su{4\mathsf{T}_\infty}\frac{\s}{3} \right)\\
	 &<\frac{\s}{6}\\
	 &<\su2.
\end{align*}
Hence $\dpr_{y'}^2K_1(y_1)$ is invertible with
\[\dpr_{y'}^2K_1(y_1)^{-1}=(\uno_d+\vae A)^{-1}T(y_0)=T(y_0)+\dst\sum_{k\geq 1}(-\vae)^k A^k T(y_0)\eqqcolon T(y_0)+\vae \wt T(y_1),\]
and
\[|\vae|\|\wt T(y_1)\|\leq |\vae|\frac{\|A\|}{1-|\vae|\|A\|}\|T\|_{\mathscr D_\sharp}\leq 2|\vae|\|A\| \|T\|_{\mathscr D_\sharp}
\le 2\frac{\s}{6}\mathsf{T}
= \mathsf{T}\frac{\s}{3}\,.\]
\noi
Similarly, from
$$
K_{yy}(z)=K_{yy}(y_0)\left(\uno_d+T(y_0)\dst\int_0^1 K_{yyy}(y_0+t(z-y_0))(z-y_0) dt\right)
$$
and
\begin{align*}
\left\|T(y_0)\dst\int_0^1 K_{yyy}(y_0+t(z-y_0))(z-y_0) dt\right\|_{r/\mathsf{C}_5,y_0}&\le \mathsf{T}\|K_{yyy}\|_{r/2,y_0}\frac{r}{\mathsf{C}_5}\le \mathsf{T}\frac{d\mathsf{K}}{r-r/2}\frac{r}{\mathsf{C}_5}<\su2
\end{align*}
 one has that, for any $z\in D_{r/\mathsf{C}_5}(y_0)$, 
\beq{invJacKyyRe}
K_{yy}(z)^{-1} \mbox{ exists and }\|K_{yy}(z)^{-1}\|\le \|K_{yy}(z)^{-1}-T(y_0)\|+\|T(y_0)\|\le 
2\su2 \mathsf{T}+\mathsf{T}=2\mathsf{T}
\eeq
Now, differentiating $F(G^{y_0}(z),z)=0$, we get, for any $z\in D_{\tilde{r}}(y_0)$,
$$
\dpr_{y'}^2K_1(G^{y_0}(z))\cdot\dpr_z G^{y_0}(z)=K_{yy}(z)\;.
$$
Therefore $G^{y_0}$ is a local diffeomorphism, with
\begin{align*}
\dpr_z G^{y_0}(z)&= \dpr_{y'}^2K_1(G^{y_0}(z))^{-1}K_{yy}(z)\\
                &=\left(K_{yy}(z)^{-1}\left(K_{yy}(z)+\vae \dpr_{y'}^2\wt K(g^{y_0}(z)) \right)\right)^{-1}\\
                &=\left(\uno_d+\vae K_{yy}(z)^{-1}\dpr_{y'}^2\wt K(g^{y_0}(z))\right)^{-1}
\end{align*}
and
$$
\|\vae K_{yy}^{-1}\dpr_{y'}^2\wt K\|_{\tilde{r},y_0}\le \| K_{yy}^{-1}\|_{\tilde{r},y_0}\|\vae\dpr_{y'}^2\wt K\|_{\tilde{r},\mathscr D_\sharp}\le 2\mathsf{T}\frac{|\vae|\mathsf{L}}{4\mathsf{T}_\infty}\le\su2|\vae|\mathsf{L}\le\frac{\s}{6}<\su2
$$
so that
\beq{Jacgy0}
\|\dpr_z G^{y_0}-\uno_d\|_{\tilde{r},y_0}\le 2 \|\vae K_{yy}^{-1}\dpr_{y'}^2\wt K\|_{\tilde{r},y_0}\le |\vae|\mathsf{L}.
\eeq
\noi
Now, we show that the family $\{G^{y_0}\}_{y_0\in\mathscr D_\sharp}$ is compatible so that, together, they define a global map on $D_{\tilde{r}}(\mathscr D_\sharp)$, say $G$ and that, in fact, $G$ is a real--analytic 
 diffeomorphism. For, assume that $z\in D_{\tilde{r}}(y_0)\bigcap D_{\tilde{r}}(\hat y_0)$, for some $y_0,\hat{y}_0\in \mathscr D_\sharp$. Then, we need to show that $G^{y_0}(z)=G^{\hat y_0}(z)$. But,  we have
$$
|G^{\hat y_0}(z)-y_0|\le |G^{\hat y_0}(z)-\hat y_0|+|\hat y_0-z|+|z-y_0|\leby{EcarY1Y0ReExt}\frac{\bar{r}}{2}+\tilde{r}+\tilde{r}<\bar{r}.
$$
Hence, $z\in D_{\tilde{r}}(y_0),\ G^{\hat y_0}(z)\in D_{\bar{r}}(y_0)$ and, by definitions, $F(G^{\hat y_0}(z),z)=0=F(G^{ y_0}(z),z)$. Then, by unicity, we get $G^{y_0}(z)=G^{\hat y_0}(z)$. Thus, the map
$$
G\colon D_{\tilde{r}}(\mathscr D_\sharp)\to\cn\quad \mbox{such that}\quad G_{|D_{\tilde{r}}(y_0)}\coloneqq G^{y_0},\quad \forall\; y_0\in \mathscr D_\sharp\;,
$$
is well--defined and, therefore, is a real--analytic local diffeomorphism. It remains only to check that $G$ is injective 
to conclude that it 
is a global diffeomorphism. Let then $z\in D_{\tilde{r}}(y_0),\hat{z}\in  D_{\tilde{r}}(\hat{y}_0)$ such that $G(z)=G(\hat{z})$, for some $y_0,\hat{y}_0\in\mathscr D_\sharp$. Then, we have 
$$
|z-\hat{z}|< \frac{r}{\mathsf{C}_{5}}-\tilde{r}.
$$
Indeed, if not then
\begin{align*}
0=|G(z)-G(\hat{z})|&\ge -|G(z)-z|+|z-\hat{z}|-|\hat{z}- G(\hat{z})|\\
                   &\geby{gy0zRe}-\bar{r}+\frac{r}{\mathsf{C}_{5}}-\tilde{r}-\bar{r}\\
                   &\ge \frac{r}{\mathsf{C}_{5}}-3\bar{r}\\
                   &\geby{rrbarAsExt}\frac{r}{\mathsf{C}_{5}}-3\frac{r}{4\mathsf{C}_{5}}\\
                   &>0\,,
\end{align*}
contradiction. Therefore,
\beq{zminhatzRe}
|z-\hat{z}|< \frac{r}{\mathsf{C}_{5}}-\tilde{r}\;.
\eeq
Thus,
$$
|\hat{z}-y_0|\le |\hat{z}-z|+|z-y_0|< \frac{r}{\mathsf{C}_{5}}-\tilde{r}+\tilde{r}=\frac{r}{\mathsf{C}_{5}}\,.
$$
Hence, $z,\hat{z}\in D_{r/\mathsf{C}_{5}}(y_0)$. But $G(z)=G(\hat{z})$ is equivalent to $K_y(z)=K_y(\hat{z})$ and then,
$$
0=K_y(z)-K_y(\hat{z})=\dst\int_0^1 K_{yy}(\hat{z}+t(z-\hat{z}))dt(z-\hat{z})\;.
$$
Thus, it is enough to show that $\dst\int_0^1 K_{yy}(\hat{z}+t(z-\hat{z}))dt$ is invertible. But
\begin{align*}
\dst\int_0^1 K_{yy}(\hat{z}+t(z-\hat{z}))dt&= K_{yy}(\hat{z})+\dst\int_0^1\int_0^1 K_{yyy}(\hat{z}+tt'(z-\hat{z}))tdt'dt\cdot(z-\hat{z})\\
      &\eqby{invJacKyyRe} K_{yy}(\hat{z})\left(\uno_d+K_{yy}(\hat{z})^{-1}\dst\int_0^1\int_0^1 K_{yyy}(\hat{z}+tt'(z-\hat{z}))tdt'dt\cdot(z-\hat{z})\right)\\
\end{align*}
and
\begin{align*}
\left\|K_{yy}(\hat{z})^{-1}\dst\int_0^1\int_0^1 K_{yyy}(\hat{z}+tt'(z-\hat{z}))tdt'dt\cdot(z-\hat{z})\right\|&\leby{invJacKyyRe} 2\mathsf{T}\cdot\su2\|K_{yyy}\|_{r/2,y_0}|z-\hat{z}|\\
  &\leby{zminhatzRe} \mathsf{T}\frac{2d\mathsf{K}}{r}\left(\frac{r}{\mathsf{C}_{5}}-\tilde{r}\right)\\
  &<\frac{\mathsf{C}_5}{2r}\frac{r}{\mathsf{C}_{5}}\\
  &=\su2.
\end{align*}
Therefore, $\dst\int_0^1 K_{yy}(\hat{z}+t(z-\hat{z}))dt$ is invertible and then we get $z-\hat{z}=0$ \ie $G$ is injective.\\
\noi
Next, we estimate $P'$. We have, for any 
$|\vae|\leq {\vae_*}$,
\[|\vae|\|g_x\|_{\bar{r},\bar{s},\mathscr D_\sharp}\leq |\vae|\mathsf{C}_1 \frac{M}{\a} \s^{-(\t+d+1)}
 \le |\vae| \bar{r}\mathsf{L}\leby{cond1ExtExt}\frac{r}{3}\frac{\s}{3}\le \frac{r}{3}\]
so that, for any $y_0\in \mathscr D_\sharp$ and $(y',x)\in D_{\bar{r},\bar{s}}(y_0)$,
\[ |y'+\vae g_x(y',x)-y_0|\leq \bar{r}+\frac{r}{3}\leq \frac{r}{2}+\frac{r}{3}=\frac{5r}{6}<r\,,\]
and thus
\begin{align*}
\|P^\ppu\|_{\bar{r},\bar{s},\mathscr D_\sharp}&\leq d^2 \|K_{yy}\|_{r,\mathscr D_\sharp}\|g_x\|_{\bar{r},\bar{s},\mathscr D_\sharp}^2\leq d^2 \mathsf{K}\left( \mathsf{C}_1 \frac{M}{\a} \s^{-(\t+d+1)}\right)^2\\
   &=d^2\mathsf{C}_1^2 \frac{\mathsf{K}M^2}{\a^2} \s^{-2(\t+d+1)}, 
\end{align*}
\begin{align*}
\|P^\ppd\|_{\bar{r},\bar{s},\mathscr D_\sharp}&\leq d\|P_y\|_{\frac{5r}{6},\bar{s},\mathscr D_\sharp}\|g_x\|_{\bar{r},\bar{s},\mathscr D_\sharp}\leq d\frac{6M}{r}\mathsf{C}_1 \frac{M}{\a} \s^{-(\t+d+1)}\\
     &= 6d\mathsf{C}_1 \frac{M^2}{\a r}\s^{-(\t+d+1)}
\end{align*}
 and, by Lemma~$\ref{fce}$--$(i)$, we have
\begin{align*}
|\vae|\|P^\ppt\|_{\bar{r},s-\frac{\s}{2},\mathscr D_\sharp}&\leq \dst\sum_{|n|_1>\k}\|P_n\|_{\bar{r},\mathscr D_\sharp}\ex^{(s-\frac{\s}{2})|n|_1}\leq M\dst\sum_{|n|_1>\k}\ex^{-\frac{\s |n|_1}{2}}\\
  &\leq M\ex^{-\frac{ \k\s}{4}}\dst\sum_{|n|_1>\k}\ex^{-\frac{\s |n|_1}{4}}\leq M\ex^{-\frac{ \k\s}{4}}\dst\sum_{|n|_1>0}\ex^{-\frac{\s |n|_1}{4}}\\
  &= M\ex^{-\frac{ \k\s}{4}} \left(\left(\dst\sum_{k\in \integer}\ex^{-\frac{\s |k|}{4}}\right)^d-1\right)=M\ex^{-\frac{ \k\s}{4}}\left(\left(1+\frac{2\ex^{-\frac{\s }{4}}}{1-\ex^{-\frac{\s }{4}}} \right)^d-1\right)\\
  &= M\ex^{-\frac{ \k\s}{4}}\left(\left(1+\frac{2}{\ex^{\frac{\s }{4}}-1} \right)^d-1\right)\leq M\ex^{-\frac{ \k\s}{4}}\left(\left(1+\frac{2}{\frac{\s }{4}} \right)^d-1\right)\\
  &\leq \s^{-d} M\ex^{-\frac{ \k\s}{4}}\left(\left(\s +8 \right)^d-\s^d\right)\leq d 8^{d}\s^{-d} M\ex^{-\frac{ \k\s}{4}}\\
  &= \mathsf{C}_2\s^{-d} M \ex^{-\frac{\k\s}{4}}\\
  &\leby{DefNArnExt} \mathsf{C}_2 \s^{-d}M \frac{\mathsf{K}|\vae|M}{\a^2}\\
  &=\mathsf{C}_2 \frac{\mathsf{K}|\vae|M^2}{\a^2}\s^{-d}\,.
\end{align*}
Hence,\footnote{Recall that $r\le r_0$ and $\s<1$.}
\begin{align*}
\|P'\|_{\bar{r},\bar{s},\mathscr D_\sharp}&\leq \|P^\ppu\|_{\bar{r},\bar{s},\mathscr D_\sharp}+\|P^\ppd\|_{\bar{r},\bar{s},\mathscr D_\sharp}+\|P^\ppt\|_{\bar{r},\bar{s},\mathscr D_\sharp}\\
  &\leq d^2\mathsf{C}_1^2 \frac{\mathsf{K}M^2}{\a^2} \s^{-2(\t+d+1)}+6d\mathsf{C}_1 \frac{M^2}{\a r}\s^{-(\t+d+1)}+\mathsf{C}_2 \frac{\mathsf{K}M^2}{\a^2}\s^{-d}\\
  &= \left(d^2\mathsf{C}_1^2 r\mathsf{K}+6d\mathsf{C}_1 \a \s^{\t+d+1}+ \mathsf{C}_2 r\mathsf{K}\s^{2\t+d+2}\right)\frac{M^2}{\a^2 r}\s^{-2(\t+d+1)}\\
  &\le \left(d^2\mathsf{C}_1^2+6d\mathsf{C}_1 +\mathsf{C}_2\right)\max(\a,r\mathsf{K})\frac{M^2}{\a^2 r}\s^{-2(\t+d+1)}\\
  &=\mathsf{C}_3 \max(\a,r\mathsf{K})\frac{ M^2}{\a^2 r}\s^{-2(\t+d+1)}
\end{align*}
Now, we need to invert $x \mapsto x+\vae\hat g_{y'}(y',x)$. But, thanks to \equ{gext003}, \equ{gext004}, \equ{gextciest} and \equ{cond1ExtExt}, we can apply Lemma~\ref{CauTypInInv} (see Appendix~\ref{appB}),
to conclude that for any given $y'\in\cn$, the map $\psi_{\vae}(y',\cdot)\colon \torus^d_{\bar{s}}\ni x\mapsto \psi_{\vae}(y',x)$ has an inverse $\f(y',x')=x'+\vae \wt\f(y',x';\vae)$, $\ci$ on $\rn\times\tn$ and real analytic on $D_{\bar{r}/4,s'}(\mathscr D_\sharp)$ such that \equ{boundalExtExtRe} holds and \equ{fi1estextExt} as well, using the multivariate F\`aa Di Bruno formula (see \cite{CS96Mult}, Theorem~$2.1$) and the real--analyticity of $g_x$. 
\qed
\newpage
\noi
Finally, we prove the convergence of the scheme by mimicking Lemma~\ref{lem:2}.
\lem{lem:2Ext} 
Let   $H_0\coloneqq H$, $K_0\coloneqq K$, $P_0\coloneqq P$, $\phi^0=\phi_0\coloneqq\id $, and $r_0$, $s_0$, $s_*$, $\s_0$, $\l_0$,  $\mathsf{W}_0$, $M_0$, $\mathsf{K}_0$, $\mathsf{K}_\infty$, $\mathsf{T}_0$, $\mathsf{T}_\infty$, $\mathsf{d}_*$, $\mathsf{e}_*$, $\mathsf{f}_*$, $\vae_1$, $\vae_2$ and $\vae_\sharp$ be  as in $\S\ref{AssumpExtArnol}$
 and for a given $\vae\not=0$, sequence of non--negative numbers $(M_j)_j$ and $j\ge 0$,  define\footnote{Notice that $s_{j}\downarrow s_*$ and $r_{j}\downarrow 0$.}
\begin{align*}
 \dst\sigma_j&\coloneqq \frac{\sigma_0}{2^j}\,,\\
  s_{j+1}&\coloneqq s_j-\sigma_j=s_*+\frac{\sigma_0}{2^j}\,,\\
  \bar s_{j}&\coloneqq s_j-\frac{2\s_i}{3}\,,\\
 \mathsf{K}_{j+1}&\coloneqq \mathsf K_0\dst\prod_{k=0}^{j}(1+\frac{\s_k}{3})\le \mathsf K_0\ex^{\frac{2\s_0}{3}}\le\mathsf{K}_\infty\,,  \\
\mathsf{T}_{j+1}&\coloneqq \mathsf T_0\dst\prod_{k=0}^{j}(1+\frac{\s_k}{3})\le \mathsf T_0\ex^{\frac{2\s_0}{3}}\le \mathsf{T}_\infty\,,\\
 \dst\l_j&\coloneqq 2^j\l_0= 2^j \log\left(\frac{\a^2 }{\mathsf{K}_0|\vae|M_0}\right)\,,\\
 \k_j &\coloneqq 4\s_j^{-1}\l_j=4^{j}\k_0=4^j \s_0^{-1} \log\left(\frac{\a^2 }{\mathsf{K}_0|\vae|M_0}\right)\,,\\
 r_{j+1}&\coloneqq \su{16\mathsf{C}_5}\dst\min\left(\frac{\mathsf{R}}{\k_j^{\t+1}}\,,\, r_j\s_j \right)\,,\\
 \tilde{r}_{j+1}&\coloneqq \frac{r_{j+1}}{\mathsf{C}_5} \,.
\end{align*}
Assume that $\vae_*$ satisfies
\beq{condExt}
\vae_*\leq \vae_\sharp\,,\quad \vae_*\,\mathsf{f}_*\,\|P\|_{r_0,s_0,\mathscr D_{r_0,\a}}\le 1\quad\mbox{and}\quad 4\vae_*^2\, \mathsf{f}_*\, \mathsf{e}_*\,\mathsf{d}_*^2\, \s_0\, \|P\|_{r_0,s_0,\mathscr D_{r_0,\a}}^2\le 3\,. 
\eeq
	where
\begin{align*}
\mathsf{d}_* &\coloneqq  2^{2\t+2d+3}\mathsf{C}_6^2\,,\\ 
\mathsf{e}_* &\coloneqq  6\dst\max\left\{\frac{4\mathsf{T}_\infty  }{\mathsf{C}_6 r_1^2}\s_0^{-(\t+d+2)}\,,\,\frac{\mathsf{C}_7}{4} \frac{ 1}{\a r_1}\s_0^{-(2\t+2d+3)}\right\}\;,\\
\mathsf{f}_* &\coloneqq  3\dst\max\left\{\frac{4\mathsf{T}_\infty  }{r_0r_1}\s_0^{-(\t+d+2)}\,,\,
\frac{\mathsf{C}_7}{4} \max(\a,r_0\mathsf{K}_\infty)\frac{ 1}{\a^2 r_1}\s_0^{-(2\t+2d+3)}\right\}\;.
\end{align*}
Then, for any $|\vae|\leq \vae_*$, 
one can construct a sequence of  diffeomorphisms 
$$
G_{j+1}\colon  D_{\tilde{r}_{j+1}}(\mathscr D_j){\to}G_{j+1}( D_{\tilde{r}_{j+1}}(\mathscr D_j))
$$
and of $\ci$--symplectomorphisms 
$$
\phi_j\colon \mathscr D\times\tn\overset{into}{\longrightarrow} \mathscr D\times\tn
$$
such that
\begin{align}
&\dpr_{y}K_{j+1}\circ G_{j+1}=\dpr_{y}K_j \;,\label{kjkjpu}\\
&\phi_{j+1}:D_{r_{j+1},s_{j+1}}(\mathscr D_{j+1})\to D_{r_{j},s_{j}}(\mathscr D_{j})\quad\qquad\quad\mbox{is real--analytic},\label{phijExt}\\
&H_{j+1}\coloneqq H_{j}\circ\phi_{j+1}\eqqcolon K_{j+1} + \vae^{2^{j+1}} P_{j+1}\qquad\quad\  \mbox{on } D_{r_{j+1},s_{j+1}}(\mathscr D_{j+1})\label{HjExtExt}
\end{align}
and converge uniformly. More precisely, given any $|\vae|\leq \vae_*$, we have the following:
\begin{itemize}
\item[$(i)$] the sequence $G^{j+1}\coloneqq G_{j+1}\circ G_j\circ\cdots\circ G_1$ converges uniformely on $\mathscr D_{r_0,\a}$ to a lipeomorphism $G_*\colon \mathscr D_{r_0,\a}\to \mathscr D_*\coloneqq G_*(\mathscr D_{r_0,\a})\subset\mathscr D$\;;
\item[$(ii)$] $\vae^{2^j}\dpr_y^{\b} P_j$ converges uniformly on 
$\mathscr D_*\times\dst\torus^d_{s_*}$ to $0$, for any $\b\in \natural_0^d$\;;
\item[$(iii)$] $\phi^j\coloneqq \phi_0\circ\phi_1\circ\phi_2\circ \cdots\circ \phi_j$ converges uniformly on 
$\mathscr D\times\tn$ to a $\ci$--symplectomorphism $\phi_*\colon \mathscr D\times\tn\overset{into}{\longrightarrow} \mathscr D\times\tn$, with  $\phi_*(y,\cdot)\colon \torus^d_{s_*}\ni x\mapsto \phi_*(y,x)$ holomorphic, for any $y\in\mathscr D$\;;
\item[$(iv)$]  $K_j$ converges uniformly on 
$\mathscr D$ to a $\ci$--map $K_*$, with
\begin{align*}
\dpr_{y_*}K_*\circ G_*&=\dpr_{y}K \quad \mbox{on} \quad \mathscr D_{r_0,\a}\;,\\
 \dpr_{y_*}^{\b}H\circ\phi_*(y_*,x)&=\dpr_{y_*}^{\b} K_*(y_*),\ \forall(y_*,x)\in\mathscr D_*\times\tn\;, \forall\;\b\in \natural_0^d\;.
\end{align*}
\end{itemize}
Finally, the following estimates hold for any $|\vae|\leq \vae_*$ and for any $i\ge 1$:
\begin{align}
|\vae|^{2^i}M_i:=|\vae|^{2^i}\|P_i\|_{r_i,s_i,\mathscr D_i}&\le \frac{(\, |2\vae|^2 \mathsf{e}_*{\,d_*^2}M_1)^{2^{i-1}}}{\mathsf{e}_* {\,d_*}^{i+1}}\ ,\label{estfin2Ext01}\\
|\meas(\mathscr D_*)-\meas(\mathscr D_{r_0,\a})|&\le \mathsf{C}_8 \vae_2 \ex^{\vae_2}\meas(\mathscr D_{r_0,\a})\ ,\label{estfin2Ext02}\\
|\mathsf{W}(\phi_*-\id)|
&\le 
\vae_1
\qquad\qquad\qquad\qquad\ \mbox{on}\quad \mathscr D_*\times\torus^d_{s_*}\label{estfin2Ext03}\ .
\end{align}
\elem 
\proof
 For $i\ge 0$, define
\begin{align*}
 \mathsf{W}_i&\coloneqq \diag\left(\su{4r_{i+1}}\uno_d,\uno_d\right)\,,\\ 
 \ovl{\mathsf{L}}_i&\coloneqq \frac{\mathsf{C}_4}{4} \max\{\a,r_i\mathsf{K}_\infty\}\frac{M_i}{\a^2 r_{i+1}}\s_i^{-(2\t+d+2)}\,,\\
 \mathsf  L_i&\coloneqq  M_i\dst\max\left\{\frac{4\mathsf{T}_\infty  }{r_ir_{i+1}}\s_i^{-(\t+d+1)}\,,\,
 \frac{\mathsf{C}_7}{4} \max\{\a,r_i\mathsf{K}_\infty\}\frac{ 1}{\a^2 r_{i+1}}\s_i^{-2(\t+d+1)}\right\}\\
             &\ge  M_i\dst\max\left\{\frac{4\mathsf{T}_\infty  }{r_ir_{i+1}}\s_i^{-(\t+d+1)}\,,\,
             \frac{4}{\mathsf{K}_i r_i^2}\,,\,\frac{\mathsf{C}_7}{4} \max\{\a,r_i\mathsf{K}_i\}\frac{ 1}{\a^2 r_{i+1}}\s_i^{-2(\t+d+1)}\right\}
 \, .
 \end{align*}
Let us assume ({\sl inductive hypothesis}) that we can iterate $j\ge 1$ times the KAM step, obtaining $j$ diffeomorphisms
$$
G_{i+1}\colon  D_{\tilde{r}_{i+1}}(\mathscr D_{i}){\to}G_{i+1}( D_{\tilde{r}_{i+1}}(\mathscr D_{j}))
$$
and $j$ $\ci$--symplectomorphisms 
\beq{bes06Ext}
\phi_{i+1}\colon \mathscr D\times\tn\overset{into}{\longrightarrow} \mathscr D\times\tn, 
\eeq
satisfying $\equ{kjkjpu}_{j=i}\div\equ{HjExtExt}_{j=i}$, for $0\le i\le j-1$, with
\beq{bbbExt}
\left\{
\begin{array}{l}
\|\dpr_y^2 K_i\|_{r_i,\mathscr D_i}\le \mathsf{K}_i\, ,\ \\  \ \\
\|T_i\|_{\mathscr D_i}\le \mathsf{T}_i\,,\ \\  \ \\ 
 \|P_{i}\|_{r_{i},s_{i},\mathscr D_{i}}\le M_{i}\,,\ \\  \ \\ 
|\vae|^{2^i} \mathsf{L}_i \le \frac{\sigma_i}{3}\ \\  \ \\ 
\l_i\ge \log\left(\frac{\a^2}{\mathsf{K}_0|\vae|^{2^i}M_i}\right)
\ . 
\end{array}\right.
\eeq
Observe that for $j=1$, it is $i=0$ and \equ{bbbExt} is implied by the definitions of $\mathsf{K}_0,\,\mathsf{T}_0,\, M_0$  and by condition \equ{condExt}.

\noi
Because of \equ{condExt} and \equ{bbbExt}, \equ{cond1ExtExt} holds 
for $H_i$ and Lemma~\ref{lem:1Ext}  can be applied to $H_i$ and  one has, for $0\le i\le j-1$ and for any $|\vae|\le \vae_*$ (see \equ{Est1Lem1bExt}, \equ{convEstExtRe} and \equ{tesitExtExtRe}): 
\begin{align}
\| G_{i+1}-\id\|_{\tilde{r}_{i+1},\mathscr D_i}&\le 2r_{i+1}\;, \label{C.1Exteq1}\\
\|\dpr_z G_{i+1}-\uno_d\|_{\tilde{r}_{i+1},\mathscr D_i}&\le |\vae|^{2^i}\mathsf L_i\;, \label{C.1Exteq2} \\
\|\mathsf{K}_{i+1}\|_{r_{i+1},\mathscr D_{i+1}}&\le \|\mathsf{K}_i\|_{r_i,\mathscr D_{i}}+|\vae|^{2^i}M_i \;, \label{C.1Exteq3}\\
\|\dpr_y^2\mathsf{K}_{i+1}\|_{r_{i+1},\mathscr D_{i+1}}&\le \|\dpr_y^2\mathsf{K}_i\|_{r_i,\mathscr D_{i}}+\mathsf K_i|\vae|^{2^i}\mathsf L_i \;,\label{C.1Exteq4}\\
\|T_{i+1}\|_{\mathscr D_{i+1}}&\le \|T_i\|_{\mathscr D_i}+\mathsf T_i|\vae|^{2^i}\mathsf L_i \;, \label{C.1Exteq5}\\
 \|\mathsf{W}_i(\phi_{i+1}-\id)\|_{r_{i+1},s_{i+1},\mathscr D_{i+1}}&\le |\vae|^{2^i}\ovl{\mathsf L}_i\;,\label{C.1Exteq6} \\
\|P_{i+1}\|_{r_{i+1},s_{i+1},\mathscr D_{i+1}}&\le M_{i+1}\coloneqq M_i \mathsf L_i\;. \label{C.1Exteq7}
\end{align}
\noi
Let $0\le i\le j-1$. Then, by definition,
\beq{rjKal}
r_1=\frac{\mathsf{R}}{16\mathsf{C}_5\k_0^{\t+1}}\le \frac{\mathsf{R}}{16\mathsf{C}_5}< \frac{\a}{\mathsf{K}_\infty} \,.
\eeq
and, since $\s_0<2^{-2(\t-1)}\mathsf{C}_5\sqrt{2}$ \ie 
\beq{sig0rj}
\frac{16\mathsf{C}_5\sqrt{2}}{2^{2(\t+1)}\s_0}>1\;,
\eeq
we have
\begin{align*}
r_{i+1}&= \dst\min\left(\frac{\mathsf{R}}{16\mathsf{C}_5\k_i^{\t+1}}\,,\,\frac{r_i\s_i}{16\mathsf{C}_5}  \right)\\
	   &= \dst\min\left(\frac{\mathsf{R}}{16\mathsf{C}_5\k_i^{\t+1}}\,,\,\frac{\mathsf{R}\s_i}{(16\mathsf{C}_5)^2\k_{i-1}^{\t+1}}\,,\,\frac{r_{i-1}\s_{i-1}\s_i}{(16\mathsf{C}_5)^2}  \right)\\
	   &\,\ \vdots \\
	   &=\dst\min\left(\frac{\mathsf{R}}{16\mathsf{C}_5\k_i^{\t+1}}\,,\,\frac{\mathsf{R}\s_i}{(16\mathsf{C}_5)^2\k_{i-1}^{\t+1}}\,,\,\cdots\,,\,\frac{r_1\s_1\cdots\s_i}{(16\mathsf{C}_5)^i}  \right)\\
	   &= \dst\min\left(\frac{\mathsf{R}}{16\mathsf{C}_5\k_i^{\t+1}}\,,\,\frac{\mathsf{R}\s_i}{(16\mathsf{C}_5)^2\k_{i-1}^{\t+1}}\,,\,\cdots\,,\,\frac{\mathsf{R}\s_1\cdots\s_i}{(16\mathsf{C}_5)^{i+1}\k_0^{\t+1}}  \right)\\
	   &=\frac{\mathsf{R}\s_1\cdots\s_i}{\left(16\mathsf{C}_5\right)^{i+1}\k_0^{\t+1}}\dst\min\left(\left(\frac{16\mathsf{C}_5\cdot 2^{\su2 (i+1)}}{2^{2(\t+1)}\s_0}\right)^i\,,\,\left(\frac{16\mathsf{C}_5\cdot 2^{\su2 i}}{2^{2(\t+1)}\s_0}\right)^{i-1}\,,\,\cdots\,,\, \left(\frac{16\mathsf{C}_5\cdot 2^{\su2}}{2^{2(\t+1)}\s_0}\right)^0\right) \\
	   &\eqby{sig0rj} \frac{\mathsf{R}\s_1\cdots\s_i}{\left(16\mathsf{C}_5\right)^{i+1}\k_0^{\t+1}}\\
	   &=2^{-\frac{i^2}{2}}\mathsf{C}_6^{-i} r_1\,.
\end{align*}
Thus,
\begin{align*}
|\vae| \mathsf L_0 (3 \sigma_0^{-1})&= 3|\vae| M_0\dst\max\left(\frac{4\mathsf{T}_\infty  }{r_0r_1}\s_0^{-(\t+d+2)}\,,\,
\frac{\mathsf{C}_7}{4} \max(\a,r_0\mathsf{K}_\infty)\frac{ 1}{\a^2 r_1}\s_0^{-(2\t+2d+3)}\right)\\
&=\mathsf{f}_*|\vae|M_0\leby{condExt}1
\end{align*}
and for $i\ge 1$\footnote{Notice that $2^i\ge i^2-1,\,\forall\; i\in\natural_0$.},
\begin{align*}
|\vae|^{2^i} \mathsf L_i (3 \sigma_i^{-1})&\eqby{rjKal} 3|\vae|^{2^i} M_i\dst\max\left(\frac{4\mathsf{T}_\infty  }{r_ir_{i+1}}\s_i^{-(\t+d+2)}\,,\,\frac{\mathsf{C}_7}{4} \frac{ 1}{\a r_{i+1}}\s_i^{-(2\t+2d+3)}\right)\\
&=3|\vae|^{2^i} M_i\dst\max\left(\frac{4\mathsf{T}_\infty\sqrt{2}  }{\mathsf{C}_6 r_1^2}\s_0^{-(\t+d+2)}\;2^{i^2}\left(2^{\t+d+1}\mathsf{C}_6^2\right)^i\,,\,\right.\\
&\hspace{5.85cm}\left.\frac{\mathsf{C}_7}{4} \frac{ 1}{\a r_1}\s_0^{-(2\t+2d+3)}\;2^{\frac{i^2}{2}}\left(2^{2\t+2d+3}\mathsf{C}_6\right)^i\right)\\
&\le 3\left(2^{2\t+2d+3}\mathsf{C}_6^2\right)^i\;2^{i^2}|\vae|^{2^i} M_i\dst\max\left(\frac{4\mathsf{T}_\infty  }{\mathsf{C}_6 r_1^2}\s_0^{-(\t+d+2)}\,,\,\frac{\mathsf{C}_7}{4} \frac{ 1}{\a r_1}\s_0^{-(2\t+2d+3)}\right)\\
&=\mathsf{e}_*\mathsf{d}_*^i \;2^{i^2-1}|\vae|^{2^i}M_i\\
&\le \mathsf{e}_*\mathsf{d}_*^i |2\vae|^{2^i}M_i\eqqcolon \frac{\theta_i}{\mathsf d_*}\;,
\end{align*}
 so that 
 $$
 \mathsf L_i<\mathsf{e}_*\,\mathsf d_*^i  M_i\;,
 $$ 
 thus by \equ{C.1Exteq7}, for any $1\le i\le j-1$,  
 $$
 |\vae|^{2^{i+1}}M_{i+1}<\mathsf{e}_*\mathsf d_*^i |\vae|^{2^{i+1}} M_i^2
 $$ 
 \ie $\theta_{i+1}<\theta_i^2$ , which iterated, yields 
 $\theta_i\le \theta_1^{2^{i-1}}$ for $1\le i\le j$. 
Next, we  show that, thanks to \equ{condExt}, \equ{bbbExt} holds also for $i=j$. In fact, by \equ{bbbExt} and \equ{C.1Exteq5},  
 we have 
$$
\|T_{i+1}\|_{\mathscr D_{i+1}}\le \|T_i\|_{\mathscr D_{i+1}}+\mathsf T_i|\vae|^{2^i}\mathsf L_i\le \mathsf T_i+\mathsf T_i\frac{\s_i}{3}=\mathsf T_{i+1}\,.
$$
and similarly for  $\|\dpr_y^2 K_{i+1}\|_{r_{i+1},\mathscr D_{i+1}}$ . Now, we  check the last relation in \equ{bbbExt} for $i=j$. 
 But, by definitions, for any $i\ge 0$, 
$$
M_{i+1}=M_i\mathsf{L}_i\ge M_i\max\{\a,r_i\mathsf{K}_\infty\}\frac{ M_i}{\a^2 r_{i+1}}\ge \frac{M_i^2 \mathsf{K}_0}{\a^2}\,,
$$
\ie
$$
\frac{|\vae|^{2^{i+1}}M_{i+1}\mathsf{K}_0}{\a^2}\ge \left(\frac{|\vae|^{2^{i}}M_i\mathsf{K}_0}{\a^2}\right)^2,
$$
which iterated yields, for any $i\ge 0$,
$$
\frac{|\vae|^{2^{i}}M_i\mathsf{K}_0}{\a^2}\ge \left(\frac{|\vae|M_0\mathsf{K}_0}{\a^2}\right)^{2^i}\,.
$$
\ie
$$
\l_i=\log\left( \left(\frac{\a^2}{|\vae|M_0\mathsf{K}_0}\right)^{2^i}\right)\ge \log\left( \frac{\a^2}{|\vae|^{2^{i}}M_i\mathsf{K}_0}\right)\,.
$$
\nl 
Now, by $\equ{estfin2Ext01}_{i=j}$,
$$|\vae|^{2^j}\mathsf L_j (3 \sigma_j^{-1})\le\frac{\theta_j}{\mathsf{d}_*}\le \su{\mathsf{d}_*}\theta_1^{2^{j-1}}\le\su{\mathsf{d}_*}(4\mathsf{e}_*{\,d_*^2}\vae_*^2 M_1)^{2^{j-1}}
\leby{condExt} \su{\mathsf{d}_*}<1\ ,$$
which implies the fourth inequality in \equ{bbbExt} with $i=j$; the proof of the induction is finished and one can construct an {\sl infinite sequence} of diffeomorphisms $G_{i+1}\colon  D_{\tilde{r}_{i+1}}(\mathscr D_i){\to}G_{i+1}( D_{\tilde{r}_{i+1}}(\mathscr D_i))$ and symplectomorphisms $\phi_i\colon \mathscr D\times\tn \righttoleftarrow$ satisfying \equ{bbbExt}, $\equ{C.1Exteq1}\div\equ{C.1Exteq7}$, \equ{estfin2Ext01} and $\equ{kjkjpu}_{j=i}\div\equ{HjExtExt}_{j=i}$ {\sl for all $i\ge 0$}. \\
\nl
Next, we show that $G^j$ converges. For any $j\ge 1$,
\begin{align*}
\|G^{j+1}-G^j\|_{\mathscr D_{r_0,\a}}= \|G_{j+1}\circ G^j-G^j\|_{\mathscr D_{r_0,\a}}\le \|G_{j+1}-\id\|_{\mathscr D_{j}}\le \|G_{j+1}-\id\|_{\tilde{r}_{j+1},\mathscr D_{j}}\leby{C.1Exteq1} 2r_{j+1}.
\end{align*}
Thus, $G^j$ is Cauchy and therefore converges uniformly on $\mathscr D_{r_0,\a}$ to a map $G_*$.

\nl
Next, we prove that $\phi_j$ is convergent by showing that it is Cauchy as well. For any $j\ge 3$, we have, using again Cauchy's estimate,
\beqano
\|\mathsf{W}_{j-1}(\phi^j-\phi^{j-1})\|_{r_j,s_j,\mathscr D_j}&=&\|\mathsf{W}_{j-1}\phi^{j-1}\circ\phi_j-\mathsf{W}_{j-1}\phi^{j-1}\|_{r_i,s_i,\mathscr D_i}\\
           &\leby{bes06Ext} &\|\mathsf{W}_{j-1}D\phi^{j-1}\mathsf{W}_{j-1}^{-1}\|_{2r_{j-1}/3, s_{j-1},\mathscr D_{j-1}}\, \|\mathsf{W}_{j-1}(\phi_j-\id)\|_{r_j,s_j,\mathscr D_j}\\
           &\leby{C.1Exteq6} & \max\left(r_{j-1}\frac{3}{r_{j-1}},\frac{3}{2\s_{j-1}}\right)    \|\mathsf{W}_{j-1}\phi^{j-1}\|_{r_{j-1}, s_{j-1},\mathscr D_{j-1}} \cdot |\vae|^{2^j}\ovl{\mathsf{L}}_j\\
           &= & \frac{3}{2\s_{i-1}}    \|\mathsf{W}_{j-1}\phi^{j-1}\|_{r_{j-1}, s_{j-1},\mathscr D_{j-1}} \cdot |\vae|^{2^j}\ovl{\mathsf{L}}_j\\
           &\le & \frac{3}{2\s_{j-1}}    \|\mathsf{W}_{j-1}\phi^{0}\|_{r_{0}, s_0,\mathscr D_{r_0,\a}} \cdot |\vae|^{2^j}\ovl{\mathsf{L}}_j\\
           &\le & \frac{3}{2\s_{j-1}}\left(\dst\prod_{i=0}^{j-2}\|\mathsf{W}_{i+1}\mathsf{W}_{i}^{-1}\| \right)\|\mathsf{W}_{0}\phi_{0}\|_{r_{0}, s_{0},\mathscr D_{r_0,\a}} \cdot |\vae|^{2^j}\ovl{\mathsf{L}}_j\\
           &=& \frac{3}{2\s_{j-1}}\left(\dst\prod_{i=0}^{j-2}\frac{r_i}{r_{i+1}} \right)\|\mathsf{W}_{0}\phi_{0}\|_{r_{0}, s_{0},y_{0}} \cdot |\vae|^{2^j}\ovl{\mathsf{L}}_j\\
           &=& \frac{3r_0}{2r_{j-1}\s_{j-1}}\|\mathsf{W}_{0}\phi_{0}\|_{r_{0}, s_{0},\mathscr D_{r_0,\a}} \cdot |\vae|^{2^j}\ovl{\mathsf{L}}_j\;.
\eeqano
\newpage
Therefore, for any $n\geq 0,\, j\geq 1$,
\begin{align*}
\|\mathsf{W}_{0}(\phi^{n+j}-\phi^n)\|_{r_{n+j},s_{n+j},\mathscr D_{n+j}}&\leq  \sum_{i=n}^{n+j}\|\mathsf{W}_{0}(\phi^{i+1}-\phi^i)\|_{r_{i+1},s_{i+1},\mathscr D_{i+1}}\\
&\le \sum_{i=n}^{n+j}\left(\dst\prod_{k=0}^{i}\|\mathsf{W}_{k}\mathsf{W}_{k+1}^{-1}\| \right)\|\mathsf{W}_{i}(\phi^{i+1}-\phi^i)\|_{r_{i+1},s_{i+1},\mathscr D_{i+1}}\\
&= \sum_{i=n}^{n+j}\left(\dst\prod_{k=0}^{i}\frac{r_{k+1}}{r_k} \right)\|\mathsf{W}_{i}(\phi_{i+1}-\phi_i)\|_{r_{i+1},s_{i+1},\mathscr D_{i+1}}\\
&= \sum_{i=n}^{n+j}\frac{r_{i+1}}{r_0}\|\mathsf{W}_{i}(\phi^{i+1}-\phi^i)\|_{r_{i+1},s_{i+1},\mathscr D_{i+1}}\\
&\leq \su2\|\mathsf{W}_{0}\phi_{0}\|_{r_{0}, s_{0},\mathscr D_{r_0,\a}} \sum_{i=n}^{n+j}\frac{r_{i+1}}{ r_i}|\vae|^{2^{i+1}}\ovl{\mathsf{L}}_{i+1}3\s_i^{-1}\\
&\leq \su2\|\mathsf{W}_{0}\phi_{0}\|_{r_{0}, s_{0},\mathscr D_{r_0,\a}} \sum_{i=n}^{n+j}\frac{r_{i+1}}{ r_i}|\vae|^{2^{i+1}}{\mathsf{L}}_{i+1}3\s_i^{-1}\\
&\leq \su2\|\mathsf{W}_{0}\phi_{0}\|_{r_{0}, s_{0},\mathscr D_{r_0,\a}} \sum_{i=n}^{n+j}\theta_{i+1}\\
&= \su2\|\mathsf{W}_{0}\phi_{0}\|_{r_{0}, s_{0},\mathscr D_{r_0,\a}} \sum_{i=n}^{n+j}\theta_1^{2^i} .
\end{align*}
Hence $\phi_j$ converges uniformly on $\mathscr D_*\times\torus^d_{s_*}\times(-\vae_*,\vae_*)$ to some $\phi_*$, which is then real--analytic function on $\mathscr D_*\times\torus^d_{s_*}\times(-\vae_*,\vae_*)$.

\nl
To estimate $\|\mathsf{W}_0(\phi_*-\id)\|_{\mathscr D_*\times \torus^d_{s_*}}$, observe that
, for $i\ge 1$,\footnote{Notce that $2^{i-1}\ge i,\, \forall\, i\ge 0$.}
$$|\vae|^{2^i}\mathsf L_i=\frac{\sigma_0}{3 \cdot 2^i}\ \mathsf{e}_*{\,d_*}^i |\vae|^{2^i} M_i < \su{3\cdot 2^{i}d_*}(|2\vae|^2 \mathsf{e}_*{\,d_*^2}M_1)^{2^{i-1}}\le \su{3d_*} \Big(\frac{4|\vae|^2\mathsf{e}_*{\,d_*^2}M_1}2\Big)^{i}$$
and therefore 
$$\sum_{i\ge 1}  |\vae|^{2^i}\mathsf L_i\le \frac{1}{3d_*}\sum_{i\ge 1} \Big(\frac{4|\vae|^2\mathsf{e}_*{\,d_*^2}M_1}2\Big)^{i}\le \frac{1}{3} 4|\vae|^2\mathsf{e}_*\mathsf{d}_* M_1 \le \su{3\mathsf{d}_*}
\ .$$ 
Moreover,
\begin{align*}
\|\mathsf{W}_0(\phi^i-\id)\|_{r_i,s_i,\mathscr D_i}&\le \|\mathsf{W}_0(\phi^{i-1}\circ\phi_i-\phi_i)\|_{r_i,s_i,\mathscr D_i}+\|\mathsf{W}_0(\phi_i-\id)\|_{r_i,s_i,\mathscr D_i}\\
&\le \|\mathsf{W}_0(\phi^{i-1}-\id)\|_{r_{i-1},s_{i-1},\mathscr D_{i-1}}+ \left(\dst\prod_{j=0}^{i-2}\|\mathsf{W}_{j}\mathsf{W}_{j+1}^{-1}\| \right) \|\mathsf{W}_{i-1}(\phi_i-\id)\|_{r_i,s_i,\mathscr D_i}\\
&= \|\mathsf{W}_0(\phi^{i-1}-\id)\|_{r_{i-1},s_{i-1},\mathscr D_{i-1}}+ \left(\dst\prod_{j=0}^{i-2}\frac{r_{j+1}}{r_j} \right) \|\mathsf{W}_{i-1}(\phi_i-\id)\|_{r_i,s_i,\mathscr D_i}\\
&= \|\mathsf{W}_0(\phi^{i-1}-\id)\|_{r_{i-1},s_{i-1},\mathscr D_{i-1}}+ \frac{r_{i-1}}{r_0}\|\mathsf{W}_{i-1}(\phi_i-\id)\|_{r_i,s_i,\mathscr D_i}\\
&\le \|\mathsf{W}_0(\phi^{i-1}-\id)\|_{r_{i-1},s_{i-1},\mathscr D_{i-1}}+|\vae|^{2^{i-1}}\ovl{\mathsf{L}}_{i-1}\ ,
\end{align*}
which iterated yields 
\begin{align*}
\|\mathsf{W}_0(\phi^i-\id)\|_{r_i,s_i,\mathscr D_i}&\le \dst\sum_{k=0}^{i-1} |\vae|^{2^k}\ovl{\mathsf{L}}_k\\
&\le |\vae|\ovl{\mathsf{L}}_0+\dst\sum_{k\ge 1} |\vae|^{2^k}{\mathsf{L}}_k\\
&\le |\vae|\ovl{\mathsf{L}}_0 +\frac{4}{3} |\vae|^2\mathsf{e}_*\mathsf{d}_* M_1\\
&= |\vae|\mathsf{C}_0 \max(\a,r_0\mathsf{K}_\infty)\frac{M_0}{\a^2 r_0}\s_0^{-(2\t+d+2)}+\frac{4}{9} |\vae|^2\mathsf{f}_*\mathsf{e}_*\mathsf{d}_* \s_0 M_0^2\\
&\le \su3 |\vae|M_0\mathsf{f}_*\s_0^{d+1}+\frac{4}{9} |\vae|^2\mathsf{f}_*\mathsf{e}_*\mathsf{d}_* \s_0 M_0^2\\
&= \vae_1
\,.
\end{align*}
Therefore, taking the limit over $i$ completes the proof of \equ{estfin2Ext03}.

\noi
Next, we show that $\|G_*-\id\|_{L,\mathscr D_{r_0,\a}}<1$, which will imply that\footnote{See Proposition II.2. in \cite{zehnder2010lectures}.} $G_*\colon \mathscr D_{r_0,\a}\overset{onto}{\longrightarrow}\mathscr D_*$ is a lipeomorphism. Indeed, for any $j\ge 1$, there exists $\hat{r}_j>0$ such that the restricted maps $G_i\colon G^{i-1}(D_{\hat{r}_j}(\mathscr D_{r_0,\a}))\to\complex$, $1\le i\le j$ with $G^0\coloneqq \id$, are well--defined\footnote{\ie $G^{i-1}(D_{\hat{r}_j}(\mathscr D_{r_0,\a}))\subset D_{\tilde{r}_i}(\mathscr D_{i-1})=\dom(G_i)$, $\forall\;1\le i\le j$.} and, therefore,
\begin{align*}
\|\dpr_z G^j-\uno_d\|_{\hat{r}_j,\mathscr D_{r_0,\a}}&\le  \|\dpr_z G^j-\dpr_z G^{j-1}\|_{\hat{r}_j,\mathscr D_{r_0,\a}}+\|\dpr_z G^{j-1}-\uno_d\|_{\hat{r}_j,\mathscr D_{r_0,\a}}\\
  &= \|\dpr_z G_j\circ\dpr_z G^{j-1} -\dpr_z G^{j-1}\|_{\hat{r}_j,\mathscr D_{r_0,\a}}+\|\dpr_z G^{j-1}-\uno_d\|_{\hat{r}_j,\mathscr D_{r_0,\a}}\\
  &\le \|\dpr_z G_j-\uno_d\|_{\tilde{r}_j,\mathscr D_{j-1}}\|\dpr_z G^{j-1} \|_{\hat{r}_j,\mathscr D_{r_0,\a}}+\|\dpr_z G^{j-1}-\uno_d\|_{\hat{r}_j,\mathscr D_{r_0,\a}}\\
  &\le \|\dpr_z G_j-\uno_d\|_{\tilde{r}_j,\mathscr D_{j-1}}(\|\dpr_z G^{j-1}-\uno_d \|_{\hat{r}_j,\mathscr D_{r_0,\a}}+1)+\|\dpr_z G^{j-1}-\uno_d\|_{\hat{r}_j,\mathscr D_{r_0,\a}}\\
  &= (\|\dpr_z G_j-\uno_d\|_{\tilde{r}_j,\mathscr D_{j-1}}+1)(\|\dpr_z G^{j-1}-\uno_d \|_{\hat{r}_j,\mathscr D_{r_0,\a}}+1)-1\\
  &\leby{C.1Exteq2} (|\vae|^{2^{j-1}}\mathsf{L}_{j-1}+1)(\|\dpr_z G^{j-1}-\uno_d \|_{\hat{r}_j,\mathscr D_{r_0,\a}}+1)-1
\end{align*}
which iterated leads to\footnote{Recall that $\ex^x-1\le x\ex^x\;,\ \forall\; x\ge 0$.}
\begin{align*}
\|\dpr_z G^j-\uno_d\|_{\hat{r}_j,\mathscr D_{r_0,\a}} &\le -1+\dst\prod_{i=1}^\infty (|\vae|^{2^{i-1}}\mathsf{L}_{i-1}+1)\\
     &\le -1+\exp\left( \sum_{i=0}^\infty |\vae|^{2^{i}}\mathsf{L}_{i}\right)\\
     &= -1+\exp\left(|\vae|\mathsf{L}_0+ \sum_{i=1}^\infty |\vae|^{2^{i}}\mathsf{L}_{i}\right)\\
     &\le -1+\exp\left(\frac{\s_0}{3}\mathsf{f}_*|\vae|M_0+\frac{1}{3} 4|\vae|^2\mathsf{e}_*\mathsf{d}_* M_1 \right)\\
     &=-1+\ex^{\vae_2}\\
     &\le \vae_2 \ex^{\vae_2}\\
     &\leby{condExt} \left(\frac{\s_0}{3}+\su{3\mathsf{d}_*}\right)\exp\left(\frac{\s_0}{3}+\su{3\mathsf{d}_*}\right)\\
     &\le\left(\frac{1}{6}+\su{6}\right)\exp\left(\frac{1}{6}+\su{6}\right)\\
     &=\su3\exp\left(\su3\right)<1.\\
\end{align*}
Thus, $G_*$ is Lipschitz continuous, with
$$
\|G_*-\id\|_{L,\mathscr D_{r_0,\a}}\le \vae_2 \ex^{\vae_2}\le \su3\exp\left(\su3\right)<1,
$$
so that, by\footnote{With $\d\coloneqq\vae_2 \ex^{\vae_2}\le \su3\exp\left(\su3\right)$.} Lemma~\ref{LebLipLem} (see Appendix~\ref{appC}), we get
\begin{align*}
|\meas(\mathscr D_*)-\meas(\mathscr D_{r_0,\a})|&\le  \left( \left(1+\su3\exp\left(\su3\right)\right)^d-1 \right)\vae_2 \ex^{\vae_2} \meas(\mathscr D_{r_0,\a})\\
   &= \mathsf{C}_8 \vae_2 \ex^{\vae_2} \meas(\mathscr D_{r_0,\a}),
\end{align*}
which proves \equ{estfin2Ext02}, Lemma~\ref{lem:2Ext} and, whence, the extension Theorem.  \qed

%
\chapter{A ``sharp'' version of Arnold's theorem \label{chpArWthLog}}
\section{Assumptions\label{AssumpArnolv2}}
Let $r_0>0,\,\t\ge d-1,\, 0<s_*<s_0\leq 1,\, y_0\in\rn$ and consider the hamiltonian parametrized by $\vae\in\real$
\[H_0(y,x;\vae)\coloneqq K_0(y)+\vae P_0(y,x),\]
with 
$$
K_0,P_0\in \mathcal{B}_{r_0,s_0}(y_0)\,.
$$
 such that
\beq{ArnoldCondv2}
 \det (\dpr^2_y K_0(y_0))\not= 0.
\eeq
Set
\[
T\coloneqq \dpr^2_y K_0(y_0)^{-1},\quad M_0\coloneqq \|P\|_{r_0,s_0,y_0},\quad\mathsf{K}_0\coloneqq \|\dpr^2_y K_0\|_{r_0,y_0},\quad \mathsf{T}_0\coloneqq \|T\|\,.
\]
Finally, define\footnote{Recall from footnote\textsuperscript{\ref{ftnarc1}} that $\mathsf{C}_0, \mathsf{C}_1>1$.}
\begin{align*}
\n			 &\coloneqq \t+1\;,\\
\hat{s}		 &\coloneqq s_0-s_* \;,\\
\s_0         &\coloneqq \frac{\hat{s}}{2}\;,\\
\eta_0 	     &\coloneqq \mathsf{T}_0\mathsf{K}_0\;,\\
\mathsf{C}_0 &\coloneqq 4\sqrt{2}\left(\frac{3}{2}\right)^{2\n+d}\dst\int_{\rn} \left( |y|_1^{\n}+d|y|_1^{2\n}\right)\ex^{-|y|_1}dy\;,\\
\mathsf{C}_1 &\coloneqq 2\left(\frac{3}{2}\right)^{\n+d}\dst\int_{\rn} |y|_1^{\n}\ex^{-|y|_1}dy\;,\\
\mathsf{C}_2 &\coloneqq 2^{3d}d\;,\\
\mathsf{C}_3 &\coloneqq	\left(d^2\mathsf{C}_1^2+6d\mathsf{C}_1 +\mathsf{C}_2\right)\sqrt{2}\;,\\
\mathsf{C}_4 &\coloneqq \max\left\{\mathsf{C}_0,\,\mathsf{C}_3\right\}\;,\\
\mathsf{C}_6 &\coloneqq \dst{\max}\left\{2^{2\n}\,,\,\frac{3\cdot 2^5d}{5}\right\}\;,\\
\mathsf{C}_7 &\coloneqq 3d\cdot 2^{6\n+2d+3}\sqrt{2}\dst\max\left\{640d^2\,,\,\mathsf{C}_4 \right\}\;,\\
\mathsf{C}_8 &\coloneqq \left(2^{-d}\mathsf{C}_6\right)^{\su8}\;,\\
\mathsf{C}_{9} &\coloneqq \frac{\mathsf{C}_6\mathsf{C}_7\mathsf{C}_8}{2^{2\n+7}d}\;,\\
\mathsf{C}_{10} &\coloneqq 3\cdot 2^d\; {\mathsf{C}_8}\;,\\
\m_*  		    &\coloneqq \max\left\{0<\m\le \ex^{-1}\; :\;\mathsf{C}_7\;\mathsf{C}_8\;\eta_0^{\frac{17}{8}}\;\s_0^{-(4\n+2d+1)}\;\m\left(\log\m^{-1} \right)^{2\n}< 1\right\}\;.
\end{align*}
\noi
\section{Statement of the KAM Theorem}
\thm{teo4v2}
Under the assumptions in $\S\ref{AssumpArnolv2}$, the following holds.
Let
\beq{smcondwhLvep}
\a\le \frac{r_0}{\mathsf{T}_0}
\eeq
and assume that
\beq{ArnoldCondv2prim}
\o\coloneqq \dpr_y K_0(y_0)\in \D^\t_\a,\quad \ie \quad |\o\cdot k|\ge \frac{\a}{|k|_1^\t},\,\forall \; k\in\zn\setminus\{0\}\;.
\eeq
Assume
\beq{smcondwhL}
|\vae|\le \m_*\;\frac{\a^2}{\mathsf{K}_0\;M_0}\;.
\eeq
Then, there exist $y_*\in \rn$ and an embedding $\phi_*\colon \tn\to D_{r_0,s_0}(y_0),$ real--analytic on $\torus^d_{s_*}$ and close to the trivial embedding 
\[\phi_0\colon x\in \tn \to (y_*,x)\in D_{r_0,s_0}(y_0),\]
and such that the $d$--torus
\beq{KronTorArnv2}
\mathcal{T}_{\o,\vae}\coloneqq \phi_*\left(\tn\right)
\eeq
is a non-degenerate invariant Kronecker torus for $H$ \ie
\beq{KronTorArnIEv2}
\phi^t_H\circ \phi_*(x)=\phi_*(x+\o t).
\eeq
Moreover, 
\begin{align}
|y_{*}-y_0|&\le \su{\mathsf{C}_{9}}\s_0^{3\n+2d+1}\frac{\a}{\mathsf{K}_0\eta_0^{\frac{17}{8}}} \;,\label{dysty0}\\
|\mathsf{W}(\phi_*-\id)|&\le \su{\mathsf{C}_{10}\;\eta_0^{\frac{1}{8}}}\;,\label{estArnTrv2}
\end{align}
uniformly on  $\{y_*\}\times \torus^d_{s_*}$ , where
$$
\mathsf{W}   \coloneqq \diag\left(\frac{\mathsf{K}_0}{{\a}}\uno_d,\uno_d\right)\;.
$$
\ethm
%
%
%

\section{Proof of Theorem~\ref{teo4v2}}
\lemtwo{lem:1bisv2}{KAM step}
Let $r>0,\,0<2\s<s\leq 1$ and consider the hamiltonian parametrized by $\vae\in\real$
$$
H(y,x;\vae)\coloneqq K(y)+\vae P(y,x),
$$
with 
$$
K,P\in \mathcal{B}_{r,s}(\mathsf{y})\,.
$$
Assume that\footnote{In the sequel, $K$ and $P$  stand for  generic real analytic hamiltonians which, later on, will respectively play the roles of $K_j$ and $P_j$,  and $\mathsf{y},\,r$, the roles of $y_j,\,r_j$ in the iterative step.}\textsuperscript{,}\footnote{Notice that $\mathsf{T}\mathsf{K}\ge \mathsf{T}\|K_{yy}(\mathsf{y})\|\ge \|T\|\|K_{yy}(\mathsf{y})\|=\|T\|\|T^{-1}\|\ge 1 $.\label{ftnTK1v2}}
\beq{RecHypArnv2}
\begin{aligned}
&\det K_{yy}(\mathsf{y})\neq 0\;, \qquad\qquad\;\, T\coloneqq K_{yy}(\mathsf{y})^{-1}\;,\\
&\|K_{yy}\|_{r,\mathsf{y}}\le \mathsf{K}\;,\qquad\qquad\quad\ \; \|T\|\le \mathsf{T}\;,\\
& \|P\|_{r,s,\mathsf{y}}\le M \;,\qquad\quad\qquad\,\,\  \o\coloneqq K_{yy}(\mathsf{y})\in \D^\t_\a\;. 
\end{aligned}
\eeq
Fix  $\vae\neq0$ and assume that
\beq{lamsup1}
\l\ge \log\left(\s^{2\n+d}\frac{{\a}^2}{|\vae|{M}\mathsf{K}}\right)\ge 1 \;.
\eeq
Let
\beq{DefNArnv2}
\begin{aligned}
&
\k\coloneqq 4\s^{-1}\l, \quad \check{r}\le \frac{5}{24d}\frac{r}{\mathsf{T} \mathsf{K}}\;,\quad
\bar{r}\le
\dst\min\left\{\frac{\a}{2d\mathsf{K}\k^{\t+1}}\,,\, \check{r} \right\}.\\
& \bar{s}\coloneqq s-\frac{2}{3}\s,\quad s'\coloneqq s-\s \,,
\end{aligned}
\eeq
and\footnote{Notice that $\mathsf{L}\ge \s^{-d}\ovl{\mathsf{L}}\ge \ovl{\mathsf{L}}$ and $\frac{40d\mathsf{T}^2\mathsf{K} }{r^2}\s^{-(\n+d)}>\frac{4\mathsf{T} }{r^2}\ge\frac{4}{\mathsf{K} r^2}$.
}
\begin{align*}
\ovl{\mathsf{L}}&\coloneqq \frac{\mathsf{C}_0}{\sqrt{2}} \max\left\{1,\frac{\a}{r\mathsf{K}}\right\}\frac{M \mathsf{K}}{\a^2}\s^{-(2\n+d)}\;,\\
\mathsf{L}&\coloneqq M\dst\max\left\{\frac{40d\mathsf{T}^2\mathsf{K}  }{r^2}\s^{-(\n+d)}\,,\,\frac{4}{\mathsf{K} r^2}\,,\,\frac{\mathsf{C}_4}{\sqrt{2}} \max\left\{1,\frac{\a}{r\mathsf{K}}\right\}\frac{ \mathsf{K}}{\a^2}\s^{-2(\n+d)}\right\}\\
		&=M\dst\max\left\{\frac{40d\mathsf{T}^2\mathsf{K}  }{r^2}\s^{-(\n+d)}\,,\,\frac{\mathsf{C}_4}{\sqrt{2}} \max\left\{1,\frac{\a}{r\mathsf{K}}\right\}\frac{ \mathsf{K}}{\a^2}\s^{-2(\n+d)}\right\}
\;.
\end{align*}
Then, there exists a generating function $g\in \mathcal{B}_{\bar r,\bar s}(\mathsf{y})$
 with the following properties:
\beq{Est1Lem1bv2}
\left\{
\begin{aligned}
&\|g_x\|_{\bar{r},\bar{s},\mathsf{y}}\le  \mathsf{C}_1 \frac{M}{\a} \s^{-(\n+d)}\,,\\
& \|g_{y'}\|_{\bar{r},\bar{s},\mathsf{y}},\, \|\dpr_{y'x}^2 g\|_{\bar{r},\bar{s},\mathsf{y}}\le \ovl{\mathsf{L}}\,,\\
&\|\dpr_{y'}^2\wt K\|_{\check{r},\mathsf{y}}\le 
\mathsf{K}\mathsf{L}\,,
\end{aligned}
\right.
\eeq
where 
$$\wt K(y')\coloneqq \average{P(y',\cdot)}\;.$$
 If, in addition,  
\beq{cond1Bisv2}
{|\vae| }{\mathsf{L}}\le \frac{\sigma}{3}\;,
\eeq
then, there exists $\mathsf{y}'\in\rn$ such that 
\beq{convEstv2}
\left\{
\begin{aligned}
&\dpr_{y'} K'(\mathsf{y}')=\o \,,\qquad \qquad\quad\qquad\quad \det\dpr^2_{y'} K'(\mathsf{y}')\neq 0\,,\\
&|\vae|\|g_x\|_{\bar{r},\bar{s},\mathsf{y}}\le \frac{r}{3}\,,\qquad \qquad\quad\qquad \ 
|\mathsf{y}'-\mathsf{y}| 
\le \frac{8|\vae|\mathsf{T} M}{r}\,,\\
&|\vae|\|\wt T\|\le \mathsf{T}|\vae|\mathsf{L}
\,, \qquad\qquad\quad\quad\quad\ \;\|P_+\|_{\bar{r},\bar s,\mathsf{y}} \le  \mathsf{L}M\,,
\end{aligned}
\right. 
\eeq
where 
$$
K'\coloneqq K+\vae\wt K\;,\qquad \left(\dpr^2_{y'} K'(\mathsf{y}')\right)^{-1}\eqqcolon T+\vae\;\wt T\;,\qquad P_+(y',x)\coloneqq P(y'+\vae g_x(y',x),x)\;.
$$  
and the following hold. For $y'\in D_{\bar{r}}(\mathsf{y})$,    the map $\psi_\vae(x):=  x+\vae g_{y'}(y',x)$ has an analytic inverse 
 $\f(x')=x'+\vae \wt{\f}(y',x';\vae)$  such that
\beq{boundalBisv2}
\|\wt{\f}\|_{\bar{r}, s',\mathsf{y}}\le   \ovl{\mathsf{L}} \qquad {\rm and}
\quad
\f=\id + \vae \wt{\f} : D_{\bar{r}/2,s'}(\mathsf{y}')\to \torus^d_{\bar{s}} \ ;
\eeq
for any $(y',x)\in D_{\bar{r},\bar s}(\mathsf{y})$, $|y'+\vae g_x(y',x)-\mathsf{y}|<\frac{2}{3} r$; the map $\phi'$ is a symplectic diffeomorphism and
\beq{phiokBis0v2}
\phi'=\big( y'+\vae g_x(y', \f(y',x')),\f(y',x')\big): D_{\bar{r}/2,s'}(\mathsf{y}')\to D_{2r/3, \bar{s}}(\mathsf{y}),
\eeq
with
\beq{phiokBis1v2}
\|\mathsf{W}\,\tilde \phi\|_{\bar{r}/2,s',\mathsf{y}'}\le \s^d{\mathsf{L}}\,,
\eeq
where $\tilde \phi$ is defined by the relation $\phi'=:\id + \vae \tilde \phi$,
$$
\mathsf{W}\coloneqq \begin{pmatrix}
\max\{\frac{\mathsf{K}}{{\a}}\;,\frac{1}r\}\;\uno_d & 0\\ \ \\
0			& \uno_d 
\end{pmatrix}
$$
and
\beq{tesitBisv2}
\|P'\|_{\bar{r}/2, s',\mathsf{y}'}\le  \mathsf{L}M\;,
\eeq
with
$$
P'(y',x')\coloneqq P_+(y',\f(x'))=P\circ \phi'(y',x')\;.
$$
\elem
\proof\\
\Giu
{\bf Step 1: Construction of the Arnold's transformation }  We seek for a near--to--the--identity symplectic transformation 
\[\phi'\colon D_{r_1,s_1}(\mathsf{y}')\to D_{r,s}(\mathsf{y}),\]
with $D_{r_1,s_1}(\mathsf{y}')\subset D_{r,s}(\mathsf{y})$,   generated by a function of the form $y'\cdot x+\vae g(y',x)$, so that
\beq{ArnTraKamv2}
\phi'\colon \left\{\begin{aligned}
y  &=y'+\vae g_x(y',x)\\
x' &=x+\vae g_{y'}(y',x)\, ,
\end{aligned}
\right.
\eeq
such that
\beq{ArnH1v2}
\left\{
\begin{aligned}
& H':= H\circ \phi'=K'+\vae^2 P'\ ,\\
& \dpr_{y'} K'(\mathsf{y}')=\o,\quad \det \dpr^2_{y'} K'(\mathsf{y}')\neq 0\,.
\end{aligned}
\right.
\eeq
By Taylor's formula, we get\footnote{Recall that $\average{\cdot}$ stands for the average over $\tn$.}
\beq{Arneq11v2}
\begin{aligned}
H(y'+\vae g_x(y',x),x)=&K(y')+\vae \wt K(y') +\vae \left[K'(y')\cdot g_x +T_{\k} P(y',\cdot)-\wt K(y') \right]+\\
						&+\vae^2 \left( P^\ppu+P^\ppd+ P^\ppt\right)(y',x) \\
			= & K'(y')+\vae \left[K'(y')\cdot g_x +T_{\k} P(y',\cdot)-\wt K(y') \right]+ \vae^2 P_+(y',x),
\end{aligned}
\eeq
with $\k\in\natural$, which will be chosen large enough so that $P^\ppt=O(\vae)$ 
and 
\beq{ArnDefPsv2}
\left\{
\begin{aligned}
P_+&\coloneqq P^\ppu+P^\ppd+ P^\ppt\\
P^\ppu &\coloneqq \su{\vae^2}\left[K(y'+\vae g_x)-K(y')-\vae K'(y')\cdot g_x \right]=\dst\int^1_0(1-t)K_{yy}(\vae t g_x)\cdot g_x\cdot g_x dt\\
P^\ppd &\coloneqq \su\vae \left[P(y'+\vae g_x,x)-P(y',x)\right]=\dst\int_0^1P_y(y'+\vae t g_x,x)\cdot g_x dt\\
P^\ppt &\coloneqq \su\vae \left[ P(y',x)-T_{\k} P(y',\cdot)\right]=\su\vae \dst\sum_{|n|_1>\k} P_n(y')\ex^{in\cdot x}\; .
\end{aligned}
\right.
\eeq
By the non--degeneracy condition in \eqref{RecHypArnv2}, for $\vae$ small enough (to be made precised below), $\det\dpr_{y'}^2 K'(\mathsf{y})\neq0$
 and, therefore, by Lemma~\ref{IFTLem}, there exists a unique $\mathsf{y}'\in D_r(\mathsf{y})$ such that the second part of \eqref{ArnH1v2} holds. 
In view of \eqref{Arneq11v2}, in order to get the first part of \eqref{ArnH1v2}, we need to find $g$ such that  $K_y(y')\cdot g_x +T_{\k} P(y',\cdot)-\wt K(y')$ vanishes; such a $g$ is indeed given by
 \beq{HomEqArnv2}
 g\coloneqq \dst\sum_{0<|n|_1\leq \k} \frac{-P_n(y')}{iK_y(y')\cdot n}\ex^{in\cdot x},
 \eeq
provided that 
\beq{CondHomEqArnv2}
K_y(y')\cdot n\neq 0, \quad \forall\; 0<|n|_1\leq \k,\quad \forall\; y'\in D_{r_1}(\mathsf{y}')\quad  \left(\subset D_{r}(\mathsf{y})\right).
\eeq
But, in fact, since $K_y(\mathsf{y})$ is rationally independent, then, given any $\k\in\natural$, there exists $\bar{r}\leq r$ such that
\beq{CondHomEqArnBisv2}
K_y(y')\cdot n\neq0,\quad \forall\; 0<|n|_1\leq \k, \quad\forall\; y'\in D_{\bar{r}}(\mathsf{y}).
\eeq
The last step is to invert the function $x\mapsto x+\vae g_{y'}(y',x)$ in order to define $P'$. But, by Lemma~\ref{IFTLem}, for $\vae$ small enough, the map $x\mapsto x+\vae g_{y'}(y',x)$ admits an real--analytic inverse of the form
\beq{InvComp2Fiv2}
\f(y',x';\vae)\coloneqq x'+\vae \wt{\f}(y',x';\vae),
\eeq
so that the Arnod's symplectic transformation is given by
\beq{ArnTrans0v2}
\phi'\colon (y',x')\mapsto \left\{
\begin{aligned}
y &= y'+\vae g_x(y',\f(y',x'))\\
x &= \f(y',x';\vae)= x'+\vae \wt{\f}(y',x';\vae) .
\end{aligned}
\right.
\eeq
Hence, \eqref{ArnH1v2} holds with
\beq{DefP1Arv2}
P'(y',x')\coloneqq P_+(y', \f(y',x')).
\eeq
{\bf Step 2: Quantitative estimates}\\
First of all, notice that
\footnote{Recall footnote \textsuperscript{\ref{ftnTK1}}.}
\beq{rrbarAsv2}
\bar{r}\le \check r\le \frac{5r}{24d}<\frac{r}{2}\;.
\eeq
\noi
We begin by extending the ``diophantine condition w.r.t. $K_y$'' uniformly to $D_{\bar{r}}(\mathsf{y})$ up to the order $\k$. Indeed, by the Mean Value Inequality and $K_y(\mathsf{y})=\o\in\D^\t_\a$,  we get, for any $0<|n|_1\leq \k$ and any $y'\in D_{\bar{r}}(\mathsf{y})$,
\begin{align}
|K_y(y')\cdot n|&=|\o\cdot n +(K_y(y')-K_y(\mathsf{y}))\cdot n|\geq |\o\cdot n|\left(1-d\frac{\|K_{yy}\|_{\bar{r},\mathsf{y}}}{|\o\cdot n|}|n|_1\bar{r}\right) \nonumber\\
         &\geq \frac{\a}{|n|_1^\t}\left(1-\frac{d\mathsf{K}}{\a }|n|_1^{\t+1}\bar{r} \right)\geq \frac{\a}{|n|_1^\t}\left(1-\frac{d\mathsf{K}}{\a }\k^{\t+1}\bar{r} \right)\ge \frac{\a}{2|n|_1^\t},\label{ArnExtDiopCondv2}
\end{align}
so that, by Lemma~\ref{fce}--$(i)$, we have
\begin{align*}
\|g_x\|_{\bar{r},\bar{s},\mathsf{y}} &\overset{def}{=}\dst\sup_{D_{\bar{r},\bar{s}}(\mathsf{y})}\left|\dst\sum_{0<|n|_1\leq \k}\frac{nP_n(y')}{K_y(y')\cdot n}\ex^{in\cdot x} \right|\leq \dst\sum_{0<|n|_1\leq \k}\frac{\|P_n\|_{\bar{r},\bar{s}, \mathsf{y}}}{|K_y(y')\cdot n|}|n|_1\ex^{\left(s-\frac{2}{3}\s\right)|n|_1}\\
   &\leq \dst\sum_{0<|n|_1\leq \k} M\ex^{-s|n|_1}\frac{2|n|_1^{\n}}{\a}\ex^{\left(s-\frac{2}{3}\s\right)|n|_1}\leq \frac{2M}{\a}\dst\sum_{n\in\zn} |n|_1^{\n}\ex^{-\frac{2}{3}\s|n|_1}\\
   &\leq \frac{2M}{\a}\dst\int_{\rn} |y|_1^{\n}\ex^{-\frac{2}{3}\s|y|_1}dy\\
   &= \left(\frac{3}{2\s}\right)^{\n+d}\frac{2M}{\a}\dst\int_{\rn} |y|_1^{\n}\ex^{-|y|_1}dy\\
   &= \mathsf{C}_1 \frac{M}{\a} \s^{-(\n+d)}\,,
\end{align*}
\begin{align*}
\|\dpr_{y'}g\|_{\bar{r},\bar{s},\mathsf{y}} &\overset{def}{=}\dst\sup_{D_{\bar{r},\bar{s}}(\mathsf{y})}\left|\dst\sum_{0<|n|_1\leq \k}\left(\frac{ \dpr_yP_n(y')}{K_y(y')\cdot n}-P_n(y')\frac{ K_{yy}(y')n}{(K_y(y')\cdot n)^2}\right)\ex^{in\cdot x} \right|\\
   &\leq \dst\sum_{0<|n|_1\leq \k}\dst\sup_{D_{\bar{r}}(\mathsf{y})}\left(\frac{\|(P_y)_n\|_{\bar{r},s, \mathsf{y}}}{|K_y(y')\cdot n|}+d\|P_n\|_{r,s, \mathsf{y}}\frac{\|K_{yy}\|_{r,\mathsf{y}}|n|_1}{|K_y(y')\cdot n|^2}\right)\ex^{\left(s-\frac{2}{3}\s\right)|n|_1}\\
   &\stackrel{\equ{RecHypArnv2}+\equ{ArnExtDiopCondv2}}{\le} \dst\sum_{0<|n|_1\leq \k}\left( \frac{M}{r-\bar{r}}\ex^{-s|n|_1}\frac{2|n|_1^{\t}}{\a}+dM\ex^{-s|n|_1}\mathsf{K}|n|_1\left(\frac{2|n|_1^{\t}}{\a}\right)^2\right)\ex^{\left(s-\frac{2}{3}\s\right)|n|_1}\\
   &\leby{rrbarAsv2} \frac{4M}{\a^2 r}\dst\sum_{0<|n|_1\leq \k}\left( |n|_1^{\t}\a +dr\mathsf{K}|n|_1^{2\t+1}\right)\ex^{-\frac{2}{3}\s|n|_1}\\
   &\le \max\left\{\a,r\mathsf{K}\right\}\frac{4M}{\a^2 r}\dst\sum_{0<|n|_1\leq \k}\left( |n|_1^{\t}+d|n|_1^{2\t+1}\right)\ex^{-\frac{2}{3}\s|n|_1}\\
   &\leq \max\left\{1,\frac{\a}{r\mathsf{K}}\right\}\frac{4M \mathsf{K}}{\a^2 }\dst\int_{\rn} \left( |y|_1^{\t}+d|y|_1^{2\t+1}\right)\ex^{-\frac{2}{3}\s|y|_1}dy \\
   &= \left(\frac{3}{2\s}\right)^{2\t+d+1}\max\left\{1,\frac{\a}{r\mathsf{K}}\right\}\frac{4M \mathsf{K}}{\a^2 }\dst\int_{\rn} \left( |y|_1^{\t}+d|y|_1^{2\t+1}\right)\ex^{-|y|_1}dy\\
   &\le \frac{\mathsf{C}_0}{\sqrt{2}} \max\left\{1,\frac{\a}{r\mathsf{K}}\right\}\frac{M \mathsf{K}}{\a^2 }\s^{-(2\t+d+1)}\\
   &\le \ovl{\mathsf{L}} \;,
\end{align*}
and, analogously,
\begin{align*}
\|\dpr^2_{y'x}g\|_{\bar{r},\bar{s},\mathsf{y}} &\overset{def}{=}\dst\sup_{D_{\bar{r},\bar{s}}(\mathsf{y})}\left|\dst\sum_{0<|n|_1\leq \k}\left(\frac{ \dpr_yP_n(y')}{K_y(y')\cdot n}-P_n(y')\frac{ K_{yy}(y')n}{(K_y(y')\cdot n)^2}\right)\cdot n\ex^{in\cdot x} \right|\\
   &\leq \dst\sum_{0<|n|_1\leq \k}\dst\sup_{D_{\bar{r}}(\mathsf{y})}\left(\frac{\|(P_y)_n\|_{\bar{r},s, \mathsf{y}}}{|K_y(y')\cdot n|}+d\|P_n\|_{r,s, \mathsf{y}}\frac{\|K_{yy}\|_{r,\mathsf{y}}|n|_1}{|K_y(y')\cdot n|^2}\right)|n|_1\ex^{\left(s-\frac{2}{3}\s\right)|n|_1}\\
   &\le \max\{\a,r\mathsf{K}\}\frac{4M}{\a^2 r}\dst\sum_{0<|n|_1\leq \k}\left( |n|_1^{\t}+d|n|_1^{2\t+1}\right)|n|_1\ex^{-\frac{2}{3}\s|n|_1}\\
   &\leq \max\left\{1,\frac{\a}{r\mathsf{K}}\right\}\frac{4M \mathsf{K}}{\a^2 }\dst\int_{\rn} \left( |y|_1^{\t}+d|y|_1^{2\t+1}\right)|y|_1\ex^{-\frac{2}{3}\s|y|_1}dy \\
   &= \left(\frac{3}{2\s}\right)^{2\t+d+2}\max\left\{1,\frac{\a}{r\mathsf{K}}\right\}\frac{4M \mathsf{K}}{\a^2 }\dst\int_{\rn} \left( |y|_1^{\t+1}+d|y|_1^{2\t+2}\right)\ex^{-|y|_1}dy\\
   &= \frac{\mathsf{C}_0}{\sqrt{2}} \max\left\{1,\frac{\a}{r\mathsf{K}}\right\}\frac{M \mathsf{K}}{\a^2 }\s^{-(2\n+d)} \\
   &=\ovl{\mathsf{L}}\;,
\end{align*}
and, for $|\vae|< {\vae_*}$,
\[\|\wt K_y\|_{r/2,\mathsf{y}}=\| \left[P_y\right]\|_{r/2,\mathsf{y}}\leq \|P_y\|_{r/2,\bar{s}, \mathsf{y}}\leq  \frac{M}{r-\frac{r}{2}}\leq \frac{2M}{r} \;,\]
\[\|\dpr_{y'}^2\wt K\|_{r/2,\mathsf{y}}=\| \left[P_{yy}\right]\|_{r/2,\mathsf{y}}\leq \|P_{yy}\|_{r/2,\bar{s}, \mathsf{y}}\leq  \frac{M}{(r-\frac{r}{2})^2}\leq \frac{4M}{r^2}\le \mathsf{K}\mathsf{L} 
\;.\]
Next, we prove the existence and uniqueness of $\mathsf{y}'$ in \eqref{ArnH1v2}. 
 Consider then
\begin{align*}
F\colon D_{\check{r}}(\mathsf{y})\times D^1_{2|\vae|}(0) &\longrightarrow \qquad \cn\\
		(y,\eta)\quad &\longmapsto K_y(y)+\eta \wt K_{y'}(y)-K_y(\mathsf{y}).
\end{align*}
Then
\begin{itemize}
\item $F(\mathsf{y},0)=0,\quad F_y(\mathsf{y},0)^{-1}=K_{yy}(\mathsf{y})^{-1}=T$;
\item For any $(y,\eta)\in D_{\check{r}}(\mathsf{y})\times D^1_{2|\vae|}(0)$,
\begin{align*}
\|\uno_d-TF_y(y,\eta)\|&\leq \|\uno_d-TK_{yy}\|+|\eta|\;\|T\|\;\|\dpr_{y'}^2\wt K\|_{r/2,\mathsf{y}}\\
	  &\leq d\|T\|\|K_{yyy}\|_{\check{r},\mathsf{y}}\check{r}+ 2|\vae|\mathsf{T}\frac{4M}{r^2}\\
      &\leq d\mathsf{T} \mathsf{K}\frac{\check{r}}{r-\check{r}}+8\mathsf{T}\frac{|\vae| M}{ r^2}\\
      &\leby{rrbarAsv2}d\mathsf{T} \mathsf{K} \frac{2\check{r}}{ r}+|\vae|\frac{8\mathsf{T} M}{r^2}\\
      &\le 2d\mathsf{T} \mathsf{K}\frac{\bar{r}}{ r}+\su2{|\vae|}\mathsf{L}\\
      &\overset{\equ{rrbarAsv2}+\equ{cond1Bisv2}}{\leq}\frac{5}{12}+\frac{\s}{6}\\
      &\le \frac{5}{12}+\su{12}=\su2\;;
\end{align*}
\item Recalling $\s\le\su2$, we have
\begin{align}
2\|T\|\|F(\mathsf{y},\cdot)\|_{2|\vae|,0}&=2\|T\|\dst\sup_{B^1_{2|\vae|}(0)}|\eta \wt K_{y'}(\mathsf{y})|\nonumber\\
		&\leq 2\mathsf{T} \frac{4|\vae| M}{r}\nonumber\\
		&\le \frac{5\cdot 2^{\n+d}}{8d}\frac{r}{\mathsf{T}\mathsf{K}}\s^{\n+d}{|\vae|}\mathsf{L}\nonumber\\
		&= 3\cdot 2^d\;(2\s)^{\n}\;\check r\;\s^{d}{|\vae|}\mathsf{L}\nonumber\\
		&\le 3\cdot 2^d\;\check r\;\s^{d}{|\vae|}\mathsf{L}\label{distyy1I}\\
		&\leby{cond1Bisv2} 3\;\check r\;(2\s)^{d}\;\frac{\s}{3}\nonumber\\
		&\le \frac{\check{r}}{2}\;.\nonumber
\end{align}
\end{itemize}
Therefore, Lemma~\ref{IFTLem} applies. Hence, there exists a function $g\colon D^1_{2|\vae|}(0)\to D_{\bar{r}}(\mathsf{y})$ such that its graph coincides with $F^{-1}(\{0\})$. In particular, $\mathsf{y}'\coloneqq g(\vae)$ is the unique $y\in D_{\bar{r}}(\mathsf{y})$ satisfying $0=F(y,\vae)=\dpr_y K'(y)-\o$ \ie the second part of \eqref{ArnH1v2}. Moreover, 
\beq{EcarY1Y0v2}
|\mathsf{y}'-\mathsf{y}|\leq 2\|T\|\|F(\mathsf{y},\cdot)\|_{2|\vae|,0}\leq \frac{8|\vae|\mathsf{T} M}{r}\leby{distyy1I} 3\cdot 2^d\;\check r\;\s^{d}{|\vae|}\mathsf{L}\leq \frac{\check{r}}{2}\;,
\eeq
so that
\beq{NextSetArnv2}
D_{\frac{\check{r}}{2}}(\mathsf{y}')\subset D_{\check{r}}(\mathsf{y}).
\eeq
Next, we prove that $\dpr^2_y K'(\mathsf{y}')$ is invertible. Indeed, by Taylor' formula, we have
\begin{align*}
\dpr^2_y K'(\mathsf{y}')&= K_{yy}(\mathsf{y})+ \dst\int_0^1 K_{yyy}(\mathsf{y}+t\vae \wt y)\cdot\vae\wt y dt+\vae \wt K_{yy}(\mathsf{y}')\\
           &= T^{-1}\left(\uno_d+\vae T\left(\dst\int_0^1 K_{yyy}(\mathsf{y}+t\vae \wt y)\cdot\wt y dt+ \wt K_{yy}(\mathsf{y}')\right)\right)\\
           &\eqqcolon T^{-1}(\uno_d+\vae A),
\end{align*}
and, by Cauchy's estimate, 
\begin{align*}
|\vae|\|A\|&\leq \|T\|\left(d\|K_{yyy}\|_{r/2,\mathsf{y}}|\vae||\mathsf{y}'-\mathsf{y}|+ |\vae|\|\dpr_{y'}^2\wt K\|_{r/2,\mathsf{y}}\right)\\
     &\leq \|T\|\left(\frac{d\|K_{yy}\|_{r,\mathsf{y}}}{r-\frac{r}{2}}|\vae||\mathsf{y}'-\mathsf{y}|+|\vae|\|\wt K_{yy}\|_{r/2,\mathsf{y}}\right)\\
	 &\leby{EcarY1Y0v2} \mathsf{T}\left(\frac{2d\mathsf{K}}{r}\frac{8|\vae|\mathsf{T} M}{r}+\frac{4|\vae|M}{r^2} \right)\\
	 &\leq \frac{4|\vae|\mathsf{T}M}{r^2}(4d\mathsf{T}\mathsf{K}+1)\\
	 &\leq\frac{20d|\vae|\mathsf{T}^2\mathsf{K} M}{r^2}\\
	 &\le \su 2|\vae|\mathsf{L}\\
	 &\leby{cond1Bisv2}\frac{\s}{6}\\
	 &\le\su2.
\end{align*}
Hence $\dpr_{y'}^2 K'(\mathsf{y}')$ is invertible with
\[\dpr_{y'}^2 K'(\mathsf{y}')^{-1}=(\uno_d+\vae A)^{-1}T=T+\dst\sum_{k\geq 1}(-\vae)^k A^k T\eqqcolon T+\vae \wt T,\]
and
\[|\vae|\|\wt T\|\leq |\vae|\frac{\|A\|}{1-|\vae|\|A\|}\|T\|\leq 2|\vae|\|A\| \|T\|
\le |\vae|\mathsf{L}\mathsf{T}
\le 2\frac{\s}{6}\mathsf{T}
= \mathsf{T}\frac{\s}{3}\,.\]
Next, we prove estimate on $P_+$. We have,
\[|\vae|\|g_x\|_{\bar{r},\bar{s},\mathsf{y}}\leq |\vae|\mathsf{C}_1 \frac{M}{\a} \s^{-(\t+d+1)} \le |\vae| \frac{r}{3}\mathsf{L}\leby{cond1Bisv2}\frac{r}{3}\frac{\s}{3}\le \frac{r}{3}\]
so that, for any $(y',x)\in D_{\bar{r},\bar{s}}(\mathsf{y})$,
\[ |y'+\vae g_x(y',x)-\mathsf{y}|\leq \bar{r}+\frac{r}{3}< \frac{r}{8d}+\frac{r}{3}<\frac{2r}{3}<r\,,\]
and thus
\begin{align*}
\|P^\ppu\|_{\bar{r},\bar{s},\mathsf{y}}&\leq d^2 \|K_{yy}\|_{r,\mathsf{y}}\|g_x\|_{\bar{r},\bar{s},\mathsf{y}}^2\leq d^2 \mathsf{K}\left( \mathsf{C}_1 \frac{M}{\a} \s^{-(\n+d)}\right)^2\\
   &=d^2\mathsf{C}_1^2 \frac{\mathsf{K}M^2}{\a^2} \s^{-2(\n+d)}, 
\end{align*}
\begin{align*}
\|P^\ppd\|_{\bar{r},\bar{s},\mathsf{y}}&\leq d\|P_y\|_{\frac{5r}{6},\bar{s},\mathsf{y}}\|g_x\|_{\bar{r},\bar{s},\mathsf{y}}\leq d\frac{6M}{r}\mathsf{C}_1 \frac{M}{\a} \s^{-(\n+d)}\\
     &= 6d\mathsf{C}_1 \frac{M^2}{\a r}\s^{-(\n+d)}
\end{align*}
and by Lemma~\ref{fce}--$(i)$, we have,
\begin{align*}
|\vae|\|P^\ppt\|_{\bar{r},s-\frac{\s}{2},\mathsf{y}}&\leq \dst\sum_{|n|_1>\k}\|P_n\|_{\bar{r},\mathsf{y}}\ex^{(s-\frac{\s}{2})|n|_1}\leq M\dst\sum_{|n|_1>\k}\ex^{-\frac{\s |n|_1}{2}}\\
  &\leq M\ex^{-\frac{ \k\s}{4}}\dst\sum_{|n|_1>\k}\ex^{-\frac{\s |n|_1}{4}}\leq M\ex^{-\frac{ \k\s}{4}}\dst\sum_{|n|_1>0}\ex^{-\frac{\s |n|_1}{4}}\\
  &= M\ex^{-\frac{ \k\s}{4}} \left(\left(\dst\sum_{k\in \integer}\ex^{-\frac{\s |k|}{4}}\right)^d-1\right)=M\ex^{-\frac{ \k\s}{4}}\left(\left(1+\frac{2\ex^{-\frac{\s }{4}}}{1-\ex^{-\frac{\s }{4}}} \right)^d-1\right)\\
  &= M\ex^{-\frac{ \k\s}{4}}\left(\left(1+\frac{2}{\ex^{\frac{\s }{4}}-1} \right)^d-1\right)\leq M\ex^{-\frac{ \k\s}{4}}\left(\left(1+\frac{2}{\frac{\s }{4}} \right)^d-1\right)\\
  &\leq \s^{-d} M\ex^{-\frac{ \k\s}{4}}\left(\left(\s +8 \right)^d-\s^d\right)\leq d 8^{d}\s^{-d} M\ex^{-\frac{ \k\s}{4}}\\
  &= \mathsf{C}_2\s^{-d} M\ex^{-\l}\\
  &\leby{lamsup1} \mathsf{C}_2\s^{-d} M\s^{-(2\n+d)}\frac{|\vae|{M}\mathsf{K}}{{\a}^2}\\
  &= \mathsf{C}_2 M\frac{|\vae|{M}\mathsf{K}}{{\a}^2}\s^{-2(\n+d)}\,.
\end{align*}
Hence\footnote{Recall that $\s<1$.},
\begin{align*}
\|P_+\|_{\bar{r},\bar{s},\mathsf{y}}&\leq \|P^\ppu\|_{\bar{r},\bar{s},\mathsf{y}}+\|P^\ppd\|_{\bar{r},\bar{s},\mathsf{y}}+\|P^\ppt\|_{\bar{r},\bar{s},\mathsf{y}}\\
  &\leq d^2\mathsf{C}_1^2 \frac{\mathsf{K}M^2}{\a^2} \s^{-2(\n+d)}+6d\mathsf{C}_1 \frac{M^2}{\a r}\s^{-(\n+d)}+\mathsf{C}_2 M\frac{|\vae|{M}\mathsf{K}}{{\a}^2}\s^{-2(\n+d)}\\
  &= \left(d^2\mathsf{C}_1^2 r\mathsf{K}+6d\mathsf{C}_1 \a \s^{\n+d}+\mathsf{C}_2 r\mathsf{K}\right)\frac{M^2}{\a^2 r}\s^{-2(\t+d+1)}\\
  &\le \left(d^2\mathsf{C}_1^2+6d\mathsf{C}_1 +\mathsf{C}_2\right)\max\left\{\a,r\mathsf{K}\right\}\frac{M^2}{\a^2 r}\s^{-2(\t+d+1)}\\
  &\leby{RecHypArnv2} \frac{\mathsf{C}_3}{\sqrt{2}} \max\left\{1,\frac{\a}{r\mathsf{K}} \right\}\frac{M^2\mathsf{K}}{\a^2 }\s^{-2(\n+d)}\\
  &\le \mathsf{L}M\;.
\end{align*}
The proof of the claims on $\phi'$ and $P'$ are proven in a similar way as in Lemma~\ref{lem:1}.
\qed
\noi
Finally, we prove the convergence of the scheme by mimicking Lemma~\ref{lem:2}.\\
\noi
Let   $H_0\coloneqq H$, $K_0\coloneqq K$, $P_0\coloneqq P$ and  $r_0$, $s_0$, $s_*$, $\s_0$,  $\mathsf{W}_0$, $M_0$, $\mathsf{K}_0$, $\mathsf{T}_0$, $\m_0$ be  as in $\S\ref{AssumpArnolv2}$
 and for a given $\vae\neq0$ and $j\ge 0$,  define\footnote{Notice that $s_{j}\downarrow s_*$ and $r_{j}\downarrow 0$.}
\begin{align*}
 \dst\s_j&\coloneqq \frac{\s_0}{2^j}\,,\\
  s_{j+1}&\coloneqq s_j-\s_j=s_*+\frac{\s_0}{2^j}\,,\\
  \bar s_{j}&\coloneqq s_j-\frac{2\s_i}{3}\,,\\
 \mathsf{K}_{j+1}&\coloneqq \mathsf K_0\dst\prod_{k=0}^{j}(1+\frac{\s_k}{3})\le \mathsf K_0\ex^{\frac{2\s_0}{3}}\le \mathsf{K}_0\sqrt{2}\,,  \\
\mathsf{T}_{j+1}&\coloneqq \mathsf T_0\dst\prod_{k=0}^{j}(1+\frac{\s_k}{3})\le \mathsf T_0\ex^{\frac{2\s_0}{3}}\le \mathsf{T}_0\sqrt{2}\,,\\
 \l_0&\coloneqq \log\m_0^{-1}\,,\\
 \mathsf{e}_*&\coloneqq \mathsf{C}_7\; \s_0^{-(4\n+2d+1)}\l_0^{2\n}\;\eta_0^2\;,\\
\mathsf{d}_*&\coloneqq 2^{2\n+2d+1}\;\mathsf{C}_6^2\;\eta_0^2\,,\\
 \k_0 &\coloneqq 4\s_0^{-1}\l_0\,,\\
 \k_j &\coloneqq 4^j\k_0\,,\\
 \wh\a  &\coloneqq \frac{\a}{\sqrt{|\vae|}}\,,\\
 \hat{r}_0 &\coloneqq \frac{r_0}{\sqrt{|\vae|}}\,,\\
 \hat{r}_{j+1}&\coloneqq \su2\min\left\{\frac{\wh{\a}}{2d\sqrt{2}\mathsf{K}_0\k_j^{\n}}\,,\, \frac{5}{48d}\frac{\wh r_j}{\eta_0} \right\}\,,\\
 r_{j+1}&\coloneqq \hat{r}_{j+1}\sqrt{|\vae|} \,,\\
 \check{r}_j&\coloneqq \frac{5}{48d}\frac{r_j}{\eta_0}\,,
\end{align*}
\begin{align*}
 \wh M_0&\coloneqq M_0\,,\\
 \wh M_{j+1}&\coloneqq 	
 			\mathsf{e}_*\mathsf d_*^{j-1} \frac{\mathsf{K}_0 {\wh M_j}^2}{\wh\a^2}\,,\\
 \m_j&\coloneqq \frac{\mathsf{K}_0\wh M_j}{\wh \a^2}\,,\\
 \th_j     &\coloneqq \mathsf{e}_*\;\mathsf d_*^j\;\m_j\,,\\
 \mathsf{W}_j&\coloneqq \diag\left(\max\left\{\frac{\mathsf{K_j}}{\wh{\a}}\;,\frac{\sqrt{|\vae|}}{r_j}\right\}\;\uno_d\,,\sqrt{|\vae|}\;\uno_d\right)\,,\\ 
\mathsf{L}_j&\coloneqq M_i\dst\max\left\{\frac{80d\sqrt{2}\;\mathsf{T}_0\;\eta_0  }{r_j^2}\s_j^{-(\n+d)}\,,\,\mathsf{C}_4 \max\left\{1,\frac{\a}{r_j\mathsf{K}_j}\right\}\frac{ \mathsf{K}_0}{\a^2}\s_j^{-2(\n+d)}\right\}\\
		&=M_j\dst\max\left\{\frac{80d\sqrt{2}\;\mathsf{T}_0\;\eta_0  }{r_j^2}\s_j^{-(\n+d)}\,,\,\frac{4}{\mathsf{K}_j r_j^2}\,,\,\mathsf{C}_4 \max\left\{1,\frac{\a}{r_j\mathsf{K}_j}\right\}\frac{ \mathsf{K}_0}{\a^2}\s_j^{-2(\n+d)}\right\} \,.
\end{align*}
Thus, for any $j\ge0$,
$$
\th_{j+1}= \mathsf{e}_*\;\mathsf d_*^{j+1}\;\m_{j+1}=\mathsf{e}_*\;\mathsf d_*^{j+1}\frac{\mathsf{K}_0\wh M_{j+1}}{\wh \a^2}=\mathsf{e}_*\;\mathsf d_*^{j+1}\frac{\mathsf{K}_0}{\wh \a^2}\;\mathsf{e}_*\mathsf d_*^{j-1} \frac{\mathsf{K}_0 {\wh M_j}^2}{\wh\a^2}= \left(\mathsf{e}_*\;\mathsf d_*^{j}\;\m_j\right)^2=\th_j^2
$$
\ie
$$
\th_j=\th_0^{2^j} \;.
$$
The very first step being quite different from all the others, it has to be done separately. Hence,\\
\lem{frstStep}
Under the above assumptions and notations, if
\beq{condBisv2}
|\vae|\le \left(\frac{r_0}{\wh{\a}\mathsf{T}_0}\right)^2\qquad\mbox{and}\qquad \max\left\{\ex\;\m_0\;,\, \th_0\right\}\le 1\;,
\eeq
then, there exist $y_1\in \mathscr D$ 
 and a real--analytic symplectomorphism
\beq{phijBis0v2}
\phi_0:D_{r_1,s_{1}}(y_1)\to D_{r_{0},s_{0}}(y_{0})\ ,
\eeq
such that, for $H_1\coloneqq H_{0}\circ\phi_0$ , we have
\beq{HjBis0v2}
\left\{
\begin{aligned}
& H_1\eqqcolon K_1 + \vae^{2} P_1\;,\\
& \dpr_{y_1} K_1(y_1)=\o \;,\quad \dpr_{y_1}^2 K_1(y_1)\neq0  \;
\end{aligned}
\right.
\eeq\
and
\begin{align}
&|y_1-y_0| \le 
\frac{8|\vae|\mathsf{T}_0 M_0}{r_0}\,,\\
& \|K_1\|_{\check r_1,y_1}\le \mathsf{K}_1\;,\qquad \|T_1\|\le  \mathsf{T}_1\;,\qquad T_1\coloneqq \dpr_{y_1}^2 K_1(y_1)^{-1}\;,\label{estfin2Bis0000v2}\\
& \vae^{2}M_1:=\vae^{2}\|P_1\|_{r_1,s_1,y_1}\le |\vae|\wh M_1\;,\label{estfin2Bis00011v2}\\
&
\|\mathsf{W}_0(\phi_0-\id)\|_{r_1,s_1,y_1}
\le \s_0^d\;|\vae|{\mathsf L}_0\cdot\sqrt{|\vae|}\;.\label{estfin2Bis010v2}
\end{align}
\elem
\proof
By
\beq{kp08}
\k_0\geby{condBisv2}4\s_0^{-1}\ge 8
\eeq
and
$$
\frac{\wh{\a}}{2d\sqrt{2}\mathsf{K}_0k_0^{\n}}\leby{kp08} \frac{1}{2d\cdot 8^{\n}\sqrt{2}\mathsf{K}_0}\frac{r_0}{\mathsf{T}_0\sqrt{|\vae|}}<\frac{\hat {r}_0}{4\mathsf{C}_5}\,,
$$
we get 
\beq{r1cal}
\hat r_1=\su2\min\left\{\frac{\wh{\a}}{2d\sqrt{2}\mathsf{K}_0\k_0^{\n}}\,,\, \frac{\hat {r}_0}{4\mathsf{C}_5} \right\}=\frac{\wh{\a}}{4d\sqrt{2}\mathsf{K}_0\k_0^{\n}}
\eeq
and, thus
\begin{align}
|\vae| \mathsf L_0 (3 \sigma_0^{-1})&\le 3|\vae|M_0\dst\max\left\{\frac{80d\sqrt{2}\;\mathsf{T}_0\;\eta_0 }{r_0^2}\s_0^{-(\n+d)}\,,\,\mathsf{C}_4 \max\left\{1,\frac{\a}{r_0\mathsf{K}_0}\right\}\frac{ \mathsf{K}_0}{\a^2}\s_0^{-2(\n+d)}\right\}\s_0^{-1}\nonumber\\
                                    &\le 3\dst\max\left\{80d\sqrt{2}\;\eta_0\frac{ \a \;\mathsf{T}_0}{{r}_{0}}\frac{ \a }{r_0 \mathsf{K}_0}\,,\,\mathsf{C}_4 \max\left\{1,\frac{\a}{r_0\mathsf{K}_0}\right\}\right\}\s_0^{-2(\n+d)-1}\frac{ \mathsf{K}_0 M_0}{\wh\a^2}\nonumber\\
                                    &\overset{\equ{condBisv2}}{\le}3\dst\max\left\{80d\sqrt{2}\,,\,\mathsf{C}_4\right\}\s_0^{-2(\n+d)-1}\m_0\;\eta_0\nonumber\\
                                    &\le\mathsf{e}_*\;\m_0\nonumber\\ 
                                    &=\th_0\leby{condBisv2} 1.\label{L0ifsg3}
\end{align}
Therefore, Lemma~\ref{frstStep} is a straightforward consequence of Lemma~\ref{lem:1bisv2}.
\qed
Once the first step is completed, all the following steps do not need any other condition. Actually, they are ``completely'' independent upon $\vae$, and, therefore, the first condition in \equ{condBisv2} is useless. To be precise, the following holds.
\lem{lem:2Bisv2} 
Assume $\equ{HjBis0v2}\div\equ{estfin2Bis00011v2}$ with some $\vae\neq 0$ and
\beq{condBisv2Prt}
 \max\left\{\ex\;\m_0\;,\, \mathsf{C}_8\;\eta_0^{\su8}\;\th_0\right\}< 1\;.
\eeq
Then, one can construct a sequence of symplectic transformations 
\beq{phijBisv2}
\phi_{j-1}:D_{r_j,s_{j}}(y_j)\to D_{r_{j-1},s_{j-1}}(y_{j-1})\;,\qquad j\ge 2
\eeq
so that
\beq{HjBisv2}
H_j:=H_{j-1}\circ\phi_{j-1}=: K_j + \vae^{2^j} P_j
\eeq
converges uniformly. More precisely, 
$\vae^{2^{j-1}} P_{j-1}$, 
$\phi^{j-1}\coloneqq \phi_1\circ\phi_2\circ \cdots\circ \phi_{j-1}$, 
 $K_{j-1}$, $y_{j-1}$ converge uniformly on 
$\{y_*\}\times\dst\torus^d_{s_*}$ to, respectively, $0$, $\phi^*$, $K_*$, $y_*$ which are real--analytic on $\dst\torus^d_{s_*}$ and  $H_1\circ\phi^*=K_*$ with $\det\dpr^2_yK_*(y_*)\neq0$. 
Finally, the following estimates hold for any $i\ge 1$:
\begin{align}
& |\vae|^{2^i}M_i\coloneqq|\vae|^{2^i}\|P_i\|_{r_i,s_i,y_i}\le |\vae|\wh M_i\ ,\label{estfin2Bis00v2}\\
&|y_{i+1}-y_i| \le \frac{8\sqrt{2}\mathsf{T}_0|\vae|^{2^i} M_i}{r_i}\ ,\label{estfin2Bis0011v2}\\
&
|\mathsf{W}(\phi^*-\id)|
\le \frac{\th_0^2}{3\cdot {2^{2d+1}}\;\mathsf{d}_*}\sqrt{|\vae|}\quad \mbox{on}\quad \{y_*\}\times\torus^d_{s_*}\label{estfin2Bis01v2}\ .
\end{align}
\elem 
\proof
 First of all, notice that, for any $i\ge 1$,
\begin{align*}
\hat r_{i+1}&= \min\left\{\frac{\wh{\a}}{4d\sqrt{2}\mathsf{K}_0\k_i^{\n}}\,,\, \frac{5}{96d}\frac{\wh r_i}{\eta_0} \right\}\\
	   &=\min\left\{\frac{\hat{r}_1}{4^{i\n}}\,,\, \frac{5}{96d\eta_0}{\hat{r}_i} \right\}\\
	   &= \min\left\{\frac{\hat{r}_1}{4^{\n i}}\,,\, \frac{5}{96d\eta_0}\frac{\hat{r}_1}{ 4^{\n(i-1)}}\,,\, \left(\frac{5}{96d\eta_0}\right)^2 {\hat{r}_{i-1}} \right\}\\
	   &\,\ \vdots \\
	   &=\min\left\{\frac{\hat{r}_1}{4^{\n i}}\,,\, \frac{5}{96d\eta_0}\frac{\hat{r}_1}{ 4^{\n(i-1)}}\,,\cdots,\, \left(\frac{5}{96d\eta_0}\right)^i {\hat{r}_{1}} \right\}\\
	   &=\frac{\hat{r}_{1}}{4^{\n i}}\dst\min\left\{\left(\frac{5\cdot 4^\n}{96d\eta_0}\right)^0\,,\,\cdots\,,\, \left(\frac{5\cdot 4^\n}{96d\eta_0}\right)^i\right\} \\
	   &= \frac{\hat{r}_{1}}{4^{\n i}}\dst{\min}^i\left\{\frac{5\cdot 4^\n}{96d\eta_0}\,,\,1\right\} \\
	   &= \hat r_1 \dst{\min}^i\left\{\frac{1}{2^{2\n}}\,,\,\frac{5}{96d\eta_0}\right\}\\
	   &=\frac{\hat r_1}{\mathbf{a}_1^{i}}\,,
\end{align*}
where
\beq{eqa1maj}
\mathbf{a}_1\coloneqq \dst{\max}\left\{2^{2\n}\,,\,\frac{96d\eta_0}{5}\right\}\le \dst{\max}\left\{2^{2\n}\,,\,\frac{96d}{5}\right\}\cdot\eta_0=\mathsf{C}_6\;\eta_0\;.
\eeq
\noi
For a given $j\ge 2$, let $(\mathscr{P}^j)$ be the following assertion:  there exist
$j-1$ symplectic transformations\footnote{Compare \equ{phiokBis0v2}.} 
\beq{bes06v2}
\phi_{i}:D_{r_{i+1},s_{i+1}}(y_{i+1})\to D_{2r_i/3, s_i}(y_i),\quad \mbox{for}\quad 1\le i\le j-1,
\eeq
 and $j-1$  Hamiltonians $H_{i+1}=H_i\circ\phi_{i}=K_{i+1}+\vae^{2^{i+1}} P_{i+1}$ real--analytic on $D_{r_{i+1},s_{i+1}}(y_{i+1})$ such that, for any $1\le i\le j-1$,
\beq{bbbBisv2}
\left\{
\begin{array}{l}
\|\dpr_y^2 K_i\|_{r_i,y_i}\le \mathsf{K}_i\, ,\ \\  \ \\
\|T_i\|\le \mathsf{T}_i\,,\ \\  \ \\ 
\dpr_{y} K_i(y_i)=\o \;,\quad \dpr_{y}^2 K_i(y_i)\neq0\,,\ \\  \ \\ 
|\vae|^{2^i}\|P_{i}\|_{r_{i},s_{i},y_{i}}\le |\vae| \wh M_{i}\,,\ \\  \ \\
 \k_i\ge 4\s_i^{-1} \log\left(\s_i^{2\n+d}\m_i^{-1}\right)\,,\ \\  \ \\
|\vae|^{2^i} \mathsf{L}_i \le \frac{\sigma_i}{3}
\  
\end{array}\right.
\eeq
\noi
and
\beq{C.1Bisv2}
\left\{
\begin{aligned}
&\dpr_{y} K_{i+1}(y_{i+1})=\o \;,\quad \dpr_{y}^2 K_{i+1}(y_{i+1})\neq0\;,\\ \ \\
&|y_{i+1}-y_i|\le \frac{8\sqrt{2}\mathsf{T}_0|\vae|^{2^i} M_i}{r_i}\;,\\ \ \\
&\|T_{i+1}\|\le \|T_i\|+\mathsf T_i|\vae|^{2^i}\mathsf L_i\;, 
 \\ \ \\
&\|K_{i+1}\|_{r_{i+1},y_{i+1}}\le \|K_i\|_{r_i,y_{i}}+|\vae|^{2^i}M_i \;,\\ \ \\ 
&\|\dpr_y^2K_{i+1}\|_{r_{i+1},y_{i+1}}\le \|\dpr_y^2K_i\|_{r_i,y_{i}}+\mathsf K_i|\vae|^{2^i}\mathsf L_i \;,\\ \ \\
&\|\mathsf{W}_i(\phi_{i}-\id)\|_{r_{i+1},s_{i+1},y_{i+1}}\le \s_i^d\;|\vae|^{2^i}{\mathsf L}_i\cdot\sqrt{|\vae|} \;, \\ \ \\
& M_{i+1}\coloneqq\|P_{i+1}\|_{r_{i+1},s_{i+1},y_{i+1}}\le  M_i \mathsf L_i\;.
\end{aligned}
\right.
\eeq
Assume $(\mathscr P^j)$, for some $j\ge 2$ and let's check $(\mathscr P^{j+1})$. Fix then $1\le i\le j-1$. Thus
$$
\|\dpr_y^2K_{i+1}\|_{r_{i+1},y_{i+1}}\leby{C.1Bisv2} \|\dpr_y^2K_i\|_{r_i,y_{i}}+\mathsf K_i|\vae|^{2^i}\mathsf L_i\leby{bbbBisv2} \mathsf K_i+\mathsf K_i\frac{\s_i}{3}=\mathsf K_{i+1}<\mathsf K_0\sqrt{2}
$$
and, similarly,
$$
\|T_{i+1}\|\le \mathsf{T}_{i+1},
$$
which prove the two first relations in \equ{bbbBisv2} for $i=j$. Also
\beq{alfhtrikpi}
\frac{\a}{r_{i}\mathsf{K}_{i}}> \frac{\a}{r_1\mathsf{K}_0\sqrt{2}}=\frac{\wh\a}{\hat r_1\mathsf{K}_0\sqrt{2}}=4d\k_0^{\n}\gtby{kp08}1\;,
\eeq
so that
\begin{align*}
|\vae|^{2^i} \mathsf L_i (3 \sigma_i^{-1})&= 3 |\vae|^{2^i}M_i\dst\max\left\{\frac{80d\sqrt{2}\mathsf{T}_0\eta_0  }{r_i^2}\s_i^{-(\n+d)}\,,\,\mathsf{C}_4 \max\left\{1,\frac{\a}{r_i\mathsf{K}_i}\right\}\frac{ \mathsf{K}_0}{\a^2}\s_i^{-2(\n+d)}\right\}\sigma_i^{-1}\\
                                          &\leby{alfhtrikpi} 3|\vae|^{2^i}M_i\dst\max\left\{\frac{80d\sqrt{2}\mathsf{T}_0\eta_0  }{r_i^2}\,,\,\mathsf{C}_4 \frac{1}{\a r_i}\right\}\s_i^{-2(\n+d)-1}\\
                                          &= 3\dst\max\left\{80d\sqrt{2}\mathsf{T}_0\eta_0 \frac{ \wh\a }{ \hat{r}_{i}}\,,\,\mathsf{C}_4 \right\}\s_i^{-2(\n+d)-1}\frac{|\vae|^{2^i}M_i}{|\vae|\wh \a \hat r_i}\\
                                          &= 3\dst\max\left\{640d^2\eta_0^2\mathbf{a}^{i-1}\k_0^{\n}\,,\,\mathsf{C}_4 \right\}\s_i^{-2(\n+d)-1}\frac{|\vae|^{2^i}M_i}{|\vae|\wh \a^2}4d\sqrt{2}\mathsf{K}_0\k_0^{\n}\mathbf{a}^{i-1}\\
                                          &\leby{kp08} 12d\sqrt{2}\dst\max\left\{640d^2\,,\,\mathsf{C}_4 \right\}\mathsf{K}_0\s_i^{-2(\n+d)-1}\frac{|\vae|^{2^i}M_i}{|\vae|\wh \a^2}\eta_0^2\mathbf{a}^{2(i-1)}\k_0^{2\n}\\
                                          &\leby{eqa1maj} 12d\sqrt{2}\dst\max\left\{640d^2\,,\,\mathsf{C}_4 \right\}\mathsf{K}_0\s_i^{-2(\n+d)-1}\frac{|\vae|^{2^i}M_i}{|\vae|\wh \a^2}\eta_0^{2i}\mathsf{C}_6^{2(i-1)}\k_0^{\n}\\
                                          &= 3d\cdot 2^{6\n+2d+3}\sqrt{2}\dst\max\left\{640d^2\,,\,\mathsf{C}_4 \right\}\mathsf{K}_0\s_0^{-(4\n+2d+1)}\left(2^{2\n+2d+1}\mathsf{C}_6^2\eta_0^2\right)^{i-1}\frac{|\vae|^{2^i}M_i}{|\vae|\wh \a^2}\left(\log\m_0^{-1}\right)^{2\n}\eta_0^2\\
                                          &\leby{bbbBisv2} \mathsf{C}_7 \s_0^{-(4\n+2d+1)}\left(\log\m_0^{-1}\right)^{2\n}\eta_0^2\:\mathsf{d}_*^{i-1}\frac{\mathsf{K}_0\wh M_i}{\wh \a^2}\\
                                          &=\mathsf{e}_*\; \mathsf{d}_*^{i-1}\m_i\\
                                          & = \frac{\theta_i}{\mathsf{d}_*}\\
                                          &= \frac{\theta_0^{2^{i}}}{\mathsf{d}_*}\\
                                          &\leby{condBisv2Prt}\su{\mathsf{d}_*}<1\;.
\end{align*}
 Moreover,
 $$
 |\vae|^{2^{i}}\mathsf L_i<\mathsf{e}_*\; \mathsf{d}_*^{i-1}\m_i\;,
 $$ 
 thus by last relation in \equ{C.1Bisv2}, for any $1\le i\le j-1$,  
 $$
 |\vae|^{2^{i+1}}M_{i+1}\le |\vae|^{2^{i}}\mathsf{L}_i\;|\vae|^{2^{i}} M_i<\mathsf{e}_*\mathsf d_*^{i-1} \; \m_i \;|\vae|^{2^{i}} M_i\leby{bbbBisv2} \mathsf{e}_*\mathsf d_*^{i-1} \; \m_i \;|\vae|\; \wh M_i = |\vae|\; \wh M_{i+1}\;,
 $$ 
 which proves the fourth relation in \equ{bbbBisv2} for $i=j$. Hence, by exactly the same computation as above, one gets
 $$
 |\vae|^{2^{i+1}} \mathsf L_{i+1} (3 \sigma_{i+1}^{-1})\le \frac{\theta_{i+1}}{\mathsf{d}_*}=\frac{\theta_0^{2^{i+1}}}{\mathsf{d}_*}<1\ ,
 $$
 which proves the last relation in \equ{bbbBisv2} for $i=j$.  It remains only to check that
the fifth relation in \equ{bbbBisv2} holds as well for $i=j$ in order to apply Lemma~\ref{lem:1bisv2} to $H_i$, $1\le i\le j$ and get \equ{C.1Bisv2} and, consequently, $(\mathscr P^{j+1})$. But in fact, 
 we have\footnote{Notice that, since $\s_0<1$ then $\mathsf{e}_*\geby{condBisv2Prt}\mathsf{C}_7>1$.}
 \begin{align*}
 4\s_j^{-1} \log\left(\s_j^{2\n+d}\m_j^{-1}\right)&\le 4\s_j^{-1} \log\left(\m_j^{-1}\right)\\
         &= 4\s_j^{-1} \log\left(\su{\mathsf{e}_*\; \mathsf{d}_*^{j}}\th_0^{-2^j}\right)\\
         &\le 4\s_j^{-1} \log\left(\left(\frac{\th_0}{\mathsf{e}_*}\right)^{-2^j}\right)\\
         &= 4\s_j^{-1} \log\left(\m_0^{-2^j}\right)\\
         &= 4^j\cdot 4\s_0^{-1} \log\left(\m_0^{-1}\right)\\
         &= \k_j \;.
 \end{align*}
\noi
To finish the proof of the induction \ie one can construct an {\sl infinite sequence} of Arnold's transformations satisfying \equ{bbbBisv2} and  \equ{C.1Bisv2} {\sl for all $i\ge 1$}, one needs only to check $(\mathscr P^{2})$. Thanks to\footnote{Observe that for $j=2$, $i=1$.} $\equ{HjBis0v2}\div\equ{estfin2Bis00011v2}$, we just need to check the two last inequalities in $\equ{bbbBisv2}_{i=1}$. But, in fact, one proves those two relations by exactly the same computation as above. Then, we apply Lemma~\ref{lem:1bisv2} to $H_1$ to get $\equ{bes06v2}_{i=1}$ and  $\equ{C.1Bisv2}_{i=1}$, which achieves the proof of $(\mathscr P^{2})$.\\
\nl
Next, we prove that $\phi^j$ is convergent by proving that it is Cauchy. For any $j\ge 3$, we have, using again Cauchy's estimate,\footnote{Recall that $2^{i-1}\ge i,\,\forall\; i\ge 0$.}
\beqano
\|\mathsf{W}_{j-1}(\phi^{j-1}-\phi^{j-2})\|_{r_j,s_j,y_j}&=&\|\mathsf{W}_{j-1}\phi^{j-2}\circ\phi_{j-1}-\mathsf{W}_{j-1}\phi^{j-2}\|_{r_j,s_j,y_j}\\
           &\leby{bes06v2}& \|\mathsf{W}_{j-1}D\phi^{j-2}\mathsf{W}_{j-1}^{-1}\|_{2r_{j-1}/3, s_{j-1},y_{j-1}}\, \|\mathsf{W}_{j-1}(\phi_{j-1}-\id)\|_{r_j,s_j,y_j}\\
           &\leby{C.1Bisv2}&  \max\left(r_{j-1}\frac{3}{r_{j-1}},\frac{3}{2\s_{j-1}}\right)    \|\mathsf{W}_{j-1}\phi^{j-2}\|_{r_{j-1}, s_{j-1},y_{j-1}} \times\\
           &&\qquad \times \|\mathsf{W}_{j-1}(\phi_{j-1}-\id)\|_{r_j,s_j,y_j}\\
           &=&  \frac{3}{2\s_{j-1}}   \|\mathsf{W}_{j-1}\phi^{j-2}\|_{r_{j-1}, s_{j-1},y_{j-1}} \, \|\mathsf{W}_{j-1}(\phi_{j-1}-\id)\|_{r_j,s_j,y_j}\\
           &\le & \frac{1}{2}    \|\mathsf{W}_{j-1}\phi^{j-2}\|_{r_{j-1}, s_{j-1},y_{j-1}} \cdot \s_{j-1}^d\left(|\vae|^{2^{j-1}}{\mathsf{L}}_{j-1}3\s_{i-1}^{-1}\right)\sqrt{|\vae|}\\
           &\le & \frac{1}{2}    \|\mathsf{W}_{j-1}\phi_1\|_{r_{2}, s_2,y_{2}} \cdot \s_{j-1}^d\;\th_{j-1}\cdot\sqrt{|\vae|}\\
           &\le &  \frac{1}{2}\left(\dst\prod_{i=1}^{j-2}\|\mathsf{W}_{i+1}\mathsf{W}_{i}^{-1}\| \right)\|\mathsf{W}_{1}\phi_1\|_{r_{2}, s_2,y_{2}} \cdot \s_{j-1}^d\;\th_{j-1}\cdot\sqrt{|\vae|}\\
           &\eqby{alfhtrikpi}& \frac{1}{2}\left(\dst\prod_{i=1}^{j-2}\frac{r_i}{r_{i+1}} \right)\|\mathsf{W}_{1}\phi_1\|_{r_{2}, s_2,y_{2}} \cdot \s_{j-1}^d\;\th_{j-1}\cdot\sqrt{|\vae|}\\
           &=& \frac{r_1}{2r_{j-1}}\|\mathsf{W}_{1}\phi_1\|_{r_{2}, s_2,y_{2}} \cdot \s_{j-1}^d\;\th_{j-1}\cdot\sqrt{|\vae|}\\
           &\leby{eqa1maj}& \su2\s_{2}^d\;\mathsf{C}_6\;\eta_0\;\|\mathsf{W}_{1}\phi_1\|_{r_{2}, s_2,y_{2}} \cdot \left(2^{-d}\mathsf{C}_6\;\eta_0\right)^{j-3}\cdot \;\th_0^{2^{j-1}}\cdot\sqrt{|\vae|}\\
           &\le& \su2\s_{2}^d\;\mathsf{C}_6\;\eta_0\;\|\mathsf{W}_{1}\phi_1\|_{r_{2}, s_2,y_{2}} \cdot \left(2^{-d}\mathsf{C}_6\;\eta_0\right)^{2^{j-4}}\cdot \;\th_0^{2^{j-1}}\cdot\sqrt{|\vae|}\\
           &=& \su2\s_{2}^d\;\mathsf{C}_6\;\eta_0\;\|\mathsf{W}_{1}\phi_1\|_{r_{2}, s_2,y_{2}} \cdot \left(\left(2^{-d}\mathsf{C}_6\;\eta_0\right)^{\su8} \th_{0}\right)^{2^{j-1}}\cdot\sqrt{|\vae|}\\
           &=& \su2\s_{2}^d\;\mathsf{C}_6\;\eta_0\;\|\mathsf{W}_{1}\phi_1\|_{r_{2}, s_2,y_{2}} \cdot \left(\mathsf{C}_8\;\eta_0^{\su8}\; \th_{0}\right)^{2^{j-1}}\cdot\sqrt{|\vae|}  \;.
\eeqano
\noi
Therefore, for any $n\ge 1,\, j\geq 0$,
\begin{align*}
\|\mathsf{W}_{1}(\phi^{n+j+1}-\phi^n)\|_{r_{n+j+2},s_{n+j+2},y_{n+j+2}}&\leq  \sum_{i=n}^{n+j}\|\mathsf{W}_{1}(\phi^{i+1}-\phi^i)\|_{r_{i+2},s_{i+2},y_{i+2}}\\
&\le \sum_{i=n}^{n+j}\left(\dst\prod_{k=1}^{i}\|\mathsf{W}_{k}\mathsf{W}_{k+1}^{-1}\| \right)\|\mathsf{W}_{i+1}(\phi^{i+1}-\phi^i)\|_{r_{i+2},s_{i+2},y_{i+2}}\\
&\eqby{alfhtrikpi} \sum_{i=n}^{n+j}\dst\prod_{k=1}^{i}\max\left\{1\;,\frac{r_{k+1}}{r_k} \right\}\|\mathsf{W}_{i+1}(\phi^{i+1}-\phi^i)\|_{r_{i+2},s_{i+2},y_{i+2}}\\
&= \sum_{i=n}^{n+j}\|\mathsf{W}_{i+1}(\phi^{i+1}-\phi^i)\|_{r_{i+2},s_{i+2},y_{i+2}}\\
&\le \su2\s_{2}^d\;\mathsf{C}_6\;\eta_0\;\|\mathsf{W}_{1}\phi_1\|_{r_{2}, s_2,y_{2}}\cdot\sqrt{|\vae|} \dst\sum_{i=n}^{n+j} \left(\mathsf{C}_8\;\eta_0^{\su8}\; \th_{0}\right)^{2^{i+1}}\;.
\end{align*}
Hence, by \equ{condBisv2Prt}, $\phi^j$ converges uniformly on $\{y_*\}\times\torus^d_{s_*}$ to some $\phi^*$, which is then real--analytic map in $x\in\torus^d_{s_*}$.

\nl
To estimate $|\mathsf{W}_0(\phi^*-\id)|$ on $\{y_*\}\times\torus^d_{s_*}$, observe that
, for $i\ge 1$,\footnote{Recall that $2^{i}\ge i+1,\, \forall\, i\ge 0$ and $\s_0\le \su2$.}
$$\s_{i}^d\;|\vae|^{2^i}\mathsf L_i\le \frac{\s_0^{d+1}}{3 \cdot 2^{i(d+1)}}\ \frac{\th_0^{2^i}}{\mathsf{d}_*} \le \su{3 \cdot 2^{(d+1)(i+1)} \mathsf{d}_*}\th_0^{{i+1}}= \su{3\mathsf{d}_*} \Big(\frac{\th_0}{2^{d+1}}\Big)^{i+1}$$
and therefore 
$$\dst\sum_{i\ge 1}  \s_{i}^d\;|\vae|^{2^i}\mathsf L_i\le \frac{1}{3\mathsf{d}_*}\sum_{i\ge 1}\Big(\frac{\th_0}{2^{d+1}}\Big)^{i+1}\le \frac{\th_0^2}{3\cdot {2^{2d+1}}\;\mathsf{d}_*}
\ .$$ 
Moreover, for any $i\ge 1$,
\begin{align*}
\|\mathsf{W}_1(\phi^i-\id)\|_{r_{i+1},s_{i+1},y_{i+1}}&\le \|\mathsf{W}_1(\phi^{i-1}\circ\phi_i-\phi_i)\|_{r_{i+1},s_{i+1},y_{i+1}}+\|\mathsf{W}_1(\phi_i-\id)\|_{r_{i+1},s_{i+1},y_{i+1}}\\
&\le \|\mathsf{W}_1(\phi^{i-1}-\id)\|_{r_{i},s_{i},y_{i}}+ \left(\dst\prod_{j=0}^{i-1}\|\mathsf{W}_{j}\mathsf{W}_{j+1}^{-1}\| \right) \|\mathsf{W}_{i}(\phi_i-\id)\|_{r_{i+1},s_{i+1},y_{i+1}}\\
&= \|\mathsf{W}_1(\phi^{i-1}-\id)\|_{r_{i},s_{i},y_{i}}+ \|\mathsf{W}_{i}(\phi_i-\id)\|_{r_{i+1},s_{i+1},y_{i+1}}\\
&= \|\mathsf{W}_1(\phi^{i-1}-\id)\|_{r_{i},s_{i},y_{i}}+ \|\mathsf{W}_{i}(\phi_i-\id)\|_{r_{i+1},s_{i+1},y_{i+1}}\\
&\le \|\mathsf{W}_1(\phi^{i-1}-\id)\|_{r_{i},s_{i},y_{i}}+\s_{i}^d\;|\vae|^{2^{i}}{\mathsf{L}}_{i}\sqrt{|\vae|}\ ,
\end{align*}
which iterated yields
\begin{align*}
\|\mathsf{W}_1(\phi^i-\id)\|_{r_i,s_i,y_i}&\le \sqrt{|\vae|}\dst\sum_{k=1}^{i-1}\s_{k}^d\; |\vae|^{2^k}{\mathsf{L}}_k\\
&\le  \sqrt{|\vae|}\dst\sum_{k\ge 1}\s_{k}^d\;|\vae|^{2^k}{\mathsf{L}}_k\\
&\le \frac{\th_0^2}{3\cdot {2^{2d+1}}\;\mathsf{d}_*}\sqrt{|\vae|}
\,.
\end{align*}
Therefore, taking the limit over $i$ completes the proof of \equ{estfin2Bis01v2}, Lemma~\ref{lem:2Bisv2}.\\
Now, to complete the proof of the Theorem, just set $\phi_*\coloneqq \phi_0\circ \phi^*$ and observe that, uniformely on $\{y_*\}\times \torus^d_{s_*}$,
\begin{align*}
|\mathsf{W}_0(\phi_*-\id)|&\le |\mathsf{W}_0(\phi_0\circ \phi^*-\phi^*)|+|\mathsf{W}_0(\phi^*-\id)|\\
&\le \|\mathsf{W}_0(\phi_0-\id)\|_{r_1,s_1,y_1}+\|\mathsf{W}_0\mathsf{W}_1^{-1}\|\;|\mathsf{W}_1(\phi^*-\id)|\\
&\le \s_{0}^d\;|\vae|{\mathsf{L}}_0\sqrt{|\vae|}+\frac{\th_0^2}{3\cdot {2^{2d+1}}\;\mathsf{d}_*}\sqrt{|\vae|}\\
&\le \left(\frac{\s_0^{d+1}}{3 }\;\th_0 +\frac{\th_0^2}{3\cdot {2^{2d+1}}\;\mathsf{d}_*}\right)\sqrt{|\vae|}\\
&\le \left(\frac{1}{3\cdot {2^{d+1}} }\;\th_0 +\frac{\th_0^2}{3\cdot {2^{2d+1}}\;\mathsf{d}_*}\right)\sqrt{|\vae|}\\
&\le \frac{\th_0}{3\cdot {2^{d}}}\sqrt{|\vae|}
\,.
\end{align*}
Moreover, for any $i\ge 1$,
\begin{align*}
|y_{i}-y_0|&\le \dst\sum_{j=0}^{i-1}|y_{j-1}-y_j|\\
 	&\leby{estfin2Bis0011v2}8\sqrt{2}\mathsf{T}_0\dst\sum_{j=0}^{i-1}\frac{|\vae|^{2^i} M_i}{r_i}\\
 	&\leby{estfin2Bis00v2} \frac{8\sqrt{2}\mathsf{T}_0}{r_1}\dst\sum_{j=0}^{\infty}\mathbf{a}_1^{i-1}|\vae|\wh M_i\\
 	&\leby{eqa1maj} \frac{64d\eta_0\k_0^{\n}|\vae|}{\a}\dst\sum_{j=0}^{\infty}(\mathsf{C}_6\;\eta_0)^{i-1}\wh M_i\\
 	&= \frac{64d\k_0^{\n}}{\mathsf{C}_6\mathsf{e}_*\mathsf{K}_0}\frac{\wh \a^2|\vae|}{\a}\dst\sum_{j=0}^{\infty}\left(\frac{\mathsf{C}_6\;\eta_0}{\mathsf{d}_*}\right)^{i}\th_0^{2^i}\\
 	&\le \frac{64d\k_0^{\n}}{\mathsf{C}_6\mathsf{e}_*\mathsf{K}_0}{\a}\dst\sum_{j=0}^{\infty}\th_0^{i+1}\\
 	&\le \frac{64d\k_0^{\n}}{\mathsf{C}_6\mathsf{e}_*\mathsf{K}_0}{\a}\dst\cdot 2\th_0\\
 	&= \frac{2^{2\n+7}d}{\mathsf{C}_6}\s_0^{-\n}\m_0\l_0^{\n}\frac{\a}{\mathsf{K}_0}\\
 	&\le \frac{2^{2\n+7}d}{\mathsf{C}_6}\s_0^{-\n}\su{\mathsf{C}_7\eta_0^2}\s_0^{4\n+2d+1}\mathsf{e}_*\m_0\frac{\a}{\mathsf{K}_0}\\
 	&= \frac{2^{2\n+7}d}{\mathsf{C}_6\mathsf{C}_7}\s_0^{3\n+2d+1}\theta_0\frac{\a}{\mathsf{K}_0\eta_0^2}\\
 	&\leby{condBisv2Prt} \frac{2^{2\n+7}d}{\mathsf{C}_6\mathsf{C}_7\mathsf{C}_8}\s_0^{3\n+2d+1}\frac{\a}{\mathsf{K}_0\eta_0^{\frac{17}{8}}}\\
 	&= \su{\mathsf{C}_{9}}\s_0^{3\n+2d+1}\frac{\a}{\mathsf{K}_0\eta_0^{\frac{17}{8}}}\;,
%
\end{align*}
and then passing to the limit, we get
$$
|y_{*}-y_0|\le  
\su{\mathsf{C}_{9}}\s_0^{3\n+2d+1}\frac{\a}{\mathsf{K}_0\eta_0^{\frac{17}{8}}}\;.
$$
  \qed


%
\part{``Sharp'' measure estimates of Kolmogorov's sets\label{part2}}
\chapter{``Explicit'' integrability on a Cantor--like set and a ``sharp'' measure estimate}
\section{Assumptions\label{AssumpExtArnolv2}}
Let $\t\ge d-1\ge 1$ and set\footnote{Notice that each $\mathsf{C}_i$ is greater than $1$ and depends only upon $d$ and $\t$. 
}
\begin{align*}
\n			 &\coloneqq \t+1\;,\\
\mathsf{C}_0 &\coloneqq 4\sqrt{2}\left(\frac{3}{2}\right)^{2\n+d}\dst\int_{\rn} \left( |y|_1^{\n}+d|y|_1^{2\n}\right)\ex^{-|y|_1}dy\;,\\
\mathsf{C}_1 &\coloneqq 2\left(\frac{3}{2}\right)^{\n+d}\dst\int_{\rn} |y|_1^{\n}\ex^{-|y|_1}dy\;,\\
\mathsf{C}_2 &\coloneqq 2^{3d}d\;,\\
\mathsf{C}_3 &\coloneqq	\left(d^2\mathsf{C}_1^2+6d\mathsf{C}_1 +\mathsf{C}_2\right)\sqrt{2}\;,\\
\mathsf{C}_4 &\coloneqq \max\left\{\mathsf{C}_0,\,\mathsf{C}_3\right\}\;,\\
\mathsf{C}_5   &\coloneqq 3d^2\cdot 2^{6\n+2d+11}\dst\max\left\{2^7d\sqrt{2}\,,\,8^{-\n}\mathsf{C}_4 \right\}\;,\\
\mathsf{C}_6 &\coloneqq 2^{\n+\frac{3}{4}d+\frac{53}{8}}d^{\frac{5}{4}}\;,\\
\mathsf{C}_7&\coloneqq 2\ex\;d\left(\frac{3}{2}\right)^{d-1}\;,\\
\mathsf{C}_8 &\coloneqq \frac{\mathsf{C}_5}{3\cdot 2^d} \;,\\
\mathsf{C}_{9} &\coloneqq  2^{3(\n+1)}\;d\;\sqrt{2}\;\mathsf{C}_6\;.\\
\end{align*}
\section{Statement of the extension Theorem}
\thm{Extteo4v2}
Under the assumptions and notations in $\S\ref{AssumpExtArnolv2}$, we have the following. Let $\mathscr D\subset\rn$ be a non--empty, bounded domain.\footnote{\ie open and connected.} 
 Consider the Hamiltonian parametrized by $\vae\in\real$
\[H(y,x;\vae)\coloneqq K(y)+\vae P(y,x),\]
where $K,P$ are real--analytic functions defined on $\mathscr D\times\tn$ with bounded holomorphic extensions to\footnote{Recall the notations in $\S\ref{parassnot}$}
$$
D_{\mathsf{r}_0,s_0}(\mathscr D)\coloneqq \dst\bigcup_{y_0\in\mathscr D}D_{\mathsf{r}_0,s_0}(y_0)\,,
$$
for some $\mathsf{r}_0>0$ and $0<s_0\le 1$, the norm being
$$
\|\cdot\|_{\mathsf{r}_0,s_0,\mathscr D}\coloneqq \dst\sup_{D_{\mathsf{r}_0,s_0}(\mathscr D)}|\cdot|\,.
$$
 Let $\a>0$, $\d>0$ and\footnote{Notice that $\mathscr D_\d$ is closed, connected, with non--empty interior of $\mathscr D_\d$ provided that $\d$ is small enough.}
\begin{align*}
\D_\a^\t &\coloneqq \left\{\o\in\rn:\quad |\o\cdot k|\ge \frac{\a}{|k|_1^\t}\,,\quad\forall\;k\in\zn\setminus\{0\} \right\}\,,\\
\mathscr D_{\d}&\coloneqq \left\{y\in \mathscr D:  B_{\d}(y)\subseteq\mathscr D\right\}\,,\\
\mathscr D_{\d,\a} &\coloneqq \left\{y_0\in\mathscr D_{\d}: \ K_y(y_0)\in \D_\a^\t\right\}\,.
\end{align*}
Assume that 
\beq{ArnoldCondExtv2}
|\det K_{yy}(y)|\neq 0\;,\quad\forall\;y\in \mathscr D_{\d,\a}\;.
\eeq
Fix 
$$0\le \s_0<\min\left\{\su2\;,\frac{d}{2^{2\n-7}}\right\}s_0
$$
 and define\footnote{Notice that $\eta_0\ge 1$.}
\begin{align*}
&s_* \coloneqq s_0-\max\left\{2\;,\frac{2^{2\n-7}}{d}\right\}\s_0\;,\\
&r_0 \coloneqq \dst\min\left\{\mathsf{r}_0,\,{32d}{\d} \right\}\;,\\
&M_0 \coloneqq \|P\|_{\mathsf{r}_0,s_0,\mathscr D}\;,\\
&\mathsf{K}_0 \coloneqq \|K_{yy}\|_{\mathsf{r}_0,\mathscr D}\;,\\
&\mathsf{T}_0 \coloneqq \|T\|_{\mathscr D}\coloneqq \sup_{y_0\in\mathscr D}\|T(y_0)\|\;,\\
&\eta_0		 \coloneqq \mathsf{T}_0\mathsf{K}_0\;,\\
&r_*		   \coloneqq \frac{\s_0^\n}{\mathsf{C}_{9}}\left(\frac{\s_0}{\eta_0}\right)^{\frac{5}{4}}\frac{\a}{\mathsf{K}_0}\;,\\
&\m_*		    \coloneqq \sup\left\{\m\le \ex^{-1}: 2\;\mathsf{C}_5\;\mathsf{C}_6\;\s_0^{4\n+2d+\frac{13}{4}}\;\eta_0^{\frac{13}{4}}\;\m\;\left(\log\m^{-1}\right)^{2\n}\le 1 \right\}\;,
\end{align*}
where $T(y)\coloneqq K_{yy}(y)^{-1}$.
 Finally, assume
\beq{smcEAr0v2}
\boxed{\a\le \frac{r_0\s_0}{\mathsf{T}_0}\qquad\mbox{and}\qquad  |\vae|\le \m_*\frac{\a^2}{\mathsf{K}_0\;M_0}\;.}
\eeq
Then, there exist $\mathscr D_*\subset {\mathscr D_{\d-r_*}}$ having the same cardinality  as $\mathscr D_{\d,\a}$, a lipeomorphism $G^*\colon \mathscr D_{\d,\a}\overset{onto}{\longrightarrow}\mathscr D_*$,  
a function $K_*\in C_W^\infty(\mathscr D_*,\real)$
and  a $C_W^\infty $--symplectic transformation\footnote{Which means that the Whitney--gradient $\nabla \phi_*=\dpr\phi_*/\dpr(y_*,x)$ satisfies $(\nabla \phi_*)\mathbb{J}(\nabla \phi_*)^T=\mathbb{J}$ uniformly on $\mathscr D_*\times \tn$, where $\mathbb{J}=\begin{pmatrix}0 & -\uno_d\\
\uno_d & 0\end{pmatrix}$.} $\phi_*\colon \mathscr D_*\times \tn\to \mathscr K\coloneqq \phi_*(\mathscr D_*\times \tn)\subset \mathscr D\times \tn$ and real--analytic in $x\in\torus^d_{s_*}$,
 such that\footnote{Notice that the derivatives are taken in the sense of Whitney.} 
\begin{align}
\dpr_{y_*}K_*\circ G^*&=\dpr_{y}K  \qquad\quad\quad\qquad\ \quad \mbox{on} \quad \mathscr D_{\d,\a}\;,\label{conjCaneq00v2}\\
\dpr^\b_{y_*}(H\circ \phi_*)(y_*,x)&=\dpr^\b_{y_*}K_*(y_*),\qquad \forall\;(y_*,x)\in \mathscr D_*\times\tn,\quad \forall\; \b\in\natural_0^d \label{conjCaneq0v2}
\end{align}
and 
\begin{align}
&\|G^*-\id\|_{\mathscr D_{\d,\a}}
\le 
r_*\,,\label{NormGstrThtv2}\\
&\|G^*-\id\|_{L,\mathscr D_{\d,\a}}\le \frac{\ex\;\s_0^{\n+d}}{\mathsf{C}_6}\,,\label{LipGstrThtv2}\\
&\meas(\mathscr D\times\tn\setminus \mathscr K)\le \left(1+\frac{d\;\ex\;\s_0^{\n+d}}{\mathsf{C}_6}\right)^d\bigg(\meas\big(({ B}_{\d \s_0}(\mathscr D)\setminus \mathscr D)\times\tn\big)+\nonumber\\
&\hspace{5cm}+\meas\big((\mathscr D\setminus \mathscr D_\d)\times\tn\big)+\meas((\mathscr D_\d\setminus \mathscr D_{\d,\a})\times\tn\big)\bigg)\;.\label{MesChPin}
\end{align}
\ethm

\noi
Now, by applying Theorem~\ref{arnolMeas1} (see Appendix~\ref{appF}) to \equ{MesChPin}, we get the following measure estimate of the unstable set $\mathscr D\times\tn\setminus \mathscr K$.
\thm{Extteo5v2}
Let the notations and assumptions in Theorem~\ref{Extteo4v2} hold, with
\beq{EqStnAlf}
0< \a\le \frac{32d\s_0}{\mathsf{T}_0}\min\left\{\frac{\mathsf{r}_0}{32d}\,,\,\frac{R(\mathscr D)}{3}\,,\,\minfoc(\dpr \mathscr D)\right\}\;,\qquad |\vae|\le \m_*\frac{\a^2}{\mathsf{K}_0\;M_0}\;.
\eeq
in place of \equ{smcEAr0v2} and
$$
\d\coloneqq\frac{\a\mathsf{T}_0}{32d\s_0}\;,
$$
where\footnote{Notice that the first condition in \equ{EqStnAlf} then reads $0<\d\le \min\left\{\frac{\mathsf{r}_0}{32d},\frac{R(\mathscr D)}{3},\minfoc(\dpr \mathscr D)\right\}$. The condition $\d\le \frac{R(\mathscr D)}{3}$ ensures that the interior of $\mathscr D_\d$ is non--empty.}
$$
R(\mathscr D)\coloneqq \sup\{R>0: B_R(y)\subseteq \mathscr D\;, \mbox{ for some } y\in\mathscr D\}\;.
$$
Futhermore, assume that the boundary $\dpr \mathscr D$ of $\mathscr D$ is a smooth hypersurface of $\rn$. Then, the conclusions in Theorem~\ref{Extteo4v2} still hold. Moreover,
\beq{MesChPinBiso}
\meas(\mathscr D\times\tn\setminus \mathscr K)\le (3\pi)^d \frac{\mathsf{T}_0}{32d\s_0}\bigg(2\mathcal{H}^{d-1}(\dpr \mathscr D)\;\a+ C(d,\s_0,\mathsf{T}_0,\mathbf{R}^{\dpr\mathscr D})\;\a^2+\meas(\mathscr D_\d\setminus \mathscr D_{\d,\a})\bigg)\,,
\eeq
where\footnote{See Appendix~\ref{appF} for the definitions.
} $\mathbf{R}^{\dpr\mathscr D}$ denotes the curvature tensor of $\dpr\mathscr D$, $\mathbf{k}_{2j}(\mathbf{R}^{\dpr\mathscr D})$, the $(2j)$--th integrated mean curvature of $\dpr\mathscr D$ in $\rn$ and 
$$
C(d,\s_0,\mathsf{T}_0,\mathbf{R}^{\dpr\mathscr D})\coloneqq \frac{\mathsf{T}_0}{16d\s_0}\dst\sum_{j=1}^{\left[\frac{d-1}{2}\right]}\frac{\mathbf{k}_{2j}(\mathbf{R}^{\dpr\mathscr D})\;}{1\cdot3\cdots (2j+1)} \left(\frac{\a\mathsf{T}_0}{32d\s_0}\right)^{2j-1}\;.
$$
\ethm

\rem{DstD0Fi0}
{\bf (i)} Notice that \equ{MesChPin} is mainly a consequence of \equ{NormGstrThtv2}; the crucial part of the proof is that one can actually extend a Lipschitz continuous function to a global Lipschitz continuous function without increasing  neither the Lipschitz constant nor the sup--norm (see Theorem~\ref{MintyExt} in Appendix~\ref{appE}).\\
{\bf (ii)}The following estimates hold as well:
\begin{align}
&|\meas(\mathscr D_*)-\meas(\mathscr D_{\d,\a})|\le \frac{\mathsf{C}_7}{2\mathsf{C}_6}\;\s_{0}^{\n+d+\frac{5}{4}}\;\eta_{0}^{-\frac{5}{4}}\meas(\mathscr D_{\d,\a})\;,\label{measKKstv2}\\
&|\mathsf{W}_0(\phi_*-\id)|\le \frac{1}{3\cdot {2^{d+1}}\mathsf{C}_6}\left(\frac{\s_0}{\eta_0}\right)^{\frac{5}{4}}\qquad\ \mbox{on}\quad \mathscr D_*\times\torus^d_{s_*}
\;,\label{estArnTrExtv2}
\end{align}
where
$$
\mathsf{W}_0   \coloneqq \diag\left(\frac{\mathsf{K}_0}{{\a}}\uno_d,\uno_d\right)\;.
$$
Notice that the constant in \equ{measKKstv2} is of order 1 and not $\a$; that is why we need Minty's Theorem (see {\bf (i)} above).\\
{\bf (iii)} Notice that the  Theorem is consistent for $\s_0=0$. In fact, in that case
$$
\vae=\a\eqby{smcEAr0v2}0.
$$
Hence, the Hamiltonian $H$ is integrable. Moreover,
$$
\mathscr D_{\d,\a}=\mathscr D_\d\;,\qquad G^*\eqby{NormGstrThtv2} \id\;,\qquad \phi_*\eqby{estArnTrExtv2} \id\;.
$$
Thus,
$$
\mathscr D_*=G^*(\mathscr D_{\d,\a})=\mathscr D_{\d,\a}=\mathscr D_\d\;.
$$
Therefore, we get $\mathscr K=\phi_*(\mathscr D_*\times \tn)=\mathscr D_\d\times \tn$, for any $\d>0$, 
 as expected.
\erem
\section{Proof of Theorem~\ref{Extteo4v2}}

\lemtwo{lem:1Extv5}{KAM step}
Let $r>0,\; 0<2\s< s\le 1$ and consider the Hamiltonian parametrized by $\vae\in \real$
$$
H(y,x;\vae)\coloneqq K(y)+\vae P(y,x)\;,
$$
where $K,P$ are real--analytic functions with bounded holomorphic extensions to $D_{r,s}(\mathscr D_\sharp)$.\\
Assume that\footnote{In the sequel, $K$ and $P$  stand for  generic real analytic Hamiltonians which, later on, will respectively play the roles of $K_j$ and $P_j$,  and $y_0,\,r$, the roles of $y_j,\,r_j$ in the iterative step.} 
\beq{RecHypArnExtv5}
\begin{aligned}
&
\det K_{yy}(y)\not=0
\;,\qquad\qquad\qquad T(y)\coloneqq K_{yy}(y)^{-1}\;,\quad \forall\;y\in\mathscr D_\sharp\;,\\
&\|K_{yy}\|_{r,\mathscr D_\sharp}\le \mathsf{K}\;,\qquad\qquad\qquad\ \, \|T\|_{\mathscr D_\sharp}\le \mathsf{T}\;,\\
& \|P\|_{r,s,\mathscr D_\sharp}\le M \;,\qquad\qquad\qquad\,\,\, K_y(\mathscr D_\sharp)\subset \D^\t_\a\;. 
\end{aligned}
\eeq
\noi
Fix  $\vae\neq0$ and assume that
\beq{DefNArnExt1v5}
\l\ge \log\left(\s^{2\n+d}\frac{{\a}^2}{|\vae|{M}\mathsf{K}}\right)\ge 1 \;.
\eeq
Let
\beq{DefNArn2v5}
\begin{aligned}
&
\k\coloneqq 4\s^{-1}\l, \quad 
\bar{r}\le
\dst\min\left\{\frac{\a}{2d\mathsf{K}\k^{\n}}\,,\, \frac{r\s}{16d\mathsf{T} \mathsf{K}} \right\},\\
&\tilde r\coloneqq \frac{\bar{r}}{16d\mathsf{T} \mathsf{K}},\quad\quad \bar{s}\coloneqq s-\frac{2}{3}\s,\quad s'\coloneqq s-\s \,,
\end{aligned}
\eeq
and\footnote{Notice that $\mathsf{L}\ge \s^{-d}\ovl{\mathsf{L}}\ge \ovl{\mathsf{L}}$ since $\s\le 1$. Notice also that $\mathsf{T}\mathsf{K}\ge 1$, so that $\frac{16\mathsf{T} }{r\bar{r}}\s^{-(\n+d)}\ge \frac{16\mathsf{T} }{r^2}\ge \frac{4}{\mathsf{K}r^2}$.
}
\begin{align*}
\ovl{\mathsf{L}}&\coloneqq\frac{\mathsf{C}_0}{\sqrt{2}} \max\left\{1,\frac{\a}{r\mathsf{K}}\right\}\frac{\mathsf{K}_0 M }{\a^2 }\s^{-(2\n+d)}\;,\\
\mathsf{L}&\coloneqq M\dst\max\left\{\frac{16\mathsf{T}  }{r\bar{r}}\s^{-(\n+d)}\,,\,
\frac{\mathsf{C}_4}{\sqrt{2}} \max\left\{1,\frac{\a}{r\mathsf{K}}\right\}\frac{\mathsf{K} }{\a^2 }\s^{-2(\n+d)}\right\}\\
          &=M\dst\max\left\{\frac{16\mathsf{T}  }{r\bar{r}}\s^{-(\n+d)}
\,,\,\frac{4}{\mathsf{K} r^2}\,,\,\frac{\mathsf{C}_4}{\sqrt{2}} \max\left\{1,\frac{\a}{r\mathsf{K}}\right\}\frac{\mathsf{K} }{\a^2 }\s^{-2(\n+d)}\right\}
\;.
\end{align*}
Then, there exists a generating function $(y',x)\mapsto y'\cdot x+g(y',x)$, with $g\in \mathcal{B}_{\bar r,\bar s}(\mathscr D_\sharp)$
 and satisfying the following inequalities:
\beq{Est1Lem1bExtv5}
\left\{
\begin{aligned}
&\|g_x\|_{\bar{r},\bar{s},\mathscr D_\sharp}\le  \mathsf{C}_1 \frac{M}{\a} \s^{-(\n+d)}\,,\\
& \|g_{y'}\|_{\bar{r},\bar{s},\mathscr D_\sharp},\, \|\dpr_{y'x}^2 g\|_{\bar{r},\bar{s},\mathscr D_\sharp}\le \ovl{\mathsf{L}}\,,\\
&\|\dpr_{y'}^2\wt K\|_{\bar{r},\mathscr D_\sharp}\le 
\mathsf{K}\mathsf{L}\,,
\end{aligned}
\right.
\eeq
where
$$\wt K(y')\coloneqq \average{P(y',\cdot)}\;.$$
 If, in addition,   
\beq{cond1ExtExtv5}
|\vae|{\mathsf{L}}\le \frac{\sigma}{3}
\ ,
\eeq
then, there exists 
a  diffeomorphism $G\colon  D_{\tilde{r}}(\mathscr D_\sharp){\to}G( D_{\tilde{r}}(\mathscr D_\sharp))$, 
such that ,
\beq{convEstExt}
\left\{
\begin{aligned}
&\dpr_{y'}K'\circ G=\dpr_{y}K\,,\qquad\qquad\qquad\quad\quad\ \; (\dpr_{y'}^2K')\circ G\neq 0\quad\mbox{ on } \mathscr D_\sharp\,,\\
&|\vae|\|g_x\|_{\bar{r},\bar{s},\mathscr D_\sharp}\le \frac{r}{3}\,,\qquad \qquad\quad\qquad\quad\ \;\,\,
\|G-\id\|_{\tilde{r},\mathscr D_\sharp}\le \s^{\n+d}\bar{r}|\vae|\mathsf{L}\;,\\
&|\vae|\|\wt T\|_{\mathscr D_\sharp'}\le \mathsf{T}|\vae|\mathsf{L}
\,, \qquad\quad\quad\quad\qquad\quad\,\;\; \|\dpr_z G-\uno_d\|_{\tilde{r},\mathscr D_\sharp}\le \s^{\n+d}|\vae|\mathsf{L}\,,\\
&\|P_+\|_{\bar{r},\bar s,\mathscr D_\sharp} \le  \mathsf{L}M\,, \qquad\qquad\quad\qquad\quad\ \,\,\, B_{\bar{r}/2}(\mathscr D_\sharp')\subset B_{\bar{r}}(\mathscr D_\sharp)\subset\mathscr D\,,
\end{aligned}
\right. 
\eeq
where
\begin{equation*}
\begin{aligned}
&\mathscr D_\sharp'\coloneqq G(\mathscr D_\sharp)\;,   \\
&\left(\dpr^2_{y'} K'(\mathsf{y}')\right)^{-1}\eqqcolon T\circ G^{-1}(\mathsf{y}')+\vae\;\wt T(\mathsf{y}')\;, \\ 
&P_+(y',x)\coloneqq P(y'+\vae g_x(y',x),x)\;.
\end{aligned}
\qquad\qquad
\begin{aligned}
&K'\coloneqq K+\vae\wt K\;,\\
 &\forall\; \mathsf{y}'\in \mathscr D_\sharp'\\ 
&
\end{aligned}
\end{equation*}
and the following hold. For any $y'\in D_{\bar{r}}(\mathscr D_\sharp)$, the map $\psi_\vae(x):=  x+\vae g_{y'}(y',x)$ has an analytic inverse 
 $\f(x')=x'+\vae \wt{\f}(y',x';\vae)$  such that 
\beq{boundalExtExt}
\|\wt{\f}\|_{\bar{r}, s',\mathscr D_\sharp}\le   \ovl{\mathsf{L}} \qquad {\rm and}
\quad
\f=\id + \vae \wt{\f} : D_{\bar{r}/2,s'}(\mathscr D_\sharp')\to \torus^d_{\bar{s}} \ ;
\eeq
for any $y_0\in\mathscr D_\sharp$ and $(y',x)\in D_{\bar{r},\bar s}(y_0)$, $|y'+\vae g_x(y',x)-y_0|<\frac{2}{3} r$; the map $\phi'$ is a symplectic diffeomorphism and
\beq{phiokExt0Ext}
\phi'=\big( y'+\vae g_x(y', \f(y',x')),\f(y',x')\big): D_{\bar{r}/2,s'}(\mathscr D_\sharp')\to D_{2r/3, \bar{s}}(\mathscr D_\sharp),
\eeq
with
\beq{phiokExt1Ext}
\|\mathsf{W}\,\tilde \phi\|_{\bar{r}/2,s',\mathscr D_\sharp'}\le \s^d{\mathsf{L}}\,,
\eeq
where $\tilde \phi$ is defined by the relation $\phi'=:\id + \vae \tilde \phi$,
$$
\mathsf{W}\coloneqq \begin{pmatrix}
\max\{\frac{\mathsf{K}}{{\a}}\;,\frac{1}r\}\;\uno_d & 0\\ \ \\
0			& \uno_d 
\end{pmatrix}
$$
and
\beq{tesitExtExt}
\begin{aligned}
&\|P'\|_{\bar{r}/2, s',\mathscr D_\sharp'}\le  \mathsf{L}M\;,
\end{aligned}
\eeq
with
$$
P'(y',x')\coloneqq P_+(y',\f(x'))=P\circ \phi'(y',x')\;.
$$
\elem
\proof\\
\noi
\Giu
{\bf Step 1: Construction of the Arnold's transformation} We seek for $r_1<r/2,\,s_1<s$, a set $\mathscr D_\sharp'\subset D_{2r_1}(\mathscr D_\sharp)$ having the same cardinality as $\mathscr D_\sharp$ and  a near--to--the--identity real--analytic symplectic transformation $\phi_1: \mathscr D\times\tn\righttoleftarrow$ satisfying
\[\phi'\colon D_{r_1,s_1}(\mathscr D_\sharp')\to D_{r,s}(\mathscr D_\sharp),\]
with $D_{r_1,s_1}(\mathscr D_\sharp')\subset D_{r,s}(\mathscr D_\sharp)$ and $\phi'$  generated by $y'\cdot x+\vae g(y',x)$ \ie
\beq{ArnTraKamExtv5}
\phi'\colon \left\{\begin{aligned}
y  &=y'+\vae\ g_x(y',x)\\
x' &=x+\vae g_{y'}(y',x)\, ,
\end{aligned}
\right.
\eeq
such that
\beq{ArnH1Extv5}
\left\{
\begin{aligned}
& H':= H\circ \phi'=K'+\vae^2 P'\quad\mbox{on }D_{r_1,s_1}(\mathscr D_\sharp')\,,\\
& \det\dpr_{y'}^2 K'(y')\not=0\,,\quad\qquad\quad\;\forall\; y'\in \mathscr D_\sharp'\,,\\
& \dpr_{y'} K'(\mathscr D_\sharp')= \dpr_{y} K(\mathscr D_\sharp)\,.
\end{aligned}
\right.
\eeq
By Taylor's formula, we get\footnote{Recall that $\average{\cdot}$ stands for the average over $\tn$.}
\beq{Arneq11Extv5}
\begin{aligned}
H(y'+\vae g_x(y',x),x)=&K(y')+\vae \wt K(y') +\vae \left[K_y(y')\cdot g_x +T_{\k} P(y',\cdot)-\wt K(y') \right]+\\
						&+\vae^2 \left( P^\ppu+P^\ppd+ P^\ppt\right)(y',x) \\
			= & K'(y')+\vae \left[K_y(y')\cdot g_x +T_{\k} P(y',\cdot)-\wt K(y') \right]+ \vae^2 P'(y',x),
\end{aligned}
\eeq
with $\k\in\natural$, which will be chosen large enough so that $P^\ppt=O(\vae)$ and 
\beq{ArnDefPsExtv5}
\left\{
\begin{aligned}
P_+&\coloneqq P^\ppu+P^\ppd+ P^\ppt\\
P^\ppu &\coloneqq \su{\vae^2}\left[K(y'+\vae g_x)-K(y')-\vae K_y(y')\cdot g_x \right]=\dst\int^1_0(1-t)K_{yy}(\vae t g_x)\cdot g_x\cdot g_x dt\\
P^\ppd &\coloneqq \su\vae \left[P(y'+\vae g_x,x)-P(y',x)\right]=\dst\int_0^1P_y(y'+\vae t g_x,x)\cdot g_x dt\\
P^\ppt &\coloneqq \su\vae \left[ P(y',x)-T_{\k} P(y',\cdot)\right]=\su\vae \dst\sum_{|n|_1>\k} P_n(y')\ex^{in\cdot x}\; .
\end{aligned}
\right.
\eeq
By the non--degeneracy condition in \eqref{ArnoldCondExt} and Lemma~\ref{inv1}, 
 for $\vae$ small enough (to be made precised below), there exists $\bar{r}\le r$ such that 
for each $\mathsf{y}\in \mathscr D_\sharp$, there exists a unique $\mathsf{y}'\in D_{\bar{r}}(\mathsf{y})$ satisfying $\dpr_{y'}K'(\mathsf{y}')=\dpr_y K(\mathsf{y})$ and $\det\dpr_{y'}^2 K'(\mathsf{y}')\not=0$; $\mathscr D_\sharp'$ is precisely the set of these $\mathsf{y}'$ when $\mathsf{y}$ runs in $\mathscr D_\sharp$.  
More precisely, $\mathscr D_\sharp'$ and $\mathscr D_\sharp$ are 
``diffeomorphic''\footnote{\ie there a exits a bijection from $\mathscr D_\sharp$ onto $\mathscr D_\sharp'$ which extends to a diffeomorphism on some neighborhood of $\mathscr D_\sharp$.}, say via $G$,
and, for each $\mathsf{y}'\in\mathscr D_\sharp'$, the matrix $\dpr^2_{y'}K_1(y_1)$ is invertible with inverse of the form
$$
T'(\mathsf{y}')\coloneqq \dpr^2_{y'}K'(\mathsf{y}')^{-1}\eqqcolon T(y_0)+\vae\wt T(\mathsf{y}'),\quad \mathsf{y}'=G(\mathsf{y}).
$$
In view of \eqref{Arneq11Extv5}, in order to get the first part of \eqref{ArnH1Extv5}, we need to find $g$ such that  $K_y(y')\cdot g_x +T_{\k} P(y',\cdot)-\wt K(y')$ vanishes; such a $g$ is indeed given by
 \beq{HomEqArn1v5}
 g(y',x)\coloneqq \dst\sum_{0<|n|_1\leq \k} \frac{-P_n(y')}{iK_y(y')\cdot n}\ex^{in\cdot x},
 \eeq
provided that 
\beq{CondHomEqArn1v5}
K_y(y')\cdot n\not= 0, \quad \forall\; 0<|n|_1\leq \k,\quad \forall\; y'\in D_{r_1}(\mathscr D_\sharp')\quad  \left(\subset D_{r}(\mathscr D_\sharp)\right).
\eeq
But, in fact, since $K_y(\mathsf{y})$ is rationally independent, for each $\mathsf{
y}\in\mathscr D_\sharp$, then, given any $\k\in\natural$, there exists $r'\leq r$ such that
\beq{CondHomEqArnExtv5}
K_y(y')\cdot n\not=0,\quad \forall\; 0<|n|_1\leq \k, \quad\forall\; y'\in D_{{r}'}(\mathscr D_\sharp).
\eeq
Then we invert the function $x\mapsto x+\vae g_{y'}(y',x)$ in order to define $P'$. But, by Lemma~\ref{IFTLem}, for $\vae$ small enough, the map $x\mapsto x+\vae  g_{y'}(y',x)$ admits an real--analytic inverse of the form
\beq{InvComp2Fiv5}
\f(y',x';\vae)\coloneqq x'+\vae \wt{\f}(y',x';\vae),
\eeq
so that the Arnod's symplectic transformation is given by
\beq{ArnTrans0v5}
\phi_1\colon (y',x')\mapsto \left\{
\begin{aligned}
y &= y'+\vae g_x(y',\f(y',x'))\\
x &= \f(y',x';\vae)= x'+\vae \wt{\f}(y',x';\vae) .
\end{aligned}
\right.
\eeq
Hence, \eqref{ArnH1Extv5} holds with
\beq{DefP1ArExtv5}
P'(y',x')\coloneqq P'(y', \f(y',x')).
\eeq
{\bf Step 2}
Above all, notice that\footnote{Recall footnote \textsuperscript{\ref{ftnTK1}}.}
\beq{rrbarAsExtv5}
\bar{r}\le\frac{r\s}{16d\mathsf{T} \mathsf{K}}\le \frac{r}{32d}< \frac{r}{2}
\,.\\
\eeq
We begin by extending the ``diophantine condition w.r.t. $K_y$'' uniformly to $D_{\bar{r}}(\mathscr D_\sharp)$ up to the order $\k$. Indeed, 
for any $\mathsf{y}\in\mathscr D_\sharp$, $0<|n|_1\leq \k$ and $y'\in D_{\bar{r}}(\mathsf{y})$,
\begin{align}
|K_y(y')\cdot n|&=|\o\cdot n +(K_y(y')-K_y(\mathsf{y}))\cdot n|\geq |\o\cdot n|\left(1-d\frac{\|K_{yy}\|_{\bar{r},\mathscr D_\sharp}}{|\o\cdot n|}|n|_1\bar{r}\right) \nonumber\\
         &\geq \frac{\a}{|n|_1^\t}\left(1-\frac{d\mathsf{K}}{\a }|n|_1^{\t+1}\bar{r} \right)\geq \frac{\a}{|n|_1^\t}\left(1-\frac{d\mathsf{K}}{\a }\k^{\t+1}\bar{r} \right)
\nonumber\\
         &\ge \frac{\a}{2|n|_1^\t},\label{ArnExtDiopCondExtv5}
\end{align}
so that, by Lemma~\ref{fce}--$(i)$, we have 

\begin{align*}
\|g\|_{\bar{r},\bar{s},\mathscr D_\sharp} &\overset{def}{=}\dst\sup_{D_{\bar{r},\bar{s}}(\mathscr D_\sharp)}\left|\dst\sum_{0<|n|_1\leq \k}\frac{P_n(y')}{K_y(y')\cdot n}\ex^{in\cdot x} \right|\leq \dst\sum_{0<|n|_1\leq \k}\frac{\|P_n\|_{\bar{r},\bar{s}, \mathscr D_\sharp}}{|K_y(y')\cdot n|}\ex^{\left(s-\frac{2}{3}\s\right)|n|_1}\\
   &\leq \dst\sum_{0<|n|_1\leq \k} M\ex^{-s|n|_1}\frac{2|n|_1^{\t}}{\a}\ex^{\left(s-\frac{2}{3}\s\right)|n|_1}\leq \frac{2M}{\a}\dst\sum_{n\in\zn} |n|_1^{\t}\ex^{-\frac{2}{3}\s|n|_1}\\
   &\leq \frac{2M}{\a}\dst\int_{\rn} |y|_1^{\t}\ex^{-\frac{2}{3}\s|y|_1}dy\\
   &= \left(\frac{3}{2\s}\right)^{\t+d}\frac{2M}{\a}\dst\int_{\rn} |y|_1^{\t}\ex^{-|y|_1}dy\\
   &\le \mathsf{C}_1 \frac{M}{\a} \s^{-(\t+d)}
\end{align*}
and analogously,
\begin{align*}
\|g_x\|_{\bar{r},\bar{s},\mathscr D_\sharp} &\overset{def}{=}\dst\sup_{D_{\bar{r},\bar{s}}(\mathscr D_\sharp)}\left|\dst\sum_{0<|n|_1\leq \k}\frac{nP_n(y')}{K_y(y')\cdot n}\ex^{in\cdot x} \right|\leq \dst\sum_{0<|n|_1\leq \k}\frac{\|P_n\|_{\bar{r},\bar{s}, \mathscr D_\sharp}}{|K_y(y')\cdot n|}|n|_1\ex^{\left(s-\frac{2}{3}\s\right)|n|_1}\\
   &\leq \dst\sum_{0<|n|_1\leq \k} M\ex^{-s|n|_1}\frac{2|n|_1^{\t+1}}{\a}\ex^{\left(s-\frac{2}{3}\s\right)|n|_1}\leq \frac{2M}{\a}\dst\sum_{n\in\zn} |n|_1^{\t+1}\ex^{-\frac{2}{3}\s|n|_1}\\
   &\leq \frac{2M}{\a}\dst\int_{\rn} |y|_1^{\t+1}\ex^{-\frac{2}{3}\s|y|_1}dy\\
   &= \left(\frac{3}{2\s}\right)^{\t+d+1}\frac{2M}{\a}\dst\int_{\rn} |y|_1^{\t+1}\ex^{-|y|_1}dy\\
   &\le \mathsf{C}_1 \frac{M}{\a} \s^{-(\t+d+1)}\,,
\end{align*}
\begin{align*}
\|\dpr_{y'}g\|_{\bar{r},\bar{s},\mathscr D_\sharp} &\overset{def}{=}\dst\sup_{D_{\bar{r},\bar{s}}(\mathscr D_\sharp)}\left|\dst\sum_{0<|n|_1\leq \k}\left(\frac{ \dpr_yP_n(y')}{K_y(y')\cdot n}-P_n(y')\frac{ K_{yy}(y')n}{(K_y(y')\cdot n)^2}\right)\ex^{in\cdot x} \right|\\
   &\leq \dst\sum_{0<|n|_1\leq \k}\dst\sup_{D_{\bar{r}}(\mathscr D_\sharp)}\left(\frac{\|(P_y)_n\|_{\bar{r},s,\mathscr D_\sharp}}{|K_y(y')\cdot n|}+d\|P_n\|_{r,s,\mathscr D_\sharp}\frac{\|K_{yy}\|_{r,\mathscr D_\sharp}|n|_1}{|K_y(y')\cdot n|^2}\right)\ex^{\left(s-\frac{2}{3}\s\right)|n|_1}\\
   &\stackrel{\equ{RecHypArnExtv5}+\equ{ArnExtDiopCondExtv5}}{\le} \dst\sum_{0<|n|_1\leq \k}\left( \frac{M}{r-\bar{r}}\ex^{-s|n|_1}\frac{2|n|_1^{\t}}{\a}+dM\ex^{-s|n|_1}\mathsf{K}|n|_1\left(\frac{2|n|_1^{\t}}{\a}\right)^2\right)\ex^{\left(s-\frac{2}{3}\s\right)|n|_1}\\
   &\leby{rrbarAsExtv5} \frac{4M}{\a^2 r}\dst\sum_{0<|n|_1\leq \k}\left( |n|_1^{\t}\a +dr\mathsf{K}|n|_1^{2\t+1}\right)\ex^{-\frac{2}{3}\s|n|_1}\\
   &\le \max(\a,r\mathsf{K})\frac{4M}{\a^2 r}\dst\sum_{0<|n|_1\leq \k}\left( |n|_1^{\t}+d|n|_1^{2\t+1}\right)\ex^{-\frac{2}{3}\s|n|_1}\\
   &\leq \max(\a,r\mathsf{K})\frac{4M}{\a^2 r}\dst\int_{\rn} \left( |y|_1^{\t}+d|y|_1^{2\t+1}\right)\ex^{-\frac{2}{3}\s|y|_1}dy \\
   &= \left(\frac{3}{2\s}\right)^{2\t+d+1}\max(\a,r\mathsf{K})\frac{4M}{\a^2 r}\dst\int_{\rn} \left( |y|_1^{\t}+d|y|_1^{2\t+1}\right)\ex^{-|y|_1}dy\\
   &\le \frac{\mathsf{C}_0}{\sqrt{2}} \max\left\{1,\frac{\a}{r\mathsf{K}}\right\}\frac{\mathsf{K} M }{\a^2 }\s^{-(2\t+d+1)}\\
   &\le \ovl{\mathsf{L}} \;,
\end{align*}
\begin{align*}
\|\dpr^2_{y'x}g\|_{\bar{r},\bar{s},\mathscr D_\sharp} &\overset{def}{=}\dst\sup_{D_{\bar{r},\bar{s}}(\mathscr D_\sharp)}\left|\dst\sum_{0<|n|_1\leq \k}\left(\frac{ \dpr_yP_n(y')}{K_y(y')\cdot n}-P_n(y')\frac{ K_{yy}(y')n}{(K_y(y')\cdot n)^2}\right)\cdot n\ex^{in\cdot x} \right|\\
   &\leq \dst\sum_{0<|n|_1\leq \k}\dst\sup_{D_{\bar{r}}(\mathscr D_\sharp)}\left(\frac{\|(P_y)_n\|_{\bar{r},s, \mathscr D_\sharp}}{|K_y(y')\cdot n|}+d\|P_n\|_{r,s,\mathscr D_\sharp}\frac{\|K_{yy}\|_{r,\mathscr D_\sharp}|n|_1}{|K_y(y')\cdot n|^2}\right)|n|_1\ex^{\left(s-\frac{2}{3}\s\right)|n|_1}\\
   &\le \max(\a,r\mathsf{K})\frac{4M}{\a^2 r}\dst\sum_{0<|n|_1\leq \k}\left( |n|_1^{\t}+d|n|_1^{2\t+1}\right)|n|_1\ex^{-\frac{2}{3}\s|n|_1}\\
   &\leq \max(\a,r\mathsf{K})\frac{4M}{\a^2 r}\dst\int_{\rn} \left( |y|_1^{\t}+d|y|_1^{2\t+1}\right)|y|_1\ex^{-\frac{2}{3}\s|y|_1}dy \\
   &= \left(\frac{3}{2\s}\right)^{2\t+d+2}\max(\a,r\mathsf{K})\frac{4M}{\a^2 r}\dst\int_{\rn} \left( |y|_1^{\t+1}+d|y|_1^{2\t+2}\right)\ex^{-|y|_1}dy\\
   &= \frac{\mathsf{C}_0}{\sqrt{2}} \max\left\{1,\frac{\a}{r\mathsf{K}}\right\}\frac{\mathsf{K} M }{\a^2 }\s^{-(2\n+d)} \\
   &=\ovl{\mathsf{L}}\;,
\end{align*}
and 
\begin{align*}
&\|\dpr_{y'}\wt K\|_{\bar{r},\mathscr D_\sharp}=\| \average{P_y}\|_{\bar{r},\mathscr D_\sharp}\leq \|P_y\|_{\bar{r},\bar{s}, \mathscr D_\sharp}\leq  \frac{M}{r-\bar{r}}\leq \frac{2M}{r} \;,\\
&\|\dpr_{y'}^2\wt K\|_{\bar{r},\mathscr D_\sharp}=\| \average{P_{yy}}\|_{\bar{r},\mathscr D_\sharp}\leq \|P_{yy}\|_{\bar{r},\bar{s}, \mathscr D_\sharp}\leq  \frac{M}{(r-\bar{r})^2}\leq \frac{4M}{r^2}\le \mathsf{K}\mathsf{L} 
\end{align*}
\noi
Next, we construct $\mathscr D_\sharp'$ 
 in \eqref{ArnH1Extv5}. 
 For, fix  $\mathsf{y}\in\mathscr D_\sharp$ 
  and consider
\begin{align*}
F\colon D_{\bar{r}}(\mathsf{y})\times D_{\tilde{r}}(\mathsf{y}) &\longrightarrow \qquad \cn\\
		(y,z)\quad &\longmapsto K_y(y)+\vae \wt K_{y'}(y)-K_y(z).
\end{align*}
Then
\begin{itemize}
\item $F_y(\mathsf{y},\mathsf{y})=\dpr^2_{y}K(\mathsf{y})+\vae\dpr^2_{y'}\wt K(\mathsf{y}) = T(\mathsf{y})^{-1}\left(\uno_d+\vae T(\mathsf{y})\dpr^2_{y'}\wt K(\mathsf{y}) \right)\eqqcolon T(\mathsf{y})^{-1}(\uno_d+\vae A_0)$ and
\begin{align*}
\|\vae A_0\|&\le \|T(\mathsf{y})\|\|\vae \dpr^2_{y'}\wt K(\mathsf{y})\|\le \mathsf T \frac{4|\vae|  M}{r^2}\leby{rrbarAsExtv5} |\vae|\frac{2\mathsf T M}{r\bar{r}}
\le \su2|\vae|\mathsf L\leby{cond1ExtExtv5} \frac{\s}{6}<\su2.
\end{align*}
Hence, $F_y(\mathsf{y},\mathsf{y})$ is invertible, with inverse 
$$
T_0\coloneqq (\uno_d+\vae A_0)^{-1}T(\mathsf{y})=\left(\uno_d+\dst\sum_{k\geq 1}(-\vae)^k A_0^k\right) T(\mathsf{y})
$$
 satisfying
\beq{T1estExt}
\|T_0\|\le \frac{\|T(\mathsf{y})\|}{1-\|\vae A_0\|}\le 2\mathsf T.
\eeq
\item For any $(y,z)\in D_{\bar{r}}(\mathsf{y})\times D_{\tilde{r}}(\mathsf{y})$,
\begin{align*}
	\|\uno_d-T_0F_y(y,z)\|&\leq \|T_0\|\|\dpr^2_{y}K(\mathsf{y})- K_{yy}\|_{\bar{r},\mathscr D_\sharp}+|\vae|\;\|T_0\|\;|\dpr^2_{y'}\wt K(\mathsf{y})|+|\vae|\;\|T_0\|\;\|\dpr_{y'}^2\wt K\|_{\bar{r},\mathscr D_\sharp}\\
		  &\leq d\cdot 2\mathsf T\|K_{yyy}\|_{\bar{r},\mathscr D_\sharp}\cdot\bar{r}+ 4|\vae|\mathsf{T}\frac{4M}{r^2}\\
	      &\leq 2d\mathsf{T} \mathsf{K}\frac{\bar{r}}{r-\bar{r}}+16\mathsf{T}\frac{|\vae| M}{ r^2}\\
	      &\leby{rrbarAsExtv5}2d\mathsf{T} \mathsf{K} \frac{2\bar{r}}{ r}+|\vae|\frac{16\mathsf{T} M}{r\bar{r}}\\
	      &\overset{\equ{rrbarAsExtv5}+\equ{cond1ExtExtv5}}{\leq}\su4+\frac{\s}{3}\\
	      &\le \su4+\su4=\su2\;;
\end{align*}
\item For any $z\in D_{\tilde{r}}(\mathsf{y})$,
\begin{align*}
2\|T_0\||F(\mathsf{y},z)|&\le 4\mathsf{T}|K_y(z)-K_y(\mathsf{y})|+ 4\mathsf{T}|\vae|\| \wt K_{y'}\|_{\bar{r},\mathscr D_\sharp}\\
&\leq 4\mathsf{T}\|K_{yy}\|_{\bar{r},\mathscr D_\sharp}\cdot \tilde{r}+ 4\mathsf{T} \frac{2|\vae| M}{r}\\
&\le 4\mathsf{T}\mathsf{K}\tilde{r}+\frac{8|\vae|\mathsf{T}_\infty M}{r}\\
&\le \frac{\bar{r}}{4d}+{|\vae|}\frac{\bar{r}}{2}\mathsf{L}\\
&\leby{cond1ExtExtv5}\frac{\bar{r}}{8}+ \frac{\bar{r}}{12}\\
&<\frac{\bar{r}}{4}\;,
\end{align*}
\ie 
$$
2\|T_0\|\|F(\mathsf{y},\cdot)\|_{\tilde{r},\mathsf{y}}\le \frac{\bar{r}}{4}\;.
$$
\end{itemize}
\noi
Therefore, Lemma~\ref{IFTLem} applies. Hence, there exists a real--analytic map $G^{\mathsf{y}}\colon D_{\tilde{r}}(\mathsf{y})\to D_{\bar{r}}(\mathsf{y})$ such that its graph coincides with $F^{-1}(\{0\})$ \ie  $y'=y'(z,\mathsf{y},\vae)\coloneqq G^{\mathsf{y}}(z)$ is the unique $y\in D_{\bar{r}}(\mathsf{y})$ satisfying $0=F(y,z)=\dpr_y K'(y)-K_y(z)$, for any $z\in D_{\tilde{r}}(\mathsf{y})$
. Moreover, $\forall\;z\in D_{\tilde{r}}(\mathsf{y})$,
\beq{EcarY1Y0Ext}
|G^{\mathsf{y}}(z)-\mathsf{y}|\leq 2\|T_0\|\|F(\mathsf{y},\cdot)\|_{\tilde{r},\mathsf{y}}
\leq \frac{\bar{r}}{4}\;,
\eeq
\beq{gy0z}
|G^{\mathsf{y}}(z)-z|\le |G^{\mathsf{y}}(z)-\mathsf{y}|+|\mathsf{y}-z|\le \frac{\bar{r}}{4}+\tilde{r}< \frac{\bar{r}}{2},
\eeq
so that
\beq{NextSetArnExt}
D_{{\bar{r}}/{4}}(G^{\mathsf{y}}(z))\subset D_{\bar{r}/2}(\mathsf{y}).
\eeq
\noi
Next, we prove that $\dpr^2_{y'}K'(y')$ is invertible, where $y'=G^{\mathsf{y}}(z)$ for some given $z\in D_{\tilde{r}}(\mathsf{y})$. Indeed, by Taylor's formula, we have,
\begin{align*}
\dpr^2_{y'} K'(y')&= K_{yy}(\mathsf{y})+ \dst\int_0^1 K_{yyy}(\mathsf{y}+t(y'-\mathsf{y}))(y'-\mathsf{y}) dt+\vae \dpr^2_{y'}\wt K(y')\\
           &= T(\mathsf{y})^{-1}\left(\uno_d+T(\mathsf{y})\left(\dst\int_0^1 K_{yyy}(\mathsf{y}+t(y'-\mathsf{y}))(y'-\mathsf{y}) dt+ \dpr^2_{y'}\wt K(y')\right)\right)\\
           &\eqqcolon T(\mathsf{y})^{-1}(\uno_d+\vae A),
\end{align*}
and, by Cauchy's estimate,\footnote{Recall footnote \textsuperscript{\ref{ftnTK1}} 
.} 
\begin{align*}
|\vae|\|A\|&\leq \|T(y_0)\|\left(d\|K_{yyy}\|_{\bar{r},\mathscr D_\sharp}|y'-\mathsf{y}|+ |\vae|\|\dpr_{y'}^2\wt K\|_{\bar{r},\mathscr D_\sharp}\right)\\
     &\leq \|T\|_{\mathscr D_\sharp}\left(\frac{d\|K_{yy}\|_{r,\mathscr D_\sharp}}{r-\bar{r}}|y'-\mathsf{y}|+|\vae|\|\dpr^2_{y'}\wt K\|_{\bar{r},\mathscr D_\sharp}\right)\\
     &\leq \mathsf{T}\left(\frac{2d\mathsf{K}}{r}\frac{\bar{r}}{2}+\frac{4M}{r^2} \right)\\
	 &\leby{rrbarAsExtv5} \mathsf{T}\left(\frac{\s}{16\mathsf{T}}+\su{4\mathsf{T}}|\vae|\mathsf{L} \right)\\
	 &\leby{cond1ExtExtv5} \mathsf{T}\left(\frac{\s}{16\mathsf{T}}+\su{4\mathsf{T}}\frac{\s}{3} \right)\\
	 &=\frac{\s}{6}\\
	 &<\su2.
\end{align*}
Hence $\dpr_{y'}^2K'(y')$ is invertible with
\[\dpr_{y'}^2K'(y')^{-1}=(\uno_d+\vae A)^{-1}T(\mathsf{y})=T(\mathsf{y})+\dst\sum_{k\geq 1}(-\vae)^k A^k T(\mathsf{y})\eqqcolon T(\mathsf{y})+\vae \wt T(y'),\]
and
\[|\vae|\|\wt T({y}')\|\leq |\vae|\frac{\|A\|}{1-|\vae|\|A\|}\|T\|_{\mathscr D_\sharp}\leq 2|\vae|\|A\| \|T\|_{\mathscr D_\sharp}
\le 2\frac{\s}{6}\mathsf{T}
= \mathsf{T}\frac{\s}{3}\,.\]
\noi
Similarly, from
$$
K_{yy}(z)=K_{yy}(\mathsf{y})\left(\uno_d+T(\mathsf{y})\dst\int_0^1 K_{yyy}(\mathsf{y}+t(z-\mathsf{y}))(z-\mathsf{y}) dt\right)
$$
and
\begin{align*}
\left\|T(\mathsf{y})\dst\int_0^1 K_{yyy}(\mathsf{y}+t(z-\mathsf{y}))(z-\mathsf{y}) dt\right\|_{r/(4d\mathsf{T}\mathsf{K}),\mathsf{y}}&\le \mathsf{T}\|K_{yyy}\|_{r/2,\mathsf{y}}\frac{r}{4d\mathsf{T}\mathsf{K}}\le \mathsf{T}\frac{d\mathsf{K}}{r-r/2}\frac{r}{4d\mathsf{T}\mathsf{K}}=\su2
\end{align*}
 one has that, for any $z\in D_{r/(4d\mathsf{T}\mathsf{K})}(\mathsf{y})$, 
\beq{invJacKyy}
K_{yy}(z)^{-1} \mbox{ exists and }\|K_{yy}(z)^{-1}\|\le \|K_{yy}(z)^{-1}-T(\mathsf{y})\|+\|T(\mathsf{y})\|\le 
2\su2 \mathsf{T}+\mathsf{T}=2\mathsf{T}\;.
\eeq
Now, differentiating $F(G^{\mathsf{y}}(z),z)=0$, we get, for any $z\in D_{\tilde{r}}(\mathsf{y})$,
$$
\dpr_{y'}^2K'(G^{\mathsf{y}}(z))\cdot\dpr_z G^{\mathsf{y}}(z)=K_{yy}(z)\;.
$$
Therefore $G^{\mathsf{y}}$ is a local diffeomorphism, with
\begin{align*}
\dpr_z G^{\mathsf{y}}(z)&= \dpr_{y'}^2K'(G^{\mathsf{y}}(z))^{-1}K_{yy}(z)\\
                &=\left(K_{yy}(z)^{-1}\left(K_{yy}(z)+\vae \dpr_{y'}^2\wt K(G^{\mathsf{y}}(z)) \right)\right)^{-1}\\
                &=\left(\uno_d+\vae K_{yy}(z)^{-1}\dpr_{y'}^2\wt K(G^{\mathsf{y}}(z))\right)^{-1}
\end{align*}
and
$$
\|\vae K_{yy}^{-1}\dpr_{y'}^2\wt K\|_{\tilde{r},\mathsf{y}}\le \| K_{yy}^{-1}\|_{\tilde{r},\mathsf{y}}\|\vae\dpr_{y'}^2\wt K\|_{\tilde{r},\mathscr D_\sharp}\le 2\mathsf{T}\frac{|\vae|\mathsf{L}}{4\mathsf{T}}\s^{\n+d}\le\su2\s^{\n+d}|\vae|\mathsf{L}<\frac{\s}{6}<\su2
$$
so that
\beq{Jacgy0}
\|\dpr_z G^{\mathsf{y}}-\uno_d\|_{\tilde{r},\mathsf{y}}\le 2 \|\vae K_{yy}^{-1}\dpr_{y'}^2\wt K\|_{\tilde{r},\mathsf{y}}\le \s^{\n+d}|\vae|\mathsf{L}.
\eeq
\noi
Now, we show that the family $\{G^{\mathsf{y}}\}_{\mathsf{y}\in\mathscr D_\sharp}$ is compatible so that, together, they define a global map on $D_{\tilde{r}}(\mathscr D_\sharp)$, say $G$ and that, in fact, $G$ is a real--analytic 
 diffeomorphism. For, assume that $z\in D_{\tilde{r}}(\mathsf{y})\bigcap D_{\tilde{r}}(\hat{\mathsf{y}})$, for some $\mathsf{y},\hat{\mathsf{y}}\in \mathscr D_\sharp$. Then, we need to show that $G^{\mathsf{y}}(z)=G^{\hat{\mathsf{y}}}(z)$. But,  we have
$$
|G^{\hat{\mathsf{y}}}(z)-\mathsf{y}|\le |G^{\hat{\mathsf{y}}}(z)-\hat {\mathsf{y}}|+|\hat{\mathsf{y}}-z|+|z-\mathsf{y}|\leby{EcarY1Y0Ext}\frac{\bar{r}}{2}+\tilde{r}+\tilde{r}<\bar{r}.
$$
Hence, $z\in D_{\tilde{r}}(\mathsf{y}),\ G^{\hat{\mathsf{y}}}(z)\in D_{\bar{r}}(\mathsf{y})$ and, by definitions, $F(G^{\hat{\mathsf{y}}}(z),z)=0=F(G^{\mathsf{y}}(z),z)$. Then, by unicity, we get $G^{\mathsf{y}}(z)=G^{\hat{\mathsf{y}}}(z)$. Thus, the map
$$
G\colon D_{\tilde{r}}(\mathscr D_\sharp)\to\cn\quad \mbox{such that}\quad G_{|D_{\tilde{r}}(\mathsf{y})}\coloneqq G^{\mathsf{y}},\quad \forall\; \mathsf{y}\in \mathscr D_\sharp\;,
$$
is well--defined and, therefore, is a real--analytic local diffeomorphism. It remains only to check that $G$ is injective 
to conclude that it 
is a global diffeomorphism. Let then $z\in D_{\tilde{r}}(\mathsf{y}),\hat{z}\in  D_{\tilde{r}}(\hat{\mathsf{y}})$ such that $G(z)=G(\hat{z})$, for some $\mathsf{y},\hat{\mathsf{y}}\in\mathscr D_\sharp$. Then, we have 
$$
|z-\hat{z}|< \frac{r}{4d\mathsf{T}\mathsf{K}}-\tilde{r}.
$$
Indeed, if not then
\begin{align*}
0=|G(z)-G(\hat{z})|&\ge -|G(z)-z|+|z-\hat{z}|-|\hat{z}- G(\hat{z})|\\
                   &\geby{gy0z}-\bar{r}+\frac{r}{4d\mathsf{T}\mathsf{K}}-\tilde{r}-\bar{r}\\
                   &\ge \frac{r}{4d\mathsf{T}\mathsf{K}}-3\bar{r}\\
                   &\geby{rrbarAsExtv5}\frac{r}{4d\mathsf{T}\mathsf{K}}-3\frac{r}{16d\mathsf{T}\mathsf{K}}\\
                   &>0\,,
\end{align*}
contradiction. Therefore,
\beq{zminhatz}
|z-\hat{z}|< \frac{r}{4d\mathsf{T}\mathsf{K}}-\tilde{r}\;.
\eeq
Thus,
$$
|\hat{z}-\mathsf{y}|\le |\hat{z}-z|+|z-\mathsf{y}|< \frac{r}{4d\mathsf{T}\mathsf{K}}-\tilde{r}+\tilde{r}=\frac{r}{4d\mathsf{T}\mathsf{K}}\,.
$$
Hence, $z,\hat{z}\in D_{r/(4d\mathsf{T}\mathsf{K})}(\mathsf{y})$. But $G(z)=G(\hat{z})$ is equivalent to $K_y(z)=K_y(\hat{z})$ and then,
$$
0=K_y(z)-K_y(\hat{z})=\dst\int_0^1 K_{yy}(\hat{z}+t(z-\hat{z}))dt(z-\hat{z})\;.
$$
Thus, it is enough to show that $\dst\int_0^1 K_{yy}(\hat{z}+t(z-\hat{z}))dt$ is invertible. But
\begin{align*}
\dst\int_0^1 K_{yy}(\hat{z}+t(z-\hat{z}))dt&= K_{yy}(\hat{z})+\dst\int_0^1\int_0^1 K_{yyy}(\hat{z}+tt'(z-\hat{z}))tdt'dt\cdot(z-\hat{z})\\
      &\eqby{invJacKyy} K_{yy}(\hat{z})\left(\uno_d+K_{yy}(\hat{z})^{-1}\dst\int_0^1\int_0^1 K_{yyy}(\hat{z}+tt'(z-\hat{z}))tdt'dt\cdot(z-\hat{z})\right)\\
\end{align*}
and
\begin{align*}
\left\|K_{yy}(\hat{z})^{-1}\dst\int_0^1\int_0^1 K_{yyy}(\hat{z}+tt'(z-\hat{z}))tdt'dt\cdot(z-\hat{z})\right\|&\leby{invJacKyy} 2\mathsf{T}\cdot\su2\|K_{yyy}\|_{r/2,\mathsf{y}}|z-\hat{z}|\\
  &\leby{zminhatz} \mathsf{T}\frac{2d\mathsf{K}}{r}\left(\frac{r}{4d\mathsf{T}\mathsf{K}}-\tilde{r}\right)\\
  &<\frac{2d\mathsf{T}\mathsf{K}}{r}\frac{r}{4d\mathsf{T}\mathsf{K}}\\
  &=\su2.
\end{align*}
Therefore, $\dst\int_0^1 K_{yy}(\hat{z}+t(z-\hat{z}))dt$ is invertible and then we get $z-\hat{z}=0$ \ie $G$ is injective.\\
\noi
Next, we estimate $\|G-\id\|_{\tilde{r},\mathscr D_\sharp}$. The strategy is to show that the expression $(K_y+\vae\wt K_{y'})^{-1}\circ K_y$ defines a map on $D_{\tilde{r}}(\mathsf{y})$ by means of the Inversion Function Lemma~\ref{inv1}; hence, we will get an explicit formula for $G$:
\beq{explForG}
G=(K_y+\vae\wt K_{y'})^{-1}\circ K_y\qquad \mbox{on}\quad D_{\tilde{r}}(\mathsf{y})\;.
\eeq
But, the proof is part of the above computation: for any $y\in D_{\bar{r}}(\mathsf{y})$,
$$
\|\uno_d-T_0 F_{y}(y,\mathsf{y})\|\le \su 2 
$$
implies, using Lemma~\ref{inv1}, that $K_y+\vae\wt K_{y'}$ admits an inverse defined on $D_{r_\sharp}(K_y(\mathsf{y})+\vae\wt K_{y'}(\mathsf{y}))$, where
$$
r_\sharp\coloneqq \frac{\bar{r}}{4\mathsf{T}}< \frac{\bar{r}}{2\|T_0\|}\;.
$$
Moreover, for any $y\in D_{\tilde{r}}(\mathsf{y})$,
\begin{align*}
|K_y(y)-(K_y(\mathsf{y})+\vae\wt K_{y'}(\mathsf{y}))|&\le \|K_{yy}\|_{\tilde{r},\mathscr D_\sharp}\cdot \tilde{r}+\|\vae\wt K_{y'}\|_{\bar{r},\mathscr D_\sharp}\\
 									   &\le \mathsf{K}\frac{\bar{r}}{16d\mathsf{T}\mathsf{K}}+ \frac{2|\vae|M}{r}\\
 									   &\le\frac{\bar{r}}{16d\mathsf{T}}+\frac{\bar{r}}{8\mathsf{T}}|\vae|\mathsf{L}\\
 									   &<\frac{\bar{r}}{4\mathsf{T}}=r_\sharp\;,
\end{align*}
and thus, \equ{explForG} is proven. Hence, for any $y\in D_{\tilde{r}}(\mathsf{y})$,
\begin{align*}
|G(y)-y|&=|(\dpr_{y'}K')^{-1}( K_y(y))-(\dpr_{y'}K')^{-1}(K_y(y)+\vae\wt K_{y'}(y))|\\
	    &\le \dst\int_0^1\|\dpr_{y'}((\dpr_{y'}K')^{-1})(K_y(y)+t\vae\wt K_{y'}(y))\|dt\;\|\vae\wt K_{y'}\|_{\bar{r},\mathscr D_\sharp}\\
	    &\le \|(\dpr_{y'}^2K')^{-1}\|_{\bar{r},\mathscr D_\sharp}\frac{2|\vae|M}{r}\\
	    &<\frac{16\mathsf{T} M}{r}|\vae|\\
	    &<\s^{\n+d}\bar{r}|\vae|\mathsf{L}\;.
\end{align*}
\noi
Now, we estimate $P_+$. We have,
\[|\vae|\|g_x\|_{\bar{r},\bar{s},\mathscr D_\sharp}\leq |\vae|\mathsf{C}_1 \frac{M}{\a} \s^{-(\n+d)}
\le |\vae| \frac{r}{3}\mathsf{L}\leby{cond1ExtExtv5}\frac{r}{3}\frac{\s}{3}\le \frac{r}{12}\]
so that, for any $\mathsf{y}\in \mathscr D_\sharp$ and $(y',x)\in D_{\bar{r},\bar{s}}(\mathsf{y})$,
\[ |y'+\vae g_x(y',x)-\mathsf{y}|\leq \bar{r}+\frac{r}{3}\le \frac{r}{32d}+\frac{r}{12}<\frac{2r}{3}<r\,,\]
and thus
\begin{align*}
\|P^\ppu\|_{\bar{r},\bar{s},\mathscr D_\sharp}&\leq d^2 \|K_{yy}\|_{r,\mathscr D_\sharp}\|g_x\|_{\bar{r},\bar{s},\mathscr D_\sharp}^2\leq d^2 \mathsf{K}\left( \mathsf{C}_1 \frac{M}{\a} \s^{-(\n+d)}\right)^2\\
   &=d^2\mathsf{C}_1^2 \frac{\mathsf{K}M^2}{\a^2} \s^{-2(\n+d)}, 
\end{align*}
\begin{align*}
\|P^\ppd\|_{\bar{r},\bar{s},\mathscr D_\sharp}&\leq d\|P_y\|_{\frac{5r}{6},\bar{s},\mathscr D_\sharp}\|g_x\|_{\bar{r},\bar{s},\mathscr D_\sharp}\leq d\frac{6M}{r}\mathsf{C}_1 \frac{M}{\a} \s^{-(\n+d)}\\
     &= 6d\mathsf{C}_1 \frac{M^2}{\a r}\s^{-(\n+d)}
\end{align*}
 and, by Lemma~$\ref{fce}$--$(i)$, we have
\begin{align*}
|\vae|\|P^\ppt\|_{\bar{r},s-\frac{\s}{2},\mathscr D_\sharp}&\leq \dst\sum_{|n|_1>\k}\|P_n\|_{\bar{r},\mathscr D_\sharp}\ex^{(s-\frac{\s}{2})|n|_1}\leq M\dst\sum_{|n|_1>\k}\ex^{-\frac{\s |n|_1}{2}}\\
  &\leq M\ex^{-\frac{ \k\s}{4}}\dst\sum_{|n|_1>\k}\ex^{-\frac{\s |n|_1}{4}}\leq M\ex^{-\frac{ \k\s}{4}}\dst\sum_{|n|_1>0}\ex^{-\frac{\s |n|_1}{4}}\\
  &= M\ex^{-\frac{ \k\s}{4}} \left(\left(\dst\sum_{k\in \integer}\ex^{-\frac{\s |k|}{4}}\right)^d-1\right)=M\ex^{-\frac{ \k\s}{4}}\left(\left(1+\frac{2\ex^{-\frac{\s }{4}}}{1-\ex^{-\frac{\s }{4}}} \right)^d-1\right)\\
  &= M\ex^{-\frac{ \k\s}{4}}\left(\left(1+\frac{2}{\ex^{\frac{\s }{4}}-1} \right)^d-1\right)\leq M\ex^{-\frac{ \k\s}{4}}\left(\left(1+\frac{2}{\frac{\s }{4}} \right)^d-1\right)\\
  &\leq \s^{-d} M\ex^{-\frac{ \k\s}{4}}\left(\left(\s +8 \right)^d-\s^d\right)\leq d 8^{d}\s^{-d} M\ex^{-\frac{ \k\s}{4}}\\
  &= \mathsf{C}_2\s^{-d} M \ex^{-\l}\\
  &\leby{DefNArnExt1v5} \mathsf{C}_2 \s^{-d}M \s^{-(2\n+d)}\frac{\mathsf{K}|\vae|M}{\a^2}\\
  &=\mathsf{C}_2 M\frac{\mathsf{K}|\vae|M}{\a^2}\s^{-2(\n+d)}\,.
\end{align*}
Hence,\footnote{Recall that $r\le r_0$ and $\s<1$.}
\begin{align*}
\|P_+\|_{\bar{r},\bar{s},\mathscr D_\sharp}&\leq \|P^\ppu\|_{\bar{r},\bar{s},\mathscr D_\sharp}+\|P^\ppd\|_{\bar{r},\bar{s},\mathscr D_\sharp}+\|P^\ppt\|_{\bar{r},\bar{s},\mathscr D_\sharp}\\
  &\leq d^2\mathsf{C}_1^2 \frac{\mathsf{K}M^2}{\a^2} \s^{-2(\n+d)}+6d\mathsf{C}_1 \frac{M^2}{\a r}\s^{-(\n+d)}+\mathsf{C}_2 M\frac{|\vae|{M}\mathsf{K}}{{\a}^2}\s^{-2(\n+d)}\\
  &= \left(d^2\mathsf{C}_1^2 r\mathsf{K}+6d\mathsf{C}_1 \a \s^{\n+d}+\mathsf{C}_2 r\mathsf{K}\right)\frac{M^2}{\a^2 r}\s^{-2(\n+d)}\\
  &\le \left(d^2\mathsf{C}_1^2+6d\mathsf{C}_1 +\mathsf{C}_2\right)\max\left\{\a,r\mathsf{K}\right\}\frac{M^2}{\a^2 r}\s^{-2(\n+d)}\\
  &\leby{RecHypArnExtv5} \frac{\mathsf{C}_3}{\sqrt{2}} \max\left\{1,\frac{\a}{r\mathsf{K}} \right\}\frac{M^2\mathsf{K}}{\a^2 }\s^{-2(\n+d)}\\
  &\le \mathsf{L}M\;.
\end{align*} 
The proof of the claims on $\phi'$ and $P'$ are proven in a similar way as in Lemma~\ref{lem:1}.
\qed
\newpage
\noi
Let   $H_0\coloneqq H$, $K_0\coloneqq K$, $P_0\coloneqq P$, $\phi^0=\phi_0\coloneqq\id $ and  
$r_0$, $s_0$, $s_*$, $\s_0$, $M_0$, $\mathsf{K}_0$, $\mathsf{T}_0$ and $\eta_0$ be  as in $\S\ref{AssumpExtArnolv2}$. For a given $\vae\neq0$ and $j\ge 0$,  define\footnote{Notice that $s_{j}\downarrow s_*$ and $r_{j}\downarrow 0$.}
\begin{align*}
 \dst\s_j&\coloneqq \frac{\s_0}{2^j}\,,\\
  s_{j+1}&\coloneqq s_j-\s_j=s_*+\frac{\s_0}{2^j}\,,\\
  \bar s_{j}&\coloneqq s_j-\frac{2\s_j}{3}\,,\\
 \mathsf{K}_{j+1}&\coloneqq \mathsf K_0\dst\prod_{k=0}^{j}(1+\frac{\s_k}{3})\le \mathsf K_0\ex^{\frac{2\s_0}{3}}<\mathsf{K}_0\sqrt{2}\,,  \\
\mathsf{T}_{j+1}&\coloneqq \mathsf T_0\dst\prod_{k=0}^{j}(1+\frac{\s_k}{3})\le \mathsf T_0\ex^{\frac{2\s_0}{3}}<\mathsf{T}_0\sqrt{2}\,,\\
 \l_0&\coloneqq \log\m_0^{-1}\,,\\
 \mathsf{e}_*&\coloneqq \mathsf{C}_5 \s_0^{-(4\n+2d+2)}\eta_0^2\l_0^{2\n}\;,\\
\mathsf{d}_*&\coloneqq \frac{2^{2\n+2d+13}d^2\eta_0^2}{\s_0^2}\,,\\
 \k_0 &\coloneqq 4\s_0^{-1}\l_0\,,\\
 \k_j &\coloneqq 4^j\k_0\,,\\
 \hat{r}_0 &\coloneqq \frac{r_0}{\sqrt{|\vae|}}\,,\\
 \hat{r}_{j+1}&\coloneqq \su2\min\left\{\frac{\wh{\a}}{2d\sqrt{2}\mathsf{K}_0\k_{j}^{\n}}\,,\, \frac{\hat{r}_{j}\s_{j}}{32d\eta_0} \right\}\,,\\
 r_{j+1}&\coloneqq \hat{r}_{j+1}\sqrt{|\vae|} \,,\\
 \tilde r_{j+1}&\coloneqq \frac{{r}_{j+1}}{8d\mathsf{T}_j\mathsf{K}_j} \,,\\
 \wh M_0&\coloneqq M_0\,,\\
 \wh M_{1}&\coloneqq 	{2\mathsf{e}_*} \frac{\mathsf{K}_0 {\wh M_{0}}^2}{\wh\a^2}\,,
 \end{align*}
\begin{align*}
 \wh M_{j+2}&\coloneqq 	
 			8\mathsf{e}_*(4\mathsf d_*)^{j} \frac{\mathsf{K}_0 {\wh M_{j+1}}^2}{\wh\a^2}\,,\\
 \m_j&\coloneqq \frac{\mathsf{K}_0\wh M_j}{\wh \a^2}\,,\\
 \th_j     &\coloneqq 8\mathsf{e}_*\;(4\mathsf d_*)^j\;\m_j\,,\\
 \mathscr D_0&\coloneqq \mathscr D_{\d,\a}\,,\\
 \mathsf{W}_0   &\coloneqq \diag\left(\frac{\mathsf{K}_0}{{\a}}\uno_d,\uno_d\right)\,,\\
 \mathsf{W}_{j+1}&\coloneqq \diag\left(\max\left\{\frac{\mathsf{K_j}}{{\a}}\;,\frac{1}{r_j}\right\}\uno_d\,,\uno_d\right)\,,\\ 
\mathsf{L}_j&\coloneqq M_j\dst\max\left\{\frac{32\sqrt{2}\mathsf{T}_0  }{r_jr_{j+1}}\s_j^{-(\n+d)}\,,\,
{\mathsf{C}_4} \max\left\{1,\frac{\a}{r_j\mathsf{K}_j}\right\}\frac{\mathsf{K}_0 }{\a^2 }\s_j^{-2(\n+d)}\right\}\\
          &=M_j\dst\max\left\{\frac{32\sqrt{2}\mathsf{T}_0  }{r_jr_{j+1}}\s_j^{-(\n+d)}
\,,\,\frac{4}{\mathsf{K}_j r_j^2}\,,\,\mathsf{C}_4 \max\left\{1,\frac{\a}{r_j\mathsf{K}_j}\right\}\frac{\mathsf{K}_0 }{\a^2 }\s_j^{-2(\n+d)}\right\}\,.
\end{align*}
Thus,
$$
\th_{1}= 8\mathsf{e}_*\;(4\mathsf d_*)\;\m_{1}=32\mathsf{e}_*\;\mathsf d_*\frac{\mathsf{K}_0\wh M_{1}}{\wh \a^2}=32\mathsf{e}_*\;\mathsf d_*\frac{\mathsf{K}_0}{\wh \a^2}\;2\mathsf{e}_* \frac{\mathsf{K}_0 {\wh M_0}^2}{\wh\a^2}= \mathsf d_*\left(8\mathsf{e}_*\;\m_0\right)^2=\mathsf d_*\;\th_0^2
$$
and, for any $j\ge1$,
\begin{align*}
\th_{j+1}&= 8\mathsf{e}_*\;(4\mathsf d_*)^{j+1}\;\m_{j+1}=8\mathsf{e}_*(4\mathsf d_*)^{j+1}\frac{\mathsf{K}_0\wh M_{j+1}}{\wh \a^2}\\
&=8\mathsf{e}_*(4\mathsf d_*)^{j+1}\frac{\mathsf{K}_0}{\wh \a^2}\;8\mathsf{e}_*(4\mathsf d_*)^{j-1} \frac{\mathsf{K}_0 {\wh M_j}^2}{\wh\a^2}= \left(8\mathsf{e}_*(4\mathsf d_*)^{j}\;\m_j\right)^2=\th_j^2
\end{align*}
\ie
$$
\th_j=\th_1^{2^{j-1}}=(\sqrt{\mathsf{d}_*}\;\th_0)^{2^j} \;.
$$
The very first step being quite different from all the others, it has to be done separately. Hence,\\
\lem{frstStepv5}
Under the above assumptions and notations, if
\beq{condBisv2v5}
|\vae|\le \left(\frac{r_0\s_0}{\wh{\a}\mathsf{T}_0}\right)^2\qquad\mbox{and}\qquad \max\left\{\ex\;\m_0\;,\, 16d\;\eta_0\;\th_0\right\}\le 1\;,
\eeq
then, there exist $\mathscr D_1\subset \mathscr D$
, a real--analytic diffeomorphism

$$
G_{1}\colon  D_{\tilde{r}_{1}}(\mathscr D_{\d,\a}){\to}G_{1}( D_{\tilde{r}_{1}}(\mathscr D_{\d,\a}))
$$
and a real--analytic symplectomorphism 
\beq{phijExtv5}
\phi_{1}:D_{r_{1},s_{1}}(\mathscr D_1)\to D_{r_{0},s_{0}}(\mathscr D_{0})
\eeq
such that
\begin{align}
&G_1(\mathscr D_{\d,\a})=\mathscr D_1\;,\label{k1k0puv5}\\
&\dpr_{y_1}K_{1}\circ G_{1}=\dpr_{y}K_0 \;,\label{kjkjpuv5}\\
&H_{1}\coloneqq H_{0}\circ\phi_{1}\eqqcolon K_{1} + \vae^{2} P_{1}\qquad\quad\ \quad\; \mbox{on } D_{r_{1},s_{1}}(\mathscr D_{1})\label{HjExtExtv5}
\end{align}
and
\begin{align}
&\mathscr D_1\subset \mathscr D_{r_1}\;,\label{estfin2Bis0001v2}\\
& \|K_1\|_{r_1,\mathscr D_1}\le \mathsf{K}_1\;,\qquad \|T_1\|_{\mathscr D_1}\le  \mathsf{T}_1\;,\qquad T_1\coloneqq (\dpr_{y_1}^2 K_1)^{-1}\;,\label{estfin2Bis0000v2}\\
& 2\vae^{2}M_1\coloneqq 2\vae^{2}\|P_1\|_{r_1,s_1,\mathscr D_1}\le |\vae|\wh M_1\;,\label{estfin2Bis000v2}\\
&
\|G_1-\id\|_{\tilde r_1,\mathscr D_{\d,\a}}\le 2\;r_1\;\s_0^{\n+d}\;|\vae|\mathsf{L}_0\;,\label{estG1idv2}\\
&\|\dpr_z G_{1}-\uno_d\|_{\tilde{r}_{1},\mathscr D_{\d,\a}}\le \s_0^{\n+d}|\vae|\mathsf L_0 \;,\label{estG1devidv2}\\
&
\|\mathsf{W}_1(\phi_1-\id)\|_{r_1,s_1,\mathscr D_1}
\le \s_0^d\;|\vae|{\mathsf L}_0\;.\label{estfin2Bis010v2}
\end{align}
\elem
\proof
By
\beq{kp08v5}
\k_0\geby{condBisv2v5}4\s_0^{-1}\ge 8
\eeq
and
$$
\frac{\wh{\a}}{2d\mathsf{K}_0\sqrt{2}\k_0^{\n}}\overset{\equ{kp08v5}+\equ{condBisv2v5}}{\le} \frac{1}{2d\cdot 8^{\n}\mathsf{K}_0\sqrt{2}}\frac{r_0\s_0}{\mathsf{T}_0\sqrt{|\vae|}}<\frac{\hat{r}_0\s_0}{32d\eta_0}\,,
$$
we get 
\beq{r1calv5}
\hat r_1=\su2\min\left\{\frac{\wh{\a}}{2d\sqrt{2}\mathsf{K}_0\k_0^{\n}}\,,\, \frac{\hat {r}_0\s_0}{32d\eta_0} \right\}=\frac{\wh{\a}}{4d\sqrt{2}\mathsf{K}_0\k_0^{\n}}
\eeq
and, thus
\begin{align}
|\vae| \mathsf L_0 (3 \sigma_0^{-1})&\le 3|\vae|M_0\dst\max\left\{\frac{32\sqrt{2}\mathsf{T}_0  }{r_0 {r}_{1}}\s_0^{-(\n+d)}\,,\,\mathsf{C}_4 \max\left\{1,\frac{\a}{r_0\mathsf{K}_0}\right\}\frac{ \mathsf{K}_0}{\a^2}\s_0^{-2(\n+d)}\right\}\s_0^{-1}\nonumber\\
                                    &\le 3\dst\max\left\{32\sqrt{2}\mathsf{T}_0\frac{ \a }{{r}_{1}}\frac{ \a }{r_0 \mathsf{K}_0}\,,\,\mathsf{C}_4 \max\left\{1,\frac{\a}{r_0\mathsf{K}_0}\right\}\right\}\s_0^{-2(\n+d)-1}\frac{ \mathsf{K}_0 M_0}{\wh\a^2}\nonumber\\
                                    &\leby{condBisv2v5}3\dst\max\left\{32\sqrt{2}\mathsf{T}_0\frac{ \a }{{r}_{1}}\,,\,\mathsf{C}_4\right\}\s_0^{-2(\n+d)-1}\m_0\nonumber\\
                                    &=3\dst\max\left\{256d\eta_0\k_0^{\n}\,,\,\mathsf{C}_4\right\}\s_0^{-2(\n+d)-1}\m_0\nonumber\\
                                    &\leby{kp08v5} 3\dst\max\left\{256d\,,\,8^{-\n}\mathsf{C}_4\right\}\eta_0\k_0^{\n}\s_0^{-2(\n+d)-1}\m_0\nonumber\\
                                    &\le\mathsf{e}_*\;\m_0\nonumber\\ 
                                    &=\th_0\leby{condBisv2v5} 1.\label{L0ifsg3v5}
\end{align}
Therefore, Lemma~\ref{frstStepv5} is a straightforward consequence of Lemma~\ref{lem:1Extv5}.
\qed
\lem{lem:2Extv5} 
Assume $\equ{HjExtExtv5}\div\equ{estfin2Bis000v2}$ with some $\vae\neq 0$ and
\beq{condBisv2Prtv5} 
 \max\left\{\ex\;\m_0\;,\, 2\mathsf{C}_{6}\;\eta_0^{\frac{5}{4}}\;\s_0^{-\frac{5}{4}}\; \th_{0}\right\}\le 1\;.
\eeq
Then, one can construct a sequence of real--analytic diffeomorphisms 
$$
G_{j}\colon  D_{\tilde{r}_{j}}(\mathscr D_{j-1}){\to}G_{j}( D_{\tilde{r}_{j}}(\mathscr D_{j-1}))\;,\qquad j\ge 2
$$
and of real--analytic symplectic transformations 
\beq{phijBisv2v5}
\phi_{j}:D_{r_{j},s_{j}}(\mathscr D_{j})\to D_{r_{j-1},s_{j-1}}(\mathscr D_{j-1})\;,
\eeq
such that
\begin{align*}
&G_j(\mathscr D_{j-1})=\mathscr D_j\subset \mathscr D_{r_j} \;,\\
&\dpr_{y}K_{j+1}\circ G_{j+1}=\dpr_{y}K_j \;,\\
&H_{j}\coloneqq H_{j-1}\circ\phi_{j}\eqqcolon K_{j} + \vae^{2^{j}} P_{j}\qquad\quad\  \mbox{on}\quad D_{r_{j},s_{j}}(\mathscr D_{j})\;,
\end{align*}
converge uniformly. More precisely, we have the following:
\begin{itemize}
\item[$(i)$] the sequence $G^{j}\coloneqq G_{j}\circ G_{j-1}\circ\cdots\circ G_2\circ G_1$ converges uniformly on $\mathscr D_{\d,\a}$ to a lipeomorphism $G^*\colon \mathscr D_{\d,\a}\to \mathscr D_*\coloneqq G^*(\mathscr D_{\d,\a})\subset\mathscr D$ and $G^*\in C^\infty_W(\mathscr D_{\d,\a})$\;.  
\item[$(ii)$] $\vae^{2^j}\dpr_y^{\b} P_j$ converges uniformly on 
$\mathscr D_*\times\dst\torus^d_{s_*}$ to $0$, for any $\b\in \natural_0^d$\;;
\item[$(iii)$] $\phi^j\coloneqq \phi_2\circ \cdots\circ \phi_j$ converges uniformly on 
$\mathscr D_*\times\tn$ to a symplectic transformation 
$$
\phi^*\colon \mathscr D_*\times\tn\overset{into}{\longrightarrow} B_{r_1}(\mathscr D_1)\times\tn,
$$
 with $\phi^*\in C^\infty_W(\mathscr D_*\times\tn)$ and $\phi^*(y,\cdot)\colon \torus^d_{s_*}\ni x\mapsto \phi^*(y,x)$ holomorphic, for any $y\in\mathscr D_*$\;;
\item[$(iv)$]  $K_j$ converges uniformly on 
$\mathscr D_*$ to a function $K_*\in C^\infty_W(\mathscr D_*)$, with
\begin{align*}
&\dpr_{y_*}K_*\circ G^*=\dpr_{y}K_0 \quad \qquad\qquad\quad\quad\mbox{on} \quad \mathscr D_{\d,\a}\;,\\
 &\dpr_{y_*}^{\b}(H_1\circ\phi^*)(y_*,x)=\dpr_{y_*}^{\b} K_*(y_*)\;,\quad \forall(y_*,x)\in\mathscr D_*\times\tn\;, \forall\;\b\in \natural_0^d\;.
\end{align*}
\end{itemize}
Finally, the following estimates hold for any $i\ge 2$:
\begin{align}
&\|G_i-\id\|_{\tilde r_{i},\mathscr D_{i-1}}\le 2\;r_{i}\;\s_{i-1}^{\n+d}\;|\vae|^{2^{i-1}}\mathsf{L}_{i-1}\;,\label{estGiidv2}\\
&\|\dpr_z G_{i}-\uno_d\|_{\tilde{r}_{i},\mathscr D_{i-1}}\le \s_{i-1}^{\n+d}\;|\vae|^{2^{i-1}}\mathsf L_{i-1} \;,\label{estGidevidv2}\\
&2^{i^2}|\vae|^{2^i}M_i\coloneqq 2^{i^2}|\vae|^{2^i}\|P_i\|_{r_i,s_i,\mathscr D_i}\le |\vae|\wh M_i\ ,\label{estfin2Ext01v5}\\
&|\meas(\mathscr D_*)-\meas(\mathscr D_{\d,\a})|\le \mathsf{C}_7\;\s_{0}^{\n+d}\;\th_0\;\meas(\mathscr D_{\d,\a})\ ,\label{estfin2Ext02v5}\\
&|\mathsf{W}_2(\phi^*-\id)|
\le \frac{\th_0}{3\cdot {2^{d}}}
\qquad\qquad\ \mbox{on}\quad \mathscr D_*\times\torus^d_{s_*}\label{estfin2Ext03v5}\ .
\end{align}
\elem 
\proof
 First of all, notice that
 \beq{riexp}
 \hat r_{i+1}={\hat{r}_1}\frac{\s_{1}\cdots\s_i}{(64d\eta_0)^i}=2^{-\frac{i^2}{2}}\left(\frac{\s_0}{64d\sqrt{2}\eta_0}\right)^i\hat{r}_1
 \;, \qquad \forall\; i\ge 0\;.
 \eeq
 Indeed, for any $j\ge 1$, we have
 \beq{ritric}
 \frac{\s_1\cdots\s_{j}}{(64d\eta_0)^j}\le \left(\frac{\s_{1}}{64d\eta_0}\right)^j\le\left(\frac{\s_{0}}{2^7d}\right)^j\le \su{4^{j\n}}\;,
 \eeq
 so that
 $$
 \hat r_{2}= \min\left\{\frac{\wh{\a}}{4d\sqrt{2}\mathsf{K}_0\k_1^{\n}}\,,\, \frac{\hat{r}_1\s_{1}}{64d\eta_0} \right\}=\hat{r}_1\min\left\{\frac{1}{4^{\n}}\,,\, \frac{\s_{1}}{64d\eta_0} \right\}\eqby{ritric}{\hat{r}_1}
\frac{\s_{1}}{64d\eta_0}\;,
 $$
and if 
$$
\hat r_{i+1}={\hat{r}_1}\frac{\s_{1}\cdots\s_i}{(64d\eta_0)^i}\;, \qquad i\ge 1\;,
$$
then
\begin{align*}
\hat r_{i+2}&= \min\left\{\frac{\wh{\a}}{4d\sqrt{2}\mathsf{K}_0\k_{i+1}^{\n}}\,,\, \frac{\hat{r}_{i+1}\s_{i+1}}{64d\eta_0} \right\}\\
	   &=\min\left\{\frac{\hat{r}_1}{4^{\n(i+1)}}\,,\, \hat{r}_{1}\frac{\s_1\cdots\s_{i+1}}{(64d\eta_0)^{i+1}} \right\}\\
	   &\eqby{ritric}\hat{r}_{1}\frac{\s_1\cdots\s_{i+1}}{(64d\eta_0)^{i+1}}\,,
\end{align*}
and \equ{riexp} is proven.\\
 
\noi
For a given $j\ge 2$, let $(\mathscr{P}^j)$ be the following assertion:  there exist
$j-1$ real--analytic diffeomorphisms
$$
G_{i+1}\colon  D_{\tilde{r}_{i+1}}(\mathscr D_{i}){\to}G_{i+1}( D_{\tilde{r}_{i+1}}(\mathscr D_{i}))\;,\quad \mbox{for}\quad 1\le i\le j-1\;,
$$
$j-1$ real--analytic symplectic transformations 
\beq{bes06v2v5}
\phi_{i+1}:D_{r_{i+1},s_{i+1}}(\mathscr{D}_{i+1})\to D_{2r_i/3, s_i}(\mathscr{D}_i),
\eeq
 and $j-1$  Hamiltonians $H_{i+1}=H_i\circ\phi_{i+1}=K_{i+1}+\vae^{2^{i+1}} P_{i+1}$ real--analytic on $D_{r_{i+1},s_{i+1}}(\mathscr{D}_{i+1})$ such that, for any $1\le i\le j-1$,
\beq{bbbBisv2v5}
\left\{
\begin{array}{l}
G_i(\mathscr D_{i-1})=\mathscr D_{i}\subset \mathscr D_{r_{i}} \, ,\ \\  \ \\
\|\dpr_y^2 K_i\|_{r_i,\mathscr{D}_i}\le \mathsf{K}_i\, ,\ \\  \ \\
\|T_i\|_{\mathscr{D}_i}\le \mathsf{T}_i\,,\ \\  \ \\ 
2^{i^2}|\vae|^{2^i}\|P_{i}\|_{r_{i},s_{i},\mathscr{D}_{i}}\le |\vae| \wh M_{i}\,,\ \\  \ \\
 \k_i\ge 4\s_i^{-1} \log\left(\s_i^{2\n+d}\m_i^{-1}\right)\,,\ \\  \ \\
|\vae|^{2^i} \mathsf{L}_i \le \frac{\sigma_i}{3}
\  
\end{array}\right.
\eeq
\noi
and
\beq{C.1Bisv2v5}
\left\{
\begin{aligned}
&\dpr_{y}K_{{i+1}}\circ G_{i+1}=\dpr_{y}K_{i} \;,\ \\  \ \\
&\| G_{i+1}-\id\|_{\tilde{r}_{i+1},\mathscr D_i}\le 2r_{i+1}\s_i^{\n+d}|\vae|^{2^i}\mathsf L_i\;,\ \\  \ \\
&\|\dpr_z G_{i+1}-\uno_d\|_{\tilde{r}_{i+1},\mathscr D_i}\le \s_i^{\n+d}|\vae|^{2^i}\mathsf L_i \;,\ \\  \ \\
&\|T_{i+1}\|_{\mathscr{D}_{i+1}}\le \|T_i\|_{\mathscr{D}_{i}}+\mathsf T_i|\vae|^{2^i}\mathsf L_i\;, 
 \\ \ \\
&\|K_{i+1}\|_{r_{i+1},\mathscr{D}_{i+1}}\le \|K_i\|_{r_i,\mathscr{D}_{i}}+|\vae|^{2^i}M_i \;,\\ \ \\ 
&\|\dpr_y^2K_{i+1}\|_{r_{i+1},\mathscr{D}_{i+1}}\le \|\dpr_y^2K_i\|_{r_i,\mathscr{D}_{i}}+\mathsf K_i|\vae|^{2^i}\mathsf L_i \;,\\ \ \\
&\|\mathsf{W}_{i+1}(\phi_{i+1}-\id)\|_{r_{i+1},s_{i+1},\mathscr{D}_{i+1}}\le \s_i^d\;|\vae|^{2^i}{\mathsf L}_i \;, \\ \ \\
& M_{i+1}\coloneqq\|P_{i+1}\|_{r_{i+1},s_{i+1},\mathscr{D}_{i+1}}\le  M_i \mathsf L_i\;.
\end{aligned}
\right.
\eeq
Assume $(\mathscr P^j)$, for some $j\ge 2$ and let us check $(\mathscr P^{j+1})$. Fix then $1\le i\le j-1$. Thus
$$
\|\dpr_y^2K_{i+1}\|_{r_{i+1},\mathscr{D}_{i+1}}\leby{C.1Bisv2v5} \|\dpr_y^2K_i\|_{r_i,\mathscr{D}_{i}}+\mathsf K_i|\vae|^{2^i}\mathsf L_i\leby{bbbBisv2v5} \mathsf K_i+\mathsf K_i\frac{\s_i}{3}=\mathsf K_{i+1}<\sqrt{2}\mathsf K_0
$$
and, similarly,
$$
\|T_{i+1}\|_{\mathscr{D}_{i+1}}\le \mathsf{T}_{i+1},
$$
which prove the second and third relations in \equ{bbbBisv2v5} for $i=j$. Therefore
\beq{alfhtrikpiv5}
\frac{\a}{r_{i+1}\mathsf{K}_{i+1}}> \frac{\a}{r_1\mathsf{K}_0\sqrt{2}}=\frac{\wh\a}{\hat r_1\mathsf{K}_0\sqrt{2}}=4d\k_0^{\n}>1
\eeq
so that
\begin{align}
|\vae|^{2^i} \mathsf L_i (3 \sigma_i^{-1})&= 3 |\vae|^{2^i}M_i\dst\max\left\{\frac{32\sqrt{2}\mathsf{T}_0  }{r_i {r}_{i+1}}\s_i^{-(\n+d)}\,,\,\mathsf{C}_4 \max\left\{1,\frac{\a}{r_i\mathsf{K}_i}\right\}\frac{ \mathsf{K}_0}{\a^2}\s_i^{-2(\n+d)}\right\}\sigma_i^{-1}\nonumber\\
                                          &\leby{alfhtrikpiv5} 3|\vae|^{2^i}M_i\dst\max\left\{\frac{32\sqrt{2}\mathsf{T}_0  }{r_i {r}_{i+1}}\,,\,\mathsf{C}_4 \frac{1}{\a r_i}\right\}\s_i^{-2(\n+d)-1}\nonumber\\
                                          &= 3\dst\max\left\{32\sqrt{2}\mathsf{T}_0\frac{ \wh\a }{ \hat{r}_{i+1}}\,,\,\mathsf{C}_4 \right\}\s_i^{-2(\n+d)-1}\frac{|\vae|^{2^i}M_i}{\a r_i}\nonumber\\
                                          &\eqby{riexp} 3\dst\max\left\{128d\sqrt{2}\eta_0\k_0^{\n}\cdot 2^{\frac{i^2}{2}}\left(\frac{64d\sqrt{2}\eta_0}{\s_0}\right)^i\,,\,\mathsf{C}_4 \right\}\s_i^{-2(\n+d)-1}\frac{|\vae|^{2^i}M_i}{\a r_i}\nonumber\\
                                          &\leby{kp08v5} 3\dst\max\left\{128d\sqrt{2}\,,\,8^{-\n}\mathsf{C}_4 \right\}\s_i^{-2(\n+d)-1}\frac{|\vae|^{2^i}M_i}{\a r_i}2^{\frac{i^2}{2}}\left(\frac{64d\sqrt{2}\eta_0}{\s_0}\right)^i\eta_0\k_0^{\n}\nonumber\\
                                          &= 3d\cdot 2^{4\n+2d+7}\sqrt{2}\dst\max\left\{128d\sqrt{2}\,,\,8^{-\n}\mathsf{C}_4 \right\}\s_0^{-(3\n+2d+2)}\left(\frac{2^{2\n+2d+7}d\sqrt{2}\eta_0}{\s_0}\right)^{i-1}\times\nonumber\\
                                          &\qquad\times\frac{2^{\frac{i^2}{2}}|\vae|^{2^i}M_i}{\a r_i}\eta_0^2\l_0^{\n}\nonumber\\
                                          &\eqby{riexp} 3d^2\cdot 2^{6\n+2d+10}\dst\max\left\{128d\sqrt{2}\,,\,8^{-\n}\mathsf{C}_4 \right\}\mathsf{K}_0\s_0^{-(4\n+2d+2)}\left(\frac{2^{2\n+2d+14}d^2\eta_0^2}{\s_0^2}\right)^{i-1}\times\nonumber\\
                                          &\qquad \times\frac{2^{\frac{i^2}{2}+\frac{(i-1)^2}{2}}|\vae|^{2^i}M_i}{|\vae|\wh\a^2}\eta_0^2\l_0^{2\n}\nonumber\\
                                          &= \mathsf{C}_5 \s_0^{-(4\n+2d+2)}\eta_0^2\l_0^{2\n}\mathsf{d}_*^{i-1}\mathsf{K}_0\frac{2^{i^2}|\vae|^{2^i}M_i}{|\vae|\wh\a^2}\nonumber\\
                                          &\leby{bbbBisv2v5} \mathsf{C}_5 \s_0^{-(4\n+2d+2)}\eta_0^2\l_0^{2\n}\mathsf{d}_*^{i-1}\frac{\mathsf{K}_0\wh M_i}{\wh \a^2}\nonumber\\
                                          &=\mathsf{e}_*\; \mathsf{d}_*^{i-1}\m_i\nonumber\\
                                          & = \frac{\theta_i}{\mathsf{d}_*}\nonumber\\
                                          &= \frac{(\sqrt{\mathsf{d}_*}\;\th_0)^{2^{i}}}{\mathsf{d}_*}\nonumber\\
                                          &\leby{condBisv2Prtv5}\su{\mathsf{d}_*}<1\;.\label{smCondSu1}
\end{align}
 Moreover,
 $$
 |\vae|^{2^{i}}\mathsf L_i<\mathsf{e}_*\; \mathsf{d}_*^{i-1}\m_i\;,
 $$ 
 thus by last relation in \equ{C.1Bisv2v5}, for any $1\le i\le j-1$,  
 $$
 2^{(i+1)^2}|\vae|^{2^{i+1}}M_{i+1}\le (2^{2i+1}|\vae|^{2^{i}}\mathsf{L}_i)(2^{i^{2}}|\vae|^{2^{i}} M_i)
 \leby{bbbBisv2v5}( 8\mathsf{e}_*(4\mathsf d_*)^{i-1} \; \m_i)(|\vae|\; \wh M_i) = |\vae|\; \wh M_{i+1}\;,
 $$ 
 which proves the fourth relation in \equ{bbbBisv2v5} for $i=j$. Hence, by exactly the same computation as above, one gets
 $$
 |\vae|^{2^{i+1}} \mathsf L_{i+1} (3 \sigma_{i+1}^{-1})\le \frac{\theta_{i+1}}{\mathsf{d}_*}=\frac{(\sqrt{\mathsf{d}_*}\;\th_0)^{2^{i+1}}}{\mathsf{d}_*}<1\ ,
 $$
 which proves the last relation in \equ{bbbBisv2v5} for $i=j$.  It remains only to check that
the fifth relation in \equ{bbbBisv2v5} holds as well for $i=j$ in order to apply Lemma~\ref{lem:1Extv5} to $H_i$, $1\le i\le j$ and get \equ{C.1Bisv2v5} and, consequently, $(\mathscr P^{j+1})$. But in fact, 
 we have, for any $i\ge 1$,\\\ \\
  $\mathsf{d}_*^{\su2 }\le\mathsf{e}_*\implies 4^i\mathsf{d}_*^{\su2 i}\le (8\mathsf{e}_*)^i\le (8\mathsf{e}_*)^{2^i-1}\implies 8\mathsf{e}_*(4\mathsf{d}_*)^{i}\le (8\mathsf{e}_*)^{2^i}\mathsf{d}_*^{\su2 i}\le (8\mathsf{e}_*)^{2^i}\mathsf{d}_*^{2^{i-2}} \implies 8\mathsf{e}_*(4\mathsf{d}_*)^{i}\mathsf{d}_*^{-2^{i-1}}\le (8\mathsf{e}_*)^{2^i}\mathsf{d}_*^{-2^{i-2}}< (8\mathsf{e}_*)^{2^i}$,\\ \ \\
 so that
 \begin{align}
 4\s_i^{-1} \log\left(\s_i^{2\n+d}\m_i^{-1}\right)&\le 4\s_i^{-1} \log\left(\m_i^{-1}\right)\nonumber\\
         &= 4\s_i^{-1} \log\left({8\mathsf{e}_*\; (4\mathsf{d}_*)^{i}}(\sqrt{\mathsf{d}_*}\;\th_0)^{-2^i}\right)\nonumber\\
         &\le 4\s_i^{-1} \log\left(\left(\frac{\th_0}{8\mathsf{e}_*}\right)^{-2^i}\right)\nonumber\\
         &= 4\s_i^{-1} \log\left(\m_0^{-2^i}\right)\nonumber\\
         &= 4^i\cdot 4\s_0^{-1} \log\left(\m_0^{-1}\right)\nonumber\\
         &= \k_i \;.\label{LamITronc}
 \end{align}
\noi
To finish the proof of the induction \ie one can construct an {\sl infinite sequence} of Arnold's transformations satisfying \equ{bbbBisv2v5} and  \equ{C.1Bisv2v5} {\sl for all $i\ge 1$}, one needs only to check $(\mathscr P^{2})$. But,\footnote{Observe that for $j=2$, $i=1$.} $\equ{HjExtExtv5}\div\equ{estfin2Bis000v2}, \equ{smCondSu1}_{i=1}$ and $\equ{LamITronc}_{i=1}$ imply $\equ{bbbBisv2v5}_{i=1}$.
 Thus, we apply Lemma~\ref{lem:1Extv5} to $H_1$ 
 to achieve the proof of $(\mathscr P^{2})$.\\
\nl
Next, we show that $G^j$ converges. For any $j\ge 1$,
\begin{align}
\|G^{j+1}-G^j\|_{\mathscr D_{0}}&= \|G_{j+1}\circ G^j-G^j\|_{\mathscr D_{0}}\nonumber\\
		&= \|G_{j+1}-\id\|_{\mathscr D_{j}}\nonumber\\
		&\le \|G_{j+1}-\id\|_{\tilde{r}_{j+1},\mathscr D_{j}}\nonumber\\
		&\leby{C.1Bisv2v5} 2r_{j+1}\s_j^{\n+d}|\vae|^{2^j}\mathsf L_j .\label{DstrClToDj}
\end{align}
Thus, $G^j$ is Cauchy and therefore converges uniformly on $\mathscr D_{\d,\a}$ to a map $G^*$.

\nl
Next, we prove that $\phi^j$ is convergent by showing that it is Cauchy as well. For any $j\ge 4$, we have, using again Cauchy's estimate,
\beqano
\|\mathsf{W}_{j-1}(\phi^{j-1}-\phi^{j-2})\|_{r_j,s_j,\mathscr{D}_j}&=&\|\mathsf{W}_{j-1}\phi^{j-2}\circ\phi_{j-1}-\mathsf{W}_{j-1}\phi^{j-2}\|_{r_{j-1}, s_{j-1},\mathscr{D}_{j-1}}\\
           &\leby{bes06v2v5}& \|\mathsf{W}_{j-1}D\phi^{j-2}\mathsf{W}_{j-1}^{-1}\|_{2r_{j-2}/3, s_{j-2},\mathscr{D}_{j-2}}\, \|\mathsf{W}_{j-1}(\phi_{j-1}-\id)\|_{r_{j-1}, s_{j-1},\mathscr{D}_{j-1}}\\
           &\leby{C.1Bisv2v5}&  \max\left(r_{j-1}\frac{3}{r_{j-1}},\frac{3}{2\s_{j-1}}\right)    \|\mathsf{W}_{j-1}\phi^{j-2}\|_{r_{j-1}, s_{j-1},\mathscr{D}_{j-1}} \times\\
           &&\qquad \times \|\mathsf{W}_{j-1}(\phi_{j-1}-\id)\|_{r_{j-1}, s_{j-1},\mathscr{D}_{j-1}}\\
           &=&  \frac{3}{2\s_{j-1}}   \|\mathsf{W}_{j-1}\phi^{j-2}\|_{r_{j-1}, s_{j-1},\mathscr{D}_{j-1}} \, \|\mathsf{W}_{j-1}(\phi_{j-1}-\id)\|_{r_{j-1}, s_{j-1},\mathscr{D}_{j-1}}\\
           &\le & \frac{1}{2}    \|\mathsf{W}_{j-1}\phi^{j-2}\|_{r_{j-1}, s_{j-1},\mathscr{D}_{j-1}} \cdot \s_{j-2}^d\left(|\vae|^{2^{j-2}}{\mathsf{L}}_{j-2}3\s_{j-2}^{-1}\right)\\
           &\le & \frac{1}{2}    \|\mathsf{W}_{j-1}\phi_2\|_{r_{2}, s_2,\mathscr{D}_{2}} \cdot \s_{j-2}^d\;\th_{j-2}\\
           &\le &  \frac{1}{2}\left(\dst\prod_{i=2}^{j-2}\|\mathsf{W}_{i+1}\mathsf{W}_{i}^{-1}\| \right)\|\mathsf{W}_{2}\phi_2\|_{r_{2}, s_2,\mathscr{D}_{2}} \cdot \s_{j-2}^d\;\th_{j-2}\\
           &\eqby{alfhtrikpiv5}& \frac{1}{2}\left(\dst\prod_{i=2}^{j-2}\frac{r_{i-1}}{r_{i}} \right)\|\mathsf{W}_{2}\phi_2\|_{r_{2}, s_2,\mathscr{D}_{2}} \cdot \s_{j-2}^d\;\th_{j-2}\\
           &=& \frac{r_1}{2r_{j-2}}\|\mathsf{W}_{2}\phi_2\|_{r_{2}, s_2,\mathscr{D}_{2}} \cdot \s_{j-2}^d\;\th_{j-2}\\
           &=& \su2\s_{1}^d\;\|\mathsf{W}_{2}\phi_2\|_{r_{2}, s_2,\mathscr{D}_{2}} \cdot 2^{\su2(j-3)^2}\left(\frac{2^{6-d}d\sqrt{2}\eta_0}{\s_0}\right)^{j-3}\cdot (\sqrt{\mathsf{d}_*}\;\th_0)^{2^{j-2}}\\
           &\le& \su2\s_{1}^d\;\|\mathsf{W}_{2}\phi_2\|_{r_{2}, s_2,\mathscr{D}_{2}} \cdot 2^{2^{j-3}}\left(\frac{2^{6-d}d\sqrt{2}\eta_0}{\s_0}\right)^{2^{j-4}}\cdot (\sqrt{\mathsf{d}_*}\;\th_0)^{2^{j-2}}\\
           &=& \su2\s_{1}^d\;\|\mathsf{W}_{2}\phi_2\|_{r_{2}, s_2,\mathscr{D}_{2}} \cdot \left(\left(\frac{2^{8-d}d\sqrt{2}\eta_0}{\s_0}\right)^{\su4} \sqrt{\mathsf{d}_*}\;\th_0\right)^{2^{j-2}}\\
           &=& \su2\s_{1}^d\;\|\mathsf{W}_{2}\phi_2\|_{r_{2}, s_2,\mathscr{D}_{2}} \cdot \left(\mathsf{C}_{6}\eta_0^{\frac{5}{4}}\s_0^{-\frac{5}{4}} \th_{0}\right)^{2^{j-2}}  \;.
\eeqano
\noi
Therefore, for any $n\ge 2,\, j\geq 0$,
\begin{align}
\|\mathsf{W}_{2}(\phi^{n+j+1}-\phi^n)\|_{r_{n+j+1},s_{n+j+1},\mathscr{D}_{n+j+1}}&\leq  \sum_{i=n}^{n+j}\|\mathsf{W}_{2}(\phi^{i+1}-\phi^i)\|_{r_{i+1},s_{i+1},\mathscr{D}_{i+1}}\nonumber\\
&\le \sum_{i=n}^{n+j}\left(\dst\prod_{k=2}^{i}\|\mathsf{W}_{k}\mathsf{W}_{k+1}^{-1}\| \right)\|\mathsf{W}_{i+1}(\phi^{i+1}-\phi^i)\|_{r_{i+1},s_{i+1},\mathscr{D}_{i+1}}\nonumber\\
&\eqby{alfhtrikpiv5} \sum_{i=n}^{n+j}\dst\prod_{k=2}^{i}\max\left\{1\;,\frac{r_{k+1}}{r_k} \right\}\|\mathsf{W}_{i+1}(\phi^{i+1}-\phi^i)\|_{r_{i+1},s_{i+1},\mathscr{D}_{i+1}}\nonumber\\
&= \sum_{i=n}^{n+j}\|\mathsf{W}_{i+1}(\phi^{i+1}-\phi^i)\|_{r_{i+1},s_{i+1},\mathscr{D}_{i+1}}\nonumber\\
&\le \su2\s_{1}^d\;\|\mathsf{W}_{2}\phi_2\|_{r_{2}, s_2,\mathscr{D}_{2}} \dst\sum_{i=n}^{n+j} \left(\mathsf{C}_{6}\eta_0^{\frac{5}{4}}\s_0^{-\frac{5}{4}} \th_{0}\right)^{2^{i+1}}\;.\label{phini}
\end{align}
Hence $\phi^j$ converges uniformly on $\mathscr D_*\times\tn$ to some $\phi^*$, which is then real--analytic function in $x\in\torus^d_{s_*}$.

\nl
To estimate $|\mathsf{W}_2(\phi^*-\id)|$ on $\mathscr{D}_*\times\torus^d_{s_*}$, observe that
, for $i\ge 1$,\footnote{Recall that $2^{i}\ge i+1,\, \forall\, i\ge 0$ and $\s_0\le \su2$.}
$$\s_{i}^d\;|\vae|^{2^i}\mathsf L_i\le \frac{\s_0^{d+1}}{3 \cdot 2^{i(d+1)}}\ \frac{(\sqrt{\mathsf{d}_*}\;\th_0)^{2^i}}{\mathsf{d}_*} \le \su{3 \cdot 2^{(d+1)(i+1)} \mathsf{d}_*}(\sqrt{\mathsf{d}_*}\;\th_0)^{{i+1}}= \su{3\mathsf{d}_*} \Big(\frac{\sqrt{\mathsf{d}_*}\;\th_0}{2^{d+1}}\Big)^{i+1}$$
and therefore 
$$\dst\sum_{i\ge 1}  \;|\vae|^{2^i}\mathsf L_i\le \frac{1}{3\mathsf{d}_*}\sum_{i\ge 1}\Big(\frac{\sqrt{\mathsf{d}_*}\;\th_0}{2}\Big)^{i+1}\le \frac{(\sqrt{\mathsf{d}_*}\;\th_0)^2}{6\;\mathsf{d}_*}= \frac{\th_0^2}{6}
\ ,$$ 
$$\dst\sum_{i\ge 1}  \s_{i}^d\;|\vae|^{2^i}\mathsf L_i\le \frac{1}{3\mathsf{d}_*}\sum_{i\ge 1}\Big(\frac{\sqrt{\mathsf{d}_*}\;\th_0}{2^{d+1}}\Big)^{i+1}\le \frac{(\sqrt{\mathsf{d}_*}\;\th_0)^2}{3\cdot {2^{2d+1}}\;\mathsf{d}_*}= \frac{\th_0^2}{3\cdot {2^{2d+1}}}
\ .$$ 
Moreover, for any $i\ge 2$,
\begin{align*}
\|\mathsf{W}_2(\phi^i-\id)\|_{r_{i},s_{i},\mathscr{D}_{i}}&\le \|\mathsf{W}_2(\phi^{i-1}\circ\phi_i-\phi_i)\|_{r_{i},s_{i},\mathscr{D}_{i}}+\|\mathsf{W}_2(\phi_i-\id)\|_{r_{i},s_{i},\mathscr{D}_{i}}\\
&\le \|\mathsf{W}_2(\phi^{i-1}-\id)\|_{r_{i-1},s_{i-1},\mathscr{D}_{i-1}}+ \left(\dst\prod_{j=0}^{i-1}\|\mathsf{W}_{j}\mathsf{W}_{j+1}^{-1}\| \right) \|\mathsf{W}_{i}(\phi_i-\id)\|_{r_{i},s_{i},\mathscr{D}_{i}}\\
&= \|\mathsf{W}_2(\phi^{i-1}-\id)\|_{r_{i-1},s_{i-1},\mathscr{D}_{i-1}}+ \|\mathsf{W}_{i}(\phi_i-\id)\|_{r_{i},s_{i},\mathscr{D}_{i}}\\
&= \|\mathsf{W}_2(\phi^{i-1}-\id)\|_{r_{i-1},s_{i-1},\mathscr{D}_{i-1}}+ \|\mathsf{W}_{i}(\phi_i-\id)\|_{r_{i},s_{i},\mathscr{D}_{i}}\\
&\le \|\mathsf{W}_2(\phi^{i-1}-\id)\|_{r_{i-1},s_{i-1},\mathscr{D}_{i-1}}+\s_{i-1}^d\;|\vae|^{2^{i-1}}{\mathsf{L}}_{i-1}\ ,
\end{align*}
which iterated yields
\begin{align*}
\|\mathsf{W}_2(\phi^i-\id)\|_{r_i,s_i,\mathscr{D}_i}&\le \dst\sum_{k=1}^{i-1}\s_{k}^d\; |\vae|^{2^k}{\mathsf{L}}_k\\
&\le  \dst\sum_{k\ge 1}\s_{k}^d\;|\vae|^{2^k}{\mathsf{L}}_k\\
&\le \frac{\th_0^2}{3\cdot {2^{2d+1}}}
\,.
\end{align*}
Therefore, taking the limit over $i$ completes the proof of \equ{estfin2Ext03v5}.

\noi
Next, we show that $\|G^*-\id\|_{L,\mathscr D_{\d,\a}}<1$, which will imply that\footnote{See Proposition II.2. in \cite{zehnder2010lectures}.} $G^*\colon \mathscr D_{\d,\a}\overset{onto}{\longrightarrow}\mathscr D_*$ is a lipeomorphism. Indeed, for any $j\ge 2$, we have
\begin{align*}
\| G^j-\id\|_{L,\mathscr D_{\d,\a}}+1&= \| (G_j-\id)\circ G^{j-1}+(G^{j-1}-\id)\|_{L,\mathscr D_{\d,\a}}+1\\
  &\le\| G_j-\id\|_{L,G^{j-1}(\mathscr D_{\d,\a})}\| G^{j-1}\|_{L,\mathscr D_{\d,\a}}+\| G^{j-1}-\id\|_{L,\mathscr D_{\d,\a}}+1\\
  &\le\| G_j-\id\|_{L,G^{j-1}(\mathscr D_{\d,\a})}(\|G^{j-1}-\id\|_{L,\mathscr D_{\d,\a}}+1)+\| G^{j-1}-\id\|_{L,\mathscr D_{\d,\a}}+1\\
  &= (\| G_j-\id\|_{L,\mathscr D_{j-1}}+1)(\| G^{j-1}-\id\|_{L,\mathscr D_{\d,\a}}+1)\\
  &\le (\|\dpr_z G_j-\uno_d\|_{\tilde{r}_j,\mathscr D_{j-1}}+1)(\| G^{j-1}-\id\|_{L,\mathscr D_{\d,\a}}+1)\\
  &\overset{\equ{estGidevidv2}+\equ{estG1devidv2}}{\le} (\s_{j-1}^{\n+d}|\vae|^{2^{j-1}}\mathsf{L}_{j-1}+1)(\| G^{j-1}-\id\|_{L,\mathscr D_{\d,\a}}+1)
\end{align*}
which iterated leads to\footnote{Recall that $\ex^x-1\le x\ex^x\;,\ \forall\; x\ge 0$.}
\begin{align}
\| G^j-\uno_d\|_{L,\mathscr D_{\d,\a}} &\le -1+\dst\prod_{i=1}^\infty (\s_{j-1}^{\n+d}|\vae|^{2^{i-1}}\mathsf{L}_{i-1}+1)\nonumber\\
     &\le -1+\exp\left( \sum_{i=0}^\infty \s_{i}^{\n+d}|\vae|^{2^{i}}\mathsf{L}_{i}\right)\nonumber\\
     &\le -1+\exp\left(\s_{0}^{\n+d}|\vae|\mathsf{L}_0+ \s_{0}^{\n+d}\sum_{i=1}^\infty |\vae|^{2^{i}}\mathsf{L}_{i}\right)\nonumber\\
     &\le -1+\exp\left(\s_{0}^{\n+d}\th_0+\s_{0}^{\n+d}\frac{\th_0^2}{6} \right)\nonumber\\
     &\le -1+\exp\left(2\s_{0}^{\n+d}\th_0\right)\nonumber\\
     &\le 2\s_{0}^{\n+d}\th_0\exp\left(2\s_{0}^{\n+d}\th_0\right)\nonumber\\
     &\ltby{condBisv2Prtv5}\frac{\ex\;\s_0^{\n+d}}{\mathsf{C}_6}<1\;.\label{LipG^jInf1}
\end{align}
Thus, letting $n\to\infty$, we get that $G^*$ is Lipschitz continuous, with
$$
\|G^*-\id\|_{L,\mathscr D_{\d,\a}}\le 2\s_{0}^{\n+d}\th_0\exp\left(2\s_{0}^{\n+d}\th_0\right)< 2\ex\;\s_{0}^{\n+d}\th_0< \frac{\ex\;\s_0^{\n+d}}{\mathsf{C}_6}<1,
$$
so that, by\footnote{With $\d\coloneqq 2\s_{0}^{\n+d}\th_0\exp\left(2\s_{0}^{\n+d}\th_0\right)$.} Lemma~\ref{LebLipLem} (see Appendix~\ref{appC}), we get
\begin{align*}
|\meas(\mathscr D_*)-\meas(\mathscr D_{\d,\a})|&\le  \left( \bigg(1+2\s_{0}^{\n+d}\th_0\exp\left(2\s_{0}^{\n+d}\th_0\right)\bigg)^d-1 \right) \meas(\mathscr D_{\d,\a})\\
   &\le d\cdot 2\s_{0}^{\n+d}\th_0\exp\left(2\s_{0}^{\n+d}\th_0\right)\left(1+2\s_{0}^{\n+d}\th_0\exp\left(2\s_{0}^{\n+d}\th_0\right)\right)^{d-1}\meas(\mathscr D_{\d,\a})\\
   &\le 2\ex\;d\left(\frac{3}{2}\right)^{d-1}\s_{0}^{\n+d}\th_0\meas(\mathscr D_{\d,\a})\\
   &= \mathsf{C}_7\;\s_{0}^{\n+d}\;\th_0\;  \meas(\mathscr D_{\d,\a}),
\end{align*}
which proves \equ{estfin2Ext02v5}.\\
Next, we show that $\phi^*\in C^\infty_W(\mathscr D_*\times\tn)$. For any $n,j\ge1$, we have
\begin{align*}
\|G^{n+j}-G^j\|_{\mathscr D_{\d,\a}}&\le \dst\sum_{k=j}^{n+j-1}\|G^{k+1}-G^k\|_{\mathscr D_{\d,\a}}\\
		&\leby{DstrClToDj}2\dst\sum_{k=j}^{n+j-1}r_{k+1}\s_k^{\n+d}|\vae|^{2^k}\mathsf L_k\\
		&\le 2r_{j+1}\s_j^{\n}\dst\sum_{k\ge 1}\s_k^{d}|\vae|^{2^k}\mathsf L_k\\
		&\le 2r_{j+1}\s_j^{\n}\;\frac{\th_0^2}{3\cdot {2^{2d+1}}}\\
		&\ltby{condBisv2Prtv5} \s_j^{\n}\frac{r_{j+1}}{16d\eta_0}\\
		&<\s_j^{\n}\;\tilde{r}_{j+1}\;.
\end{align*}
Now, letting $n\to\infty$, we get
\beq{G^*G^jDist}
\|G^{*}-G^j\|_{\mathscr D_{\d,\a}}<\s_j^{\n}\;\tilde{r}_{j+1}<\frac{\tilde{r}_{j+1}}{4}\;.
\eeq
Hence,\footnote{Recall that, by definition, $G^j(\mathscr D_{\d,\a})=\mathscr D_j$ and $G^*(\mathscr D_{\d,\a})=\mathscr D_*$. } for any $j\ge 1$,
\beq{HaDstr}
G^j\big(D_{\frac{\tilde{r}_{j+1}}{8}}(\mathscr D_{\d,\a})\big)\stackrel{\equ{LipG^jInf1}}{\subset}D_{\frac{\tilde{r}_{j+1}}{4}}\big(G^j(\mathscr D_{\d,\a}))\stackrel{\equ{G^*G^jDist}}{\subset} D_{\frac{\tilde{r}_{j+1}}{2}}(\mathscr D_*)\stackrel{\equ{G^*G^jDist}}{\subset} D_{\tilde{r}_{j+1}}(\mathscr D_j)\subset D_{r_j}(\mathscr D_j)\;.
\eeq
 Therefore, for any $n\ge 1$, we have
\beqano
\dst\sum_{j\ge 3}\|\mathsf{W}_2(\phi^{j}-\phi^{j-1})\|_{\tilde{r}_{j+1}/2,s_j,\mathscr{D}_*}\; \left(\frac{\tilde{r}_{j+1}}{2}\right)^{-n}&\leby{HaDstr}& 2^{n+4}d\;\eta_0^n\dst\sum_{j\ge 3}\|\mathsf{W}_2(\phi^{j}-\phi^{j-1})\|_{r_j,s_j,\mathscr{D}_j}\; r_{j+1}^{-n}\\
	    &\overset{\equ{phini}+\equ{riexp}}{\le}& 2^{n+3}d\;\eta_0^n\;\s_{1}^d\;r_1^{-n}\;\|\mathsf{W}_{2}\phi_2\|_{r_{2}, s_2,\mathscr{D}_{2}}\times\\
		&&\times \dst\sum_{j\ge 3} \left(\mathsf{C}_{6}\eta_0^{\frac{5}{4}}\s_0^{-\frac{5}{4}} \th_{0}\right)^{2^{j}}{2^{n\frac{j^2}{2}}\left(\frac{64d\sqrt{2}\eta_0}{\s_0}\right)^{nj}}\\
&<& \infty\;,
\eeqano
since, for $j$ sufficiently large,
$$
\left(\mathsf{C}_{6}\eta_0^{\frac{5}{4}}\s_0^{-\frac{5}{4}} \th_{0}\right)^{2^{j}}{2^{n\frac{j^2}{2}}\left(\frac{64d\sqrt{2}\eta_0}{\s_0}\right)^{nj}}<\left(\sqrt{2}\mathsf{C}_{6}\eta_0^{\frac{5}{4}}\s_0^{-\frac{5}{4}} \th_{0}\right)^{2^{j}}\qquad\mbox{and}\qquad \sqrt{2}\mathsf{C}_{6}\;\eta_0^{\frac{5}{4}}\;\s_0^{-\frac{5}{4}}\; \th_{0}\leby{condBisv2Prtv5} \su{\sqrt{2}}<1.
$$
Thus, writing
$$
\phi^{j}=(\phi^{j}-\phi^{j-1})+\cdots+(\phi^{3}-\phi^{2})\;,\qquad j\ge 3\;,
$$
and invoking Lemma~\ref{Whit1} (see Appendix \ref{appD}), we conclude that $\phi^*\in C^\infty_W(\mathscr D_*\times\tn)$.\\

\noi
Finally, we prove $G^*\in C^\infty_W(\mathscr D_{\d,\a})$ analogously. For any $j\ge 2$ and $n\ge 1$, we have
$$
G^j=(G^{j}-G^{j-1})+\cdots+(G^{2}-G^{1})\;,
$$
and
\begin{align*}
\dst\sum_{j\ge 1}\|G^{j+1}-G^{j}\|_{\tilde{r}_{j+1}/8,\mathscr D_{\d,\a}}\;\left(\frac{\tilde{r}_{j+1}}{8}\right)^{-n}&=8^n\dst\sum_{j\ge 1}\|(G_{j+1}-\id)\circ G^{j}\|_{\tilde{r}_{j+1}/8,\mathscr D_{\d,\a}}\;\tilde{r}_{j+1}^{-n}\\
    &\leby{HaDstr}8^n\dst\sum_{j\ge 2}\|G_{j+1}-\id\|_{\tilde{r}_{j+1},\mathscr D_{j}}\;\tilde{r}_{j+1}^{-n}\\
	&\le 2^{3n+1}\dst\sum_{j\ge 1}r_{j+1}\;\tilde{r}_{j+1}^{-n}\;\s_{j}^{\n+d}|\vae|^{2^j}\mathsf{L}_j\\
	&<\infty\;,
\end{align*}
which proves that $G^*\in C^\infty_W(\mathscr D_{\d,\a})$.\\

\noi
\noi
Now, to complete the proof of Theorem~\ref{Extteo4v2}, 
observe that, uniformly on $\mathscr{D}_*\times \torus^d_{s_*}$,
\begin{align*}
|\mathsf{W}_1(\phi_*-\id)|&\le |\mathsf{W}_1(\phi_1\circ \phi^*-\phi^*)|+|\mathsf{W}_1(\phi^*-\id)|\\
&\le \|\mathsf{W}_1(\phi_1-\id)\|_{r_1,s_1,\mathscr{D}_1}+\|\mathsf{W}_1\mathsf{W}_2^{-1}\|\;|\mathsf{W}_2(\phi^*-\id)|\\
&\le \s_{0}^d\;|\vae|{\mathsf{L}}_0+\frac{\th_0^2}{3\cdot {2^{2d+1}}}\\
&\le \frac{\s_0^{d+1}}{3 }\;\th_0 +\frac{\th_0^2}{3\cdot {2^{2d+1}}}\\
&\le \frac{1}{3\cdot {2^{d+1}} }\;\th_0 +\frac{\th_0^2}{3\cdot {2^{2d+1}}}\\
&\le \frac{\th_0}{3\cdot {2^{d}}}\,.
\end{align*}
Moreover, setting $G_{0}\coloneqq \id$, we have 
for any $i\ge 3$,
\begin{align*}
|G^{i}-\id|_{\mathscr D_{\d,\a}}&\le \dst\sum_{j=0}^{i-1}|G^{j+1}-G^{j}|_{\mathscr D_{\d,\a}}\\
		&=\dst\sum_{j=0}^{i-1}|G_{j}-\id|_{\mathscr D_{j-1}}\\
		&\overset{\equ{estGiidv2}+\equ{estG1idv2}}{\le}2\dst\sum_{j=0}^{i-1}r_{j+1}\s_j^{\n+d}|\vae|^{2^j}\mathsf L_j\\
		&\le 2\s_0^{\n}\;r_1\dst\sum_{j=0}^\infty \s_j^d|\vae|^{2^j}\mathsf L_j\\
		&\le \frac{2\;\s_0^{\n}\;\th_0}{3\cdot {2^{d}}}\;r_1\;,
\end{align*}
and then passing to the limit, we get
$$
|G^{*}-\id|_{\mathscr D_{\d,\a}}=|G^{*}-\id|_{\mathscr D_{\d,\a}}\le  \frac{2\;\s_0^{\n}\;\th_0}{3\cdot {2^{d}}}\;r_1
\le \frac{\s_0^\n}{\mathsf{C}_{9}}\left(\frac{\s_0}{\eta_0}\right)^{\frac{5}{4}}\frac{\a}{\mathsf{K}_0}=r_*\;.
$$
Finally, we prove \equ{MesChPin}. By Theorem~\ref{MintyExt}, $G^*-\id$ can be extended to a global Lipschitz continuous function $f\colon \rn\righttoleftarrow$, with\footnote{Where $
|y|_2\coloneqq \sqrt{y_1^2+\cdots+y_d^2}\;,$ and recall that $|y|\le |y|_2\le \sqrt{d}\;|y|$, for any $y\in\rn$.
}
\begin{align}
\dst\sup_{\rn}|f|_{2}&=\dst\sup_{\mathscr D_{\d,\a}}|G^*-id|_{2}
\;,\label{LipgloGstr0K}\\
\dst\sup_{\substack{y,y'\in\rn\\y\neq y'}}\frac{|f(y)-f(y')|_{2}}{|y-y'|_2}&=\dst\sup_{\substack{y,y'\in\mathscr D_{\d,\a}\\y\neq y'}}\frac{|(G^*-id)(y)-(G^*-id)(y')|_{2}}{|y-y'|_2}
\;.\label{LipgloGstr1K}
\end{align}
Hence,
\begin{align}
\|f\|_{\rn}&\overset{def}{=}\dst\sup_{\rn}|f|\nonumber\\
		   &\le \dst\sup_{\rn}|f|_{2}\nonumber\\
	       &\eqby{LipgloGstr0}\dst\sup_{\mathscr D_{\d,\a}}|G^*-id|_{2}\nonumber\\
	       &\le \sqrt{d}\dst\sup_{\mathscr D_{\d,\a}}|G^*-id|\nonumber\\
	       &\leby{NormGstrThtv2}\sqrt{d}\;r_*\nonumber\\
	       &\ltby{smcEAr0v2}\su{32d}\frac{r_0\s_0}{\eta_0}\nonumber\\
	       &\le {\d \s_0} \label{LipgloGstr01K}
\end{align}
and
\begin{align}
\|f\|_{L,\rn}&\overset{def}{=}\dst\sup_{\substack{y,y'\in\rn\\y\neq y'}}\frac{|f(y)-f(y')|}{|y-y'|}\nonumber\\
		     &\le \dst\sup_{\substack{y,y'\in\rn\\y\neq y'}}\frac{|f(y)-f(y')|_{2}}{|y-y'|_2/\sqrt{d}}\nonumber\\
		     &\eqby{LipgloGstr1K}\sqrt{d}\dst\sup_{\substack{y,y'\in\mathscr D_{\d,\a}\\y\neq y'}}\frac{|(G^*-id)(y)-(G^*-id)(y')|_{2}}{|y-y'|_2}\nonumber\\
		     &\le \sqrt{d}\dst\sup_{\substack{y,y'\in\mathscr D_{\d,\a}\\y\neq y'}}\frac{\sqrt{d}|(G^*-id)(y)-(G^*-id)(y')|}{|y-y'|}\nonumber\\
		     &=d\|G^*-\id\|_{L,\mathscr D_{\d,\a}}\nonumber\\
		     &\leby{LipGstrThtv2} d\frac{\ex\;\s_0^{\n+d}}{\mathsf{C}_6}<\su 2\;.\label{LipgloGstr11K}
\end{align}
Set $g\coloneqq f+\id$. Then, by Lemma~\ref{LipRang}, 
\beq{GsDdeltMT}
\mathscr D\subset g\big(\ovl{ B}_{\d \s_0}(\mathscr D)\big)
\;.
\eeq 
\newpage
Notice also that, by \equ{LipgloGstr11K},\footnote{See  \cite[Proposition~II.2.]{zehnder2010lectures}.} $g$ is a homeomorphism of $\rn$. Consequently,
\beqano
\meas(\mathscr D\times\tn\setminus \mathscr K)&=&\meas(\mathscr D\times\tn)-\meas\big(\phi_*(\mathscr D_*\times\tn)\big)\\
 &=& \meas(\mathscr D\times\tn)-\meas(\mathscr D_*\times\tn)\\
 &\leby{GsDdeltMT}&\meas\big(g(\ovl{B}_{\d \s_0}(\mathscr D))\times\tn\big)-\meas(\mathscr D_*\times\tn)\\
 &=&(2\pi)^d \meas\big(g(\ovl{B}_{\d \s_0}(\mathscr D))\setminus g(\mathscr D_{\d,\a})\big)\\
 &=&(2\pi)^d\meas\big(g(\ovl{B}_{\d \s_0}(\mathscr D)\setminus \mathscr D_{\d,\a})\big)\qquad\qquad \mbox{(because $g$ is injective)}\\
 &\le&(2\pi)^d\|g\|_{L,\rn}^d\meas\big(\ovl{B}_{\d \s_0}(\mathscr D)\setminus \mathscr D_{\d,\a}\big)\\
 &\le&(2\pi)^d(1+ \|f\|_{L,\rn})^d\meas\big(\ovl{B}_{\d \s_0}(\mathscr D)\setminus \mathscr D_{\d,\a}\big)\\
 &\leby{LipgloGstr11K}&(2\pi)^d\left(1+\frac{d\;\ex\;\s_0^{\n+d}}{\mathsf{C}_6}\right)^d\meas\big(\ovl{B}_{\d \s_0}(\mathscr D)\setminus \mathscr D_{\d,\a}\big)\\
 &=&(2\pi)^d\left(1+\frac{d\;\ex\;\s_0^{\n+d}}{\mathsf{C}_6}\right)^d\bigg(\meas\big(\ovl{B}_{\d \s_0}(\mathscr D)\setminus \mathscr D\big)+\meas(\mathscr D\setminus \mathscr D_\d)+\\
 &&\hspace{8cm}\;+\meas(\mathscr D_\d\setminus \mathscr D_{\d,\a}\big)\bigg)\\
\eeqano
 \qed


\section{Sharp measure estimate of the complement of $\mathscr K$ in an arbitrary set}
The strategy here is to localize the {\it Kolmogorov set} and then sum them up. Thus, we start by examining the cube case.
\subsection{Local analysis: the case where $\mathscr D$ is a cube}
\thm{ApArnC1l1}
Let
$$
\mathscr D=B_{R}(y_0)\;,\qquad R>0\,
$$
and let the assumptions in Theorem~\ref{Extteo4v2} hold, with
\begin{align*}
\mathsf{C}&\coloneqq \frac{1}{32}+\frac{d\;\sqrt{d}}{\mathsf{C}_{9}}+\dst\sum_{0\neq k\in\zn}\frac{1}{|k|_1^\n}\;,\\
\d&\le 
\dst\min\left\{\frac{\mathsf{T}_0}{32d\s_0}\;\a \;,\;\frac{\mathsf{r}_0}{32d}\;,\;\frac{R}{4}\right\}\;.
\end{align*}
Furthermore, assume that
$$
K_y\colon \mathscr D\to\O\coloneqq K_y(\mathscr D) 
$$
is a diffeomorphism.
Then, 
$$
\meas\big(\mathscr D\times\tn\setminus \mathscr K\big)
\le \mathsf{C}\;(6\pi)^d\;\frac{\vth_0\mathsf{T}_0}{\s_0}\;{R^{d-1}\a}\;,
$$
with\footnote{Indeed, pick any matrix $A=[a_{ij}]_{1\le i,j\le d}$. Then $\|A\|=\dst\max_{1\le i\le d}|a_{i1}|+\cdots+|a_{id}|$ and $|\det A|=|\sum_{\x\in\;\X_d}a_{1\x(1)}\cdots a_{d\x(d)}|\le \sum_{\x\in\;\X_d}|a_{1\x(1)}|\cdots|a_{d\x(d)}|\le \prod_{1=1}^d(|a_{i1}|+\cdots+|a_{id}|)\le \|A\|^d$, where $\X_d$ is the set of permutations of $\{1,\cdots,d\}$.\label{vth0sup11}}
$$
\vth_0\coloneqq \frac{\mathsf{K}_0^d}{\varrho_0}\ge 1\;,\qquad\varrho_0 \coloneqq \dst\inf_{y\in 
\mathscr D_{\d}}|\det K_{yy}(y)|>0\;.
$$
\ethm
\proof We shall denote the {\it Euclidean norm} by\footnote{Recall that $|y|\le |y|_2\le \sqrt{d}\;|y|$, for any $y\in\rn$.}
$$
|y|_2\coloneqq \sqrt{y_1^2+\cdots+y_d^2}\;.
$$
Recall that 
$$\mathscr D_\d\overset{def}{=}B_{R-\d}(y_0)\overset{def}{=} \{y\in\rn\;:\;|y-y_0|=\dst\max_{1\le j\le d}|y_j-y_{0j}|<R-\d\}$$
 and $\phi_*(\mathscr D_*\times\tn)\subset \mathscr D\times\tn$. Therefore, 
\begin{align*}
\meas\big(\mathscr D\times\tn\setminus \phi_*(\mathscr D_*\times\tn)\big)&= \meas\big(\mathscr D\times\tn\big)-\meas\big(\phi_*(\mathscr D_*\times\tn)\big)\\
	&= \meas\big(\mathscr D\times\tn\big)-\meas\big(\mathscr D_*\times\tn\big)\\
	&= (2\pi)^d\big(\meas(\mathscr D)-\meas(\mathscr D_*)\big)\\
	&\le (2\pi)^d\big(\meas(\mathscr D\setminus\mathscr D_\d)+\meas(\mathscr D_\d\setminus \mathscr D_*)\big)\\
	&= (2\pi)^d\big((2R)^d-(2R-2\d)^d+\meas(\mathscr D_\d\setminus G_*(\mathscr D_0))\big)\\
	&\le (2\pi)^d\big(2^dd\;R^{d-1}\d+\meas(\mathscr D_\d\setminus G_*(\mathscr D_0))\big)\;.\\
%
\end{align*}
It remains to estimate $\meas(\mathscr D_\d\setminus G_*(\mathscr D_0))$. Firstly, thanks to Theorem~\ref{MintyExt} (see Appendix~\ref{appE}), 
$G_*-\id$ can be extended to a global Lipschitz continuous function $f\colon \rn\righttoleftarrow$, with
\begin{align}
\dst\sup_{\rn}|f|_{2}&=\dst\sup_{\mathscr D_0}|G_*-id|_{2}
\;,\label{LipgloGstr0}\\
\dst\sup_{\substack{y,y'\in\rn\\y\neq y'}}\frac{|f(y)-f(y')|_{2}}{|y-y'|_2}&=\dst\sup_{\substack{y,y'\in\mathscr D_0\\y\neq y'}}\frac{|(G_*-id)(y)-(G_*-id)(y')|_{2}}{|y-y'|_2}
\;.\label{LipgloGstr1}
\end{align}
Hence,
\begin{align}
\|f\|_{\rn}&\overset{def}{=}\dst\sup_{\rn}|f|\nonumber\\
		   &\le \dst\sup_{\rn}|f|_{2}\nonumber\\
	       &\eqby{LipgloGstr0}\dst\sup_{\mathscr D_0}|G_*-id|_{2}\nonumber\\
	       &\le \sqrt{d}\dst\sup_{\mathscr D_0}|G_*-id|\nonumber\\
	       &\leby{NormGstrThtv2}\sqrt{d}\;r_*\eqqcolon \wh r \label{LipgloGstr01}
\end{align}
and
\begin{align}
\|f\|_{L,\rn}&\overset{def}{=}\dst\sup_{\substack{y,y'\in\rn\\y\neq y'}}\frac{|f(y)-f(y')|}{|y-y'|}\nonumber\\
		     &\le \dst\sup_{\substack{y,y'\in\rn\\y\neq y'}}\frac{|f(y)-f(y')|_{2}}{|y-y'|_2/\sqrt{d}}\nonumber\\
		     &\eqby{LipgloGstr1}\sqrt{d}\dst\sup_{\substack{y,y'\in\mathscr D_0\\y\neq y'}}\frac{|(G_*-id)(y)-(G_*-id)(y')|_{2}}{|y-y'|_2}\nonumber\\
		     &\le \sqrt{d}\dst\sup_{\substack{y,y'\in\mathscr D_0\\y\neq y'}}\frac{\sqrt{d}|(G_*-id)(y)-(G_*-id)(y')|}{|y-y'|}\nonumber\\
		     &=d\|G_*-\id\|_{L,\mathscr D_0}\nonumber\\
		     &\leby{LipGstrThtv2} d\frac{\ex\;\s_0^{\n+d}}{\mathsf{C}_6}<\su 2\;.\label{LipgloGstr11}
\end{align}
Set $g\coloneqq f+\id$. Then, by Lemma~\ref{LipRang},
\beq{GstrDdelt}
\mathscr D_\d\subset g\left(\ovl{B}_{\wh r}(\mathscr D_\d)\right)=g\left(\ovl{B_{R-\d+\wh r}(y_0)}\right)\;.
\eeq
Notice also that, by \equ{LipgloGstr11},\footnote{See  \cite[Proposition~II.2.]{zehnder2010lectures}.} $g$ is a homeomorphism of $\rn$. Consequently,
\beqano
\meas(\mathscr D_\d\setminus G_*(\mathscr D_0))&\leby{GstrDdelt}&\meas\left(g\left(\ovl{B}_{R-\d+\wh r}(y_0)\right)\setminus G_*(\mathscr D_0)\right)\\
											   &\overset{\mbox{def}}{=}& \meas\left(g\left(\ovl{B}_{R-\d+\wh r}(y_0)\right)\setminus g(\mathscr D_0)\right)\\
											   &=& \meas\left(g\left(\ovl{B}_{R-\d+\wh r}(y_0)\setminus \mathscr D_0\right)\right)\qquad\qquad \mbox{(because $g$ is injective)}\\
											   &\le& \|g\|_{L,\rn}^d\meas\left(\ovl{B}_{R-\d+\wh r}(y_0)\setminus \mathscr D_0\right)\\
											   &\le& (1+ \|f\|_{L,\rn})^d\bigg(\meas\left(\ovl{B}_{R-\d+\wh r}(y_0)\setminus \mathscr D_\d\right)+\meas\left(\mathscr D_\d\setminus \mathscr D_0\right)\bigg)\\
											   &\leby{LipGstrThtv2}&\left(\frac{3}{2}\right)^d\bigg(2^d(R-\d+\wh r)^d-2^d(R-\d)^d+\meas\left(\mathscr D_\d\setminus \mathscr D_0\right)\bigg)\\
											   &\le& \left(\frac{3}{2}\right)^d\bigg(2^ddR^{d-1}\sqrt{d}\;\frac{\s_0^\n}{\mathsf{C}_{9}}\left(\frac{\s_0}{\eta_0}\right)^{\frac{5}{4}}\frac{\a}{\mathsf{K}_0}+\meas\left(\mathscr D_\d\setminus \mathscr D_0\right)\bigg)\\
											   &\le& R^{d-1}\frac{3^d\;d\;\sqrt{d}}{\mathsf{C}_{9}}\frac{\a}{\mathsf{K}_0}+\left(\frac{3}{2}\right)^d\meas\left(\mathscr D_\d\setminus \mathscr D_0\right)  \;.
\eeqano
Finally,
\begin{align*}
\meas(\mathscr D_\d\setminus\mathscr D_0)&=\dst\int_{\mathscr D_\d\setminus\mathscr D_0}dy\\
			&=\dst\int_{\{y\in\mathscr D_\d\;:\;K_y(y)\not\in \D_\a^\t\}}dy\\
			&=\dst\int_{\{z\in K_y(B_{R-\d}(y_0))\;:\;z\not\in \D_\a^\t\}}|\det K_{yy}|^{-1} dz\\
			&\le \su{\varrho_0}\dst\int_{\left\{z\in B_{(R-\d)\|K_{yy}\|_{r_0,\mathscr D_0}}(K_y(y_0))\;:\;z\not\in \D_\a^\t\right\}} dz\\
			&= \su{\varrho_0}\dst\int_{\dst\bigcup_{0\neq k\in\zn}\left\{z\in B_{(R-\d)\mathsf{K}_0}(K_y(y_0))\;:\;|k\cdot z|<\frac{\a}{|k|_1^\t}\right\}} dz\\
			&\le \su{\varrho_0}\dst\sum_{0\neq k\in\zn}\dst\int_{\left\{z\in B_{(R-\d)\mathsf{K}_0}(K_y(y_0))\;:\;|k\cdot z|<\frac{\a}{|k|_1^\t}\right\}} dz\\
			&\le \su{\varrho_0}\dst\sum_{0\neq k\in\zn}(2(R-\d)\mathsf{K}_0)^{d-1}\frac{2\a}{|k|_1^\n}\\
			&= \mathbf{a}_1 \;2^d\;(R-\d)^{d-1}\;\a\;,
\end{align*}
where
$$
\mathbf{a}_1\coloneqq \frac{\mathsf{K}_0^{d-1}}{\varrho_0}\dst\sum_{0\neq k\in\zn}\frac{1}{|k|_1^\n}\;.
$$
Putting all together, we get
\begin{align*}
\meas\big(\mathscr D\times\tn\setminus \phi_*(\mathscr D_*\times\tn)\big)&\le (2\pi)^d\bigg(2^dd\;R^{d-1}\d+R^{d-1}\frac{3^d\;d\;\sqrt{d}}{\mathsf{C}_{9}}\frac{\a}{\mathsf{K}_0}+\mathbf{a}_1 \;3^d\;(R-\d)^{d-1}\;\a\bigg)\\
	&= (2\pi)^d\left(2^dd\;\frac{\mathsf{K}_0\d}{\a}+\frac{3^d\;d\;\sqrt{d}}{\mathsf{C}_{9}}+3^d\mathsf K_0\;\mathbf{a}_1\right)\frac{R^{d-1}\a}{\mathsf{K}_0}\\
	&\le (6\pi)^d\left(\mathsf{K}_0\;\frac{d\d}{\a}+\frac{d\;\sqrt{d}}{\mathsf{C}_{9}}+\dst\sum_{0\neq k\in\zn}\frac{1}{|k|_1^\n}\right)\frac{\mathsf{K}_0^{d-1}}{\varrho_0}{R^{d-1}\a}\\
	&\le (6\pi)^d\left(\frac{\eta_0}{32\s_0}+\frac{d\;\sqrt{d}}{\mathsf{C}_{9}}+\dst\sum_{0\neq k\in\zn}\frac{1}{|k|_1^\n}\right)\frac{\mathsf{K}_0^{d-1}}{\varrho_0}{R^{d-1}\a}\\
	&\le (6\pi)^d\left(\frac{1}{32}+\frac{d\;\sqrt{d}}{\mathsf{C}_{9}}+\dst\sum_{0\neq k\in\zn}\frac{1}{|k|_1^\n}\right)\frac{\mathsf{T}_0\mathsf{K}_0^{d}}{\s_0\varrho_0}{R^{d-1}\a}\\
	&= \mathsf{C}\;(6\pi)^d\;\frac{\vth_0\mathsf{T}_0}{\s_0}\;{R^{d-1}\a}\;.
\end{align*}
\qed
\subsection{Global analysis}
Let $\mathfrak{D}$ be any non--empty, bounded subset of $\rn$.  Consider the Hamiltonian parametrized by $\vae\in\real$
\[H(y,x;\vae)\coloneqq K(y)+\vae P(y,x),\]
where $K,P$ are real--analytic functions defined on $\mathfrak{D}\times\tn$ with bounded holomorphic extensions to\footnote{Recall the notations in $\S\ref{parassnot}$}
$$
D_{\mathsf{r}_0,s_0}(\mathfrak{D})\coloneqq \dst\bigcup_{y_0\in\mathfrak{D}}D_{\mathsf{r}_0,s_0}(y_0)\,,
$$
for some $\mathsf{r}_0>0$ and $0<s_0\le 1$, the norm being
$$
\|\cdot\|_{\mathsf{r}_0,s_0,\mathfrak{D}}\coloneqq \dst\sup_{D_{\mathsf{r}_0,s_0}(\mathfrak{D})}|\cdot|\,.
$$
Assume that 
\beq{ArnoldCondExtv2}
\varrho\coloneqq \dst\inf_{y\in \mathfrak{D}}|\det K_{yy}(y)|> 0\;.
\eeq
Fix $0<s_*<s_0$ and define\footnote{Recall footnote\textsuperscript{\ref{vth0sup11}}.}
\begin{align*}
&\s_0         \coloneqq 2^{7-2\n}d(s_0-s_*)\;,\\
&M \coloneqq  \|P\|_{\mathsf{r}_0,s_0,\mathfrak{D}}\;,\\
&\mathsf{T} \coloneqq \|T\|_{\mathfrak{D}}\coloneqq \sup_{y\in\mathfrak{D}}\|T(y)\|\;,\\
&\mathsf{K} \coloneqq \|K_{yy}\|_{\mathsf{r}_0,\mathfrak{D}}\;,\\
&\eta		 \coloneqq \mathsf{T}\mathsf{K}\ge 1\;,\\
&\vth\coloneqq \frac{\mathsf{K}^d}{\varrho}\ge 1\;,\\
\d&\coloneqq \frac{\mathsf{T}}{32d\s_0}\a\,,\\
R_0&\coloneqq \frac{\mathsf{r}_0}{64d\eta^2}\,,\\
&r_*		   \coloneqq \frac{\s_0^\n}{\mathsf{C}_{9}}\left(\frac{\s_0}{\eta}\right)^{\frac{5}{4}}\frac{\a}{\mathsf{K}}\;,\\
&\m_*		    \coloneqq \sup\left\{\m\le \ex^{-1}: 2\;\mathsf{C}_5\;\mathsf{C}_6\;\s_0^{4\n+2d+\frac{13}{4}}\;\eta^{\frac{13}{4}}\;\m\;\left(\log\m^{-1}\right)^{2\n}\le 1 \right\}\;,
\end{align*}
where $T(y)\coloneqq K_{yy}(y)^{-1}$. Let
\begin{align*}
\mathsf{C}&\coloneqq  \frac{1}{32}+\frac{d\;\sqrt{d}}{\mathsf{C}_{9}}+\dst\sum_{0\neq k\in\zn}\frac{1}{|k|_1^\n}\;,\\
\wh{\mathsf{C}}&\coloneqq 64\;d\;\mathsf{C}\;.
\end{align*}
Given $R>0$, define the set $\mathscr C_R$ of coverings of $\mathfrak{D}$ by cubes as follows: $F\in \mathscr C_R$ {\it iff} there exists $n_F\in \natural$ and $n_F$ cubes, say $F_i,$ $1\le i\le n_F$, of equal side--length $2R$, centered at a point  $y_i\in\mathfrak{D}$ and such that 
\beq{CovIm}
F=\big\{F_i:1\le i\le n_F\big\} \qquad \mbox{and}\qquad\mathfrak{D}\subset \dst\bigcup_{i=1}^{n_F}F_i\;.
\eeq
Then define
\beq{DefR00}
\mathscr R\coloneqq \bigg\{
0<R\le R_0: \mathscr C_{R}\neq \emptyset\ \mbox{and}\ \inf_{F\in \mathscr C_{R}}n_F (2R)^d\le 2^d\meas(\mathfrak{D})\bigg\}
\eeq
and
\beq{Defn0St}
n_{\mathfrak{D}}\coloneqq \dst\min_{R\in\mathscr R}\min \bigg\{n_F: F\in \mathscr C_{R}\quad\mbox{and}\quad n_F R^d\le \meas(\mathfrak{D})\bigg\}\;.
\eeq 
Pick any $R_0'\in\mathscr R$ and $F^0\in \mathscr C_{R_0'}$ such that $n_{F^0}=n_{\mathfrak{D}}$. Then
$$
F^0=\big\{F_i^0\coloneqq B_{R_0'}(p_i')\;,\quad \mbox{for some}\quad p_i'\in \mathfrak{D}\;,\quad 1\le i\le n_{\mathfrak{D}}\big\}\;.
$$
Thus, the following holds.
\thm{LebLocGen} Let the above assumptions and notations hold. Assume
\beq{smcEAr0v2K}
\boxed{\a\le 
8d\frac{R_*\s_0}{\mathsf{T}} \qquad\mbox{and}\qquad  |\vae|\le \m_*\frac{\a^2}{\mathsf{K}\;M}\;,}
\eeq
where $R_*\in\{R_0,R_0'\}$.\\
{\bf Part I: Description of the local Kolmogorov's sets $\mathscr K^i$}\\
There exists $n_*\in \natural$ and and $\mathbf{p}_i\in\mathfrak{D}$, $1\le i\le n_*$, such that
$$
\mathfrak{D}\subset \dst\bigcup_{i=1}^{n_*}B_{R_*}(\mathbf{p}_i)\;,
$$
and the following holds. Pick any $1\le i\le n_*$. Define
\begin{align*}
H^i&\coloneqq K^i+\vae P^i\coloneqq  H_{|B_{R_*}(\mathbf{p}_i)\times\tn}\,,\\
\D_\a^\t &\coloneqq \left\{\o\in\rn:\quad |\o\cdot k|\ge \frac{\a}{|k|_1^\t}\,,\quad\forall\;k\in\zn\setminus\{0\} \right\}\,,\\
\mathscr D_0^i &\coloneqq \left\{y_0\in B_{R_*-\d}(\mathbf{p}_i): \ K_y(y_0)\in \D_\a^\t\right\}\,.
\end{align*}
 Then, there exist $\mathscr D_*^i\subset  B_{R_*-\d+r_*}(\mathbf{p}_i)$ having the same cardinality  as $\mathscr D_0^i$, a lipeomorphism $G^{*i}\colon \mathscr D_0^i\overset{onto}{\longrightarrow}\mathscr D_*^i$,  with $(G^{*i})^{-1}\in C_W^\infty(\mathscr D_*^i)$, 
a function $K_*^i\in C_W^\infty(\mathscr D_*^i,\real)$
and  a $C_W^\infty $--symplectic transformation\footnote{Which means that the Whitney--gradient $\nabla \phi_*^i=\dpr\phi_*^i/\dpr(y_*,x)$ satisfies $(\nabla \phi_*^i)\mathbb{J}(\nabla \phi_*)^T=\mathbb{J}$ uniformly on $\mathscr D_*^i\times \tn$, where $\mathbb{J}=\begin{pmatrix}0 & -\uno_d\\
\uno_d & 0\end{pmatrix}$.} $\phi_*^i\colon \mathscr D_*^i\times \tn\to \mathscr K^i\coloneqq \phi_*^i(\mathscr D_*^i\times \tn)\subset B_{R_*}(\mathbf{p}_i)\times \tn$ and real--analytic in $x\in\torus^d_{s_*}$,
 such that\footnote{Notice that the derivatives are taken in the sense of Whitney.} 
\begin{align}
\dpr_{y_*}K_*^i\circ G^{*i}&=\dpr_{y}K^i  \qquad\quad\quad\qquad\ \quad \mbox{on} \quad \mathscr D_0^i\;,\label{conjCaneq00v2K}\\
\dpr^\b_{y_*}(H^i\circ \phi_*^i)(y_*,x)&=\dpr^\b_{y_*}K_*^i(y_*),\qquad \forall\;(y_*,x)\in \mathscr D_*^i\times\tn,\quad \forall\; \b\in\natural_0^d \label{conjCaneq0v2K}
\end{align}
and 
\begin{align}
&\|G^{*i}-\id\|_{\mathscr D_0^i}
\le 
r_*\,,\label{NormGstrThtv2K}\\
&\|G^{*i}-\id\|_{L,\mathscr D_0^i}\le \frac{\ex\;\s_0^{\n+d}}{\mathsf{C}_6}\,,\label{LipGstrThtv2K}\;.
\end{align}
{\bf Part II: Sharp measure estimate of the complement of $\mathscr K$}\\
Set
$$
\mathscr K\coloneqq \dst\bigcup_{i=1}^{n_*}\mathscr K^i\subset B_{R_*-\d+r_*}(\mathfrak{D})\times\tn\;.
$$
{\bf Case 1:} $R_*=R_0$.\\
In that case,\footnote{$n_*=N$ and $\mathsf{p}_i=p_i$, where $N$ and $p_i$ are the ones appearing in Lemma~\ref{CovLem}, with $R=R_0$.}
\beq{beQNst}
1\le n_*\le \left(\left[\frac{\diam \mathfrak{D}}{R_0}\right]+1\right)^d
\eeq
and
\beq{PrMeEstDBL}
\meas\left(\left(\dst\bigcup_{i=1}^{n_*}B_{R_0}(\mathbf{p}_i)\times\tn\right)\setminus \mathscr K\right)\le \wh{\mathsf{C}}\;(6\pi)^d\;\frac{\vth\;\eta^2\;\mathsf{T}}{\s_0\;\mathsf{r}_0}\left(\diam\mathfrak{D}+\frac{\mathsf{r}_0}{2^6d\eta^2}\right)^{d}\a\;.
\eeq
{\bf Case 2:} $R_*=R_0'$. \\
In that case, $n_*=n_{\mathfrak{D}}$, $\mathsf{p}_i=p_i'$ and
\beq{PrMeEstD}
\meas\left(\left(\dst\bigcup_{i=1}^{n_{\mathfrak{D}}}B_{R_0'}(p_i')\times\tn\right)\setminus \mathscr K\right)\le \mathsf{C}\;(12\pi)^d\;n_{\mathfrak{D}}^{\su d}\;\frac{\vth\;\mathsf{T}}{\s_0}\;\meas(\mathfrak{D})^{\frac{d-1}{d}}\;{\a}\;.
\eeq
\ethm
\ \\

\noi
We shall need the following elementary Lemmas in the proof.
\lemtwo{CovLem}{Covering Lemma}
For any $R>0$, there exist\footnote{$[x]$ denotes the integer part function $\max\{n\in \integer: n\le x\}$, while $\lceil x\rceil$ denotes the ''ceiling function'' $\min\{n\in \integer: n\ge x\}$.}
$$
1\le N\le \left(\left[\frac{\diam \mathfrak{D}}{R}\right]+1\right)^d
$$
and $p_i\in\mathfrak{D}$, $1\le i\le N$ such that
$$
\mathfrak{D}\subseteq\dst\bigcup_{i=1}^{N}B_{R}({p}_i)\;.
$$
\elem
\proof Let $\r\coloneqq \diam \mathfrak{D}$ and $z_i\coloneqq \inf\{y_i:y\in E\}$. Then $\mathfrak{D}\subseteq z+B_{\r}(0)$. Now, let $0<R'<R$ close enough to $R$ so that
$$
\left\lceil\frac{\r}{R'}\right\rceil=\left[\frac{\r}{R}\right]+1\eqqcolon \wh N\;.
$$
Then, $\mathfrak{D}$ can be covered by $\wh N$ closed, contiguous cubes $\D_i$, $1\le j\le \wh N^d$\;, of equal side--length $2R'$.  
 Let $i_j$ those indices such that $\D_{i_j}\bigcap\mathfrak{D}$ and pick $p_j\in \D_{i_j}\bigcap\mathfrak{D}$; let $N$ be the number of such cubes. But then, one has $\D_{i_j}\subseteq B_{R}({p}_j)$, for each $1\le j\le N\le \wh N$. The Lemma is therefore proved.
\qed
\lem{LwBndDet}
Let\footnote{$\mathcal{S}_d(\cn)$ denotes the symmetric matrices of order $d$, with entries in $\cn$.} $A\colon D_R(y_0)\to \mathcal{S}_d(\cn)$ be a matrix--valued function. Assume that 
$$
a\coloneqq\dst\sup_{y\in D_R(y_0)}\|A(y)\|<1.
$$ 
Then, for any $y\in D_R(y_0)$, the eigenvalues $\uno_d+A(y)$ are bounded in modulus from below by $1-a$. hence, in particular,
\beq{DetBlwIdp}
|\det\big(\uno_d+A(y)\big)|\ge (1-a)^d\;,\qquad\forall\;y\in D_R(y_0)\;.
\eeq
\elem
\proof
Let $y\in D_R(y_0)$ and $v\neq 0$ an eigenvector of $\uno_d+A(y)$ with associated eigenvalue $\l$. Then
\begin{align*}
|\l|\|v\|&=\|v+A(y)v\|\\
	     &\ge \|v\|-\|A(y)v\|\\
	     &\ge \|v\|-\|A(y)\|\|v\|\\
	     &\ge (1-a)\|v\|>0\;.
\end{align*}
Thus, the Lemma is proven since the determinant is equal to the product of the eigenvalues, counted with multiplicities.
\qed
\proof {\bf of Theorem~\ref{LebLocGen}} Set
$$
r_0\coloneqq \min\{\mathsf{r_0}\;,\;32d\d\}\;.
$$
Then,
$$
\frac{\mathsf{T}}{32d\s_0}\;\a \leby{smcEAr0v2K} \frac{R_*}{4}
\le \frac{R_0}{4}<\frac{\mathsf{r}_0}{32d}\;.
$$
Hence,
\begin{align}
\d&= \dst\min\left\{\frac{\mathsf{T}}{32d\s_0}\;\a \;,\;\frac{\mathsf{r}_0}{32d}\;,\;\frac{R_*}{4}\right\}\;,\label{deltminCond}\\
r_0&= 32d\d \;,\nonumber
\end{align}
so that
\beq{Cond1Cub}
\a=\frac{r_0\s_0}{\mathsf{T}}\;.
\eeq
 Thus, thanks to \equ{deltminCond} and \equ{Cond1Cub}, we need only to check that $K_y$ is injective on $F^0_i$ in order to apply Theorem~\ref{ApArnC1l1} to $H_i$. 
But, for any $y\in D_{\mathsf{r}_0/(4\eta)}(\mathsf{p}_i)$,
\begin{align*}
\|\uno_d-T(\mathsf{p}_i)K_{yy}(y)\|&\le \mathsf{T}\|K_{yy}(y)-K_{yy}(\mathsf{p}_i)\|\\
		&\le \mathsf{T}\|K_{yyy}\|_{\mathsf{p}_i,\mathsf{r}_0/2}\;\frac{\mathsf{r}_0}{4\eta}\\
		&\le \mathsf{T} \frac{\|K_{yy}\|_{\mathsf{p}_i,\mathsf{r}_0}}{\mathsf{r}_0/2}\frac{\mathsf{r}_0}{4\eta}\\
		&\le\mathsf{T}\mathsf{K}\frac{1}{2\eta}\\
		&=\su2.
\end{align*}
Thus, by Lemma~\ref{inv1}, $g\coloneqq (K_y)^{-1}$ is a real analytic diffeomorphism on $D_{r'}(z_i)$, where 
$$
z_i\coloneqq K_y(\mathsf{p}_i) \qquad\mbox{and}\qquad r'\coloneqq \frac{\mathsf{r}_0}{8\eta\mathsf{T}}\le \su{2\|T(\mathsf{p}_i)\|}\frac{\mathsf{r}_0}{4\eta}\\;.
$$
Furthermore,
\beq{estInvKy0}
\dst\sup_{D_{r'}(z_i)}\|g_z\|\le 2\|T(\mathsf{p}_i)\|\le 2\mathsf{T}\;.
\eeq
Set $T'\coloneqq g_z(z_i)^{-1}=K_{yy}(\mathsf{p}_i)$. Then, for any $z\in D_{r'/(8\eta)}(z_i)$,
\begin{align*}
\|\uno_d-T'g_z(z)\|&\le \|T'\|\|g_z(z)-g_z(z_i)\|\\
		&\le \mathsf{K}\|g_{zz}\|_{z_i,r'/2}\;\frac{r'}{8\eta}\\
		&\le \mathsf{T} \frac{\|g_z\|_{z_i,r'}}{r'/2}\frac{r'}{8\eta}\\
		&\leby{estInvKy0}2\mathsf{T}\mathsf{K}\frac{1}{4\eta}\\
		&=\su2.
\end{align*}
Thus, again by Lemma~\ref{inv1}, the inverse of $g$, \ie $K_y$, is a real analytic diffeomorphism on $D_{r''}(\mathsf{p}_i)$ (since $g(z_i)=\mathsf{p}_i$), where 
$$
R_*\le R_0< r''\coloneqq \frac{\mathsf{r}_0}{64\eta^2}= \frac{r'}{16\eta\mathsf{K}}\le \su{2\|T'\|}\frac{r'}{8\eta}\;.
$$ 
Moreover, in exactly the same way as above, one gets
\beq{SmpKyy}
\dst\sup_{y\in D_{R_0}(\mathsf{p}_i)}\|\uno_d-T(\mathsf{p}_i)K_{yy}(y)\|\le 2\mathsf{T}\frac{\mathsf{K}}{\mathsf{r}_0/2}R_0\leby{DefR00}\su{32d\eta}<\su2\;.
\eeq
Hence,
\begin{align}
\dst\inf_{y\in D_{R_0}(\mathsf{p}_i)}|\det K_{yy}(y)|&= \dst\inf_{y\in D_{R_0}(\mathsf{p}_i)}\left|\det\bigg(K_{yy}(\mathsf{p}_i)\big\{\uno_d-(\uno_d-T(\mathsf{p}_i)K_{yy}(y))\big\}\bigg)\right|\nonumber\\
		&=\dst\inf_{y\in D_{R_0}(\mathsf{p}_i)}|\det K_{yy}(\mathsf{p}_i)|\left|\det\bigg(\uno_d-(\uno_d-T(\mathsf{p}_i)K_{yy}(y))\bigg)\right|\nonumber\\
		&\overset{\equ{SmpKyy} + \equ{DetBlwIdp} }{\ge}|\det K_{yy}(\mathsf{p}_i)|\left(1-\su2\right)^d\nonumber\\
		&\ge \frac{\varrho}{2^d}>0\;. \label{infDetKyyR0}
\end{align}
The estimates \equ{PrMeEstD} and \equ{PrMeEstDBL} then follow easily. For instance, let us treat the second case \ie $R_*=R_0'$. The case $R_*=R_0$ is proved in a similar way by firstly using Lemma~\ref{CovLem}, with $R=R_0$; then setting $\mathsf{p}_i=p_i$, $n_*=N$ and thus applying Theorem~\ref{ApArnC1l1} to each $H^i$.\\
Let then $R_*=R_0'$.  Thus, we can apply Theorem~\ref{ApArnC1l1} to 
  $H^i$. Hence, there exists a {\it Kolmogorov set} 
\beq{LocKolSet}
\mathscr{K}^i\subset  F^0_i\times\tn\;,
\eeq
 associated to 
 $H^i$, with all the desired properties and satisfying\footnote{Where, $\varrho_0$ is replaced by ${\varrho}/{2^d}$, thanks to \equ{infDetKyyR0}; $\mathsf{T}_0$ and $\mathsf{K_0}$ by $\mathsf{T}$ and $\mathsf{K}$ respectively.}
\beq{EqFi0K}
\meas(F_i^0\times\tn\setminus \mathscr{K}_i)\le \mathsf{C}\;(12\pi)^d\;\frac{\vth\;\mathsf{T}}{\s_0}\;{R_0^{d-1}\a}\;.
\eeq
Therefore
\begin{align*}
\meas(\mathfrak{D}\times\tn\setminus \mathscr{K})&\leby{CovIm} \meas\left(\left(\dst\bigcup_{i=1}^{n_0}F_i^0\times\tn\right)\setminus \left(\dst\bigcup_{i=1}^{n_0}\mathscr{K}^i\right)\right)\\
		&\leby{LocKolSet} \dst\sum_{i=1}^{n_0}\meas\left(F_i^0\times\tn\setminus \mathscr{K}^i\right)\\
		&\leby{EqFi0K} \dst\sum_{i=1}^{n_0}\mathsf{C}\;(6\pi)^d\;\frac{\vth\;\mathsf{T}}{\s_0}\;{R_0'^{d-1}\a}\\
		&=\mathsf{C}\;(12\pi)^d\;n_0^{\su d}\;\frac{\vth\;\mathsf{T}}{\s_0}\;\bigg(n_0\;R_0'^d\bigg)^{\frac{d-1}{d}}\;{\a}\\
		&\leby{Defn0St} \mathsf{C}\;(12\pi)^d\;n_0^{\su d}\;\frac{\vth\;\mathsf{T}}{\s_0}\;\big(\meas\mathfrak{D}\big)^{\frac{d-1}{d}}\;{\a}\;.
\end{align*}
\qed



\cleardoublepage
\appendix
\chapter*{Appendices}
\addcontentsline{toc}{chapter}{Appendices}
\setcounter{subsection}{0}
\renewcommand{\thesubsection}{\Alph{subsection}}

\appA
{On the initial order of truncation $\k_0$ of the Fourier series in Theorem~\ref{teo1}}
Let 
\[\Th>0,\, 0<\vth<1,\,0<\s\leq \su{20},\, \n>\bar{\n}>d\geq 2,\, \b\coloneqq 1-\frac{1}{\bar{\n}}+\frac{1}{\n},\,  0<\wt c\leq (1-\b)\ex,\]
with\footnote{Notice that 
\[\frac{\bar{\n}}{\wt c(d-1)^{\b}}\geq \frac{\n\bar{\n}}{(\n-\bar{\n})\ex}\cdot \frac{\bar{\n}}{(d-1)^\b}>\frac{\bar{\n}^2}{(\bar{\n}-1)^\b \ex}>\frac{\bar{\n}}{\ex}>\frac{2}{\ex}> \su2,\]
 so that one can choose $\vth=\su2$ and in that case, if one chooses in addition $\wt c=\bar{c}$, then $\wt \k_0=\k_0,\, \wt C_{11}=C_{11},\, \wt C_{12}=C_{12}$ and $ \wt C_{13}=C_{13}$, with $\bar{c},\, C_{11},\, C_{12}$ and $C_{13}$ as in $\S\ref{AssumpPosc}$ and $\S\ref{AssumpExtPosc}$.}
\[\vth\leq \frac{\bar{\n}}{\wt c(d-1)^{\b}},\]
\beqano
\wt \k_0 &= & \left[-\frac{\log\Th}{(1-\vth)\s}\right],\\
\wt C_{11} &= & \exp\left((1-\vth)\left(\left(\frac{\bar{\n}}{\wt c\vth}\right)^{1/\b}+\frac{1}{20}\right) \right),\\
\wt C_{12}&= & \left( \frac{\ex^{-\frac{1-\vth}{20}}}{2C_6}\left((1-\vth)\wt c\right)^{\bar{\n}}\right)^{-\n/{\bar{\n}}},\\
\wt C_{13}&= & \exp\left((1-\vth)\left(\left(\frac{2C_0C_4}{ C_5}\right)^{1/\bar{\n}}+\frac{1}{20}\right)\right).
\eeqano
Then
\lem{lemK0}

\noi
\begin{itemize}
\item[$(i)$] If\quad $\Th< \min(20^{\n}\wt C_{11}^{-1},\wt C_{12}^{-1})\s^{\n}$ \quad then
\beq{eqApp101}
 \k_0^{\bar{\n}}\s^{\bar{\n}} \ex^{-\k_0\s}\leq \Th\ex^{\frac{1-\vth}{20}} <\frac{1}{2C_6\k_0^{\bar{\n}}},\quad \k_0\s>d-1.
\eeq
\item[$(ii)$] If \quad $\Th\leq 20^{\n}\wt C_{13}^{-1}\s^{\n}$\quad and \quad $2h_0\k_0^{\bar{\n}}\leq \a$ \quad then
\beq{eqApp103}
\frac{4C_4h_0}{\a \s^{\bar{\n}}}\leq \frac{C_5}{C_0}.
\eeq
\end{itemize}
\elem

\proof 
Above all, notice that (for any $0<\b<1$)
\beq{eqApp104}
\forall\, t>1,\quad \frac{t}{\log t}\geq(1-\b)\ex t^\b\geq \wt c t^\b.
\eeq
Let \quad $t\coloneqq \k_0\s$.\\
Let's prove $(i)$. Assume that  $\Th< \min(20^{\n}\wt C_{11}^{-1},\wt C_{12}^{-1})\s^{\n}$. Then
\begin{align*}
t > -\frac{\log \Th }{1-\vth}-\s > -\frac{\log\left( (20\s)^{\n} \wt C_{11}^{-1} \right)}{1-\vth}-\frac{1}{20}
 &\geq -\frac{\log( \wt C_{11}^{-1} )}{1-\vth}-\frac{1}{20}\\
 &=\left(\frac{\bar{\n}}{\wt c\vth}\right)^{1/\b}\geq d-1\geq 1.
\end{align*}
Therefore $t>1$ and $t>\left(\frac{\bar{\n}}{\wt c\vth}\right)^{1/\b}$, so that
\beqano
\frac{t}{\log t}\overset{\eqref{eqApp104}}{\geq} \wt c t^\b\geq \frac{\bar{\n}}{\vth}&\Longrightarrow & t^{\bar{\n}}\leq \ex^{\vth t}\\
		&\Longrightarrow & t^{\bar{\n}} \ex^{-t}\leq \ex^{-(1-\vth)t}\leq \ex^{\log\Th+(1-\vth)\s}\leq \Th \ex^{\frac{1-\vth}{20}}.
\eeqano
On the other hand, since\quad $\Th\leq (20\s)^{\n}\wt C_{11}^{-1}\leq \wt C_{11}^{-1}< 1$\quad then
\beqano
\Th \k_0^{\bar{\n}} & \leq & \Th \left(\frac{\log \Th^{-1} }{(1-\vth)\s}\right)^{\bar{\n}} \overset{\eqref{eqApp104}}{\leq} \frac{ \Th^{1-{\bar{\n}}(1-\b)} }{\left((1-\vth)\wt c\s\right)^{\bar{\n}}}=\frac{ \Th^{{\bar{\n}}/\n} }{\left((1-\vth)\wt c\s\right)^{\bar{\n}}}\\
           &< & \frac{ (\wt C_{12}^{-1}\s^{\n})^{\bar{\n}/\n} }{\left((1-\vth)\wt c\s\right)^{\bar{\n}}} = \frac{\ex^{-\frac{1-\vth}{20}}}{2C_6}.
\eeqano
Finally we prove $(ii)$. If \quad $\Th\leq 20^{\n}\wt C_{13}^{-1}\s^{\n}$\quad and \quad $2h_0\k_0^{\bar{\n}}\leq \a$ \quad then
\beqano
\frac{4C_4h_0}{\a \s^{\bar{\n}}}& \leq & \frac{2C_4}{ t^{\bar{\n}}}\\
					  & <    & 2C_4\left(-\frac{\log \Th}{1-\vth}-\s\right)^{-{\bar{\n}}}\\
					  &< & 2C_4\left(-\frac{\log\left( (20\s)^{\n} \wt C_{13}^{-1}\right)}{1-\vth}-\frac{1}{20}\right)^{-{\bar{\n}}}\\
					  &\leq & 2C_4\left(-\frac{\log\left( \wt C_{13}^{-1}\right)}{1-\vth}-\frac{1}{20}\right)^{-{\bar{\n}}}= \frac{C_5}{C_0}.\\
\eeqano
\qed
\appB
{Smooth contraction mapping Lemma \label{appB}}
Let $r,\, s,\,\s,\,\d,\,\mathrm{L}>0$. 
Let $u\in C^\infty(\rn\times\tn,\tn)$, with $x\mapsto x+u(y,x)$ holomorphic on $\torus^d_{s+\d}$ for any given $y\in\rn$. 
   Assume,
\beq{cauTIn20}
\su{\s}\|u\|_{0,s+\d},\ \|u_x\|_{0,s+\d}\le\mathrm{L}\le \d\le \su2\;, 
\eeq
where
$$
\|\cdot\|_{0,s+\d}\coloneqq \dst\sup_{\rn\times\torus^d_{s+\d}}|\cdot|\;.
$$
Assume also that for any $n\in\natural$ there exists a  constant  
 $C_{n}>0$ with the following property: for any $\b_1,\b_2\in\natural_0^d$ with $|\b_1|_1+|\b_2|_1\le n$,
\beq{cauTIn}
{r}^{|\b_1|_1}\s^{|\b_2|_1-1}\|\dpr_{y}^{\b_1}\dpr_{x}^{\b_2}u\|_{0,s}\le C_{u,n}{\mathrm L}\;,.
\eeq
\lem{CauTypInInv}
Under the above assumptions,  there exists a unique map $v\in C^\infty(\rn\times\tn,\tn)$, with $x\mapsto x+u(y,x)$ holomorphic on $\torus^d_{s}$ such that for any given $y\in\rn$, the map $x\mapsto x+v(y,x)$ is the inverse of $x\mapsto x+u(y,x)$. Moreover, for any $n\in\natural$ there a constant $\wh C_{n}>0$ such that for any $\b_1,\b_2\in\natural_0^d$ with $|\b_1|_1+|\b_2|_1\le n$,
\beq{cauTInvInv}
{r}^{|\b_1|_1}\s^{|\b_2|_1-1}\|\dpr_{y}^{\b_1}\dpr_{x}^{\b_2}v\|_{0,s}\le \wh C_{n}{\mathrm L}
\eeq
 and 
\beq{cauTInEstv}
\|v\|_{0,s}\le \|u\|_{0,s+\d}\;,\qquad \|v_x\|_{0,s}\le \frac{\|u_x\|_{0,s+\d}}{1-\d}\;.
\eeq
Furthermore, $\wh C_{n}$ can be expressed in term of $C_{n}$, for any $n\in\natural_0$.
\elem
\proof
Let $\mathcal{F}$ be the set of $w\in C^0(\rn\times\torus^d_{s},\cn)$ such that 
$$
\|w\|_{0,s}\le \d\;.
$$
Then, $(\mathcal{F},\; \|\cdot\|_{0,s})$ is a Banach space and for any $w\in\mathcal{F}$,
$$
\|\Im(x+w(y,x))\|_{0,s}< s+\|w\|_{0,s}\le s+\d.
$$
Hence the map
$$F\colon \mathcal{F}\ni w \mapsto -u(\pi_1,\pi_2+w)\in \mathcal{F}$$
is well--defined.
Notice that 
$$x+v(y,x)+u(y,x+v(y,x))=x \Longleftrightarrow v(y,x)=-u(y,x+v(y,x))\Longleftrightarrow v=F(v).$$ 
Hence, we have to show that $F$ admits a unique fixed point. 
But
$$
\| F(w_1)- F(w_2)\|_{0,s}\le\|u_x\|_{0,s+\d}\|w_1-w_2\|_{0,s}\le  \d \|w_1-w_2\|_{0,s}\;, \quad\forall\;w_1,w_2\in {\mathcal{F}}
$$
\ie, $F$ is a contraction. Therefore, by the Banach's Fixed Point (or contraction mapping) Theorem, $F$ admits a unique fixed point, say $v$, and $v$ is obtained as the uniform limit of the sequence $(F^n(0))_n$. Thus, by Weierstrass's Theorem, $x\mapsto x+v(y,x)$ is holomorphic on $\torus^d_{s}$, for each $y\in\rn$. Moreover
$$
\|v\|_{0,s}=\|F(v)\|_{0,s}\le\|u\|_{0,s+\d}
$$
and, by differentiating $v=F(v)$ \wrt $x$, we get
$$
v_x=-(\uno_d+u_x)^{-1}u_x=-\left(\dst\sum_{n=0}^\infty (-u_x)^n\right)u_x
$$
so that
$$
\|v_x\|_{0,s}\le \left(\dst\sum_{n=0}^\infty\|u_x\|_{0,s+\d}^n\right)\|u_x\|_{0,s+\d}\le \frac{\|u_x\|_{0,s+\d}}{1-\d}\;,
$$
which conclude the proof of \equ{cauTInEstv}. 
Next, we shall proceed inductively for the remainder of the proof. We have
$$
\su{\s}\|v\|_{0,s}\le \su{\s}\|u\|_{0,s+\d} \le \mathrm{L}\;, 
$$
which proofs $\equ{cauTInvInv}_{n=0}$. Set $w(y,x)=(x+v(y,x)$ \ie $w=\pi_2+v$.  
Now, fix $m\in\natural_0$ and assume that $v\in C^m(\rn\times\torus^d_{s},\cn)$ and $\equ{cauTInvInv}_n$ holds, for any $0\le n\le m$. Then, using the multivariate F\`aa Di Bruno's formula (see \cite{CS96Mult}, Theorem~$2.1$) to differentiate $v=F(v)$, for any $\b=(\b_1,\b_2)\in\natural_0^{d}\times \natural_0^{d}$, with $|\b|_1= m+1$,
 we have\footnote{With the convention $0^0=1$.}
\begin{align*}
 -\su{\s}\L^{\b}\dpr^{\b}v&=-\su{\s}\L^{\b}\dpr_{y}^{\b_1}\dpr_{x}^{\b_2} v\\
   &=-\su{\s}\L^{\b}\dpr_{y}^{\b_1}\dpr_{x}^{\b_2} F(v)\\
   &=  \su{\s}\dst\sum_{\substack{\l_1\in\natural^{d}_0\\ \l_1\preccurlyeq \b_1}}\L^{(\l_1,0)}\;\dst\sum_{\substack{\l_2\in\natural^{d}_0\\1\le |\l_2|_1\le |(\b_1-\l_1,\b_2)|_1}}\L^{(0,\l_2)}\;\dpr^{\l_2}_x\dpr^{\l_1}_y u\sum_{j=1}^{|(\b_1-\l_1,\b_2)|_1}\,\sum_{(k,l)\in\mathscr S(j,(\b_1-\l_1,\b_2),\l_2)}\\
   &\quad\L^{(\b_1-\l_1,\b_2-\l_2)}\;((\b_1-\l_1,\b_2))!\,\prod_{i=1}^{j}\frac{(\dpr^{l_i}w)^{k_i}}{k_i!(l_i!)^{|k_i|_1}} \\
   &=\dst\sum_{\substack{\l=(\l_1,\l_2)\in\natural^{d}_0\times\natural^{d}_0\\ \l_1\preccurlyeq \b_1\\1\le |\l|_1\le m+1}}\su{\s}\L^{\l}\;\dpr^{\l} u\sum_{j=1}^{m+1-|\l_1|_1}\sum_{(k,l)\in\mathscr S(j,\b-(\l_1,0),\l_2)}(\b_1-\l_1)!\b_2!\prod_{i=1}^{j}\frac{(\su{\s}\L^{l_i}\dpr^{l_i}w)^{k_i}}{k_i!(l_i!)^{|k_i|_1}} \\
   &=\dst\sum_{\substack{\l=(\l_1,\l_2)\in\natural^{d}_0\times\natural^{d}_0\\ \l_1\preccurlyeq \b_1\\\l_1\neq0\;\mbox{ \scriptsize or}\;|\l_2|_1\neq1\\1\le |\l|_1\le m+1}}\su{\s}\L^{\l}\;\dpr^{\l} u\sum_{j=1}^{m+1-|\l_1|_1}\sum_{(k,l)\in\mathscr S(j,\b-(\l_1,0),\l_2)}(\b_1-\l_1)!\b_2!\prod_{i=1}^{j}\frac{(\su{\s}\L^{l_i}\dpr^{l_i}w)^{k_i}}{k_i!(l_i!)^{|k_i|_1}} +\\
   &\quad+\dst\sum_{\substack{(\l_1,\l_2)\in\natural^{d}_0\times\natural^{d}_0\\ \l_1=0, \;|\l_2|_1=1}}\dpr^{\l} u\sum_{j=1}^{1}\sum_{(k,l)=(\l_2,\b)}\b!\prod_{i=1}^{1}\frac{(\su{\s}\L^{l_i}\dpr^{l_i}w)^{k_i}}{k_i!(l_i!)^{|k_i|_1}} \\
   &=\dst\sum_{\substack{\l=(\l_1,\l_2)\in\natural^{d}_0\times\natural^{d}_0\\ \l_1\preccurlyeq \b_1\\\l_1\neq0\;\mbox{ \scriptsize or}\;|\l_2|_1\neq1\\1\le |\l|_1\le m+1}}\su{\s}\L^{\l}\;\dpr^{\l} u\sum_{j=1}^{m+1-|\l_1|_1}\sum_{(k,l)\in\mathscr S(j,\b-(\l_1,0),\l_2)}(\b_1-\l_1)!\b_2!\prod_{i=1}^{j}\frac{(\su{\s}\L^{l_i}\dpr^{l_i}w)^{k_i}}{k_i!(l_i!)^{|k_i|_1}}\; +\\
   &\quad+u_x\cdot \su{\s}\L^{\b}\dpr^{\b}(\pi_2+v) \\
%
\end{align*}
\ie,
\begin{align*}
 \su{\s}\L^{\b}\dpr^{\b}v &= -(\uno_d+u_x)^{-1}\left( \su{\s}\L^{\b} u_x\; \dpr^{\b}\pi_2 +\dst\sum_{\substack{\l=(\l_1,\l_2)\in\natural^{d}_0\times\natural^{d}_0\\ \l_1\preccurlyeq \b_1\\\l_1\neq0\;\mbox{ \scriptsize or}\;|\l_2|_1\neq1\\1\le |\l|_1\le m+1}}\su{\s}\L^{\l}\;\dpr^{\l} u\sum_{j=1}^{m+1-|\l_1|_1}\sum_{(k,l)\in\mathscr S(j,\b-(\l_1,0),\l_2)}\right.\\
 &\quad\left.(\b_1-\l_1)!\b_2!\prod_{\substack{i=1\\ \\ \\ \\ \\}}^{j}\frac{(\su{\s}\L^{l_i}\dpr^{l_i}(\pi_2+v))^{k_i}}{k_i!(l_i!)^{|k_i|_1}} \right)\;,\\
\end{align*}
where
\begin{align*}
&\bullet\L\coloneqq (\underbrace{r,\cdots,r}_{d},\underbrace{\s,\cdots,\s}_{d})\;,\\
&\bullet\mathscr S(j,\b-(\l_1,0),\l_2)\coloneqq \left\{(k,l)=(k_1,\cdots,k_j,l_1,\cdots,l_j)\in \left(\natural^{d}_0\right)^j\times\left(\natural^{2d}_0\right)^j : \,\dst\prod_{i=1}^j|k_i|_1>0,\right.\nonumber\\
                   &\qquad\qquad\qquad\qquad\qquad\quad\left. 0\prec l_1\prec\cdots\prec l_j\,,\quad\dst\sum_{i=1}^jk_i=\l_2\quad\mbox{and}\quad \dst\sum_{i=1}^j|k_i|_1 l_i=\b-(\l_1,0) \right\}\;,\\
&\bullet\forall\; k\in\natural,\; (a,b)\in \natural_0^k\times\natural_0^k\;, \left(a\preccurlyeq b\Longleftrightarrow a_j\le b_j\;,\forall 1\le j\le k\right) \;,
\end{align*}
and, for all $k\in\natural,\; (a,b)\in \natural_0^k\times\natural_0^k\;,  a\prec b$ if and only if one of the following holds
\begin{itemize}
\item[$(i)$] $|a|_1< |b|_1$ or
\item[$(ii)$] $|a|_1=|b|_1$ and there exists $1<j\le k$ such that $a_i=b_i$ for all $1\le i<j-1$ and $a_j<b_j$ .
\end{itemize}

Therefore, $v\in C^{m+1}(\rn\times\torus^d_{s},\cn)$. Moreover, since $\|(\uno_d+u_x)^{-1}\|\le \frac{1}{1-\d}\le 2$, by the inductive hypothesis, \equ{cauTIn20} and \equ{cauTIn}, we have 
\begin{align*}
\su{\s}\L^{\b}\|\dpr^{\b}v\|_{0,s}&\le 2\left( \|u_x\|_{0,s+\d}+\dst\sum_{\substack{\l=(\l_1,\l_2)\in\natural^{d}_0\times\natural^{d}_0\\ \l_1\preccurlyeq \b_1\\\l_1\neq0\;\mbox{ \scriptsize or}\;|\l_2|_1\neq1\\1\le |\l|_1\le m+1}}C_{|\l|_1}\mathrm{L}\sum_{j=1}^{m+1-|\l_1|_1}\sum_{(k,l)\in\mathscr S(j,\b-(\l_1,0),\l_2)}\right.\\
						&\quad\left.(\b_1-\l_1)!\b_2!\prod_{\substack{i=1\\ \\ \\ \\ \\}}^{j}\frac{(1+\wh C_{|l_i|_1}\mathrm{L})^{|k_i|_1}}{k_i!(l_i!)^{|k_i|_1}}\right)\\
						&\le 2\left( \mathrm{L}+\dst\sum_{\substack{\l=(\l_1,\l_2)\in\natural^{d}_0\times\natural^{d}_0\\ \l_1\preccurlyeq \b_1\\\l_1\neq0\;\mbox{ \scriptsize or}\;|\l_2|_1\neq1\\1\le |\l|_1\le m+1}}C_{|\l|_1}\mathrm{L}\sum_{j=1}^{m+1-|\l_1|_1}\sum_{(k,l)\in\mathscr S(j,\b-(\l_1,0),\l_2)}(\b_1-\l_1)!\b_2!(1+\wh C_{m})^{|\l_2|_1}\right)\\
						&\le \wh C_{m+1}\mathrm{L}\;,
\end{align*}
where $\wh C_{m+1}>0$ is an universal constant, independent upon $\b$
. Finally, notice that $\wh C_{m+1}$ can be expressed in term of $C_{m+1}$ if $\wh C_{m}$ can be expressed in term of $C_{m}$. These concludes the proof of the Lemma.
\qed
\newpage
\appE{Extension of Lipschitz--H\"older continuous functions with control on the sup--norm\label{appE}}
We aim to recall here a very deep Extension Theorem for Lipschitz--H\"older continuous function, following closely \cite{minty1970extension}.\footnote{Recall that, Kirszbraun’s Theorem (see \cite{federergeometric}, $\S2.10.43$) asserts only that one can extend a Lipschitz continuous function without increasing the Lipschitz constant.} 
\thmtwo{MintyExt}{G.~J. Minty\cite{minty1970extension}}
Let $(V,\average{\cdot\;,\cdot})$ be a separable inner product space, 
$\emptyset\neq A\subseteq V$, $b>0,\;0<a\le 1$ and $g\colon A\to \rn$ a $(a,b)$--Lipschitz--H\"older continuous function on $A$ \ie\footnote{Recall that $|\cdot|_2$ denotes the Euclidean norm on $\rn$.}
\beq{aLipHol}
|g(x_1)-g(x_2)|_2\le b\; \|x_1-x_2\|^a\;, \qquad\forall\; x_1,x_2\in A\;.
\eeq
Then, there exists a global $(a,b)$--Lipschitz--H\"older continuous function\footnote{\ie satisfying \equ{aLipHol} on $V$.} $G\colon V\to \rn$ such that $G|_A=g$. Futhermore, $G$ can be chosen in such away that $G(V)$ is contained in the closed convex hull of $g(A)$. Hence, in particular,
\beq{MinTextMint}
\dst\sup_{x\in V} \|G(x)\|=\dst\sup_{x\in A} \|g(x)\|\quad \mbox{and}\quad \dst\sup_{x_1\neq x_2\in V} \frac{\|G(x_1)-G(x_2)\|}{\|x_1-x_2\|}=\dst\sup_{x_1\neq x_2\in A} \frac{\|g(x_1)-g(x_2)\|}{\|x_1-x_2\|}\;.
\eeq
\ethm
We need some preliminaries to prove the Theorem. Given $n\in \natural$, we shall denote
$$
\Upsilon_n\coloneqq \left\{\l=(\l_1,\cdots,\l_n)\in [0,1]^n\;:\ \l_1+\cdots+\l_n=1 \right\}\;.
$$
\dfntwo{defKF}{Kirszbraun function}
Let $V_1$ be a $\real$--vector space and $X$ a non--empty set. A function $f\colon V_1\times X\times X\to \real$ is called Kirszbraun function (K--function) if:
\begin{itemize}
\item[$(i)$] $f$ is convex and for any $x_1,x_2\in X$ and for any finite--dimensional subspace $S$ of $V_1$, the function $f(\cdot,x_1,x_2)\colon S\ni y\to f(y,x_1,x_2)$ is Lower semicontinuous\footnote{\ie for any $t\in\real$, the sublevel set $\{y\in S: f(y,x_1,x_2)\le t\}$ is closed in $S$ endowed with the canonical topology.};
\item[$(ii)$] for any $n\in\natural$, for any $(y_1,x_1),\cdots,(y_n,x_n)\in V_1\times X$, for any $x\in X$ and for any $(\l_1,\cdots,\l_n)\in \Upsilon_n$, the inequality
\beq{inKFun}
\dst\sum_{1\le i,j\le n}\l_i\l_jf(y_i-y_j,x_i,x_j)\ge 2\dst\sum_{i=1}^n\l_if(y_i-y,x_i,x)\;,\qquad y\coloneqq \dst\sum_{j=1}^n\l_j y_j
\eeq
holds.
\end{itemize}
\edfn
Then, the following holds.
\thmtwo{MintyKfunc}{G.~J. Minty\cite{minty1970extension}}
Let $f\colon \rn\times V\times V\to \real$  be a K--function, $n\in\natural$, $(y_1,x_1),\cdots,(y_n,x_n)\in \rn\times V$. Assume that $f$ is continuous and  for any $1\le i,j\le n$,
\beq{fyiyjlt0}
f(y_i-y_j,x_i,x_j)\le 0\;.
\eeq
Then, given any $x\in V$, there exists $y$ in the convex hull of $\{ y_1,\cdots,y_n\}$ such that $f(y_i-y,x_i,x)\le 0$, for any $1\le i\le n$.
\ethm
\proof
Consider the function
$$
F\colon \Upsilon_n\times \Upsilon_n\ni (\l,\m)\mapsto \dst\sum_{i=1}^n\l_if\left( y_i-\dst\sum_{j=1}^n\m_j y_j,x_i,x\right)\;.
$$
Then, it is clear that $F$ is convex and lower semicontinuous in $\m$, concave and upper semicontinuous in $\l$. Thus, since $\Upsilon_n$ is compact and thanks to the von Neumann's Minimax Theorem, there exists $(\l^0,\m^0)\in \Upsilon_n\times\Upsilon_n$ such that
$$
F(\l^0,\m^0)\le\dst\max_{\l\in\Upsilon_n}F(\l,\m^0)=\min_{\m\in\Upsilon_n}\max_{\l\in\Upsilon_n}F(\l,\m)=\max_{\l\in\Upsilon_n}\min_{\m\in\Upsilon_n}F(\l,\m)=\dst\min_{\m\in\Upsilon_n}F(\l^0,\m)\le  F(\l^0,\m^0)\;.
$$
Hence,
\beq{Flmlm0}
F(\l,\m^0)\le F(\l^0,\m^0)\le F(\l^0,\m)\;,\qquad \forall\; \l,\m\in\Upsilon_n\;.
\eeq
But,
$$
2F(\l^0,\l^0)=2\dst\sum_{i=1}^n\l^0_if\left( y_i-\sum_{j=1}^n\l^0_j y_j,x_i,x\right)\leby{inKFun}\sum_{1\le i,j\le n}\l_i\l_jf(y_i-y_j,x_i,x_j)\leby{fyiyjlt0}0\;.
$$
Set
$$
y^0\coloneqq \dst\sum_{j=1}^n\m^0_j y_j\;.
$$
Therefore, for any $1\le i\le n$,
$$
f(y_i-y^0,x_i,x)=F(\d^i_i,\m^0)\leby{Flmlm0}F(\l^0,\l^0)\le 0\;,
$$
where $\d^i_j$ is the Kronecker delta: $\d^i_j=1$ if $i=j$ and $0$ otherwise.
\qed
We shall need also the following. 
\lem{MoMin}
Let $x_1,\cdots,x_n\in \rn$. Then,
\begin{itemize}
\item[$(i)$]  given any $\b,a_1,\cdots,a_n>0$, we have\footnote{$\G$ being the Euler's Gamma function.}
\beq{MoMinId}
\dst\sum_{1\le i,j\le n}\frac{\average{x_i\;,\;x_j}}{(a_i+a_j)^\b}=\frac{1}{\G(\b)}\int_0^\infty\left\|\sum_{i=1}^n\ex^{-a_it}x_i\right\|^2t^{\b-1}dt\ge 0\;.
\eeq
\item[$(ii)$]given any  $(\l_1,\cdots,\l_n)\in\Upsilon_n$ and any $0<a\le1$,
\beq{FoxIn}
\dst\sum_{1\le i,j\le n}\l_i\l_j|x_i-x_j|_2^{2a}\le \sum_{1\le i,j\le n}\l_i\l_j(|x_i|_2^2+|x_j|_2^2)^{a}\le 2\sum_{i=1}^n\l_i\left|x_i\right|_2^{2a}\;.
\eeq
\end{itemize}
\elem
\proof $(i)$ is trivial. Let us prove $(ii)$. Above all, we recall the Bernoulli inequality:
\beq{BernIneq}
(1+x)^r\le 1+rx\;,\qquad \forall\; x\ge -1\;,\ \forall\;0\le r\le 1\;.
\eeq
The case $a=1$ is obvious. Let then $0<a<1$. By the continuity of the norm, up to approximating the zero vector by a sequence of non--zero vectors, we can assume that each $x_i\neq 0$, $i=1,\cdots,n$. Thus, we have
\begin{align*}
\dst\sum_{1\le i,j\le n}\l_i\l_j|x_i-x_j|_2^{2a}&= \dst\sum_{1\le i,j\le n}\l_i\l_j\average{x_i-x_j,x_i-x_j}^{a}\\
		&=\sum_{1\le i,j\le n}\l_i\l_j(|x_i|_2^2+|x_j|_2^2)^{a}\left(1-\frac{2\average{x_i\;,\;x_j}}{|x_i|^2+|x_j|_2^2}\right)^a\\
		&\leby{BernIneq} \sum_{1\le i,j\le n}\l_i\l_j(|x_i|_2^2+|x_j|_2^2)^{a}\left(1-\frac{2a\average{x_i\;,\;x_j}}{|x_i|_2^2+|x_j|_2^2}\right)\\
		&=\sum_{1\le i,j\le n}\l_i\l_j(|x_i|_2^2+|x_j|_2^2)^{a}-2a\sum_{1\le i,j\le n}\frac{\average{\l_ix_i\;,\l_jx_j}}{(|x_i|_2^2+|x_j|_2^2)^{1-a}}\\
		&\leby{MoMinId}\sum_{1\le i,j\le n}\l_i\l_j(|x_i|_2^2+|x_j|_2^2)^{a}\\
		&= \sum_{1\le i,j\le n}\l_i\l_j\max\{|x_i|_2,|x_j|_2\}^{2a}\left(1+\frac{\max\{|x_i|_2,|x_j|_2\}^{2}}{\max\{|x_i|_2,|x_j|_2\}^2}\right)^{a}\\
		&\leby{BernIneq} \sum_{1\le i,j\le n}\l_i\l_j\max\{|x_i|_2,|x_j|_2\}^{2a}\left(1+a\left(\frac{\min\{|x_i|_2,|x_j|_2\}}{\max\{|x_i|_2,|x_j|_2\}}\right)^2\right)\\
		&\le \sum_{1\le i,j\le n}\l_i\l_j\max\{|x_i|_2,|x_j|_2\}^{2a}\left(1+\left(\frac{\min\{|x_i|_2,|x_j|_2\}}{\max\{|x_i|_2,|x_j|_2\}}\right)^{2a}\right)\\
		&= \sum_{1\le i,j\le n}\l_i\l_j\left(|x_i|_2^{2a}+|x_j|_2^{2a} \right)\\
		&=2\sum_{i=1}^n\l_i\left|x_i\right|_2^{2a}\;.
\end{align*}
\qed
Now, we are in position to prove Theorem~\ref{MintyExt}.
\proof {\bf{of Theorem~\ref{MintyExt} }} The proof is divided into three steps.\\
{\bf{Step 1 }}We show that
$$
f\colon \rn\times V\times V\ni (y,x_1,x_2)\mapsto |y|_2^{2}- b^2\|x_1-x_2\|^{2a}
$$
is a K--function. $f$ is obviously continuous (actually, $\ci$) and convex in $y$. Now, let $n\in\natural$, $(y_1,x_1),\cdots,(y_n,x_n)\in \rn\times V$,  $x\in V$ and  $(\l_1,\cdots,\l_n)\in \Upsilon_n$, and set $y\coloneqq \dst\sum_{j=1}^n\l_j y_j$. Then
\begin{align*}
\dst\sum_{1\le i,j\le n}\l_i\l_jf(y_i-y_j,x_i,x_j)&= \dst\sum_{1\le i,j\le n}\l_i\l_j|(y_i-y)-(y_j-y)|_2^2-b^2\sum_{1\le i,j\le n}\l_i\l_j\|(x_i-x)-(x_j-x)\|^{2a}\\
	&=\dst\sum_{1\le i,j\le n}\l_i\l_j(|y_i-y|_2^2+\|y-y_j\|^2+2\average{y_i-y,y-y_j})-\\
	&\quad -b^2\sum_{1\le i,j\le n}\l_i\l_j\|x_i-x_j\|^{2a}\\
	&=2\dst\sum_{i=1}^n\l_i|y_i-y|_2^2+2b^2\average{\dst\sum_{i=1}^n \l_i(y_i-y),\dst\sum_{i=1}^n\l_j(y-y_j)}-\\
	&\quad -b^2\sum_{1\le i,j\le n}\l_i\l_j|x_i-x_j|^{2a}\\
	&=2\dst\sum_{i=1}^n\l_i|y_i-y|_2^2-b^2\sum_{1\le i,j\le n}\l_i\l_j\|(x_i-x)-(x_j-x)\|^{2a}\\
	&\geby{FoxIn} 2\dst\sum_{i=1}^n\l_i|y_i-y|_2^2-2b^2\dst\sum_{i=1}^n\l_i\|x_i-x\|^{2a}\\
	&= 2\dst\sum_{i=1}^n\l_if(y_i-y,x_i,x)\;.
\end{align*}
{\bf{Step 2 }} We want to show that we can extend $g$ to $A\bigcup \{x_0\}$ in such away that the image of $x_0$ by the extension lies in closed convex hull $\ovl{\conv(g(A))}$ of $g(A)$, for any $x_0\in V$. 
If $x_0\in A$, there is nothing to do. Let then $x_0\in V\setminus A$. Set\footnote{Notice that, for any $x\in A$, $b\|x-x_0\|^a>0$ and  $\{y\in \rn: f(g(x)-y,x,x_0)\le 0 \}$ is the closed ball (with respect to the Euclidean norm) centered at $g(x)$ with radius $b\|x-x_0\|^a$.}
$$
\mathscr{C}(x)\coloneqq \ovl{\conv(g(A))}\bigcap \left\{y\in \rn: f(g(x)-y,x,x_0)\le 0 \right\}\;,\qquad x\in A\;.
$$
Then, for any $x\in A$, $\mathscr{C}(x)$ is a compact convex subset of $\rn$. Now pick any $x_1,\cdots,x_{d+1}\in A$ and set $y_i\coloneqq g(x_i)\;,\ 1\le i\le d+1$. Thus, \equ{aLipHol} implies
$$
f(y_i-y_j,x_i,x_j)\le 0\;, \qquad\forall\; 1\le i\le d+1\;.
$$
Thanks to {\bf{Step 1}}, we can apply Theorem~\ref{MintyKfunc}. Therefore, there exists $y_0$ in the convex hull of $\{ y_1,\cdots,y_n\}$ such that $f(y_i-y_0,x_i,x_0)\le 0$, for any $1\le i\le d+1$. Hence,
$$
\dst\bigcap_{i=i}^{d+1}\mathscr C(x_i)\neq \emptyset\qquad (\mbox{since it contains } y_0).
$$
Thus, by Helly's Theorem\footnote{See \cite{danzer1921helly}}, there exists
$$
y_{x_0}\in\dst\bigcap_{x\in A}\mathscr C(x)\;.
$$
Consequently, the extension $g_{x_0}$ of $g$ to $A\bigcup \{x_0\}$ is obtained by setting $g_{x_0}(x_0)\coloneqq y_{x_0}$.\\
{\bf{Step 3 }} Pick any countable dense subset $D$ of $V$. Then, by {\bf{Step 2}}, we can extend $g$ inductively to $A\bigcup D$. Denote by $g_D$ such an extension and notice that $g_D(A\bigcup D)\subset \ovl{\conv(g(A))}$ and satisfies \equ{aLipHol} on $A\bigcup D$,  by construction. Now, pick any $x^0,x^1\in V\setminus A$ and  sequences $\{x^i_n\}\subset D$ converging to $x^i$, $i=0,1$. Fix $i=0,1$.  Then, for any $n,m\in \natural$,
$$
|g_D(x^i_n)-g_D(x^i_m)|_2\le b\|x^i_n-x^i_m\|^a\;.
$$
Hence, the sequence $\{g_D(x^i_n)\}\subset \rn$ is Cauchy and, therefore, converges to a $y^i\in\rn$ and $y^i$ does not depend upon the sequence chosen but only upon $x^i$. Now, by 
\begin{align*}
&g_D(x^0_n),\;g_D(x^1_n)\in g_D( D)\subset \ovl{\conv(g(A))}\;,\\
&|g_D(x^i_n)-g(x)|_2\le b\|x^i_n-x\|^a\;, \qquad 
\forall\;x\in A
\end{align*}
 and
$$
|g_D(x^0_n)-g_D(x^1_n)|_2\le b\|x^0_n-x^1_n\|^a\;, 
$$
for all $n\ge 0$, we get, by passing to the limit, 
\begin{align*}
&y^0,y^1\in \ovl{\conv(g(A))}\;,\\
&|y^i-g(x)|_2\le b\|x^i-x\|^a\;, \qquad 
\forall\;x\in A
\end{align*}
 and
$$
|y^0-y^1|_2\le b\|x^0-x^1\|^a\;.
$$
Then, a desired extension is obtained by just setting
$$
f(x)\coloneqq  g(x)\;,
$$
for $x\in A$ and
$$
f(x)\coloneqq \dst\lim g_D(x_n)\;,
$$
for $x\in V\setminus A$ and $\{x_n\}\subset D$ any sequence converging to $x$.
\qed
\newpage
\appC{Lebesgue measure and Lipschitz continuous map \label{appC}}
\lem{LebLipLem}
Let $\emptyset\not=A\subset\rn$ be a Lebesgue--measurable set and $f\colon A\to \rn$ be Lipschitz continuous with
\beq{apceq1}
\|f-\id\|_{L,A}\coloneqq \dst\sup_{\substack{x,y\in A\\x\not=y}}\frac{|f(x)-f(y)|}{|x-y|}\le \d\;.
\eeq
Then
\beq{EstApd}
|\meas(f(A))-\meas(A)|\le ((1+\d)^d-1)\meas(A)\;.
\eeq
\rem{appD0}
Notice that the inequality \equ{EstApd} is sharp as shown by the example $f=(1+\d)\;\id$.
\erem
\elem
\proof
By Theorem~\ref{MintyExt} (see Appendix~\ref{appE}), $f-\id$ can be extended to a Lipschitz continuous $g\colon \rn\righttoleftarrow$ with
$$
\|g\|_{L,\rn}=\|f-\id\|_{L,A}\le \d\;. 
$$
Now, by Rademacher’s Theorem, there exists a set $N\subset\rn$ with $\meas(N)=0$ and such that $g$ is differentiable on $\rn\setminus N$. Then\footnote{Let's point out that $\rn\setminus N$ is non--convex if $N$ is non--empty. Netherless, one can just approximate a segment by curves contained in $\rn\setminus N$ and with length arbitrarily close to the length of the segment.}
$$
\|g_y\|_{\rn\setminus N}=\|g\|_{L,\rn\setminus N}\le \|g\|_{L,\rn}\le \d\;. 
$$
Now pick $y\in \rn\setminus N$. Then, 
\begin{align*}
|\det(\uno_d+g_y(y))-1|&=\left|\dst\int_0^1\frac{d}{dt}\det(\uno_d+tg_y)dt \right|\\
					&=\left|\dst\int_0^1\tr\left(\adj(\uno_d+tg_y)g_y\right) dt \right|\\
					&\le \dst\int_0^1d\|\uno_d+tg_y\|^{d-1}\|g_y\|dt\\
					&\le \dst\int_0^1 d\left(1+\d t\right)^{d-1}\d dt\\
					&= (1+\d)^d-1\\
\end{align*}
Thus, by the change of variable (or area or coarea) formula\footnote{See \cite{evans2015measure}, $\S3.3$}, we have
\begin{align*}
|\meas(f(A))-\meas(A)|&=|\meas(g(A))-\meas(A)|\\
					  &=\left|\dst\int_{(\id+g)(A)}dy-\int_A dy \right|\\
					  &= \left|\dst\int_{(\id+g)(A\setminus N)}dy-\int_{A\setminus N} dy \right|\\
					  &= \left|\dst\int_{A\setminus N}|\det (\uno_d+g_y)| dy-\int_{A\setminus N} dy \right|\\
					  &\le \dst\int_{A\setminus N}|\det(\uno_d+ g_y)-1| dy\\
					  &\le ((1+\d)^d-1)\meas(A)\;.
\end{align*}
\qed
\newpage
\appD{Whitney's smoothness\label{appD}}
\dfn{WhitDef}
Let $A\subset \rn$ be non--empty and $n\in\natural_0$, $m\in\natural$. A function $f\colon A\to \real^m$ is said $C^n$ on $A$ in the Whitney sense, with Whitney derivatives $(f_\n)_{\n\in \natural_0^d,{|\n|_1\le n}}$ , $f_0=f$, and we write $f\in C^n_W(A,\real^m)$, if for any $\vae>0$ and $y_0\in A$, there exists $\d>0$ such that, for any $y,y'\in A\cap B_\d(y_0)$ and $\n\in \natural_0^d$, with  ${|\n|_1\le n}$,
\beq{whitdef}
\left|f_\n(y')-\dst\sum_{\substack{\m\in \natural_0^d\\ {|\m|_1\le n-|\n|_1}}} \frac{1}{\m!}f_{\n+\m}(y)(y'-y)^\m\right|\le \vae |y'-y|^{n-|\n|_1}\;.
\eeq
\edfn
The following is proven in \cite[$\S2.7$, pg.~58]{chierchia1986quasi} for $d=1$. 
\lem{Whit1}
Let $A\subset \rn$ be non--empty and $n\in\natural_0$. For $m\in\natural$, let $f_m$ be  a real--analytic function with holomorphic extension to $D_{r_m}(A)$, with $r_m\downarrow0$ as $m\rightarrow\infty$. Assume that
\beq{hypWhit}
a\coloneqq\dst\sum_{m=1}^\infty \|f_m\|_{r_m,A}\;r_m^{-n}<\infty,\qquad \|f_m\|_{r_m,A}\coloneqq \dst\sup_{B^d_{r_m}(A)}|f_m|\;.
\eeq
Then $f\coloneqq\dst\sum_{m=1}^\infty f_m\in C^n_W(A,\real)$ with Whitney's derivatives $f_\n\coloneqq\dst\sum_{m=1}^\infty\dpr_y^\n f_m$.
\elem
\proof
Let $\n\in \natural_0^d$, with ${|\n|_1\le n}$. We start showing that 
$$
f_\n= \dst\sum_{m=1}^\infty \dpr_y^\n f_m
$$
is well--defined on $A$. Indeed, for any $m\ge 1$, $f_m\in\ci(A)$ and, by Cauchy's estimate,
$$
\|\dst\sum_{m=1}^\infty \dpr_y^\n f_m\|_A\le \dst\sum_{m=1}^\infty \|\dpr_y^\n f_m\|_{\frac{r_m}{2}, A}\le 2^n\sum_{m=1}^\infty \| f_m\|_{r_m, A}\; r_m^{-|\n|_1}\le 2^n r_1^{n-|\n|_1}\sum_{m=1}^\infty \| f_m\|_{r_m, A}\; r_m^{-n}\ltby{hypWhit}\infty \;,
$$
where
$$
\|\cdot\|_A\coloneqq  \dst\sup_A|\cdot|\;.
$$
Finally, we show that $(\dpr_y^\n f)_{\n\in \natural_0^d,{|\n|_1\le n}}$ are the Whitney's derivatives of $f$. Fix then $y_0\in A$, $0<\vae<a$ and $\n\in \natural_0^d$, with ${|\n|_1\le n}$. Set
$$
b\coloneqq 2^n\; a\dst\sum_{\substack{\m\in \natural_0^d\\ {|\m|_1= n+1}}} 1\;.
$$
Let $m_1\in \natural$ such that\footnote{Such a $m_1$ exists by \equ{hypWhit}.} 
\beq{n0choi}
2^n\dst\sum_{m=m_1}^\infty \| f_m\|_{r_m, A}\; r_m^{-n}\le \frac{\vae^{n+1}}{(2b)^n}\qquad (<\vae)\;. 
\eeq
	%
Let\footnote{Let us point out that $\d$ does not depend upon $\n$. These is crucial! Actually, $\d$ does not even depend upon $y_0$.}
$$
\d\coloneqq \frac{\vae}{4b}\;r_{m_1}
$$ 
and
$$
f^{[k]}\coloneqq \dst\sum_{m=1}^{k} f_m \;, \qquad k\ge 1\;.
$$
Now, pick $y,y'\in A\cap B_{\d}(y_0)$, with $y\neq y'$. Let then\footnote{Notice that such a $m_2$ exists since $|y'-y|\le |y'-y_0|+|y_0-y|<2\d=\frac{\vae}{2b}r_{m_1}$ and the sequence $(r_m)_m$ is strictly decreasing.} $m_2\ge m_1$ such that
\beq{rm1yy}
\frac{\vae}{2b}\;r_{m_2+1}\le |y'-y| < \frac{\vae}{2b}\;r_{m_2}\;.
\eeq
Notice that $f^{[m_2]}$ is holomorphic on $D_{r_{m_2}}(A)$ and 
$$
0<r\coloneqq |y'-y| < \frac{\vae}{2b}\;r_{m_2}\le \frac{\vae}{2a}\;r_{m_2}<\frac{r_{m_2}}{2}\;.
$$
Moreover, for any $1\le m\le m_2$,
\beq{rmyyp0}
\left\{
\begin{aligned}
&r_m-|y'-y|\ge r_{m_2}-|y'-y|\gtby{rm1yy}\frac{b}{\vae}|y'-y| +\left(\frac{b}{\vae}-1\right)|y'-y|> \frac{b}{\vae}|y'-y| \;, \\
& r_m-|y'-y|=\left(\frac{r_m}{2}-|y'-y| \right)+\frac{r_m}{2}\ge \left(\frac{\vae}{2b}\;r_{m_2}-|y'-y| \right)+\frac{r_m}{2}\gtby{rm1yy} \frac{r_m}{2}\;.
\end{aligned}
\right.
\eeq
 Therefore, by Taylor--Lagrange's formula and Cauchy's estimates, we have (for some $0<t<1$)
\begin{align}
\left|f^{[m_2]}_\n(y')-\dst\sum_{\substack{\m\in \natural_0^d\\ {|\m|_1\le n-|\n|_1}}} \frac{1}{\m!}f_{\n+\m}^{[m_2]}(y)(y'-y)^\m\right|&= \left|\dst\sum_{\substack{\m\in \natural_0^d\\ {|\m|_1= n-|\n|_1+1}}} \frac{1}{\m!}f_{\n+\m}^{[m_2]}(y+t(y'-y))(y'-y)^\m\right|  \nonumber\\
                &\le \dst\sum_{\substack{\m\in \natural_0^d\\ {|\m|_1= n-|\n|_1+1}}} \sum_{m=1}^{m_2} \|\dpr^{\n+\m}_y f_m\|_{r,A}\; r^{|\m|_1}  \nonumber\\
                &\le r^{n-|\n|_1+1}\dst\sum_{\substack{\m\in \natural_0^d\\ {|\m|_1= n-|\n|_1+1}}} \sum_{m=1}^{m_2} \| f_m\|_{r_m,A}\;(r_m-r)^{-(|\n|_1+|\m|_1)} \nonumber \\
                &\le \frac{b}{2^n a} r^{n-|\n|_1+1}\dst \sum_{m=1}^{m_2} \| f_m\|_{r_m,A}\;(r_m-r)^{-(n+1)} \nonumber \\
                &\leby{rmyyp0}\frac{b}{2^n a} r^{n-|\n|_1+1}\dst \sum_{m=1}^{m_2} \| f_m\|_{r_m,A} \left(\frac{r_m}{2}\right)^{-n}\left(\frac{b}{\vae}|y'-y|\right)^{-1}\nonumber\\
                &=\vae\;  r^{n-|\n|_1}\;\su a\dst \sum_{m=1}^{m_2} \| f_m\|_{r_m,A} {r_m}^{-n}\nonumber\\
                &\le \vae\;  r^{n-|\n|_1}
 \;.\label{gpetio}
\end{align}
Furthermore, for any $\m\in \natural_0^d$, with $|\m|_1\le n$,
\begin{align}
\dst\sum_{m> m_2} \|\dpr_\m f_m\|_{A}&\le \dst\sum_{m> m_2} \|\dpr_\m f_m\|_{\frac{\vae}{4b}r_m,A}\nonumber\\
	&\le \dst\sum_{m> m_2} \| f_m\|_{\frac{\vae}{2b}r_m,A}\left(\frac{\vae}{4b}r_m\right)^{-|\m|_1}\nonumber\\
	&= \dst\sum_{m> m_2} \| f_m\|_{\frac{\vae}{2b}r_m,A}\left(\frac{\vae}{4b}r_m\right)^{-n}\left(\frac{\vae}{4b}r_m\right)^{n-|\m|_1}\nonumber\\
	&\le \left(\frac{\vae}{4b}r_{m_2+1}\right)^{n-|\m|_1}\left(\frac{2b}{\vae}\right)^{n}\;2^n\dst\sum_{m> m_1} \| f_m\|_{r_m,A}\;r_m^{-n}\nonumber\\
	&\overset{\equ{rm1yy}+\equ{n0choi}}{\le} r^{n-|\m|_1}\;\left(\frac{2b}{\vae}\right)^{n}\;\frac{\vae^{n+1}}{(2b)^n}\nonumber\\
	&= {\vae}\;r^{n-|\m|_1}\;.\label{delchoi}
\end{align}
Thus, 
\beqano
\left|f_\n(y')-\dst\sum_{\substack{\m\in \natural_0^d\\ {|\m|_1\le n-|\n|_1}}} \frac{1}{\m!}f_{\n+\m}(y)(y'-y)^\m\right| &\le& |f_\n(y')-f^{[m_2]}_\n(y')|+\\
 &+&\left|f^{[m_2]}_\n(y')-\dst\sum_{\substack{\m\in \natural_0^d\\ {|\m|_1\le n-|\n|_1}}} \frac{1}{\m!}f_{\n+\m}^{[m_2]}(y)(y'-y)^\m\right|+\\
 &+&\left|\dst\sum_{\substack{\m\in \natural_0^d\\ {|\m|_1\le n-|\n|_1}}} \frac{1}{\m!}\left( f_{\n+\m}^{[m_2]}(y)-f_{\n+\m}(y)\right)(y'-y)^\m\right|\\
 &\leby{gpetio}& \dst\sum_{m> m_2} \|\dpr_\n f_m\|_{A}+\\
 &+&\vae\;  r^{n-|\n|_1}+\\
 &+&\dst\sum_{\substack{\m\in \natural_0^d\\ {|\m|_1\le n-|\n|_1}}} \dst\sum_{m> m_2} \|\dpr_{\n+\m} f_m\|_{ A}\;r^{|\m|_1}\\
 &\overset{\equ{delchoi}}{\le}& \vae\;  r^{n-|\n|_1}+\\
 &+&\vae\;  r^{n-|\n|_1}+\\
 &+&\vae\;  r^{n-|\n|_1}\dst\sum_{\substack{\m\in \natural_0^d\\ {|\m|_1\le n-|\n|_1}}} 1\\
 &\le& (2+(n+1)^d)\;\vae\;|y'-y|^{n-|\n|_1}\;,
\eeqano
which concludes the proof, by the arbitrariness of $0<\vae<a$.
\qed
\rem{appDRe}
\begin{enumerate}
\item 
Actually, we proved something stronger. Namely, for any $\vae>0$, there exists $\d>0$ such that, for any $y,y'\in A $ 
and $\n\in \natural_0^d$, with $|y'-y|<\d$ and ${|\n|_1\le n}$,
\beq{whitfr}
\left|f_\n(y')-\dst\sum_{\substack{\m\in \natural_0^d\\ {|\m|_1\le n-|\n|_1}}} \frac{1}{\m!}f_{\n+\m}(y)(y'-y)^\m\right|\le \vae\; a\; |y'-y|^{n-|\n|_1}\;.
\eeq
\item 
In fact, $f$ satisfies the following ``uniform'' Whitney's condition, provided $n\ge 1$: for any $y,y'\in A $ and $\n\in \natural_0^d$, with $|y'-y|\le r_0$ and ${0\le |\n|_1\le n-1}$,
\beq{whitfr2}
\left|f_\n(y')-\dst\sum_{\substack{\m\in \natural_0^d\\ {|\m|_1\le n-1-|\n|_1}}} \frac{1}{\m!}f_{\n+\m}(y)(y'-y)^\m\right|\le a\left(2^n+2\ex^d\right)|y'-y|^{n-|\n|_1}\;.
\eeq
Indeed, for such given $y,y',\n$, let $m\in\natural_0$ such that $r_{m_3+1}<|y'-y|\le r_{m_3}.$ Then, by similar computations as above, we have
\beqano
\left|f_\n(y')-\dst\sum_{\substack{\m\in \natural_0^d\\ {|\m|_1\le n-1-|\n|_1}}} \frac{1}{\m!}f_{\n+\m}(y)(y'-y)^\m\right| &\le& |f_\n(y')-f^{[m_3]}_\n(y')|+\\
 &+&\left|f^{[m_3]}_\n(y')-\dst\sum_{\substack{\m\in \natural_0^d\\ {|\m|_1\le n-1-|\n|_1}}} \frac{1}{\m!}f_{\n+\m}^{[m_3]}(y)(y'-y)^\m\right|+\\
 &+&\left|\dst\sum_{\substack{\m\in \natural_0^d\\ {|\m|_1\le n-1-|\n|_1}}} \frac{1}{\m!}\left( f_{\n+\m}^{[m_3]}(y)-f_{\n+\m}(y)\right)(y'-y)^\m\right|\\
 &\le& \frac{2^n\;a}{(n-m)!}\;  |y'-y|^{n-|\n|_1}+\\
 &+&2\dst\sum_{\substack{\m\in \natural_0^d\\ {|\m|_1\le n-1-|\n|_1}}}\frac{1}{\m!} \dst\sum_{m> m_3} \|\dpr_{\n+\m} f_m\|_{ A}\;|y'-y|^{|\m|_1}\\
 &\le& \left(\frac{2^n\;a}{(n-m)!}+2\;a \dst\sum_{\m\in \natural_0^d}\su{\m!}\right)\;|y'-y|^{n-|\n|_1}\\
 &=& a\;\left({2^n}+2\ex^2\right)\;|y'-y|^{n-|\n|_1}\;.
\eeqano
\end{enumerate}
\erem
Finally, we recall the very deep Whitney's extension theorem.
\thmtwo{Whi34}{Whitney \cite{whitney1934analytic}}
Let $A\subset \rn$ and $f\in C^n_W(A,\real)$, $n\in\natural_0$. If $A$ is closed, then there exists $\bar{f}\in C^n(\rn,\real)$, real--analytic on $\rn\setminus A$ and such that $\bar{f}_\n=f_\n$ on $A$, for any $\n\in \natural_0^d$, with ${|\n|_1\le n}$.
\ethm
\newpage
\appF{Generalized Steiner's formula\label{appF}}
We aim here to recall the generalized Steiner's formula to compute the volume of the two halves tubes that composes a (uniform) tubular neighborhood of an embedded hypersurface without boundary of $\rn$.\footnote{For a genralization to a non--uniform tubular neighborhood, see \cite{roccaforte2013volume}.}\\
Let $\mathfrak{S}$ be a smooth, bounded, orientable hypersurface without boundary of $\rn$ (equipped with the Euclidean metric). Fix the orientation given by a smooth unit normal vector field of $\mathfrak{S}$
$$
\mathbf{n}=(\mathbf{n}_1,\cdots,\mathbf{n}_d)\colon \mathfrak{S}\to N\mathfrak{S}=(T\mathfrak{S})^{\perp}\;, \quad |\mathbf{n}|_2=
\mathbf{n}_1^2+\cdots+\mathbf{n}_d^2=1\;.
$$
Let $dy\coloneqq dy_1\wedge\cdots\wedge dy_d$ be the volume form on $\rn$ (which induces the Lebesgue--measure $\meas$ on $\rn$) and $\nabla$ be the Levi--Civita connection on $\rn$. Then, let $d{\mathfrak{S}}$ be the induced area--form on $\mathfrak{S}$, defined by
$$
d{\mathfrak{S}}(X_1,\cdots,X_{d-1})=dy(X_1,\cdots,X_{d-1},\mathbf{n})\;,
$$
for any $X_1,\cdots,X_{d-1}\in \G(\mathfrak{S})$, where $\G(\mathfrak{S})$ denotes the Lie algebra of smooth vector fields on $\mathfrak{S}$. Define the shape operator $\mathrm{S}\colon \G(\mathfrak{S})\to \G(\mathfrak{S})$ by
$$
\mathrm{S}X=-\nabla_X\mathbf{n}\;.
$$
Define the map $e_c\colon \{(y,u): y\in \mathfrak{S}\,, u=\pm \mathbf{n}(y)\}\to [0,\infty]$ by\footnote{$e_c(y,u)$ is the distance from $y\in\mathfrak{S}$ to its {\it cut-focal point} in the direction $u$ if such a cut--focal point exists; otherwise $e_c(y,u)=\infty$.}
$$
e_c(y,u)\coloneqq \dst\sup\{t>0:\dist(y+tu,\mathfrak{S})=t\}\;.
$$
Then, define the {\it minimal focal distance}
$$
\minfoc(\mathfrak{S})\coloneqq \dst\inf\{e_c(y,u):y\in \mathfrak{S}\,, u=\pm \mathbf{n}(y)\}\,.
$$
 Given $\d>0$, define the two half--tubes about $\mathfrak{S}$
 $$
\mathscr T^{\pm}(\mathfrak{S},\d)\coloneqq \{y\pm t\mathbf{n}(y): y\in\mathfrak{S},\, 0\le t\le \d\}
 $$
 and the $\d$--tubular neighborhood of $\mathfrak{S}$
 $$
 \mathscr T(\mathfrak{S},\d)\coloneqq \mathscr T^{+}(\mathfrak{S},\d)\bigcup \mathscr T^{-}(\mathfrak{S},\d)\;.
 $$
The contraction operators $\mathrm{C}^j$ on the space of double forms of type\footnote{A double form of type $(p,q)$ is a $\mathfrak{F}(\mathfrak{S})$--linear map $\L\colon \G(\mathfrak{S})^p\times \G(\mathfrak{S})^q\to \mathfrak{F}(\mathfrak{S})$, which is antisymmetric in the first $p$ variables and in the last $q$ as well, where $\mathfrak{F}(\mathfrak{S})$ denotes the space of smooth functions on $\mathfrak{S}$.} $(p,q)$ are defined inductively as follows: $\mathrm{C}^0(\L)=\L$ and, for $j\ge1$,
$$
\mathrm{C}^j(\L)(X_1,\cdots,X_{p-j})(Y_1,\cdots,Y_{q-j})=\dst\sum_{i=0}^{d-1}\mathrm{C}^{j-1}(\L)(X_1,\cdots,X_{p-j},E_i)(Y_1,\cdots,Y_{q-j},E_i)\;,
$$
where $\{E_1,\cdots,E_{d-1}\}$ is any orthonormal frame field of $\mathfrak{S}$. Let $\mathbf{R}^{\mathfrak{S}}$ be the {\it curvature tensor} of $\mathfrak{S}$. $\mathbf{R}^{\mathfrak{S}}$ being a double form of type $(2,2)$, one can then take the wedge product of $\mathbf{R}^{\mathfrak{S}}$ with itself $j$ times to get the double form $(\mathbf{R}^{\mathfrak{S}})^j$ of type $(2j,2j)$. Set $\mathrm{C}^0((\mathbf{R}^{\mathfrak{S}})^0))=1$ and define the $(2j)$--th and $(2j+1)$--th integrated mean curvatures of $\mathfrak{S}$ in $\rn$ as follows ($j\ge 0$):
\begin{align*}
\mathbf{k}_{2j}(\mathbf{R}^{\mathfrak{S}})&\coloneqq \frac{1}{j!(2j)!}\dst\int_\mathfrak{S}\mathrm{C}^{2j}((\mathbf{R}^{\mathfrak{S}})^j)\;d\mathfrak{S}\;,\\
\mathbf{k}_{2j+1}(\mathbf{R}^{\mathfrak{S}},\mathrm{S})&\coloneqq \frac{1}{j!(2j)!}\dst\int_\mathfrak{S}\left\{\tr(\mathrm{S})\,\mathrm{C}^{2j}((\mathbf{R}^{\mathfrak{S}})^j)-2j\tr\left(\mathrm{S}\mathrm{C}^{2j-1}((\mathbf{R}^{\mathfrak{S}})^j)\right)\right\} d\mathfrak{S}\;.
\end{align*}
Thus, the foolowing holds.
\thmtwo{arnolMeas1}{\cite{gray2012tubes}, pg.~224}
\beq{ArnMeq1}
\meas\left(\mathscr T^{\pm}(\mathfrak{S},\d)\right)= \dst\sum_{j=0}^{\left[\frac{d-1}{2}\right]}\frac{\mathbf{k}_{2j}(\mathbf{R}^{\mathfrak{S}})\;\d^{2j+1}}{1\cdot3\cdots (2j+1)}\mp\sum_{j=0}^{\left[\frac{d}{2}-1\right]}\frac{\mathbf{k}_{2j+1}(\mathbf{R}^{\mathfrak{S}},\mathrm{S})\;\d^{2j+2}}{1\cdot3\cdots (2j+1)(2j+2)}\;,
\eeq
for any $0\le \d\le \minfoc(\mathfrak{S})$.
\ethm
\newpage
\appG{Some others facts on Lipschitz continuous functions\label{appG}}
In the following, we prove that any set is contained in some enlargement of itself through any contracting mapping which is bounded on the former set.
\lem{LipRang}
Let $g\colon \cn\to\cn$ be Lipschitz continuous function. Assume that
\begin{align}
&\d\coloneqq \sup_{\rn}|g-\id|<\infty\;,\label{EqApG1}\\
&\|g-\id\|_{L,\rn}<1.\label{EqApG2}
\end{align}
Then, for any 
$\emptyset\neq A\subset\cn$,\footnote{Where $\ovl {D_{\d}( A)}$ denotes the closed $\d$--neighborhood of $A$ in $\cn$.} 
$$
A\subset g\left(\ovl {D_{\d}( A)}\right)\;.
$$
\elem
\proof Set $f\coloneqq g-\id$  and let $\bar y\in A$. It is enough to show that there exists $|y|\le \d$ such that $\bar y=g(y+\bar y)$ \ie $y=-f(y+\bar{y})$ \ie $y$ is a fixed point of the map 
$$h\colon \ovl{D_{\d}(0)}\ni y\mapsto -f(y+\bar{y}).$$
But, for any $y\in \ovl{D_{\d}(0)}$,
$$
|h(y)|=|f(y+\bar{y})|\le \|f\|_{\rn}\leby{EqApG1} \d\;,
$$
\ie $h\colon \ovl{D_{\d}(0)}\to \ovl{D_{\d}(0)}$. Moreover, $h$ is a contraction since $\|h\|_{L,\ovl{D_{\d}(0)}}\le \|f\|_{L,\rn}\ltby{EqApG2}1$. Thus, we can apply the Banach's fixed point Theorem to complete the proof.
\qed
\bibliographystyle{alpha}
\bibliography{BibtexDatabase}



\end{document}